\documentclass[11.5pt]{article}

\usepackage[noadjust]{cite}
\usepackage{amsfonts}
\usepackage{comment}
\usepackage{tabularx}
\usepackage[protrusion=true,expansion=true]{microtype}
\usepackage{enumerate}
\usepackage{bbm} 
\usepackage{paralist}  
\usepackage{mathrsfs}  
\usepackage{graphicx} 
\usepackage{tikz}
\usepackage{latexsym}
\usepackage{url} 
\usepackage[normalem]{ulem} 
\usepackage{amsmath, amsthm, amssymb, color, graphicx} 
\usepackage{hyperref}
\usepackage[margin=1in]{geometry}

\makeatletter
\DeclareFontFamily{U}{tipa}{}
\DeclareFontShape{U}{tipa}{m}{n}{<->tipa10}{}
\newcommand{\arc@char}{{\usefont{U}{tipa}{m}{n}\symbol{62}}}%

\newcommand{\arc@arc}[2]{%
	\sbox0{$\m@th#1#2$}%
	\vbox{
		\hbox{\resizebox{\wd0}{\height}{\arc@char}}
		\nointerlineskip
		\box0
	}%
}
\makeatother

\DeclareGraphicsExtensions{.jpg,.pdf,.eps}
\graphicspath{{./Figures/}}
\newtheorem{theorem}{Theorem}[section]

\newtheorem{corollary}[theorem]{Corollary}
\newtheorem{lemma}[theorem]{Lemma}
\newtheorem{definition}[theorem]{Definition}
\newtheorem{proposition}[theorem]{Proposition}

\newtheorem{conv}[theorem]{Convention}
\newtheorem{remark}[theorem]{Remark}

\newcommand{\arm}{\alpha}
\newcommand{\Mink}{\mathfrak{m}}

\newcommand{\markl}{\ell}
\newcommand{\Eb}{\mathrm{Eb}}

\newcommand{\Bol}{\op{Bol}}

\newcommand{\dsk}{\Map}

\newcommand{\dlp}{\Upsilon}
\newcommand{\BD}{\op{BD}}

\def\pv{\mathrm{piv}}

\def\constp{\mathfrak{c}}% quantum-pivotal-constant
\def\sle{\eta}

\newcommand{\fh}{\mathfrak{h}}

\newcommand{\notion}[1]{{\bf\textit{#1}}}
\newcommand{\xncomment}[1]{{#1}}

\newcommand{\old}[1]{{}} 

\newcommand\epsrhodelta{\omega^{\eps,\rho,\mesh}}
\newcommand\orhodelta{\omega^{0,\rho,\mesh}}
\newcommand\oodelta{\omega^{\mesh}}
\newcommand\epsodelta{\omega^{\eps,\mesh}}   
\newcommand\xor{\mathbin{\triangle}}
\newcommand{\bz}{{\mathbf z}}
\newcommand{\mpivep}{{\nu^\eps_{h,\Gamma}}}
\newcommand{\Pd}{{\P^{\op d}}}
\newcommand{\hd}{{h^{\op d}}}
\newcommand{\hw}{{h^{\op w}}}
\newcommand{\hs}{{h^{\op s}}}
\newcommand{\Dh}{{\mathcal{DH}}}
\newcommand{\CD}{{\mathcal{D}_{*,*}}}
\newcommand{\dbl}{{\op {dbl}}}
\newcommand{\area}{{\op {area}}}
\newcommand{\reg}{{\op {reg}}}

\newcommand{\GHPU}{{\mathrm{GHPU}}}
\newcommand{\GHPUL}{{\mathrm{GHPUL}}}
\def\BM{\BB M}

\newcommand{\Map}{{\mathsf{M}}}

\newcommand{\bh}{\mathbf{h}}

\newcommand{\C}{\mathbb{C}}

\newcommand{\1}{\mathbf{1}}
\newcommand{\D}{\mathbb{D}}
\newcommand{\E}{\mathbb{E}}
\newcommand{\N}{\mathbb{N}}
\newcommand{\Q}{\mathbb{Q}}
\newcommand{\Z}{\mathbb{Z}}
\newcommand{\R}{\mathbb{R}}
\renewcommand{\P}{\mathbb{P}}
\newcommand{\bbH}{\mathbb{H}}
\newcommand{\ol}{\overline}
\newcommand{\wt}{\widetilde}

\newcommand{\eps}{\varepsilon}
\newcommand{\mesh}{\delta}
\def\defeq{:=}
\def\eqd{\overset{d}{=}}
\newcommand{\p}{\parital}
\def\P{\mathbb{P}}
\def\E{\mathbb{E}}
\DeclareMathOperator{\Cov}{Cov}
\DeclareMathOperator{\Cdy}{Cdy}
\DeclareMathOperator{\SLE}{SLE}
\DeclareMathOperator{\CLE}{CLE}
\DeclareMathOperator{\LQG}{LQG}  
\DeclareMathOperator{\Var}{Var} 
\def\p{\partial}

\def\bT{\mathbb{T}}

\def\cW{\mathcal{W}}
\def\cV{\mathcal{V}}

\def\cT{\mathcal{T}}
\def\cS{\mathcal{S}}
\def\cR{\mathcal{R}}
\def\cQ{\mathcal{Q}}
\def\cP{\mathcal{P}}

\def\cN{\mathcal{N}}
\def\cM{\mathcal{M}}
\def\cL{\mathcal{L}}

\def\cI{\mathcal{I}}
\def\cH{\mathcal{H}}

\def\cF{\mathcal{F}}
\def\cE{\mathcal{E}}
\def\cD{\mathcal{D}}
\def\cC{\mathcal{C}}
\def\cB{\mathcal{B}}
\def\cA{\mathcal{A}}

\def\Tg{\mathbb{T}} %%Triangular grid

\def\@rst #1 #2other{#1}
\newcommand{\arxiv}[1]{\href{http://arxiv.org/abs/#1}{#1}}
\def\MR#1{}

\newcommand{\aryb}{\begin{eqnarray*}}
	\newcommand{\arye}{\end{eqnarray*}}
\def\alb#1\ale{\begin{align*}#1\end{align*}}
\newcommand{\eqb}{\begin{equation}}
\newcommand{\eqe}{\end{equation}}
\newcommand{\eqbn}{\begin{equation*}}
\newcommand{\eqen}{\end{equation*}}

\newcommand{\Piv}{{\op{PV}}}
\newcommand{\bv}{{\mathbf v}}
\newcommand{\BB}{\mathbbm}
\newcommand{\op}{\operatorname}
\newcommand{\bd}{\mathbf}

\newcommand{\frk}{\mathfrak}
\newcommand{\eqD}{\overset{d}{=}}

\newcommand{\ep}{\varepsilon}
\newcommand{\rta}{\rightarrow}

\newcommand{\wh}{\widehat} 
\newcommand{\mcl}{\mathcal}

\newcommand{\bdy}{\partial}
\newcommand{\Ber}{\operatorname{Ber}}
\newcommand{\pivm}{\cM_{\bh,\Gamma}}

\begin {document}
\title{Convergence of uniform triangulations under the Cardy embedding}
\author{
	\begin{tabular}{c} Nina Holden\\ [-5pt] \small ETH Z\"urich \end{tabular}
	\begin{tabular}{c} \\[-5pt]\small  \end{tabular}
	\begin{tabular}{c} \\[-5pt]\small  \end{tabular}
	\begin{tabular}{c} \\[-5pt]\small  \end{tabular}
	\begin{tabular}{c} Xin Sun\\ [-5pt] \small University of Pennsylvania \end{tabular}
}
\date{}
\maketitle

\begin{abstract}
	We consider an embedding of  planar maps into an equilateral triangle $\Delta$ which we call the Cardy embedding. 
	The embedding is a discrete approximation of a conformal map based on percolation observables that are used in Smirnov's proof of Cardy's formula.
	Under the Cardy embedding, the planar map induces a metric and an area measure on $\Delta$ and a boundary measure on $\partial \Delta$. 
	We prove that for uniformly sampled triangulations, 
	the metric and the measures converge jointly in the scaling limit to the Brownian disk conformally embedded into $\Delta$ (i.e., to the $\sqrt{8/3}$-Liouville quantum gravity disk). 
	As part of our proof, we prove scaling limit results for critical site percolation on the uniform triangulations,   in a quenched sense. 
	 In particular, we  establish the scaling limit of the percolation crossing probability for a uniformly sampled triangulation with four boundary marked points.
\end{abstract}

\section{Introduction}\label{emb-label:intro}

Random planar geometry has been a central topic in probability in the last two decades.
\xncomment{One of the main goals} is to construct and study random surfaces. A natural approach is to consider the scaling limit of random planar maps.  
Inspired by Riemannian geometry, a natural point of view is to consider a  planar map as an abstract metric measure space. In this regards, Le Gall \cite{legall-uniqueness}, Miermont \cite{miermont-brownian-map}, and others (e.g.\ \cite{bjm-uniform,abraham-bipartite,ab-simple,beltran-legall-pendant}) proved that a large class of uniformly sampled random planar maps converge in the scaling limit 
to a random  metric measure space with  the topology of the sphere,  known as \notion{the Brownian map}.  In the case where the random planar map has a macroscopic boundary, the scaling limit is  the \notion{Brownian disk} \cite{bet-mier-disk}, which is a metric measure space with the topology of a disk. 

\notion{Liouville quantum gravity} (LQG) is another approach for constructing a random surface, which takes the perspective of conformal geometry.
Since the foundational work of Polyakov \cite{polyakov-qg1}, LQG has been an active research area in theoretical physics. 
The mathematical study of LQG was initiated by Duplantier and Sheffield \cite{shef-kpz}.
The idea is to  consider an instance $h$ of the Gaussian free field (GFF) on a planar domain $D$ and study the surface with volume measure $e^{\gamma h}d^2z$.
This definition does not make rigorous sense since $h$ is a distribution and not a function. However, by
first regularizing $h$ and then taking a limit, for each $\gamma\in (0,2)$, the random area measure $\mu_h\defeq e^{\gamma h}\,d^2z$ on $D$ is well defined and nontrivial. If $D$ has a nontrivial boundary, the measure $\xi_h\defeq e^{\gamma h/2}dz$  on $\p D$ can also be defined.
Very recently, Ding, Dubédat, Dunlap, and Falconet \cite{DDDF19} and Gwynne and Miller \cite{GM-metric}  proved that one may construct a unique metric (i.e., a distance function) $d_h$ by regularizing the metric tensor $e^{2\gamma h/\dim_\gamma}(dx^2+dy^2)$, where $\dim_\gamma$ is the Hausdorff dimension of the surface \cite{gp-kpz}.
For $\gamma=\sqrt{8/3}$, this metric agrees with the metric constructed earlier
by Miller and Sheffield \cite{lqg-tbm1,lqg-tbm2,lqg-tbm3}, which gives a metric space with the law of a Brownian surface, \xncomment{namely the random metric measure spaces which describe the scaling limits of uniform planar maps as mentioned in the previous paragraph}. 
There is a coordinate change rule depending on $\gamma$ 
that relates fields on two conformally equivalent domains such that $(d_h, \mu_h,\xi_h)$ is invariant under conformal maps. The random geometry defined by $(h,d_h,\mu_h,\xi_h)$  is called $\gamma$-LQG. \xncomment{We refer to Section~\ref{subsec:lqg} for more details.}

A fundamental belief in random planar geometry which has been guiding  its development is the following. 
Given any $\gamma\in(0,2)$, there is a family of random planar maps whose scaling limit under discrete conformal embeddings is $\gamma$-LQG.
In particular, uniform  random planar maps converge to $\sqrt{8/3}$-LQG in this sense. Here a discrete conformal embedding means a discrete approximation of the Riemann mapping.
Notable examples include the circle packing and the Tutte embedding. See e.g.\ \cite{shef-kpz,legall-icm,dkrv-lqg-sphere} for precise conjectures. Before the current paper, this convergence had not been verified for any natural combinatorial random planar maps under any discrete conformal embedding.  See Section~\ref{subsec:outlook} for results on planar maps obtained from coarse graining of a $\gamma$-LQG surface.

\xncomment{As pointed out in in~\cite{lsp94}, it was conjectured by Aizenman that critical planar percolation is conformally invariant. This conjecture was checked numerically for the crossing probability in~\cite{lsp94}.}
%Based on Aizenman's \xncomment{suggestion}, it  was conjectured \cite{lsp94} that the crossing probability for critical planar percolation is conformally invariant.
Cardy \cite{cardy-formula} then predicted an explicit formula for the left/right crossing probability for rectangles of any aspect ratio.
Cardy's formula was proved by Smirnov \cite{smirnov-cardy} in the case of site percolation on the triangular
lattice. 
A by-product of Smirnov's proof is a discrete conformal embedding based on percolation observables, which we call the \notion{Cardy embedding} (see Definition~\ref{emb-def:cardy}). 
In this paper, we prove that  large uniform triangulations converge to $\sqrt{8/3}$-LQG under the Cardy embedding (see Theorem~\ref{emb-thm:Cardy}).

\begin{comment}
Smirnov's proof of Cardy's formula and 
Schramm's discovery of the Schramm-Loewner evolution (SLE) \cite{schramm0}
mark the beginning of a range of works which greatly improved our understanding of the scaling limit of critical percolation on the triangular lattice
\cite{camia-newman-sle6,smirnov-werner-percolation,ss-planar-perc,  gps-pivotal,hlls-pivot}.
Smirnov's proof is famously difficult to adapt to percolation in other settings \cite{beffara-easy}, even for bond percolation on $\Z^2$.
In this paper we prove that in  the random environment defined by uniform triangulations, critical site percolation has a \emph{quenched} scaling limit (see Theorems~\ref{thm:crossings} and~\ref{emb-thm:main}). To our knowledge this is the first full quenched scaling limit result for critical percolation beyond site percolation on the triangular lattice.
The only other quenched scaling limit result we are aware of is for the crossing probability of squares for  Poisson Voronoi  percolation \cite{agmt16}.
\end{comment}

This paper is the culmination of a seven-paper research program which also includes \cite{hlls-cut,hlls-pivot,bhs-site-perc,aasw-type2, ghs-metric-peano,ghss18}. 
Other papers that are important to this program include \cite{gps-fourier,gps-pivotal,gps-dynamic,wedges,gwynne-miller-perc}. 
See Section~\ref{subsec:outline} for an overview of the program  and an outline of this paper.

\subsection{The Cardy embedding as a discrete conformal embedding}\label{subsec:Cardy}
The Riemann mapping theorem asserts that any  two simply connected planar domains with boundary  are related by a conformal map.
The Riemann mapping admits natural discrete approximations which we call \emph{discrete conformal embeddings}. 
As a notable example, Thurston conjectured that the circle packing gives an approximation of the Riemann mapping from a simply connected domain to the unit disk.
This conjecture was proved  by Rodin and Sullivan\cite{rs87}.

\begin{figure}
	\centering
	\includegraphics[scale=1]{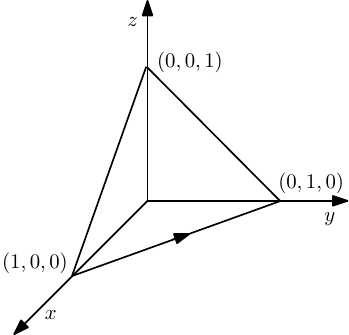}\qquad\qquad\includegraphics[]{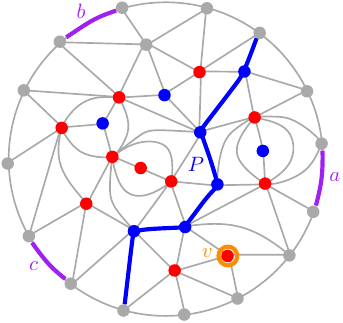}
	\caption{{\bf Left}: Illustration of $\Delta$ as an oriented surface with disk topology. The arrow indicates the counterclockwise orientation of $\p\Delta$.
		{\bf Right}: Illustration of the event $E_a(v)$.}
	\label{fig:crossing-event-and-delta}
\end{figure}

Consider the equilateral triangle $\Delta \defeq\{(x,y,z):x+y+z=1,\,x,y,z>0 \}$.
We view $\Delta$ as an oriented surface with disk topology and boundary $\p\Delta$ where the orientation is such that $(1,0,0)$, $(0,1,0)$, and $(0,0,1)$ are ordered counterclockwise. 
See Figure~\ref{fig:crossing-event-and-delta} for an illustration.
Given a Jordan domain $D$ with three distinct boundary points $a,b,c$ in counterclockwise order,
there exists  a unique Riemann mapping from $D$ to $\Delta$ that maps $a$, $b$, and $c$ to $(1,0,0)$, $(0,1,0)$, and $(0,0,1)$, respectively.  
We  denote this mapping  by $\Cdy_D$. The dependence on $(a,b,c)$ is dropped to lighten the notation.
Smirnov's elegant proof of Cardy's formula  provides   an approximation scheme for $\Cdy_D$
based on percolation observables.  This gives another example of a discrete conformal embedding
which we call  the \emph{Cardy embedding}.

We now define the Cardy embedding in the general setting of  triangulations of polygons.
Recall that a planar map is a planar graph (multiple edges and self-loops allowed) embedded into the
sphere, viewed modulo orientation-preserving homeomorphisms. For a planar map $M$, we write $\mcl V(M)$, $\mcl E(M)$, and $\mcl F(M)$ for the set of vertices, edges, and faces, respectively. 
A map is \emph{rooted} if one of its edges, called the \emph{root edge}, is 
distinguished and oriented.  The face to the right of the root edge is called the \emph{root face}. Given an integer $\ell\ge 2$, a rooted planar map $M$ is called \emph{a triangulation with boundary length $\ell$} if every face in $\cF(M)$  has degree 3, except the root face, which has degree $\ell$. We write $\partial M$ for the graph consisting of the edges and vertices on the root face of $M$.  
A vertex on  $M$ is called a \emph{boundary vertex} if it is on $\bdy M$. Otherwise, it is called an \emph{inner vertex}. We similarly define \emph{boundary edges} and \emph{inner edges}. 
If $\p M$ is simple, namely, consists of $\ell$ \emph{distinct}  boundary vertices, we say that $M$ is a \emph{triangulation of an $\ell$-gon}.
Let $\frk T(\ell)$ be the set of triangulations of an $\ell$-gon and define $\frk T\defeq\cup_{\ell\ge 2} \frk T(\ell)$. 
We call an element in $\frk T$ a \notion{triangulation of a polygon}. 

Given $M\in \frk T$, a \emph{site percolation} on $M$ is a coloring of $\cV(M)$ in two colors, say, red and blue. 
The \emph{Bernoulli-$\frac12$ site percolation} on $M$ is the random site percolation $\omega$ on $M$ such that each inner vertex is independently colored red or blue with equal probability. 
The coloring of the boundary vertices is called the \emph{boundary condition} of $\omega$ and can be prescribed arbitrarily.

Given a triangulation of a polygon $M$ with  three distinct boundary edges $a,b,c$ ordered counterclockwise,
we denote by $(a,b)$ the set of boundary  vertices of $M$ situated between $a$ and $b$ in counterclockwise order  (including one endpoint of $a$ and one endpoint of $b$).
Define $(b,c)$ and $(c,a)$ similarly. For a vertex $v\in \cV(M)$, let
$E_a(v)$ be the event that there exists a simple path 
(i.e., a sequence of distinct vertices on $M$ where any two  consecutive vertices are adjacent) $P$ on $M$ such that 
\begin{itemize}
	\item[(a)] $P$ contains one endpoint in $(c,a)$ and one endpoint in $(a,b)$, while all other vertices of $P$ are inner blue vertices;  
	\item[(b)] either $v\in P$ or $v$ is on the same side of $P$ as the edge $a$.
\end{itemize}
See Figure \ref{fig:crossing-event-and-delta} for an illustration.
We  define the events $E_b(v)$ and $E_c(v)$ similarly. Note that $E_a(v)$, $E_b(v)$, and $E_c(v)$ do not depend on the boundary condition of $\omega$.

%We  define the events $E_b(v)$ and $E_c(v)$ similarly. Note that $E_a(v)$, $E_b(v)$, and $E_c(v)$ do not depend on the boundary condition of $\omega$.

Given any nonnegative vector $(x,y,z)\in [0,\infty)^3$, let $(x,y,z)_{\Delta}\defeq (x+y+z)^{-1}(x,y,z)$, with the convention that $(0,0,0)_\Delta\defeq(1/3,1/3,1/3)$.
In other words, $(x,y,z)_\Delta$ is the projection of  $(x,y,z)$ onto the equilateral triangle  $\Delta$ along its own direction.   
The Cardy embedding is a mapping from the vertex set of a triangulation of a polygon to the closed triangle $\ol\Delta\defeq\Delta\cup\p \Delta$, defined using 
observables of site percolation on top of it.
\begin{definition}[Cardy embedding]\label{emb-def:cardy}
	Given a triangulation of a polygon $M$ with  three distinct boundary edges $a,b,c$ ordered counterclockwise,
	let $\Ber_M$ be the probability measure corresponding to the Bernoulli-$\frac12$ site percolation on $M$.
	The  Cardy embedding $\Cdy_M$ of $(M,a,b,c)$ is the function from $\cV(M)$ to  $\ol\Delta$ given by 
	$$\Cdy_M(v)= (\Ber_M [E_a(v)],\Ber_M[E_b(v)],\Ber_M[E_c(v)])_\Delta\qquad \textrm{for all }v\in\cV(M).$$ 
\end{definition}

Smirnov's theorem \cite{smirnov-cardy} can be phrased in terms of the Cardy embedding as follows.
Suppose $D$ is a Jordan domain with three distinct marked boundary points $a,b,c$ ordered counterclockwise. 
Let $\Tg$ denote the triangular lattice. Given a small mesh size  $\mesh>0$, let $D^\mesh$ be a lattice approximation of $D$ via $\mesh\Tg$ such that $D^\delta$ is a triangulation of a polygon (see Section~\ref{subsec:notation} for a precise definition).
Let $a^\mesh,b^\mesh,c^\mesh$ be points on $\p D^\mesh$ that approximate $a,b,c$, respectively. Let $\Cdy^\delta$ be the Cardy embedding of $(D^\mesh,a^\mesh,b^\mesh,c^\mesh)$ and recall the Riemann mapping $\Cdy_D$ from $D$ to $\Delta$ defined above.
\begin{theorem}[Smirnov] \label{thm:Smirnov}
	In the setting above,\footnote{Smirnov's definition of crossing probabilities is slightly different from ours, but the difference between the definitions is negligible in the scaling limit.}
	$$
	\lim_{\delta\rta 0}\sup_{v\in D^\mesh}\big|\Ber_{D^\mesh} [E_{a^\delta}(v)]+\Ber_{D^\mesh}[E_{b^\delta}(v)]+\Ber_{D^\mesh}[E_{c^\delta}(v)]-1\big| = 0
	$$
	and
	$$
	\lim_{\delta\rta 0}\sup_{v\in\cV (D^\mesh)} | \Cdy^\delta(v)-\Cdy_D(v)|=0.
	$$
\end{theorem}
In Definition~\ref{emb-def:cardy},  let $e$ be an edge lying on the arc $(c,a)$ and let $v$ be the endpoint of $e$ closer to $a$.
Then $\Ber_{M}[E_a(v)]$ is the  so-called \emph{crossing probability} between $(c,e)$ and $(a,b)$. 
Let $D=[0,R]\times [0,1]$ for some $R>0$ and let the marked boundary points of $D$ be $(R,0)$, $(R,1)$, and $(0,1)$. 
By Theorem \ref{thm:Smirnov}, the $x$ coordinate of $\Cdy_D(0,0)$ is the $\mesh\to0$ limit of the crossing probability between the left and right sides of $D^\mesh$. By the Schwarz-Christoffel formula, the value of $\Cdy_D(0,0)$ can be expressed explicitly as a function of $R$, which agrees with Cardy's formula for this crossing probability in~\cite{cardy-formula}. 
Therefore Theorem~\ref{thm:Smirnov} gives a rigorous proof of Cardy's formula, which explains why we call our embedding the Cardy embedding.

\subsection{Main result}\label{emb-subsec:Cdy-conv}
\subsubsection{Scaling limit of uniform triangulations under the Cardy embedding}\label{subsub:scaling} 
Our main result is that large uniform triangulations of polygons converge to $\sqrt{8/3}$-LQG under the Cardy embedding.
We will focus on a particular variant where self-loops are not allowed  while  multiple-edges are allowed; these are often called type II triangulations of polygons. See Remark~\ref{rmk:type} for extensions to other variants. 
We consider the critical Boltzmann measure, which is defined as follows. For $\ell\ge 3$, let $\frk T_2(\ell)$ be the set of maps in $\frk T(\ell)$ with no self-loops (but multiple-edge are allowed).
Given  $\ell\ge 3$, it is well-known that  if each element $M\in \frk T_2(\ell)$ is assigned weight $(2/27)^n$, where $n$ is the number of vertices of $M$,
then the resulting measure on $\frk T_2(\ell)$ is  finite.  Let $\Bol_2(\ell)$ be the probability measure obtained by normalizing this measure. 
Following \cite{angel-schramm-uipt}, we call a map with law  $\Bol_2(\ell)$ a \notion{Boltzmann triangulation} of type $\op{II}$ with boundary length $\ell$.

Fix a sequence of integers $\{\ell^n\}_{n\in\BB N}$ such that $\ell^n\ge 3$ for all $n\in\BB N$ and $(3n)^{-1/2} \ell^n \to 1$  as $n\rta\infty$. 
Let $\dsk^n$  be sampled from $\Bol_2(\ell^n)$.
Denote the root edge of $\dsk^n$ by $a^n$ and sample two other boundary edges $b^n$ and $c^n$ uniformly at random,
conditioning on $a^n,b^n,c^n$ being distinct and ordered counterclockwise.  
Let $d^{\op{gr}}_{\dsk^n}:\cV(\dsk^n)\times\cV(\dsk^n)\to\N\cup\{0 \}$ be the graph distance  of $\dsk^n$ and define $d^n\defeq(3n/4)^{-1/4}d^{\op{gr}}_{\dsk^n}$.
Let $\mu^n$ be $(2n)^{-1}$ times the counting measure on $\cV(\dsk^n)$. Let $\xi^n$ be $1/\ell^n$ times the counting measure on $\cV(\p \dsk^n)$.
We obtain a random compact metric space  endowed with two measures, which we denote by $\cM^n=(\Map^n,d^n,\mu^n,\xi^n)$.
In collaboration with Albenque \cite{aasw-type2}, we proved that $\cM^n$ converge in law  to a variant of the Brownian disk called the \emph{free Brownian disk  with unit  perimeter}, which we  denote by $\BD_1$ (see Theorem~\ref{thm:AHS}).  Moreover, the marked edges $(a^n,b^n,c^n)$ converge to three marked points on the boundary of $\BD_1$.
By works of Miller and Sheffield  \cite{lqg-tbm1,lqg-tbm2,lqg-tbm3} \xncomment{(see Corollary 1.5 of \cite{lqg-tbm2})}, 
there exists a variant $h_\Delta$ of the Gaussian free field on $\Delta$ such that 
\xncomment{$$(\ol\Delta,d_{\Delta},\mu_{\Delta},\xi_{\Delta})\defeq(\ol\Delta,c_{\op d}d_{h_\Delta}, c_{\op m}\mu_{h_\Delta},\xi_{h_\Delta})$$} 
has the law of $\BD_1$ with the three marked points  being $(1,0,0)$, $(0,1,0)$, and $(0,0,1)$.  
Here $(d_{h_\Delta}, \mu_{h_\Delta}, \xi_{h_\Delta})$  is the metric/measure triple in $\sqrt{8/3}$-LQG corresponding to $h_\Delta$ as mentioned above Section~\ref{subsec:Cardy},
and  $c_{\op d},c_{\op m}$ are implicit positive constants coming from Miller and Sheffield's theorem. See Theorem~\ref{thm:LQG-TBM}  and Definition~\ref{def:hdelta} for precise definitions.

Let $\Cdy^n$ be the Cardy embedding of $(\dsk^n,a^n,b^n,c^n)$.  
Now we define a triple $(d^n_\Delta,\mu^n_\Delta,\xi^n_\Delta)$ which is the pushforward of $\cM^n$ onto $\ol\Delta$ 
under $\Cdy^n$. 
To be precise, for $x\in \ol\Delta$, let $\frk v(x)$ be the vertex of $\dsk^n$ which is closest to $x$ under the Cardy embedding, i.e., we let $\frk v(x)$ be the vertex $v\in\cV(\dsk^n)$ such that $|\Cdy_{\Map^n}(v)-x |$ is minimized over $v\in\cV(\Map^n)$; if there is a tie we resolve it in some arbitrary way.
Let\footnote{By Theorem~\ref{emb-thm:CLE} and~\eqref{eq:max0}, the measure $\xi^n_\Delta$ concentrates near $\p \Delta$, although we view it as a measure on $\ol\Delta$.}
\begin{align} 
	d^n_\Delta(x,y) &\defeq d^n(\frk v(x),\frk v(y)),&& \textrm{for }x,y\in\ol\Delta, \label{eq:graph-dist}\\
	\mu^n_\Delta(U) &\defeq \mu^n\left(\{ v\in \cV(\dsk^n)\,:\,\Cdy_{\Map^n}(v)\in U \}\right), && \textrm{for each Borel set }U\subset\ol\Delta,\nonumber\\
	\xi^n_\Delta(U) &\defeq   \xi^n\left(\{ v\in \cV(\p\dsk^n)\,:\,\Cdy_{\Map^n}(v)\in U \}\right), &&\textrm{for each Borel set }U\subset \ol\Delta.\nonumber
\end{align}
Our main result can be stated as follows.
\begin{theorem}\label{emb-thm:Cardy}
	In the setting above, $(d^n_\Delta,\mu^n_\Delta,\xi^n_\Delta)$ converge jointly in law to $(d_\Delta,\mu_\Delta,\xi_\Delta)$ as $n\rta\infty$, where we equip the first coordinate with the uniform topology  and the latter two coordinates with the Prokhorov topology on Borel measures on $\ol\Delta$.  
\end{theorem}
To draw an analogy with Theorem~\ref{thm:Smirnov},  Theorem~\ref{emb-thm:Cardy} asserts that the Cardy embedding of $\dsk^n$  provides a discretization of the conformal embedding of the Brownian disk onto $\ol\Delta$.      

\begin{figure}
	\centering
	\includegraphics[scale=1]{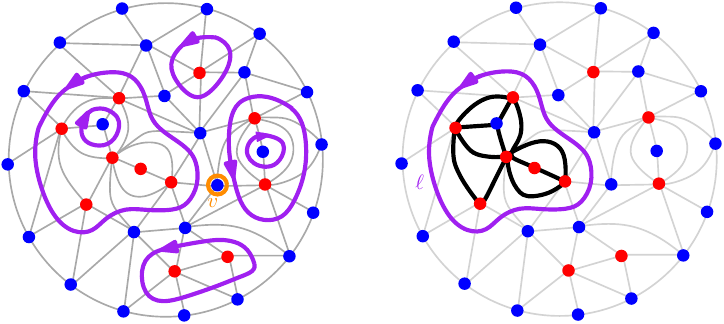}
	\caption{{\bf Left}: The loop ensemble $\Gamma(M,\omega)$ of the percolation $\omega$ is shown in purple. The vertex $v$ marked in orange is a pivotal point since two of the loops will be joined together to form a longer loop if the color of $v$ is flipped.
	{\bf Right}: The subgraph with bold edges is reg$(\ell)$ as defined in Section \ref{subsub:DP}.
	}
	\label{fig:delta}
\end{figure}

Theorem~\ref{emb-thm:main} still holds under slight modifications to the definition of the Cardy embedding in Definition~\ref{emb-def:cardy}. 
For example, by Proposition~\ref{prop:max2}, we have the  following analogue of the first equation of Theorem \ref{thm:Smirnov}: 
\begin{equation}\label{eq:max0}
	\max_{v\in \cV(\Map^n)}|\Ber_{\Map^n}(E_{a^n}(v))+\Ber_{\Map^n}(E_{b^n}(v))+\Ber_{\Map^n}(E_{c^n}(v))-1|=o_n(1). 
\end{equation}
Therefore the projection $(\cdot,\cdot,\cdot)_\Delta$  in Definition~\ref{emb-def:cardy} is not essential.
We can also modify some details in the definition of $E_a(v)$, such as letting $a,b,c$ be vertices instead of edges, or requiring that $v$ does not lie on $P$.
Using ideas from a recent alternative proof of Cardy's formula on the triangular lattice \cite{khristoforov-thesis}, it is possible to modify in such a way that the three crossing probabilities in~\eqref{eq:max0} always sum to exactly 1.

\subsection{Quenched scaling limits of site percolation}\label{subsub:quenched}

We prove  Theorem~\ref{emb-thm:Cardy} by establishing \emph{quenched} scaling limit results for site percolation on uniform triangulations.
To explain what we mean by quenched, let us start by considering the simplest percolation observable, namely the crossing probability between two boundary arcs. 
Let $(\Map^n,a^n,b^n,c^n)$ and $h_\Delta$ be as in Theorem~\ref{emb-thm:Cardy}.
Conditioning on  $(\Map^n,a^n,b^n,c^n)$, uniformly sample an edge $e^n$ on the arc $(c^n,a^n)$ and let $v^n$ be the endpoint of $e^n$ which is closer to $a^n$.
By the discussion below Theorem~\ref{thm:Smirnov}, $\Ber_{\Map^n} [E_{a^n}(v^n)]$ is the crossing probability between the arcs $(c^n,e^n)$ and $(a^n,b^n)$.
In the continuum, let $\mathbf v$ be a point on the counterclockwise arc on $\p \Delta$ from $(0,0,1)$ to $(1,0,0)$ sampled according to the measure 
$\xi_\Delta$ on $\p \Delta$ restricted to this arc. In other words,
$\mathbf v$ is a random point on this arc such that conditioning on $h_\Delta$, the ratio between the $\xi_\Delta$-masses of the counterclockwise arcs from  $(0,0,1)$ to $v$ and the one from $(0,0,1)$ to $(1,0,0)$
is uniformly distributed between $0$ and $1$. Let $x(\mathbf v)$ be the $x$-coordinate of $\mathbf v$.  Then we have the following.
\begin{theorem}
	In the setting described above,  $\Ber_{\Map^n} [E_{a^n}(v^n)]$ converge in law to $x(\mathbf v)$. 
	\label{thm:crossings}
\end{theorem}
It is clear from Theorem~\ref{thm:crossings} that the following more symmetric looking variant  holds. Let  $(e^n_1,e^n_2,e^3_n,e^n_4)$ be four uniformly sampled edges on $\p\Map^n$, conditioning on the edges being distinct and ordered counterclockwise. Then the crossing probability between the arcs $(e^n_1,e^n_2)$ and $(e^n_3,e^n_4)$ converge in law to a random variable, whose law is straightforward to describe in terms of the measure $\xi_\Delta$. We skip a more formal statement to avoid extra notations.

Earlier scaling limit results for percolation on random planar maps have considered observables involving 
\emph{both} the randomness of the planar map and the percolation. This includes for example \cite{gwynne-miller-perc,bhs-site-perc,curien-kortchemski-looptree-perc,angel-uihpq-perc} and Theorem~\ref{emb-thm:upper} below. In the context of random processes in random environment, \xncomment{these types} of statements are referred as \notion{annealed} scaling limit results.   
Alternatively, we can consider percolation observables which are  functions only of the environment, in our case, the underlying planar map. The crossing probability $\Ber_{\Map^n} [E_{a^n}(v^n)]$ in Theorem \ref{thm:crossings} is an example of such an observable. Convergence of such observables are referred to as \notion{quenched} scaling limit results.

Smirnov's proof of Cardy's formula is famously difficult to adapt to percolation in other settings \cite{beffara-easy}, even for bond percolation on $\Z^2$. To our best knowledge, this paper is the first work where quenched scaling limit results for percolation on random planar maps are established. Even for general environments beyond the triangular lattice, the only other quenched scaling limit result we are aware of is for the crossing probability of squares
for Poisson Voronoi  percolation~\cite{agmt16}. We also note that a variant of Theorem~\ref{thm:crossings} with $\SLE_6$ in place of percolation is stated in \cite{curien-glimpse}  as a theorem conditional on an unproven assertion.

There is a  close relationship between quenched scaling limit results and the convergence of certain embeddings, which is well known in the context of random walk in random environment. There the embedding is the so-called  Tutte embedding. 
See \cite{berger-biskup-perc-rw,gms-tutte} and reference therein. Our proof of Theorem~\ref{emb-thm:Cardy} is also based on this connection. More precisely, 
by the disk variant of Le Gall~\cite{legall-uniqueness}, Miermont~\cite{miermont-brownian-map} (see Theorem~\ref{thm:AHS}), and Miller-Sheffield~\cite{lqg-tbm1,lqg-tbm2}, there exists a sequence of embeddings $\{\Eb_n\}$ of $\Map^n$ to $\Delta$ such that Theorem~\ref{emb-thm:Cardy} holds with $\Cdy_{\Map^n}$ replaced by $\Eb_n$. One example of $\{\Eb_n\}$ can be obtained from the framework of mating of trees~\cite{bhs-site-perc,ghs-metric-peano}. However, the embeddings $\{\Eb_n\}$ are rather implicit and a priori do not carry any information about the conformal structure of $\Map^n$. 
Our approach to Theorem~\ref{emb-thm:Cardy} can be understood as first proving that
under the random environment obtained by embedding $\Map^n$ via $\Eb_n$, the critical site percolation has a quenched scaling limit as if the  environment is just the regular triangular lattice.
Then since the Cardy embedding is defined via percolation observables,   the difference between $\Eb_n$ and $\Cdy_{\Map^n}$ must vanish as $n\to \infty$, hence Theorem~\ref{emb-thm:Cardy} follows. In Section~\ref{subsub:multi}, we formulate a variant of this approach without introducing the extra embeddings $\{\Eb_n\}$.

\subsubsection{Scaling limit of multiple site percolations  on uniform triangulations}\label{subsub:multi}
Recall $\Map^n$ in Theorem~\ref{emb-thm:Cardy}.
Conditioning on $\Map^n$, let $\{\omega^n_i\}_{i\in\N}$ be a sequence of independent samples from $\Ber_{\Map^n}$.
In this section we  formulate a convergence result for $\{(\Map^n,\omega^n_i)\}_{i\in\N}$ (Theorem~\ref{emb-thm:main}) which is sufficient for the proof of Theorem~\ref{emb-thm:Cardy}.

Recall that $\Map^n$  is sampled from $\Bol_2(\ell^n)$ and has a root edge denoted by $a^n$. %Also recall that $\cM^n=(\Map^n,d^n,\mu^n,\xi^n)$. In Section~\ref{subsub:scaling},  $\xi^n$ is viewed as the uniform measure on $\cV(\p \Map^n)$. In this section, instead of a measure, we think of $\xi^n$ as a  curve of duration $[0,1]$, tracing $\p \Map^n$ clockwise starting and ending at $a^n$. 
Also recall that $\cM^n=(\Map^n,d^n,\mu^n,\xi^n)$, where $\xi^n$ is the uniform measure on $\cV(\p \Map^n)$. In this section, instead of a measure $\xi^n$, we consider a curve $\vec\xi^n$ of duration $[0,1]$, tracing $\p \Map^n$ clockwise starting and ending at $a^n$ such that each boundary edge is traced $1/\ell^n$ units of time.
This way, we view $\cM^n$ as a compact metric measure space decorated by a curve.  The natural topology for such objects is the so-called  \emph{Gromov-Hausdorff-Prokhorov-uniform} (GHPU) topology, which is introduced in \cite{gwynne-miller-uihpq}.  It is the natural variant of the Gromov-Hausdorff topology for spaces which are also equipped with a measure and a curve. In the continuum, the free Brownian disk with unit perimeter $\BD_1$ can also be naturally viewed as a compact metric measure space decorated by a curve. See Section~\ref{subsec:pre} for more details on the GHPU topology and the Brownian disk.

With Albenque, we proved the following.
\begin{theorem}[\cite{aasw-type2}] \label{thm:AHS}
	$\cM^n$ converge in law  to $\BD_1$ in the GHPU topology as $n\to\infty$. 
\end{theorem}

In order  to capture the full information of the percolation,  we consider  the  loop   ensemble observable  \cite{camia-newman-sle6}, which is defined as follows.  Given a triangulation of a polygon $M$,  let $\omega$ be a site percolation on  $M$ with \emph{monochromatic blue boundary condition}. Namely, the color of each boundary vertex is blue.
Removing all edges on $M$ whose endpoints have opposite colors, we call each connected component in the remaining graph a \emph{percolation cluster}, or simply a cluster, of $\omega$. 
By definition, vertices in each cluster have the same color. 
Moreover, each pair of neighboring vertices that are on different clusters must have opposite colors. 
We call the cluster containing $\bdy M$ the \emph{boundary cluster}. If $\cC$ is a non-boundary cluster of $\omega$, one can canonically define a loop on $M$ surrounding $\cC$ as a path of vertices in the dual map of $M$. We orient the path such that the vertices to the left (resp., right) of the path are red (resp., blue). 
The collections of such loops is called the \emph{loop ensemble} of $\omega$, and we denote it by $\Gamma(M,\omega)$. See  Figure \ref{fig:delta} for an illustration. 
Note that $\omega$ is uniquely determined by $\Gamma(M,\omega)$.

Given a Jordan domain $D$,  a loop ensemble in $D$ is a collection of oriented loops, each viewed as a curve in $D\cup\p D$ modulo monotone reparametrization and rerooting.
Let $\cL(D)$ denote the space of loop ensembles in $D$.
Recall the lattice approximation $D^\mesh$ to $D$ in Theorem~\ref{thm:Smirnov}. Let $\omega^\mesh$ be sampled from $\Ber_{D^\mesh}$ with monochromatic blue boundary  condition.
It was proved in \cite{camia-newman-sle6} that $\Gamma(D^\mesh,\omega^\mesh)$  converge in law as $\mesh\to 0$ to a random variable $\Gamma$ taking values in $\cL(D)$
which is called a \emph{conformal loop ensemble} with parameter $\kappa=6$ ($\CLE_6$) on  $D$.\footnote{In Section~\ref{subsec:SLE}, $\Gamma$ is called a $\CLE_6$ with monochromatic blue boundary condition.} See Theorem~\ref{emb-thm:CLE}  for a precise statement of this  result  including the topology of convergence.

Given $\Map^n$ and  $\{\omega^n_i\}_{i\in\N}$ as above, let $\dlp^n_i\defeq\Gamma(\Map^n,\omega^n_i)$ be the loop ensemble associated with $\omega^n_i$ as defined  in Section~\ref{emb-subsec:Cdy-conv}. Then $(\cM^n,\dlp^n)$ can be viewed  as a compact metric measure space decorated by a (boundary) curve and a loop ensemble.
The natural topology for such objects is the so-called  \emph{Gromov-Hausdorff-Prokhorov-uniform-loop} (GHPUL) topology, which was first introduced in~\cite{ghs-metric-peano}. This is the natural variant of the GHPU topology for cases where the metric space is further decorated by a loop ensemble; see Section~\ref{subsec:pre}.

In the continuum, there exists a variant of the GFF on the unit disk $\D$, denoted by $\bh$, such that $(\D\cup\p\D,c_{\op d}d_\bh, c_{\op m}\mu_\bh,\vec{\xi}_\bh )$ has the law of 
$\BD_1$ as a metric measure space decorated by a curve \cite{lqg-tbm1,lqg-tbm2,lqg-tbm3}, where the constants $c_{\op d},c_{\op m}$ are as in the definition of $(d_\Delta,\mu_\Delta)$ in Theorem~\ref{emb-thm:Cardy}.
The curve $\vec{\xi}_\bh$ is defined by tracing $\p \D$ clockwise,  starting and ending at 1, with the speed prescribed by the boundary measure $\xi_\bh$.
Since $(\Delta,h_\Delta)$ in Theorem~\ref{emb-thm:Cardy} and $(\D,\bh)$  both correspond to $\BD_1$, the two fields are related (in law) by a conformal map between $\D$ and $\Delta$ and the change of coordinates formula for $\sqrt{8/3}$-LQG (see \eqref{emb-eqn-lqg-coord} below).
Let $\{\Gamma_i\}_{i\in \N}$ be a sequence of independent samples of $\CLE_6$ on $\D$ which are also independent of $\bh$. 
Then $(\D\cup\p\D,c_{\op d}d_\bh, c_{\op m}\mu_\bh,\vec{\xi}_\bh,\Gamma_i)$ can be viewed as a compact metric measure space decorated by a curve and a loop ensemble; see Section~\ref{subsec:SLE}.  
For simplicity, we write $(\D\cup\p\D, c_{\op d}d_\bh, c_{\op m}\mu_\bh,\vec{\xi}_\bh,\Gamma_i)$ as $(\D,\bh,\Gamma_i)$.
The following theorem is a precise formulation of the aforementioned convergence of $\{(\Map^n,\omega^n_i)\}_{i\in\N}$. 
\begin{theorem}\label{emb-thm:main}
	In the setting of the paragraph  above, for each $k\in\N$,  $\{(\cM^n,\dlp^n_i)\}_{1\le i\le k}$ jointly converge in law to $\{(\D,\bh, \Gamma_i)\}_{1\le i\le k}$ in the GHPUL topology.
\end{theorem}
\xncomment{We point out that Theorem~\ref{emb-thm:main} for $k>1$ does not easily  follow the $k=1$ case.
	The reason is that the macroscopic behavior of the percolation could depend on microscopic details of the map which disappear in the
	limit. This way, two independent copies of percolation on the same map could have some correlation in the scaling limit. 
	In general,  conditional distributions do not behave well  under convergence in law.}

Theorems~\ref{emb-thm:Cardy} and~\ref{thm:crossings} are easy consequences of  Theorem~\ref{emb-thm:main}. We briefly explain the idea here and refer to Section~\ref{emb-subsec:Cardy-conv} for details.

For Theorem~\ref{thm:crossings}, recall $v^n$ defined there.  For $i\in \N$, let $E^{i}_{a^n}(v^n)$ be defined as $E_{a^n}(v^n)$, with $\omega^n_i$ being the site percolation on $\Map^n$.
Our proof of Theorem~\ref{emb-thm:main} implies that $\{\1_{E^i_{a^n}(v^n)}\}_{1\le i\le k}$ also converge jointly to their continuum counterparts.
By the law of large numbers,  $\Ber_{\Map^n}[E_{a^n}(v^n)] - k^{-1}\sum_{1}^{k}\1_{E^i_{a^n}(v^n)}$ converge to $0$ in probability as $k\to\infty$ \xncomment{uniformly in $n$}.
This proves Theorem~\ref{thm:crossings}. 

Now suppose we are in the setting of Theorem~\ref{emb-thm:Cardy}. 
By the same reasoning as in the previous paragraph, if $v^n$ is sampled uniformly from $\cV(\Map^n)$, 
then $\Ber_{\Map^n}(E_{a^n}(v^n))$, $\Ber_{\Map^n}(E_{b^n}(v^n))$, and $\Ber_{\Map^n}(E_{c^n}(v^n))$ jointly converge to their continuum counterparts. 
This essentially gives the convergence of $\mu^n_\Delta$ to $\mu_\Delta$. A similar argument gives the convergence of $\xi^n_\Delta$.  For the metric $d^n_\Delta$, let $(v^n,u^n)$ be a pair of  vertices uniformly sampled from $\cV(\Map^n)\times \cV(\Map^n)$. Then by the GHPU convergence of $\cM^n$,  $d^n (v^n,u^n)$ 
converge to its continuum counterpart.  Now the uniform convergence of $d^n_\Delta$ follows from the  continuity of $d_\Delta$.
This gives Theorem~\ref{emb-thm:Cardy}.

\subsubsection{On the universality}
We now comment on the universality of our results within the realm of uniform maps and percolation observables. 
See Section \ref{subsec:outlook} for discussion of (nonuniform) planar map models decorated by other statistical physics models.

\begin{remark}[Other variants of uniform triangulations]\label{rmk:type}
	Recall that a triangulation is of type I (resp., type II; type III) if multi-edges and self-loops are allowed (resp., multi-edges are allowed but not self-loops; neither multi-edges nor self-loops are allowed). In \cite{aasw-type2} we consider natural couplings between Boltzmann triangulations of types I, II, and III, and prove that triangulations of polygons of all three types converge in the scaling limit to the Brownian disk. \xncomment{ In these couplings triangulations of different types are related by collapsing self-loops
		or multiple edges. On the other hand, whether a crossing event occurs for a site percolation does not change under these operations.
		Therefore} Theorems ~\ref{emb-thm:Cardy},~\ref{thm:crossings}, and~\ref{emb-thm:main} still hold
	for Boltzmann triangulations of types I and III.
	We also expect these results to hold for uniformly sampled planar maps with other local constraints (quadrangulations, general maps, etc). Establishing these results require nontrivial work.
	The main ingredient which is missing is convergence of the pivotal measure on the planar map. In the case of type II triangulations we obtain this via the bijection in \cite{bhs-site-perc}.
\end{remark}

\begin{remark}[Surfaces with other topologies]
	Our proof techniques can also give variants of Theorem~\ref{emb-thm:main} on uniform triangulations with other topologies. More precisely, given some surface topology (sphere, torus, etc.), if one knows that a uniformly sampled triangulation with this topology converges to a Brownian surface, then one can establish a variant of Theorem~\ref{emb-thm:main}. 
	Furthermore, we get quenched scaling limit results for macroscopic observables of Bernoulli-$\frac12$ site percolation, similar to~Theorem~\ref{thm:crossings}. For example, for uniform triangulation on the sphere with four uniformly sampled vertices $a,b,c,d$, in which case the convergence to the Brownian surface has been established, our method gives that
	the probability that $a,b$ and $c,d$ are separated by a red cycle has a scaling limit. 
	For uniform triangulation on the torus, if the convergence to Brownian torus is shown, then the probability that there exists a non-contractible red cluster has a scaling limit.
	\xncomment{Although the Cardy embedding is specific to surfaces conformally equivalent to the disk,
		for other surfaces we can use other percolation observables to define discrete conformal embeddings.}
\end{remark}

\subsection{Outline of the program}\label{subsec:outline}
\begin{figure}
	\centering
	\includegraphics{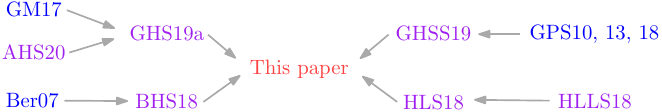}
	\caption{The figure shows other papers in our program (purple), other papers directly relevant for our program (blue), and the dependencies between these papers (gray arrows). In addition to the papers shown in the figure, several stages of our program rely heavily on the fundamental papers \cite{schramm0,smirnov-cardy,camia-newman-sle6,legall-uniqueness,miermont-brownian-map,wedges,lqg-tbm1,lqg-tbm2,lqg-tbm3}.}
	\label{fig:flow-chart}
\end{figure}

Recall that the current work is the final paper in a program also involving \cite{hlls-cut,hlls-pivot,bhs-site-perc,aasw-type2, ghs-metric-peano,ghss18}. See Figure \ref{fig:flow-chart} for an overview of the dependencies between these papers and other papers relevant for the program. The bulk of this paper (Sections~\ref{emb-sec:overview},~\ref{emb-sec:pivot}, and~\ref{sec:GPS}), as well as the bulk of the whole program, is to establish Theorem~\ref{emb-thm:main}.
In this section we give an overview of this program by giving an outline of the proof of Theorem~\ref{emb-thm:main}.

\subsubsection{Annealed scaling limit for one site percolation}\label{subsub:annealed}
The $k=1$ case of Theorem~\ref{emb-thm:main} is proved in our joint work with Gwynne.
\begin{theorem}[\cite{ghs-metric-peano}]\label{emb-thm:upper}
	Theorem~\ref{emb-thm:main} holds when $k=1$.
\end{theorem}
The single interface variant of Theorem~\ref{emb-thm:upper} was proved in \cite{gwynne-miller-perc}, conditioning on Theorem~\ref{thm:AHS}, which was proved in \cite{aasw-type2}. Based on this variant, Theorem~\ref{emb-thm:upper} was proved in  \cite{ghs-metric-peano} via an iterative construction of $\CLE_6$ with chordal $\SLE_6$ (see Lemma~\ref{lem:iteration} for this construction) and its discrete analog.

Theorem~\ref{emb-thm:upper} is an example of an annealed scaling limit result for percolated triangulations, where the convergence is in the sense of GHPUL. 
In another paper of this program \cite{bhs-site-perc}, we discovered, together with Bernardi, a bijection between  lattice walks with steps in $\{(0,1), (1,0), (-1,-1)\}$
and percolated type II triangulations. This bijection builds on an earlier bijection of Bernardi \cite{bernardi-dfs-bijection} between lattice walks in the first quadrant and trivalent maps decorated by a depth-first-search tree. Many percolation observables are encoded nicely  in this bijection. The two most relevant examples are the crossing events in Definition~\ref{emb-def:cardy}, along with the counting measure on self-intersection and mutual-intersection points of macroscopic loops in the loop ensemble. These points are called \emph{pivotal points}. See Section~\ref{subsub:DP}.

The bijection in \cite{bhs-site-perc} is an example of a \emph{mating-of-trees} bijection. 
Its continuum counterpart is an encoding of a $\CLE_6$ and an independent $\sqrt{8/3}$-LQG surface by a 2D Brownian motion. This encoding was introduced in a foundational paper by Duplantier, Miller, and Sheffield \cite{wedges}. See also \cite[Section~7]{bhs-site-perc} and~\cite{ghs-mot-survey}.
Using this bijection and the continuum theory in \cite{wedges}, the scaling limit of many percolation observables were established in \cite{bhs-site-perc}, including those concerning crossing events and pivotal points. This type of scaling limit result is sometimes referred to as convergence in the mating-of-trees sense. In \cite{ghs-metric-peano}, it was proved that the GHPUL convergence in Theorem~\ref{emb-thm:upper} holds jointly with the mating-of-trees  convergence in \cite{bhs-site-perc}. See Proposition~\ref{prop:flip} and \eqref{emb-eq:conv1} for consequences of such  joint convergence. 

The two works \cite{bhs-site-perc} and~\cite{ghs-metric-peano} give a rather complete annealed scaling limit result for percolation on triangulations.
This was achieved by employing the full strength of the continuum theory of $\SLE_6$ and $\CLE_6$  coupled with  $\sqrt{8/3}$-LQG
(including \cite{wedges,gwynne-miller-sle6} and \cite[Section~7]{bhs-site-perc}),
as well as three powerful tools in the discrete: a labeled tree encoding of the graph metric in the spirit of Schaeffer \cite{schaeffer-bijection}  (see \cite{aasw-type2}), a Markovian exploration of uniform triangulations called the peeling process~\xncomment{\cite{angel-uipt} (also see \cite{gwynne-miller-perc})}, and the mating-of-trees  bijection in \cite{bhs-site-perc}. 

When  attacking Theorem~\ref{emb-thm:main} for $k\ge 2$, the toolbox becomes quite limited.
The main methodological innovation of this paper is to supply an approach for doing so, which we explain in Sections~\ref{subsub:DP} and~\ref{subsub:LDP}.

\subsubsection{Dynamical percolation on uniform triangulations}\label{subsub:DP}
It will be apparent from  Section~\ref{emb-sec:overview} that all the difficulties with proving Theorem~\ref{emb-thm:main} for general $k\in\N$ are present already in the $k=2$ case. Therefore we focus on this case.

Our high level idea is the following. 
Let $(\D,\bh, \Gamma_i)_{i=1,2}$ be a subsequential limit of $(\cM^n,\dlp^n_i)_{i=1,2}$, whose existence is guaranteed by Theorem~\ref{emb-thm:upper}.
It suffices to show that $\Gamma_1$ and $\Gamma_2$ are independent.
Suppose we have a dynamic $(\ol\omega^n_t)_{t\ge 0}$ which is stationary conditioned on $\Map^n$ and has one-time conditional marginal law $\Ber_{\Map^n}$.
Moreover, suppose the process  $(\cM^n,\Gamma(\Map^n,\ol\omega^n_t))_{t\ge 0}$ has a GHPUL scaling limit whose one-time marginal 
law is given by $(\D,\bh,\Gamma_1)$.
We denote this process by $(\D, \bh, \ol\Gamma_t)_{t\ge 0}$.
For $t>0$, since $\omega^n_1$ and $\omega^n_2$ are completely independent while $\ol\omega^n_0$ and $\ol\omega^n_t$ may not be,
the correlation between $\Gamma_1$ and $\Gamma_2$ should be no stronger than that of  $\ol \Gamma_0$ and $\ol \Gamma_t$.
If we further know that $(\ol \Gamma_t)_{t\ge 0}$ is ergodic, then by sending $t\to \infty$ we must have that $\Gamma_1$ and $\Gamma_2$ are independent. 
See Section~\ref{emb-sec:overview} for a precise version of this reasoning.

It remains to establish the existence of a dynamic as described in the previous paragraph. The most natural candidate is the following. 
Let $\Map^n$ be as in Theorem~\ref{emb-thm:main} and let $\ol\omega^n$ be sampled from $\Ber_{\Map^n}$. 
Given $(\dsk^n,\ol\omega^n)$, put i.i.d.\ exponential  clocks of rate $n^{-1/4}$ at each interior vertex.\footnote{An exponential clock of rate $r>0$ is a clock which rings at a discrete set of times such that the time between two consecutive rings is given by independent exponential random variables with parameter $r$. In other words, the set of times at which the process rings has the law of a Poisson process on $(0,\infty)$ of intensity $r$.}
When the clock at $v$ rings, flip the color of $v$. For $t\ge 0$, let $\ol\omega^n_t$ be the site percolation at time $t$.
We call $(\ol\omega^n_t)_{t\ge 0}$ a \notion{dynamical percolation} on $\dsk^n$.

We set the clock rate to be $n^{-1/4}$ because \xncomment{the number of vertices on $\Map^n$ governing macroscopic changes is of order $n^{1/4}$ (see the discussion on pivotal points below), so that under this rate the number of updates for these vertices is of constant order.}
We expect that the scaling limit of $(\Map^n,\ol\omega^n_t)_{t\ge 0}$ satisfies the desired ergodic property described in the second paragraph.
If $\Map^n$ is replaced by $\mesh\Tg$ for $\mesh>0$, then the same dynamic was studied by Garban, Pete, and Schramm \cite{gps-pivotal,gps-dynamic}, who established the existence of a scaling limit. However, 
their proof is hard to adapt to the random triangulation case since it relies on the fact that $\Tg$ is nicely embedded into $\C$ (see \cite[Section~8]{gps-dynamic} in particular).  
We expect that proving the aforementioned convergence of $(\Map^n,\ol\omega^n_t)_{t\ge 0}$ is a technically challenging problem.

To get around
this difficulty, we introduce a \emph{cutoff} variant of $(\ol\omega^n_t)_{t\ge 0}$. In this variant of the process, we only update vertices that cause macroscopic changes. 

Let us first quantify the notion of macroscopic change. 
Let $\omega^n$ be a site percolation on $\Map^n$ with monochromatic blue boundary condition.
Given  a non-boundary cluster $\cC$ of $\omega^n$, 
let $\neg\cC$ be  the connected component of $\cV(\Map^n)\setminus \cV(\cC)$ containing $\bdy\Map^n$.
Let $\ol\cC$ be the largest subgraph of $\Map^n$ such that $v\in \cV(\cC)$ if and only if  $v\notin \neg \cC$. 
For each loop $\ell\in\Gamma(\Map^n,\omega^n)$, 
let $\reg(\ell)=\ol\cC$ where $\cC$ is the cluster of $\omega^n$ surround by $\ell$.  
We call $\area(\ell):=\mu^n(\reg(\ell))$ the \emph{area} of $\ell$.
For $v\in \cV(\Map^n)\setminus\cV(\bdy \Map^n)$,  let $\omega^n_v$ be obtained from $\omega^n$  by flipping the color of $v$, and let $\cL^n_v$ be the symmetric difference between $\Gamma(\Map^n,\omega^n)$ and $\Gamma(\Map^n,\omega^n_v)$.  
For $\ep>0$, we say that $v$ is an \notion{$\ep$-pivotal point} of $\omega^n$
if there are at least three  loops in $\cL^n_v$ with  area  at least $\ep$.
Morally speaking, $v$ is an $\ep$-pivotal point if flipping the color of $v$ results in a macroscopic  change of ``size'' at least $\ep$.

We now consider the following modification of $(\Map^n,\ol\omega^n_t)_{t\ge 0}$: when the clock of a vertex $v$ rings at time $t$, the color of $v$ is flipped if and only if $v$ is an $\eps$-pivotal point of $\ol\omega^n_t$.  We denote this modified dynamic by  $(\Map^n,\ol\omega^{\eps,n}_t)_{t\ge 0}$. 

Let $\bh$ be as in Theorem~\ref{emb-thm:main} and let $\Gamma$ be a $\CLE_6$ on $\D$ independent of $\bh$.
We can mimic the definition in the discrete to define the $\eps$-pivotal points of $(\bh,\Gamma)$ (see Definition~\ref{def:eps}).
Let $\cP_\eps$ be the set of $\eps$-pivotal points of $(\bh,\Gamma)$. Then $\cup_{\eps>0}\cP_\eps$ is simply the collection of all self-intersections and mutual intersections of loops in $\Gamma$.
We call points in $\cup_{\eps>0}\cP_\eps$  the \emph{pivotal point}s of $\Gamma$. The analogue of color flipping in the continuum is merging and splitting of loops of $\Gamma$; see Section~\ref{subsec:SLE}.

In \cite{bhs-site-perc}, a measure $\nu^\eps_{\bh,\Gamma}$ supported on the $\eps$-pivotal points of $(\bh,\Gamma)$, called  the \notion{$\sqrt{8/3}$-$\LQG$ $\eps$-pivotal measure},  was defined 
based on the theory of mating of trees \cite{wedges}. (See Definition~\ref{def:eps-piv} for a precise definition.) 
Let $\nu^{\eps,n}_{\pv}$ be $n^{-1/4}$ times the counting measure on the $\eps$-pivotal points of $\ol\omega^n_0$. 
As alluded to in Section~\ref{subsub:annealed}, it was proved in \cite{bhs-site-perc,ghs-metric-peano} that for some constant $c_{\op p}>0$,  
\begin{equation}\label{eq:piv-conv}
	\textrm{$(\cM^n, \nu^{\eps,n}_{\pv},  \Gamma(\Map^n,\ol\omega^n_0))$ converge in law to $(\D, \bh, c_{\op p}\nu^\eps_{\bh,\Gamma}, \Gamma))$.}
\end{equation}
Here the convergence is for a variant of the GHPUL topology that takes into account the additional measure $\nu^{\eps,n}_{\pv}$.

The Markovian dynamic $(\ol\omega^{\eps,n}_t)_{t\ge 0}$ can be described as follows. Starting from the configuration at time $t=0$,  
we wait for an exponential clock of rate $\nu^{\eps,n}_{\pv}(\cV(\Map^n))$ to ring.
Once the clock rings, a vertex $v$ is chosen according to $\nu^{\eps,n}_{\pv}$ and the color of $v$ is flipped. 
Then we iterate this procedure. 
In light of this description and \eqref{eq:piv-conv}, we can show that $(\cM^n,\Gamma(\Map^n,\ol\omega^{\eps,n}_t))_{t\ge 0}$ has a GHPUL scaling limit whose one-time marginal law is given by $(\D,\bh,\Gamma_1)$. We denote this process by $(\D, \bh, \ol\Gamma^\eps_t)_{t\ge 0}$. For each $\eps>0$, the process $(\ol\Gamma^\eps_t)_{t\ge 0}$ is \emph{not} ergodic.
However, we will prove in Section~\ref{subsub:LDP} that 
\begin{equation}\label{eq:eps-lim}
	\textrm{$(\ol\Gamma^\eps_t)_{t\ge 0}$ converge to an ergodic process  as $\eps\to0$.}
\end{equation}
Recall the setting of the second paragraph. 
The correlation between $\Gamma_1$ and $\Gamma_2$ should be no stronger than that of $\ol \Gamma^\eps_0$ and $\ol \Gamma^\eps_t$ for each $\eps>0$ and $t>0$.
In light of~\eqref{eq:eps-lim}, by first sending $\eps\to 0$ and then $t\to \infty$, we can still establish the $k=2$ case of Theorem~\ref{emb-thm:main}. 
Again see Section~\ref{emb-sec:overview} for how to make this reasoning rigorous.

\subsubsection{Quantum pivotal measure and Liouville dynamical percolation}\label{subsub:LDP}
The proof of \eqref{eq:eps-lim} is done in Sections~\ref{emb-sec:pivot} and~\ref{sec:GPS}, based on \cite{hlls-cut,hlls-pivot,ghss18}.

The first key step is to achieve  a good understanding of the measure $\nu^\eps_{\bh,\Gamma}$.
Recall $(\D,\bh,\Gamma)$ in \eqref{eq:piv-conv} and  the set  $\cP_\eps$ of $\eps$-pivotal points of $(\bh,\Gamma)$ in Section~\ref{subsub:DP}.
By \eqref{eq:piv-conv}, $\nu^\eps_{\bh,\Gamma}$ is the scaling limit of $\mu^n$ restricted to  the discrete analog of $\cP_\eps$ under a proper renormalization.  

Fix $\delta>0$, and suppose $\D^\mesh$ is the lattice approximation of $\D$ via $\mesh\Tg$.  Let $\omega^\mesh$ be sampled from $\Ber_{\D^\mesh}$.
In~\cite{gps-pivotal}, it was proved that the counting measure on the pivotal points of $\omega^\mesh$ under proper rescaling converge to a random measure $\Mink$; see the discussion below Definition~\ref{def:pivm} for a precise description of $\Mink$.  
The convergence is joint with  the loop ensembles. Now suppose  $\{\omega^\mesh\}_{\mesh>0}$ are coupled such that the loop ensemble convergence holds almost surely. 
Suppose $\bh$ is independent of $\{\omega^\mesh\}_{\mesh>0}$.
For each loop $\ell$ of $\omega^\mesh$
let $\mu_\bh(\reg(\ell))$ be the \emph{area} of $\ell$ and define the  $\eps$-pivotal points for $(\bh,\omega^\mesh)$ as in Section~\ref{subsub:DP} with this notion of loop area.
Let $\cP^\mesh_\eps$ be the union of all hexagons corresponding to $\eps$-pivotal points of $(\bh,\omega^\mesh)$. 
It is not hard to show that 
under proper rescaling, as $\mesh\to 0$,  
the  \xncomment{\emph{Gaussian multiplicative chaos} (GMC)} measure $e^{\bh/\sqrt6}d^2z$ restricted to $\cP^\mesh_\eps$ converge in probability to a random measure $\pivm^\eps$; see Section~\ref{sec:prop:quantum-Mink}. 
Moreover, $\pivm^\eps=(e^{\bh/\sqrt 6}\Mink)|_{\cP_\eps}$ a.s., where the right side is understood as the restriction of a GMC; see~\cite{rhodes-vargas-review,berestycki-gmt-elementary} and Definition~\ref{def:GMC}.
It is well-known that $\cP_\eps$ is a fractal of dimension $3/4$ (see e.g. \cite{smirnov-werner-percolation}). 
The  so-called \emph{Knizhnik-Polyakov-Zamolodchikov} (KPZ) relation (see e.g. \cite{shef-kpz} and Remark~\ref{rmk:KPZ}) suggests that 
\begin{equation}\label{eq:KPZ}
	\nu^\eps_{\bh,\Gamma}=\constp\pivm^\eps  \textrm{ a.s.\ for a deterministic constant }\constp.
\end{equation}
\xncomment{Here the exponent $1/\sqrt6$ in $\pivm^\eps $ is precisely related to the dimension $3/4$ of $\cP_\eps$  via the KPZ relation.}
We restate \eqref{eq:KPZ} as Proposition~\ref{prop:quantum-Mink} and prove it in Section~\ref{sec:prop:quantum-Mink}. Most of the work is carried out in Section~\ref{emb-sec:pivot}, where we prove Proposition~\ref{prop:pivot-bichordal}, a local version of Proposition~\ref{prop:quantum-Mink}. We say that it is local because we will cover $\cP_\eps$ by finitely many sets which are the scaling limits of the pivotal points of the crossing event for certain topological quads (see Lemma~\ref{lem:cover}), and Proposition~\ref{prop:pivot-bichordal} is the variant of~\eqref{eq:KPZ} for these sets.

Although  the argument is  quite technical,
the underlying idea behind Propositions~\ref{prop:quantum-Mink} and~\ref{prop:pivot-bichordal} is simply that both $\nu^\eps_{\bh,\Gamma}$ and $\pivm^\eps$ 
are canonical in the sense that they satisfy a few natural properties that uniquely determine the measure. To carry out this idea, we need an intrinsic characterization of  the aforementioned measure $\Mink$ 
that does not refer to the limiting procedure. With this in mind, we proved with Lawler and Li~\cite{hlls-cut} that 
$r^{d-2}$ times the Lebesgue measure restricted to the $r$-neighborhood of  cut points of a planar Brownian motion has a scaling limit as $r\to0$, 
which we call the \emph{$3/4$-dimensional occupation measure}.
Using a connection between Brownian cut points and the scaling limit of pivotal points of quad-crossing events (see Proposition~\ref{prop:cutpt}), 
we proved with Li~\cite{hlls-pivot} that restricting to  the scaling limit of the pivotal points of quad crossing events, the measure $\Mink$
equals the corresponding $3/4$-dimensional occupation measure on these points.

With the results from \cite{hlls-cut,hlls-pivot} at hand, we first prove the variant of~\eqref{eq:KPZ} with $\cP_\eps$ replaced by Brownian cut points (i.e.\ Lemma~\ref{lem:pre}). 
This is based on the theory of quantum zippers in \cite{shef-zipper,wedges} and the coordinate change formula for GMC over occupation measures. Then using the connection between Brownian cut points and  the scaling limits of pivotal points of quad crossing events, we conclude the proof of Proposition~\ref{prop:pivot-bichordal}. We finally prove~\eqref{eq:KPZ} (i.e.\ Proposition~\ref{prop:quantum-Mink}) via a covering argument.

Given \eqref{eq:KPZ}, we will approximate the process $(\ol\Gamma^\eps_t)_{t\ge 0}$ in \eqref{eq:eps-lim} by a variant of dynamical percolation on the triangular lattice $\Tg$. This enables us to  use powerful tools that are only available for site percolation on $\Tg$, including various scaling limit results and the sharp noise sensitivity established in \cite{gps-fourier}.

Fix $\delta>0$, and suppose
that $\omega^\mesh$ is sampled from $\Ber_{\D^\mesh}$ independently of $\bh$.
In light of \eqref{eq:KPZ}, we can consider  a 
variant of the dynamical percolation on $\D^\mesh$, where the rate of the exponential clock at a vertex $v$ is proportional to  (a regularized version of) $e^{\bh(v)/\sqrt6}$.  
This is the so-called discrete \notion{Liouville dynamical percolation} (LDP) driven by $e^{\bh/\sqrt 6}$ introduced by Garban, Sep{\'u}lveda, and us in \cite{ghss18}; see  Section~\ref{subsec:LDP}.
Now we can define an $\eps$-cutoff dynamic of the discrete LDP on the triangular lattice by mimicking the definition of $(\ol \omega^{\eps,n}_t)_{t\ge 0}$ in Section~\ref{subsub:DP}, and then use \eqref{eq:KPZ} to argue that the loop ensemble evolution  of this cutoff dynamic  converge to the process $(\ol\Gamma^\eps_t)_{t\ge 0}$ in~\eqref{eq:eps-lim}. 

Now to conclude the proof of \eqref{eq:eps-lim}, we just need to show that as $\eps\to0$, the $\eps$-cutoff  dynamic of the discrete LDP driven by $e^{\bh/\sqrt 6}$ stabilize to a limiting ergodic process. The paper \cite{ghss18} achieved this goal modulo two differences.
First, following \cite{gps-pivotal,gps-dynamic},  in \cite{ghss18} we work under a different cutoff on the pivotal points which is based on alternating four arm events. (See the notion of $\rho$-important points in Section~\ref{subsub:rho-important}.) 
Compared to the $\eps$-pivotal points, this cutoff is not so natural in the context of random planar maps because it relies on the ambient space. However, it is convenient for
fine multi-scale analysis on $\Tg$, which gives the desired stability when removing the cutoff.  
The limiting process is called the \emph{continuum Liouville dynamical percolation} driven by $e^{\bh/\sqrt{6}}$.
In Section~\ref{sec:GPS} we study the relation between the two cutoffs and show that $\lim_{\eps\to 0}(\ol\Gamma^\eps_t)_{t\ge 0}$ exists and is given by the continuum Liouville dynamical percolation  driven by $e^{\bh/\sqrt{6}}$.

The second difference from \cite{ghss18} is that there
the planar percolation is not encoded by the loop ensemble, but rather by crossing information for all topological rectangles in the plane. 
The latter is called the \emph{quad-crossing configuration}. 
Similarly as above, the quad crossing configuration is not so natural  in the context of random planar maps due to its dependence on the ambient space. On the other hand, the quad crossing perspective is crucial in our proof of the ergodicity of continuous LDP in \cite{ghss18}, which relies on Fourier analysis of Boolean functions following \cite{gps-fourier}.
This difference in observable will not be a problem if we know that 
the $\CLE_6$  and  the scaling limit of the quad-crossing configuration of $\omega^\mesh$  determine each other. 
This has long been conjectured to be true (see \cite{ss-planar-perc}). The fact that the $\CLE_6$ determines quad-crossing configuration is essentially proved in  \cite{camia-newman-sle6}, as pointed out in \cite{gps-pivotal}. We establish measurability in the reverse direction in this paper; see Theorem~\ref{thm:quad}. This concludes our proof.

\subsection{Related works and outlook}\label{subsec:outlook}
Theorem \ref{emb-thm:Cardy} solves a special case of the aforementioned  conjecture that Liouville quantum gravity describes the scaling limit of random planar maps under discrete conformal embeddings. 
The general version of the conjecture can be formulated as follows.

For the ease of discussion, assume that there are $m_1$ different ways to sample a random planar map of a given size. The map can be required to be a triangulation, quadrangulations, simple map, etc., and the probability measure can be uniform (like in our paper) or nonuniform. For example, we can reweight the uniform distribution by the partition function of a statistical physics model such as the uniform spanning tree (UST), the (critical) Ising model, or the  Fortuin-Kasteleyn (FK) random cluster model.
We also assume  that there are $m_2$ different ways to conformally embed a  planar map into $\C$.
Besides the Cardy embedding considered in this paper and the aforementioned circle packing and the Tutte embedding, one can also consider the square tiling and the embedding obtained by applying the uniformization theorem to the planar map viewed as a piecewise smooth 2D Riemannian manifold. 

The general conjecture predicts convergence of random planar maps under conformal embedding to $\gamma$-LQG in each of the $m_1m_2$ situations obtained by specifying the law of the random planar map and the embedding method, where the value of $\gamma$ depends on the law of the planar map. For example, uniformly sampled planar maps give $\gamma=\sqrt{8/3}$. 
Consider a statistical physics model on a planar map whose partition function is approximately $(\det\Delta)^{-c/2}$, where $\det \Delta$ represents the determinant of the Laplacian of the planar map and $c\in \R$ is the so-called \emph{central charge} of the model. Suppose our random planar map is sampled such that the probability of sampling a particular map is proportional to the partition function of the statistical physics model on the planar map. 
Choose $\gamma\in (0,2)$ such that $c=25-6(2/\gamma+\gamma/2)^2$. Then the scaling limit of the random planar map is conjecturally given by $\gamma$-LQG. For example, the UST has central change $c=-2$, and therefore the scaling limit of UST weighted random planar maps is $\sqrt2$-LQG. For the Ising model, we have $c=1/2$ and $\gamma=\sqrt3$. 
Our paper is the first work which solves one version of this conjecture. 

We remark that convergence to LQG under a conformal embedding (namely, the Tutte embedding) has been established earlier for a large class of random planar maps obtained from coarse-graining an LQG surface, e.g.\ the so-called mated-CRT map \cite{gms-tutte} and the Poisson Voronoi tessellation of the Brownian disk \cite{gms-voronoi} \xncomment{and its extension to general $\gamma$-LQG in~\cite{afs-volume}}, except that the convergence established there is only for the vertex counting measures, not for the measures and the graph metric jointly.

The Cardy embedding is a representative for a class of embeddings which are defined using observables  of statistical physics models on planar maps. The 
Tutte embedding is another such example, where the model is simple random walk and the observables are given by the  harmonic measure. One can define natural embeddings of planar maps in other universality classes by using observables of other statistical physics models. For example, in the case of the FK random cluster model one can use properties of the FK loops to define an embedding similarly to the case of percolation. For a UST weighted map with sphere topology one can first send three uniformly sampled vertices $v_1$, $v_2$, and $v_3$ to 0, 1, and $\infty$, respectively, and then determine the position in $\BB C$ of an arbitrary vertex $w$ by considering the topology of the tree branches connecting $w$, $v_1$, $v_2$, and $v_3$. In light of this, the  ``number''  $m_2$ of possible discrete conformal embeddings is quite large.

Using the aforementioned $m_1$ random planar map models and $m_2$ discrete conformal embeddings, we obtain $m_1m_2$ random environments in which we can consider statistical physics models, such as random walk or percolation.  We conjecture the following universality. 
If the  random process  converges to a conformally invariant process on a regular lattice, then the same convergence holds for the random process in one of these
$m_1m_2$ random environments, in a quenched sense. For example, our results in Section~\ref{subsub:quenched} imply this type of convergence where the random  process is site percolation, while the random environment is provided by the uniform triangulation under the Cardy embedding, or any other embedding for which the analog of Theorem~\ref{emb-thm:Cardy} holds.
As another example, we expect that  since random walk on regular lattices converge to planar Brownian motion, the  random walk in one of these  $m_1m_2$  environments converge to planar Brownian motion in a quenched sense.  
Our results in Section~\ref{subsub:quenched} are the only such quenched scaling limit results in the literature for natural model-decorated combinatorial random planar maps. 
The quenched scaling limit of random walk has been established in \cite{gms-random-walk} for a large class of random planar maps obtained by coarse \xncomment{graining} $\LQG$.

It may be possible to use the approach introduced in this paper to prove the conjectures above when the random planar map is weighted by a statistical physics model and  the discrete conformal embedding is defined using observables of the same  model. 
In this case, if one can establish the analogue of Theorem \ref{emb-thm:main}, then one can prove the analogue of Theorem~\ref{emb-thm:Cardy}.
Note that in our case, uniform planar maps can be thought of as percolation weighted planar maps and the Cardy embedding is defined via percolation observables.  
At a conceptual level, our dynamical approach should still work in the more general setting.
However, carrying out this approach  beyond the setting of our current paper is a challenge. 
In particular, we  use the metric convergence of uniform  triangulations to the Brownian disk and a sharp mixing property for the scaling limit of dynamical percolation on the planar map. 
Both of these ingredients are currently missing for other planar maps and statistical physics models, each of which is a major open question \xncomment{in its own right}.

Convergence of model-decorated random planar maps to $\LQG$ has been established for a much more  general class of planar map models in the so-called \emph{peanosphere sense}. 
This convergence is based on the mating-of-trees framework of \cite{wedges}. 
In the discrete,  a number of mating-of-trees type bijections have been discovered, similar in spirit as the one we discovered with Bernardi \cite{bhs-site-perc}.
With such kind of bijections and the mating-of-trees framework for LQG coupled with SLE/CLE, convergence in the peanosphere sense means convergence to Brownian motion of the random walk encoding the decorated map. 
This idea was first proposed and carried out in  \cite{shef-burger}. See~\cite[Section 5.1]{ghs-mot-survey} for a survey with further examples.
Here we point out that this convergence does not concern the metric or conformal structure of the map. Moreover, it is an annealed instead of quenched result if we view it as a convergence result for a random process in a random environment.

Dynamical percolation is an important tool in the current paper, and we prove a weak notion of convergence of dynamical percolation on the random planar map to Liouville dynamical percolation; namely, we prove convergence of the variant of the process where only $\eps$-pivotal points change color, and we prove that the limiting process stabilizes to the continuous LDP as $\eps\to0$. 
An interesting open problem is to prove convergence of true dynamical percolation on the random planar map to the continuous LDP. One can also attempt to establish similar scaling limit results for models closely related to dynamical percolation, such as the minimal spanning tree, invasion percolation, and near-critical percolation. See \cite{gps-dynamic,gps-mst} for scaling limit of results for these models on the triangular lattice.

\subsection*{Structure of the paper}\label{subsub:structure}
In Section \ref{sec:pre} we provide necessary background on $\sqrt{8/3}$-LQG, SLE$_6$, CLE$_6$, and the topological spaces relevant for the convergence results. In Section \ref{emb-sec:overview} we prove Theorem \ref{emb-thm:main}, assuming two lemmas which are proved in Section \ref{sec:GPS}. In Section \ref{emb-subsec:Cardy-conv} we conclude the proof of Theorems \ref{emb-thm:Cardy} and \ref{thm:crossings} using Theorem \ref{emb-thm:main}. In Section \ref{emb-sec:pivot} we establish a preliminary version of \eqref{eq:KPZ} via an extensive analysis of the CLE$_6$ pivotal points. In Section \ref{sec:GPS} we establish the two aforementioned lemmas using Liouville dynamical percolation, in addition to concluding the proof of \eqref{eq:KPZ}. 

\subsection*{Acknowledgements}
We are grateful for enlightening discussions with Ewain Gwynne and Scott Sheffield  at the early stage of this project. 
We also thank Vincent Tassion for pointing us to the paper \cite{agmt16}, and we thank an anonymous referee for careful reading of the paper and many helpful comments.
The research of N.H.\ is supported by the Norwegian Research Council, Dr.\ Max R\"ossler, the Walter Haefner Foundation, and the ETH Z\"urich Foundation.
The research of X.S.\ is supported by the Simons Foundation as a Junior Fellow at the Simons Society of Fellows, by NSF grant DMS-1811092, and by Minerva fund at the Department of Mathematics at Columbia University.

\section{Preliminaries} \label{sec:pre}
\subsection{Basic notations}\label{subsec:notation}
\noindent {\bf Sets.} Let $\N=\{1,2,\dots \}$ be the set of positive integers. Let $\C$ be the complex plane. 
Let $\D=\{z\in \C: |z|<1 \}$, $\bbH=\{z:\op{Re} z>0\}$, and $\cS=\R\times (0,\pi)$. 

\vspace{7pt}
\noindent 
{\bf Domains}. A (planar) domain is a connected open subset of  $\C$. 
Given a domain $D$, let $\bdy D$ denote the set of prime ends of $D$.
If $\bdy D$ is a simple closed curve, we call $D$ a Jordan domain. 
Given a simply connected domain $D$, we say $D$ is $C^0$ if any  conformal  map  $\phi:\D \to D$ 
can be extended continuously  to $\bdy \D$. (Here, if $D$ is unbounded, we use the spherical metric  on $\C\cup\{\infty\}$).
If $D$ is $C^0$ and the continuous extension of $\phi$ is smooth except for finitely many points, we say that $D$ is \notion{piecewise smooth}. 
Given two domains $D_1,D_2\subset\C$ we write $D_1\Subset D_2$ if $D_1\cup \p D_1\subset D_2$. For two distinct points $a,b$ on $\p D$, let $\p_{a,b}D$ be the counterclockwise arc on $\p D$  from $a$ to $b$.

\vspace{7pt}
\noindent{\bf Lattice}. Let $\bT$ denote the regular  triangular lattice where each face is an equilateral triangle and the points $(0,0), (1,0)$ belong to $\bT$.
For $\mesh>0$, let $\mesh \bT$  be $\bT$  rescaled by $\mesh$.  A Jordan domain $D$ is called a \notion{$\mesh$-polygon} if $\p D$ lies on  $\mesh\bT$. If $D$ is a general Jordan domain, let $D^\mesh$ be the largest $\mesh$-polygon whose set of inner vertices (namely, vertices on $\mesh\bT$ that are inside the $\mesh$-polygon) is contained in  $D$ and forms a connected set on $\mesh\bT$.\footnote{In case of a draw, we choose $D^\mesh$ arbitrarily from the set of largest $\mesh$-polygons, but note that $D^\mesh$ will be uniquely determined for all sufficiently small $\mesh$.}    Including all vertices and edges in $D^\mesh\cap \mesh\bT$, we obtain a triangulation of a polygon, which we  call the \notion{$\mesh$-approximation} of $D$ and still denote by $D^\mesh$.

\vspace{7pt}
\noindent{\bf Measures.}  Given measurable spaces $E,F$, a measure $\mu$ on $E$, and a measurable map $\phi:E\to F$, 
the pushforward of $\mu$ under $\phi$ is denoted by $\phi_*\mu$. Let $f$ be a measurable nonnegative function on $E$. We let $f\mu$ denote the measure 
whose Radon-Nikodym with respect to $\mu$ is $f$.

\vspace{7pt}
\noindent{\bf Random variables.}  Given two random variables $X$ and $Y$, we write $X\overset{d}{=}Y$ if $X$ and $Y$ have the same law.  If $Z$ and $W$ are two random variables on the same probability space, we say that \emph{$Z$ (almost surely) determines $W$}  if and only if there exists a random variable $W'$ measurable with respect to the $\sigma$-algebra generated by $Z$ such that $W=W'$ almost surely.

\subsection{Topological preliminaries}
\label{subsec:pre}
In this section we define the topologies used in   Theorems~\ref{thm:AHS} and~\ref{emb-thm:main}, following \cite{ghs-metric-peano}. We start by defining the GHPU topology in Theorem~\ref{thm:AHS}. 
Given a metric space $(X,d)$, 
for two closed sets $E_1,E_2 \subset X$, their \emph{Hausdorff distance}   is given by 
\[
\BB d_d^{\op{H}}(E_1,E_2):=\max\{ \sup_{x\in E_1}\inf_{y\in E_2} d(x,y), \sup_{y\in E_2}\inf_{x\in E_1} d(x,y) \}.
\]
For a closed set $A\subset X$ and $\eps>0$ define $A_\eps=\{z\in X: d(x,z)\le \eps\textrm{ for some } x\in A\}$ to be the $\eps$-neighborhood of $A$. Then, for two finite Borel measures $\mu_1,\mu_2$ on $X$, their \emph{Prokhorov distance}   is given by   
	\[
	\BB d^{\op{P}}_d (\mu_1,\mu_2) =\inf \{ \ep>0: \mu_1(A)\le \mu_2(A_\eps)+\ep \;\textrm{and}\; \mu_2(A)\le\mu_1(A_\eps)+\ep \; \textrm{for all closed sets }A\subset X \}.
	\]

Let $C_0(\BB R , X)$ be the space of continuous curves $\xi : \BB R\rta X$ which extend continuously to the extended real line $[-\infty,\infty]$, i.e., the limits $\lim_{t\rta+\infty}\xi(t)$ and $\lim_{t\rta-\infty}\xi(t)$ exist.
The \emph{uniform distance} between $\xi_1,\xi_2\in C_0(\R,X)$ 
is given by 
$$\BB d^{\op U}_d(\xi_1,\xi_2)\defeq\sup_{t\in \R} d(\xi_1(t) ,\xi_2(t)).$$
For a finite interval $[a,b]$, we can view a curve $\xi : [a,b] \rta X$ as an element of $C_0(\BB R ,X)$ by defining $\xi(t) = \xi(a)$ for $t < a$ and $\xi(t) = \xi(b)$ for $t> b$.  

Let $\BM^\GHPU$ be the set of quadruples $\frk X  = (X , d , \mu , \xi)$ where $(X,d)$ is a compact metric space, 
$\mu$ is a finite Borel measure on $X$, and $\xi \in C_0(\BB R,X)$. 
If we are given elements $\frk X^1 = (X^1 , d^1, \mu^1 , \xi^1) $ and $\frk X^2 =  (X^2, d^2,\mu^2,\xi^2) $ of $ \BM^\GHPU$ and isometric embeddings $\iota^1 : (X^1 , d^1) \rta (W,D)$ and $\iota^2 : (X^2 , D^2) \rta  (W,D)$ for some metric space $(W,D)$, we define the \emph{GHPU distortion} of $(\iota^1,\iota^2)$ by
\begin{align}
	\label{eqn-ghpu-var}
	\op{Dis}_{\frk X^1,\frk X^2}^\GHPU\left(W,D , \iota^1, \iota^2 \right)   
	:=  \BB d^{\op{H}}_D \left(\iota^1(X^1) , \iota^2(X^2) \right) +   
	\BB d^{\op{P}}_D \left(( (\iota^1)_*\mu^1 ,(\iota^2)_*\mu^2) \right) + 
	\BB d_D^{\op{U}}\left( \iota^1 \circ \xi^1 , \iota^2 \circ\xi^2 \right).
\end{align}
The \emph{Gromov-Hausdorff-Prokhorov-Uniform distance} between $\frk X^1$ and $\frk X^2$ is given by
\begin{align} \label{eqn-ghpu-def}
	\BB d^\GHPU\left( \frk X^1 , \frk X^2 \right) 
	= \inf_{(W, D) , \iota^1,\iota^2}  \op{Dis}_{\frk X^1,\frk X^2}^\GHPU\left(W,D , \iota^1, \iota^2 \right)      ,
\end{align}
where the infimum is over all compact metric spaces $(W,D)$ and isometric embeddings $\iota^1 : X^1 \rta W$ and $\iota^2 : X^2\rta W$.
By~\cite[Section 2]{gwynne-miller-uihpq}, $\BB d^\GHPU$ is a complete separable metric on~$\BM^\GHPU$ provided we identify any two elements of~$\BM^\GHPU$ which differ by a measure- and curve-preserving isometry.

Given a graph $G$,  identify each edge of $G$ with a copy of the unit interval $[0,1]$. We define a metric $d^{\op{gr}}_G$ on $G$ by requiring that 
this identification is an isometric  embedding  of $[0,1]$ into $(G,d_G,\mu_G)$. Let $\mu_G$ denote the counting measure on the vertex set of $G$.
For a discrete interval $[a,b]_\Z\defeq [a,b]\cap \Z$, a function $ \rho : [a,b]_\Z \rta \mcl E(G)$ is called an \emph{edge path} if $\rho(i)$ and $\rho(i+1)$ share an endpoint for each $i\in [a,b-1]_\Z$. We can extend an edge path $\rho$ from $[a,b]_{\BB Z}$ to $[a-1,b] $ in such a way that $\rho$ is continuous and $\rho([i-1,i])$ lies on the edge $\rho(i)$. Note that there are multiple ways to extend $\rho$, but any two different extensions result in curves with uniform distance  at most~$1$.

Recall  the Boltzmann triangulation  $\Map^n$ in Theorem~\ref{thm:AHS}, whose boundary length $\ell^n$ satisfies $\left(3n\right)^{-1/2}\ell^n \rta1$.
Then $\p \Map^n$ can be viewed as an edge path $\beta^n$ tracing the boundary clockwise\footnote{In contrast to some other papers \cite{aasw-type2,ghs-metric-peano},  we orient $\p \Map^n$ clockwise because  in Theorem~\ref{emb-thm:upper}, the percolation has monochromatic blue boundary condition. We want to  be consistent with the orientation induced by the percolation where blue color is on the right-hand side. Also see Section~\ref{subsec:SLE}, where we require the domain to have clockwise oriented boundary when the $\CLE_6$ has monochromatic blue boundary condition. Note that the law of $(\Map^n,d^n,\mu^n,\vec{\xi}^n)$ in $\BM^\GHPU$ is unchanged if we swap the orientation of $\p \Map^n$.}
starting and ending at the root edge.  Set
\begin{equation}\label{eq:scale}
	d^n\defeq(3n/4)^{-1/4}d^{\op{gr}}_{\Map^n},  \quad\quad 	\mu^n\defeq (2n)^{-1}\mu_{\Map^n},\quad\textrm{ and }\quad \vec{\xi}^n(t)\defeq\beta^n(t \ell^n)\; \textrm{for }t\in[0,1].
\end{equation}
Then $\cM^n\defeq \left( \Map^n , d^n , \mu^n , \vec{\xi}^n \right)$ is a random variable in $\BM^\GHPU$.  
Now the precise meaning of Theorem~\ref{thm:AHS} becomes clear. 
It states that $\cM^n$ converge in law to a random variable $\BD_1$ in the GHPU topology.
A random variable with the law of $\BD_1$ is called  a \notion{free Brownian disk with unit perimeter}. 
We refer to \cite{bet-mier-disk} for an explicit construction of $\BD_1$ using the Brownian snake. 
For the purpose of this paper, we can take Theorem~\ref{thm:AHS} as our definition of $\BD_1$. 
Alternatively, Theorem~\ref{thm:LQG-TBM} below specifies $\BD_1$ as well.

Now we define the GHPUL topology used in Theorem~\ref{emb-thm:main}.
Given a metric space $(X, d)$,  an \emph{unrooted oriented loop} on $X$  is a continuous map from the circle to $X$ identified up to reparametrization by orientation-preserving homeomorphisms of the circle.
Define the pseudo-distance  between two continuous maps from the circle $\BB R/\Z$ to $X$  by
\[
\BB d^{\op u}_{d} (\ell,\ell')= \inf_{\psi} \sup_{t\in \R /\Z} d(\ell(t), \ell'(\psi(t)),
\]
where the infimum is taken over all orientation-preserving homeomorphisms $\psi : \BB R/\BB Z\rta\BB R/\BB Z$. 
%\xncomment{Suppose $(X,d)$ is compact so that it is complete and separable.} Then $\BB d^{\op u}_{d}$ induces a complete metric, which we still denote by $\BB d^{\op u}_d$, on unrooted oriented loops. Moreover, the space of unrooted loops is separable with respect to $\BB d^{\op u}_{d}$.

A closed set of unrooted oriented loops on $X$ with respect to the $\BB d^{\op u}_d$-metric is called a \emph{loop ensemble} on $X$.  
We let  $\cL(X)$ be the space of loop ensembles on $X$ equipped with the Hausdorff metric 
\begin{equation}\label{emb-eq:dist-ensemble}
	\BB d_d^{\op L} (c,c')=\max\{ \BB d_d^{\op L,0} (c,c'),\BB d_d^{\op L,0} (c',c) \},
\end{equation}
where
\begin{equation}
	\BB d_d^{\op L,0}(c,c')=\inf\{\eps>0: \forall \ell \in c, \exists \ell' \in c'\,\textrm{such that} \;\BB d_d^{\op u} (\ell,\ell')\le \eps\}.
\end{equation}

Let $\BM^\GHPUL$ be the set of $5$-tuples $\frk X  = (X , d,  \mu , \eta, c)$ where $(X,d)$ is a compact metric space, $\mu$ is a finite Borel measure on $X$, $\eta \in C_0(\BB R,X)$, and $c\in \cL(X)$.   If we are given elements $\frk X^1  = (X^1 , d^1,  \mu^1 , \eta^1,c^1)$ and $\frk X^2  = (X^2 , d^2,  \mu^2 , \eta^2,c^2)$  in $ \BM^\GHPUL$ and isometric embeddings 
$\iota^1 : (X^1 , d^1) \rta (W,D)$ and $\iota^2 : (X^2 , d^2) \rta  (W,D)$ for some metric space $(W,D)$, we define the \emph{GHPU-Loop (GHPUL) distortion} of $(\iota^1,\iota^2)$ by
\[
\op{Dis}_{\frk X^1,\frk X^2}^\GHPUL\left(W,D , \iota^1, \iota^2 \right)   
:= \op{Dis}_{\frk X^1,\frk X^2}^\GHPU\left(W,D , \iota^1, \iota^2 \right)    +  \BB d_d^{\op L} \left(\iota^1(c^1) , \iota^2(c^2) \right),
\]
where $\op{Dis}_{\frk X^1,\frk X^2}^\GHPU(\cdot)$ is the GHPU distortion as defined in~\eqref{eqn-ghpu-var}.

The \emph{GHPUL distance} between $\frk X^1$ and $\frk X^2$ is given by
\[
\BB d^\GHPUL\left( \frk X^1 , \frk X^2 \right) 
= \inf_{(W, D) , \iota^1,\iota^2}  \op{Dis}_{\frk X^1,\frk X^2}^\GHPUL\left(W,D , \iota^1, \iota^2 \right),   
\]
where the infimum is over all compact metric spaces $(W,D)$ and isometric embeddings $\iota^1 : X^1 \rta W$ and $\iota^2 : X^2\rta W$.  
\xncomment{Following the same argument for the  completeness and separability of $(\BM^\GHPU,\BB d^\GHPU)$  in \cite[Proposition 1.3 and Section 2.2]{gwynne-miller-uihpq},} 
we see that  the space $(\BM^\GHPUL,\BB d^\GHPUL)$ is a complete separable metric space.

Recall $\Map^n$ in Theorem~\ref{emb-thm:main}.  
Let $\omega^n$ be sampled from
$\Ber_{\Map^n}$ with monochromatic  blue boundary condition and 
let $\dlp^n\defeq\Gamma(\Map^n,\omega^n)$ be the loop ensemble of $\omega^n$ defined in Section~\ref{subsub:quenched}.
Given a loop $\ell\in\dlp^n$, the edges  traversed by $\ell$ form an edge path. Therefore $\ell$ can be viewed as an unrooted oriented loop on $\Map^n$.
This way, $\dlp^n$ can be viewed as an element in $\cL(\Map^n)$ and $(\Map^n,d^n,\mu^n,\vec{\xi}^n,\dlp^n)$ is a random variable in $\BM^\GHPUL$.  
We write $(\Map^n,d^n,\mu^n,\vec{\xi}^n,\dlp^n)$ as $(\cM^n,\dlp^n)$ for simplicity. In Theorem~\ref{emb-thm:main}, $\{ (\cM^n,\dlp^n_i) \}_{i\in \N}$ 
should be understood as a sequence of identically distributed random variables in $\BM^\GHPUL$ with the law of $(\cM^n,\dlp^n)$. 

\subsection{$\sqrt{8/3}$-Liouville quantum gravity}
\label{subsec:lqg}
Let us recall the definition of the Gaussian free field (GFF). Let $D\subsetneq \C$ be a simply connected domain and let $h$ be a random distribution on $D$. 
We call $h$ a \emph{zero-boundary GFF on $D$} if for any compactly supported smooth function $f: D\to\R$,
$(h,f)$ is a centered Gaussian with variance $\iint f(x)G_{\op D}(x,y)f(y)\,d^2 x\, d^2 y$, where $G_{\op D}(\cdot,\cdot)$ is the Green's function on $D$ with Dirichlet boundary condition.
We call $h$ a \emph{free-boundary GFF on $D$} if for any smooth function $g$ on $D$ with $\int_D g(x) \,d^2 x=0$, $(h,g)$  is a centered Gaussian  with variance $\iint f(x)G_{\op N}(x,y)f(y)\, d^2x\, d^2y$,
where $G_{\op N}(\cdot,\cdot)$ is the Green's  function on $D$ with Neumann boundary condition. The law of the zero-boundary GFF is unique while the law of free-boundary GFF is only unique up to additive constant. The zero-boundary GFF and the free-boundary GFF are not pointwise defined functions, but almost surely belong to  the Sobolev space $H^{-1}(D)$. 
We refer to \cite{shef-gff,shef-zipper,wedges} for more details on the GFF. 

Let $\ol\Dh=\{(D,h): \textrm{$D\subsetneq  \C$  is a simply connected $C^0$ domain,  $h$ is a distribution on $D$} \}$. 
Fix $\gamma\in (0,2)$. Given $(D,h), (\wt D,\wt h) \in \ol\Dh$, let $\phi:\wt D\to D$ be a conformal map. We write
\begin{equation}\label{emb-eqn-lqg-coord}
	(D,h)\overset{\phi}{\sim}_{\gamma} (\wt D,\wt h)\textrm{ if and only if } \wt h = h \circ \phi + {\mathbf Q}\log |\phi'| \textrm{ for }{\mathbf Q}\defeq 2/\gamma+\gamma/2.
\end{equation}
We write  $(D,h)\sim_{\gamma} (\wt D,\wt h)$ if and only if there exists a conformal map $\phi:\wt D\to D$ such that $(D,h)\overset{\phi}{\sim}_{\gamma} (\wt D,\wt h)$. Then $\sim_\gamma$ defines an equivalence relation on $\ol\Dh$.
Let  $\Dh_\gamma:=\ol\Dh/\mathord\sim_\gamma$.  
By the Riemann mapping theorem, $\Dh_\gamma$ is in bijection with  distributions on $\bbH$  if we identify distributions $h$ and $\wt h$ on $\bbH$ satisfying $(\bbH,h){\sim}_{\gamma} (\bbH,\wt h)$.
This allows us to define a topology on $\Dh_\gamma$ from the natural topology on distributions on $\bbH$ so that we can  consider the Borel $\sigma$-algebra and probability  measures on $\Dh_\gamma$. An element in $\Dh_\gamma$ is called a \emph{generalized surface} with disk topology. 
A random variable taking values in $\Dh_\gamma$ is called
a \emph{$\gamma$-Liouville quantum gravity surface} ($\gamma$-LQG surface). More generally, we can define  generalized surfaces decorated with additional structures, such as metrics, measures, points, and curves.
\begin{definition}\label{def:mark}
	For $i=1,2$, let  $(D^i,h^i)\in \ol\Dh$. Let $d^i$, $\mu^i$, $x^i$, and $\eta^i$ be a metric, a measure, a point, and a curve on $D\cup \bdy D$, respectively.  
	Let $\phi: D^2\to D^1$ be a conformal map. If
	$(D^1,h^1)\overset{\phi}{\sim}_\gamma(D^2,h^2)$, $d^2(\cdot,\cdot)=d^1(\phi(\cdot),\phi(\cdot))$, $\mu^1=\phi_*\mu^2$, $x^1=\phi(x^2)$, and $\eta^1=\phi\circ \eta^2$,
	we write  $(D^1,h^1,d^1,\mu^1,x^1,\eta^1) \overset{\phi}{\sim}_\gamma(D^2,h^2,d^2,\mu^2,x^2,\eta^2)$.
	If there are multiple metrics, measures, points, and/or curves, define $\overset{\phi}{\sim}_\gamma$ similarly. 
	We define the equivalence relation $\sim_\gamma$ for these tuples  in the same way as we defined $(D,h)\sim_\gamma(\wt D,\wt h)$.
\end{definition}
\begin{conv}\label{conv:DF}
	In this paper we focus on $\gamma=\sqrt{8/3}$. Accordingly, ${\mathbf Q}=5/\sqrt{6}$ in \eqref{emb-eqn-lqg-coord}.
	We will simply write $\Dh$,  $\overset{\phi}{\sim}$, and $\sim$ instead of $\Dh_{\sqrt{8/3}}$,  $\overset{\phi}{\sim}_{\sqrt{8/3}}$, and $\sim_{\sqrt{8/3}}$, respectively.
	In particular, if $S$ is an element in $\ol\Dh$, possibly with decorations as in Definition~\ref{def:mark}, then we write its equivalence class under $\sim$ as $S/\mathord\sim$. 
\end{conv}

Next we introduce a general class of random distributions which covers all GFF type distributions considered in this paper,
such as the ones in Definition~\ref{def:disk} and in Section~\ref{subsub:surface}.
\begin{definition}[Free Liouville field]\label{def:free}
	A random distribution $\wh h$ on $\bbH$ is called a \emph{free Liouville field on $\bbH$} if there exists a pair $(h',g)$ such that 
	\begin{enumerate}
		\item\label{item:law} 
		$h'$ is a free-boundary GFF on $\bbH$,  $g$ is a random  function on $\bbH\cup\partial\bbH$ which is continuous except at finitely many points on $\p\bbH$;
		\item \label{item:sum} the law of $\wh h$ is absolutely continuous with respect to the law of $h'+g|_{\bbH}$.
	\end{enumerate}
	Given a simply connected domain $D$, a random distribution $h$ on $D$ is called a \emph{free Liouville field} on $D$ 
	if there exists a free Liouville field $\wh h$ on $\bbH$ such that  $(D,h)\sim(\bbH, \wh h)$.
\end{definition}

Set $\gamma=\sqrt{8/3}$ as in Convention~\ref{conv:DF}.
Let $D$ be a simply connected $C^0$ domain and let $h$ be a free Liouville field on $D$.
According to \cite{shef-kpz}, one can define the  \emph{$\sqrt{8/3}$-LQG area measure} $\mu_h=:``e^{\gamma h} \,d^2 z$'' by a regularization procedure   $\lim_{\eps\to 0} \eps^{\gamma^2/2} e^{\gamma h_\eps}$, where $h_\eps$ is the circle average modification of $h$; see Definition~\ref{def:GMC}.
Let  $\phi:\bbH \to D$ be a conformal map and $\wt h$ be such that $(D,h)\overset{\phi}{\sim} (\bbH,\wt h)$.
One can similarly define a nontrivial measure $\xi_{\wt h} :=``e^{\gamma \wt h(z)/2}\, d z$'' on $\bdy\bbH$ and then define 
$\xi_{h}:=(\phi^{-1})_*\xi_{\wt h}$.  By \cite{shef-kpz}, the definition of $\xi_h$ does a.s.\ not depend on the choice of $\phi$ (see also \cite{shef-wang-lqg-coord}).
We call $\xi_{h}$ the \emph{$\sqrt{8/3}$-LQG boundary measure} of $(D,h)$. 
By \cite{lqg-tbm1,lqg-tbm2}
a  metric $d_h$ corresponding to the metric tensor $(e^{\gamma h/4})^2(dx^2+dy^2)$ may be defined on $D\cup \p D$ using a growth process called the quantum Loewner evolution (QLE). 
Recently, \cite{GM-metric,DDDF19} constructed $d_h$ via a direct regularization procedure similar to the area case.
We list two important properties of $(d_h,\mu_h,\xi_h)$:
\begin{align}
	\label{eq:scaling}
	&\mu_{h+c}=e^{\gamma c}\mu_h, \quad \xi_{h+c}=e^{\gamma c/2}\xi_h, \quad \textrm{and}\quad d_{h+c}=e^{\gamma c/4}d_h\qquad \text{a.s.\,\,} \forall c\in\R.\\
	&(D,h,d_h,\mu_h,\xi_h)    \overset{\phi}{\sim} (\bbH,\wt h, d_{\wt h}, \mu_{\wt h},\xi_{\wt h})  \textrm{ a.s.}\label{eq:conf-inv}
\end{align}

Now we introduce the main $\sqrt{8/3}$-LQG surface that will be considered in this paper. 
It will be most convenient to introduce it on the horizontal strip $\cS=\R\times (0,\pi)$.
Let $h$ be a free-boundary GFF on $\cS$. Then $h$ can uniquely written as $h=h^{\op c}+h^{\ell}$, where 
$h^{\op c}$ is constant on vertical lines of the form $u + [0,i\pi]$ for $u\in \R$, and $h^{\ell}$  has mean zero on all such
vertical lines. Since  the law of the free-boundary GFF is unique modulo an additive constant, 
the law of $h^{\ell}$ does not depend on the choice of additive constant for $h$, and we call $h^{\ell}$ the \emph{lateral component} of the free-boundary GFF on $\cS$.

\begin{definition}[$\sqrt{8/3}$-LQG disk]
	\label{def:disk} 
	Let $\gamma=\sqrt{8/3}$,  ${\mathbf Q}=5/\sqrt{6}$, and $a={\mathbf Q}-\gamma=1/\sqrt{6}$.
	Let $(X_t)_{t\in\R}$ be  such that  $(X_t)_{t\ge 0}$ has the law of $B_{2t}-at$, where $B_t$ is a standard Brownian motion starting at the origin.
	Furthermore, $(X_{-t})_{t\ge 0}$ is independent of $(X_t)_{t\ge 0}$ and has the law of $B_{2t}-at$ conditioned on being negative.\footnote{Here we condition on a zero probability event. This can be made sense of via a limiting procedure.}
	Let  $h^1(z)=X_t$ for each $z\in \cS$ and $t\in \R$ with $\op{Re}z=t$.  
	Let $h^2$ be a random distribution on $\cS$ independent of $X_t$ which has the law of the lateral component of the free-boundary GFF on $\cS$.
	Let $\hs=h^1+h^2$ and $M:=\sup_{t\in \R} X_t$.
	Let $\hd$ be a random distribution on $\cS$, whose law is given by 
	\begin{equation}\label{eq:Boltzmann}
		\hs-2\gamma^{-1}\log\xi_\hs(\partial\cS) \qquad\textrm{reweighted  by $e^{- 2({\mathbf Q}-\gamma)M}\xi_\hs(\partial \cS)^{4/\gamma^2-1}$}.
	\end{equation}
\end{definition}
\begin{comment}
The reason we define one-point disk instead of two-point, is because in the rest of the paper we solely work on one-point disk, although the definition above is removing one point from the two-point disk. 
\end{comment}
\begin{remark}[Equivalence of definitions of $\sqrt{8/3}$-LQG disk]\label{rmk:disk}
	Various equivalent definitions of the unit boundary length $\sqrt{8/3}$-LQG disk are given in~\cite{wedges,sphere-constructions}. 
	We choose to work with Definition~\ref{def:disk} because the field is  described explicitly.
	Here we show the equivalence of Definition~\ref{def:disk} and the construction in~\cite[Section~4.5]{wedges}. In the notations of Definition~\ref{def:disk}, the construction in \cite{wedges} can be described as follows.
	Let $\P^{\op s}$ be the probability measure given by $\hs$ before the reweighting in \eqref{eq:Boltzmann} and let  $\ol\hs\defeq \hs - M$. Let $\ol\bdy\defeq \xi_{\ol\hs}(\bdy \cS)$ so that $e^{-2({\mathbf Q}-\gamma)M}\xi_\hs(\partial \cS)=\ol\bdy$.
	Let the pair $(e^*,\ol h^s)$ be sampled from the product measure $\1_{x>0}x^{4/\gamma^2}\, d x\otimes d \P^{\op s}$. 
	Then the conditional law of $(\cS, \ol \hs+2\gamma^{-1}\log e^*,+\infty)$ given the event $e^{*} \ol \bdy =1$ is the unit boundary $\sqrt{8/3}$-LQG disk as defined in \cite{wedges}.
	
	To see the equivalence with Definition~\ref{def:disk}, we first note that  when $e^{*} \ol \bdy =1$, we have $ \ol \hs+2\gamma^{-1}\log e^*= \ol \hs -2\gamma^{-1}\log \ol\bdy= \hs-2\gamma^{-1}\log\xi_\hs(\partial\cS)$.
	Moreover,  for each $\eps>0$, by Bayes' rule,  the conditional law  $\P^{\op s}\left[\cdot \mid  e^* \ol\bdy\in [1, 1+\eps] \right ]$ equals $c \ol\bdy^{4/\gamma^2-1}\, d\P^{\op s}$, where $c$ is a normalizing constant not depending on $\eps$. Sending $\eps\rta 0$ we obtain the equivalence.
\end{remark} 

We now give the precise definition of the field $\bh$ in Theorem~\ref{emb-thm:main}.
\begin{definition}\label{def:bh}
	Let $\phi:\D\to\cS$ be the conformal map satisfying $\phi(0)=\pi i/2$ and $\phi(1)=+\infty$. 
	Let $\bh$ be the free Liouville  field on $\D$ such that  $(\cS,h^{\op d})\overset{\phi}{\sim} (\D,\bh)$, where $h^{\op d}$ is as in Definition~\ref{def:disk}.
\end{definition}
By \eqref{eq:conf-inv},  Theorem~\ref{emb-thm:main} remains true  if we replace $\phi$ by another conformal map from $\D$ to $\cS$.
We choose this particular definition both for concreteness and for technical convenience in Section~\ref{sec:GPS} (see Lemma~\ref{lem:g}).

The Brownian disk $\BD_1$ can be identified with $(\D,\bh)$  in Theorem~\ref{emb-thm:main} as follows.
\begin{theorem} (\cite{lqg-tbm2}) \label{thm:LQG-TBM}
	Let  $\bh$ be as in Definition~\ref{def:bh}, let $(d_\bh,\mu_\bh,\xi_\bh)$ be as above \eqref{eq:scaling},
	and let $\vec{\xi}_\bh$ be a curve of duration 1 which traces $\p D$ clockwise starting from $1$ in the speed specified by $\xi_\bh$.
	%Identify the boundary measure $\xi_\bh$ with a curve of duration 1 which traces $\p D$ clockwise starting from $1$ in the speed specified by $\xi_\bh$.	
	Then there exist constants $c_{\op d}, c_{\op m} > 0$ such that
	$(\D\cup \p \D, c_{\op d} \,d_{\bh},  c_{\op m}\, \mu_{\bh}, \vec{\xi}_{\bh})$, viewed as a random variable  in $\BM^\GHPU$, is a free Brownian disk with unit perimeter. 
\end{theorem} 
We conclude this section by the precise description of the law of $h_\Delta$ in Theorem~\ref{emb-thm:Cardy}.
Let $\bh$ be as in Definition~\ref{def:bh}. Conditioning on $\bh$, independently sample two  points $v_1,v_2$ on $\p \D$ according to the measure $\xi_\bh$. By possibly relabeling $v_1$ and $v_2$, we assume that $1,v_1,v_2$ are ordered counterclockwise. Let $\psi:\D\to\Delta$ be the conformal map that maps $1$, $v_1$, and $v_2$ to $(1,0,0)$, $(0,1,0)$, and $(0,0,1)$, respectively.
\begin{definition}\label{def:hdelta}
	In Theorem~\ref{emb-thm:Cardy}, $h_\Delta$ denotes a random distribution with the law of $\bh\circ\psi +{\mathbf Q}\log|\psi'|$, where
	$(\bh,\psi)$ is defined as in the paragraph above.
	Moreover, $d_\Delta\defeq c_{\op d} \,d_{h_\Delta},  \mu_\Delta\defeq c_{\op m}\, \mu_{h_\Delta}$, and $\xi_\Delta\defeq \xi_{h_\Delta}$, 
	with $d_{h_\Delta},\mu_{h_\Delta}, \xi_{h_\Delta}$ as described above \eqref{eq:scaling} and constants $c_{\op d}, c_{\op m}$ as in Theorem~\ref{thm:LQG-TBM}.
\end{definition}

\subsection{Chordal SLE$_6$ and CLE$_6$}\label{subsec:SLE}  
Let $\CD=\{(D,a,b): \textrm{ $D$ is a simply connected $C^0$ domain, } a,b\in \bdy D,\, a\neq b \}$.
The clockwise (resp.,  counterclockwise) arc on $\bdy D$ from $a$ to $b$ is called the left (resp., right) boundary of $(D,a,b)$.
Suppose $\eta$ is a curve on $D\cup \bdy D$ from $a$ to $b$ for some $(D,a,b)\in \CD$.
For each $t\ge 0$ with $\eta(t)\in D\cup\p D$, let $D_t$ be the  connected component of $D\setminus \eta([0,t])$ whose boundary contains $b$. Otherwise, let $D_t=\emptyset$.
For each $(D,a,b)\in \CD$, the (chordal) $\SLE_6$  on $(D,a,b)$ is a probability measure  on  non-self-crossing curves on $D\cup \bdy D$ from $a$ to $b$ modulo increasing reparametrization.
$\SLE_6$ is uniquely characterized by the following three properties.
\begin{compactitem}
	\item {\bf Conformal invariance:} Suppose $\phi$ is a conformal map from $D$ to another simply connected $C^0$ domain $D'$. 
	Then $\eta$ has the law of an $\SLE_6$ on $(D,a,b)$ if and only if $\phi\circ\eta$ (modulo increasing parametrization)  has the law of an $\SLE_6$ on $(D',\phi(a),\phi(b))$. 
	\item {\bf Domain Markov property:} Let $\eta$ be an $\SLE_6$ on $(D,a,b)$, parametrized such that the parametrization on each initial segment is determined by the same segment modulo increasing parametrization. For each $t>0$, on the event $D_t\neq \emptyset$, we have that $D_t$ is $C^0$ a.s.\ and 
	the conditional law of $\eta$ after $t$ is that of an $\SLE_6$ on $(D_t, \eta(t),b)$.
	\item {\bf Target invariance:}  Let $\eta$ (resp., $\eta'$) be a chordal SLE$_6$ on $(D,a,b)$ (resp.,  $(D,a',b')$) such that $b\neq b'$. Let $\tau$ (resp., $\tau'$)  be the first time $\eta$ (resp., $\eta'$) hits the arc on $\bdy D$ between $b$ and $b'$ that does not contain $a$.
	Then $\eta|_{[0,\tau]}$ and $\eta|_{[0,\tau']}$ are equal in law modulo increasing reparametrization.
\end{compactitem}
It is proved by Schramm \cite{schramm0} that the first two properties define a one-parameter family of random curves called (chordal) $\SLE_\kappa$ with $\kappa\in(0,\infty)$.
The target invariance property singles out $\SLE_6$. 
By \cite{schramm-sle}, if $\eta$ is an SLE$_6$ curve on $(D,a,b)$, then $\eta$ is a.s.\ a non-simple curve which create ``bubbles'' (bounded simply connected domains) by hitting its past and the domain boundary. Furthermore, the range of $\eta$  has zero Lebesgue measure a.s. When $D_t\neq\emptyset$, let $\eta^t_{\markl}$ and $\eta^t_{\op r}$  be the left and right, respectively, boundary of $(D_t,\eta(t),b)$. 
For $t>0$,   the  laws of $\eta^t_{\markl}$ and $\eta^t_{\op r}$ away from $\p D$ are variants of $\SLE_{8/3}$ \cite{dubedat-duality}.  We refer to \cite{lawler-book} for more background on $\SLE_6$.

Given $\mesh>0$ and a Jordan domain $D$, let $D^\mesh$ be the $\mesh$-approximation of $D$ (see Section~\ref{subsec:notation}). 
Let $\omega^\mesh$ be a  Bernoulli-$\frac12$ site percolation on $D^\mesh$, namely, each inner vertex of $D^\mesh$ is colored  red or blue independently with probability $\frac12$.
Let $\Gamma^\mesh$ be the loop ensembles of $\omega^\mesh$ with monochromatic blue boundary condition. 
\xncomment{Let  $\cL(D)$ be the space of loop ensembles on $D\cup \bdy D$ 
endowed with the $\BB d^{\op L}_d$-metric (see Section~\ref{subsec:pre}), where $d$ is the Euclidean metric on $D$. 
Note that $\cL(D)$ is complete and separable; see \cite[Section 2.2]{camia-newman-sle6}.}
\begin{theorem}[\cite{camia-newman-sle6}]  \label{emb-thm:CLE}
	As  $\mesh\to 0$, $\Gamma^\mesh$ converge in law to a random variable $\Gamma$ in
	$\cL(D)$ in the $\BB d^{\op L}_d$-metric.
\end{theorem}
We take Theorem~\ref{emb-thm:CLE} as our definition of $\CLE_6$ on $D$. 
\begin{definition}[$\CLE_6$]\label{def:CLE}
	A random variable  in $\cL(D)$ with the law of $\Gamma$ is called a \emph{$\CLE_6$ on $D$ with monochromatic blue boundary condition}. 
	A random variable with the law of the  loop ensemble obtained by reversing the orientation of each loop in $\Gamma$ is called a \emph{$\CLE_6$ on $D$ with monochromatic red boundary condition}. 
\end{definition}

For $\Gamma$ in Definition~\ref{def:CLE}, with probability 1, for each  $z\in D$, the loop whose range is the single point $z$ belongs to  $\Gamma$. 
We call these loops trivial loops in $\Gamma$. There are countably many nontrivial loops in $\Gamma$ almost surely, whose $\BB d^{u}_d$-closure equals $\Gamma$. Throughout the paper when we declare a loop $\ell\in \Gamma$ we always assume that $\ell$ is a nontrivial loop.

We now explain how  to sample a $\CLE_6$ (with monochromatic boundary condition)
iteratively from chordal $\SLE_6$.  
We start by assigning an orientation to $\p D$.  
If we want the $\CLE_6$ to have  blue (resp., red) boundary condition, we assign clockwise (resp., counterclockwise) orientation to $\p D$.
Fix two distinct points $a,b\in \p D$. Let $\ol{ab}$ be the  segment on $\bdy D$ from $a$ to $b$ in the \emph{same} orientation as $\p D$.
We first sample an $\SLE_6$ $\eta^{ab}$ on $(D,a,b)$.
A connected component of $D\setminus\eta^{ab}$ is called a \emph{dichromatic bubble} if its boundary has non-empty intersection with $\ol{ab}$.
Let $\cB$ be a dichromatic bubble and let $x_\cB$ and $\wh x_\cB$ be the last and first, respectively, point on $\bdy \cB$ visited by $\eta^{ab}$, and let $\eta^\cB$ be the segment of $\eta^{ab}$ in between. For each dichromatic bubble $\cB$, conditioning on $\eta$, let $\eta_\cB$ be a chordal $\SLE_6$ on $(\cB,x_\cB,\wh x_\cB)$. Moreover, we assume that these $\eta_\cB$'s are conditionally independent given $\eta$.
Let $\ell_\cB$ be the oriented loop obtained by concatenating $\eta^\cB$ and $\eta_\cB$.
Let $\Gamma_a^b=\{\ell_\cB: \textrm{ $\cB$ is a dichromatic bubble} \}$.  
Suppose $\cB'$ is a connected component of $D\setminus\cup_{\ell\in \Gamma_a^b} \ell$. 
The orientation of loops in $\Gamma_a^b$ and  $\bdy D$ together
define an orientation on  $\bdy \cB'$, either clockwise or counterclockwise. 
If the orientation is clockwise (resp., counterclockwise), we call $\cB'$ a monochromatic blue (resp., red) bubble.
Conditioning on $\Gamma_a^b$, for each monochromatic bubble $\cB'$,  
\xncomment{we independently iterate the sampling procedure with the domain $D$ replaced by $\cB'$ with the already assigned orientation on $\bdy \cB'$.}

\begin{lemma}[\cite{camia-newman-sle6}]
	\label{lem:iteration}
	Given $\eta^{ab}$, $\Gamma_a^b$, and $\{\Gamma_{\cB'}\}$ as above,  let $\Gamma$ be the union of $\Gamma_a^b$ and the collection of nontrivial loops in $\Gamma_{\cB'}$, where $\cB'$ ranges over all monochromatic bubbles. Then if $\p D$ is oriented clockwise (resp., counterclockwise), then  $\Gamma$ has the law of the nontrivial loops of a $\CLE_6$ on $D$ with monochromatic blue (resp., red) boundary condition.  Moreover, $\Gamma$ determines $\Gamma_a^b$ and $\eta^{ab}$ almost surely. We call $\eta^{ab}$ the \emph{interface of $\Gamma$ on $(D,a,b)$}.
\end{lemma}
Both  $\Gamma_a^b$ and $\eta^{ab}$ can be defined as explicit functions of $\Gamma$. Consider all the loops in $\Gamma$ having nonempty intersection with $\ol{ab}$. There is a natural partial order $\prec$ on these loops where $\ell \prec \ell'$ if and only if $\ell$  is in a connected component of $D\setminus\ell'$ whose boundary contains neither $a$ nor $b$. Then $\Gamma_a^b$ is exactly the set of maximal elements for the partial order $\prec$. Moreover, for each loop $\ell\in\Gamma_a^b$, it is possible to recover its corresponding dichromatic bubble $\cB$, $\eta_\cB$ and $\eta^\cB$. By concatenating $\eta^\cB$ for all $\cB$ and taking a closure, we obtain $\eta^{ab}$.   

As a consequence of the iterative construction above and the conformal invariance of $\SLE_6$, 
the law of $\CLE_6$ is also conformally invariant. 
Namely, let $\Gamma$ be a $\CLE_6$ on a Jordan domain $D$. Let $D'$ be another Jordan domain  and let $\phi:D\to D'$ be a deterministic conformal map. 
Then the law of $\{\phi\circ \ell\}_{\ell\in \Gamma} $ is a $\CLE_6$ on $D'$ with the same boundary condition as $\Gamma$.

Now we record some important geometric properties of $\CLE_6$.
Suppose we are in the setting of Definition~\ref{def:CLE}. 
For each $\ell\in \Gamma$, let $\neg\ell$ be the connected component of $\C\setminus \ell$ whose closure contains $\bdy D$, where (here and below) we identify  $\ell$ with its range.
Let  $\reg (\ell)$ be the closure of the union of all connected components of $\C\setminus\ell$ other than $\neg\ell$ whose boundary is visited by  $\ell$  in the same orientation as $\ell$ is visiting $\bdy(\neg\ell)$. 
We call $\reg(\ell)$ the \emph{region enclosed by $\ell$}. 
Given $\ell\neq \ell'\in \Gamma$, we say that $\ell$ and $\ell'$ are \emph{nested} if and only if $\ell\subset \reg(\ell')$ or $\ell'\subset\reg(\ell)$.
\begin{figure}\label{fig:loop}
	\includegraphics[scale=1]{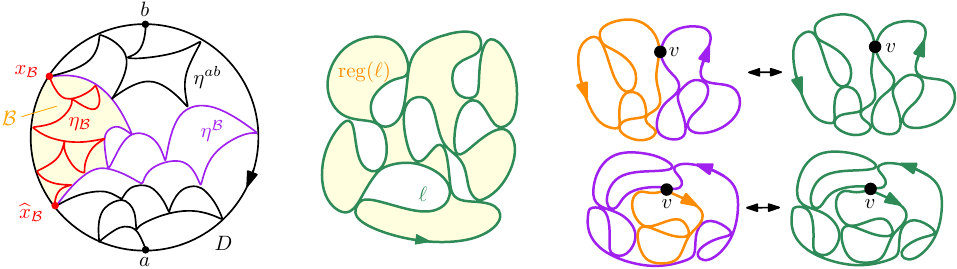}
	\caption{{\bf Left}: Illustration of the construction of a CLE$_6$ loop. The concatenation of the black curves and the purple curve is an SLE$_6$ $\eta^{ab}$ from $a$ to $b$. The domain $\cB$ (light yellow) is a dichromatic bubble. The CLE$_6$ loop $\ell_{\cB}$ is the concatenation of $\eta^{\cB}$ (purple) and $\eta_{\cB}$ (red). 
		{\bf Middle}: Illustration of the region $\reg(\ell)$ (light yellow) surrounded by the CLE$_6$ loop $\ell$.
		{\bf Right}: Illustration of the operation of flipping the color at a pivotal point $v$. In Case~\ref{item:int} of Definition~\ref{def:pivotal} we go from left to right, and in Case~\ref{item:double} of Definition~\ref{def:pivotal} we go from right to left. 
		The loops on the top (resp., bottom) left  are non-nested (resp., nested).}
\end{figure}

\begin{definition}[Pivotal point]\label{def:pivotal}
	Suppose $D$ and $\Gamma$ are as in Theorem~\ref{emb-thm:CLE}. A point $v\in \D$ is called a pivotal point of $\Gamma$ if one of the following two occurs:
	\begin{enumerate}
		\item There exist two loops $\ell,\ell'\in\Gamma$ such that $v\in\ell\cap\ell'$. 
		\label{item:int}
		\item There exists a  loop $\ell\in\Gamma$ that visits $v$ and $\ell$ visits $v$ at least twice.  
		\label{item:double}
	\end{enumerate} 
\end{definition}

The following basic properties of $\CLE_6$  are extracted from \cite{camia-newman-sle6}.
\begin{lemma}\label{lem:property}
	If $D$ and $\Gamma$ are as in Theorem~\ref{emb-thm:CLE}, then the following hold almost surely.
	\begin{itemize}
		\item {\bf local finiteness}: For each $\eps>0$, there exist finitely many loops in $\Gamma$  with diameter larger than $\eps$.
		\item {\bf finite chaining}: Given any $\ell\in \Gamma$ and $\ell' \in \Gamma\cup \{\p D \}$, there is a 
		finite set of  loops  $\ell_0=\ell,\ell_1,\dots, \ell_k=\ell'$ in $\Gamma$ such that for all $i\in \{1,\dots,k\}$, $\ell_{i-1}\cap\ell_i\neq \emptyset$. 
		\item {\bf no triple points}: $\Gamma$ has no pivotal points on $\p D$. If $v$ is a pivotal point of $\Gamma$, exactly one of the following holds: There exist exactly two loops $\ell,\ell'\in\Gamma$ that visit $v$, each of which visits $v$ exactly once; or 
		there exists a unique  loop $\ell\in\Gamma$ that visits $v$, and $\ell$ visits $v$ exactly twice.
		\item {\bf parity}:  Given any pair of loops in $\ell,\ell' \in \Gamma$ with $\ell\cap \ell' \neq\emptyset$,  $\ell,\ell'$ have opposite orientation if and only if they are nested.
		If $\ell\cap \p D\neq \emptyset$, then $\ell$ must be an outermost loop,  in the sense that  there exists no $\ell' \in \Gamma$ other than $\ell$ with $\ell \subset \reg(\ell')$.
	\end{itemize}
\end{lemma}
\xncomment{The local finiteness follows from~\cite[Lemma~6.6]{camia-newman-sle6}. The finite chaining and no-triple-point properties follow from~\cite[Theorem 2]{camia-newman-sle6}.
	The parity property follows from the no-triple-point property.}

If $v$ is a pivotal point of $\Gamma$, 
by \emph{flipping the color} at $v$, we mean merging $\ell,\ell'$ into a single loop in Case~\ref{item:int} of Definition~\ref{def:pivotal} and splitting ${\ell}$ into two loops in Case~\ref{item:double} of Definition~\ref{def:pivotal}. (See Figure~\ref{fig:loop}.)
If a loop does not visit $v$, flipping the color at $v$ keeps the loop unchanged.  
Let $\Gamma_v$ denote the set of loops obtained after flipping the color at $v$. 
By the parity  property of $\CLE_6$,   $\Gamma$ induces an orientation on each loop in $\Gamma_v$, making it an element of $\cL(D)$ (after including trivial loops). 
By the no-triple point property,  the symmetric difference  $\cL_v$ of $\Gamma$ and $\Gamma_v$ always contains exactly  three loops. Now we define the continuum $\eps$-pivotal points by mimicking the discrete definition in Section~\ref{subsub:DP}.
\begin{definition}[$\eps$-pivotal point]\label{def:eps}
	Given a Jordan domain $D$, let $\Gamma$ be a $\CLE_6$ on $D$  and let $h$ be a free Liouville field (see Definition~\ref{def:free}) on $D$ independent of $\Gamma$. 
	Given a pivotal point $v$ of $\Gamma$ and $\eps>0$, we call $v$ an $\eps$-pivotal point of $(h,\Gamma)$ if  $\mu_h(\reg(\ell))\ge \eps$ for all $\ell\in\cL_v$.
\end{definition}

\begin{remark}[$\CLE_6$ on top of $\sqrt{8/3}$-LQG]\label{rmk:CLE}
	Suppose we are in the setting of Theorem~\ref{thm:LQG-TBM}.
	Let $\Gamma$ be a $\CLE_6$ on $\D$ with monochromatic blue boundary condition. 
	Then $(\D\cup \p \D, c_{\op d} \,d_{\bh},  c_{\op m}\, \mu_{\bh}, \vec{\xi}_{\bh},\Gamma)$ is a random variable in $\BM^\GHPUL$.
	When it is clear from context, we will denote this random variable by
	$(\D,\bh,\Gamma)$. In  particular, $(\D,\bh,\Gamma_i)$ in Theorem~\ref{emb-thm:main}  should be understood in this sense. 
	Now Theorem~\ref{emb-thm:upper} asserts that $(\cM^n, \dlp^n)$ defined at the end of Section~\ref{subsec:pre} converge in law to $(\D,\bh,\Gamma)$ in the GHPUL topology.
\end{remark}

\section{A dynamical percolation on random triangulations}\label{emb-sec:overview}
In this section we prove Theorem~\ref{emb-thm:main}. The argument is ``soft''  as long as the ``hard'' input Lemmas~\ref{lem:eps-LDP} and~\ref{emb-prop:LDP}
are supplied. We postpone the proofs of these two lemmas to Section~\ref{sec:GPS}.

For $\eps>0$, recall the dynamics $(\Map^n,\ol\omega^{\eps,n}_t)_{t\ge 0}$ defined in Section~\ref{subsub:DP}.
The following elementary  observation is crucial to us. We leave the proof to the reader.
\begin{lemma}\label{lem:stationary}
	Conditioning on $\Map^n$,  the process $(\ol\omega^{\eps,n}_t)_{t\ge 0}$ is stationary.
\end{lemma}
For $t>0$, let $\ol\dlp^{\eps,n}_t\defeq \Gamma(\Map^n, \ol\omega^{\eps,n}_t)$ be the loop ensemble of $\ol\omega^{\eps,n}_t$.
Recall $\cM^n\in\BM^\GHPU$  in Section~\ref{subsec:pre}, which is obtained by rescaling $\Map^n$ according to \eqref{eq:scale}.	
We view $(\cM^n, \ol\dlp^{\eps,n}_t)_{t\ge 0}$ as a process taking values in $\BM^{\GHPUL}$ as explained at the end of Section~\ref{subsec:pre}.
In Section~\ref{sec:GPS}, we will prove the following. 
\begin{lemma}\label{lem:eps-LDP}
	For any fixed $\eps>0$,  $(\cM^n,\ol\dlp^{\eps,n}_i)_{i\in \N }$ converge in law as $n\rta\infty$ to a stationary  sequence $( Y^\eps_i )_{i\in \N }$  in the $\GHPUL$ topology  
\end{lemma}
We restrict the index set to positive integers in Lemma~\ref{lem:eps-LDP} to avoid unnecessary topological technicalities for continuous time processes.

Recall $(\D,\bh,\Gamma)$ in Remark~\ref{rmk:CLE}.
By Theorem~\ref{emb-thm:upper}, for each $i\in\N$, $Y^\eps_i$ in Lemma~\ref{lem:eps-LDP} is equal in law to $(\D,\bh,\Gamma)$ as a random variable in $\BM^\GHPUL$.
More generally, there exists a sequence of $\CLE_6$'s $(\ol\Gamma^\eps_i)_{i\in \N}$ coupled with $\bh$ such that $(Y^\eps_i)_{i\in \N}\eqd (\D,\bh,\ol\Gamma^\eps_i)_{i\in \N}$.

\begin{lemma}\label{emb-prop:LDP}
	Let $(\bh,\ol\Gamma^\eps_i)_{i\in \N}$ be defined as above.  There exists a sequence of $\CLE_6$'s $(\ol\Gamma_i)_{i\in \N}$ coupled with $\bh$ 
	such that as $\eps\rta 0$, $(\bh,\ol\Gamma^\eps_i)_{i\in \N}$ converge in law to $(\bh,\ol\Gamma_i)_{i\in \N}$ in the $H^{-1}(\D)\times \cL(\D)$ topology. Moreover, $(\ol\Gamma_i)_{i\in \N}$ is 
	stationary and ergodic.
\end{lemma}
To deduce Theorem~\ref{emb-thm:main} from the above lemmas we  use the following observation.
\begin{lemma}\label{emb-lem:basic}
	Let $X$ and $(Y_i)_{i\in \N}$ be random variables on the same probability space. Suppose $(X,Y_i)_{i\in \N} $ is stationary  and  
	$(Y_i)_{i\in \N}$ is ergodic. Then $X$ and $Y_1$ are independent.
\end{lemma}
\begin{proof}	
	Let $f$ and $g$ be two  bounded real-valued measurable functions defined on the space in which $X$ and $Y_i$, respectively, take values. By stationarity of $(X,Y_i)_{i\in\N}$,
	\[
	\Cov(g(X),f(Y_1)) = \Cov\left(g(X), \frac 1m\sum_{i=1}^{m}f(Y_i)\right).
	\]
	Now Lemma~\ref{emb-lem:basic} follows from the Birkhoff ergodic theorem. 
\end{proof}
\begin{proof}[Proof of Theorem~\ref{emb-thm:main}]
	Fix $\eps\in (0,1)$.  Consider the process $(\Map^n,\ol\omega^{\eps,n}_t)_{t\ge 0}$ in Lemma~\ref{lem:stability}.
	Conditioning on $\Map^n$, let $\omega^n$ be sampled from $\Ber_{\Map^n}$ such that $\omega^n$ is conditionally independent of $(\ol\omega^{\eps,n}_t)_{t\ge 0}$.
	Let $\dlp^n=\Gamma(\Map^n,\omega^n)$.
	By Theorem~\ref{emb-thm:upper}, $(\cM^n, \ol\dlp^{\eps, n}_i)_{i\in \N}$  and $(\cM^n,\dlp^n)$ are tight in the $\GHPUL$ topology.
	\xncomment{At this moment, we do not know that they jointly converge in law. Indeed, the joint convergence of $(\cM^n, \ol\dlp^{\eps, n}_1)$  and $(\cM^n,\dlp^n)$
		is precisely the $k=2$ case of Theorem~\ref{emb-thm:main}. We now use the previous lemmas to show that in any subsequential limit, the joint limiting law of $(\cM^n, \ol\dlp^{\eps, n}_1)$  and $(\cM^n,\dlp^n)$ is as desired.}
	
	By the Skorokhod representation theorem, given any subsequence $\cN\subset\N$,  we can  choose a further subsequence $\cN'\subset\cN$ such that  there exists  a coupling of 
	$\{(\Map^n, \omega^n,\ol\omega^{\eps,n}_i)_{i\in\N} : n\in \cN'\}$ 
	where both $(\cM^n, \ol\dlp^{\eps, n}_i)_{i\in \N}$  and $(\cM^n,\dlp^n)$ have almost sure $\GHPUL$ limits as $n\rta \infty$ along $\cN'$.
	By Lemma~\ref{thm:LQG-TBM} the  $\GHPU$ limit of $\cM^n$ can be written as $(\ol\Delta, c_{\op d} \,d_{\bh},  c_{\op m}\, \mu_{\bh}, \vec{\xi}_{\bh})$, where $\bh$ is as defined in Definition~\ref{def:bh}.
	As in Lemma~\ref{lem:eps-LDP} we denote the $\GHPUL$ limit of $(\cM^n, \ol\dlp^{\eps, n}_i)_{i\in \N}$ by $(\D,\bh,\ol \Gamma^\eps_i)_{i\in \N}$, where $(\ol \Gamma^\eps_i)_{i\in\N}$ is a sequence of $\CLE_6$'s on $\D$.  
	By Theorem~\ref{emb-thm:upper}, there exists a $\CLE_6$ $\Gamma$ on $\D$ with monochromatic blue boundary condition such that $(\D, \bh,\Gamma)$ is the $\GHPUL$ limit of $(\cM^n,\dlp^n)$.
	Moreover, $(\bh, \Gamma, \ol\Gamma^{\eps}_i)_{i\in\N}$ is stationary. 
	
	\xncomment{By Lemma~\ref{emb-prop:LDP},  the $\eps$-indexed family  $\{(\bh, \Gamma, \ol\Gamma^{\eps}_i)_{i\in\N}: \eps>0 \}$ is tight. }
	Therefore we can choose a sequence $\eps_m\downarrow 0$ such that as $m\to\infty$,
	$(\bh, \Gamma, \ol\Gamma^{\eps_m}_i)_{i\in\N}$ converge in law to a stationary sequence, which we denote by $(\wt h,\wt \Gamma,\ol\Gamma_i)_{i\in \N}$. 
	\xncomment{(Similarly as above,  due to the presence of $\Gamma$ we need to choose a sequence $\eps_m\downarrow 0$ to ensure  convergence in law instead of simply taking $\eps\downarrow 0$.)}
	Applying Lemma~\ref{emb-lem:basic} to $X=(\wt h,\wt \Gamma)$ and $Y_i=\ol\Gamma_i$, we see that $(\wt h,\wt \Gamma)$ is independent of $\ol\Gamma_1$.
	Since the law of $(\Map^n,\dlp^n,\ol\dlp^{\eps,n}_1)$ is equal to the law of $(\Map^n,\dlp^n_1,\dlp^n_2)$ in Theorem~\ref{emb-thm:main}, which does not depend on $\eps$,
	the law of $(\bh,\Gamma,\ol\Gamma^\eps_1)$ does not depend on $\eps$ either. In fact, it must equal the law of $(\wt h,\wt \Gamma,\ol\Gamma_1)$.
	Therefore $(\bh,\Gamma)$ is independent of $\ol\Gamma^\eps_1$.
	In particular, the law of $(\bh,\Gamma,\ol\Gamma^\eps_1)$ does not depend on the choice of subsequences $\cN$ and $\cN'$. 
	Therefore  $(\cM^n,\dlp^n)$ and $(\cM^n,\dlp^{\eps,n}_1)$  jointly converge in law to $(\D,\bh,\Gamma)$ and $(\D,\bh,\ol\Gamma^\eps_1)$, respectively. 
	This gives Theorem~\ref{emb-thm:main} when $k=2$.
	
	For $k\ge 3$ we assume by induction that Theorem~\ref{emb-thm:main} holds for $k-1$. 
	Now we replace $\omega^n$ above by $k-1$ independent percolations sampled from $\Ber_{\Map^n}$ 
	and apply the exact same argument as above. Then by the induction hypothesis, $\Gamma$ above becomes $k-1$ independent copies of $\CLE_6$ which are also independent of $\bh$. We again use Lemma~\ref{emb-lem:basic} to conclude the proof.
\end{proof}

\section{Convergence under the Cardy embedding}
\label{emb-subsec:Cardy-conv}

In this section we will conclude the proof of Theorems \ref{emb-thm:Cardy} and \ref{thm:crossings}.

Recall  $h_\Delta$, $d_\Delta=c_{\op d} d_{h_\Delta}$, $\mu_\Delta=c_{\op m}\mu_{h_\Delta}$,  and \xncomment{$\xi_\Delta=\xi_{h_\Delta}$} in Theorem~\ref{emb-thm:Cardy}, whose precise meaning can be found in  Definition~\ref{def:hdelta}. Let $\Gamma$ be a CLE$_6$ on $\Delta$ with monochromatic blue boundary condition independent of $h_\Delta$.
Then we can identify $(\Delta,h_\Delta, \Gamma)$ with a random variable in $\BM^{\GHPUL}$ as explained in Remark~\ref{rmk:CLE}, with $(\D,\bh)$ replaced by $(\Delta, h_\Delta)$.
We first state a basic variant of Theorem~\ref{emb-thm:main} for maps with marked points. 
Note that  elements in $\BM^\GHPUL$ with marked points can be naturally endowed a topology as in Section~\ref{subsec:pre} which  includes the convergence of the marked points.

\begin{lemma}
	Let $(\Map^n,a^n,b^n,c^n)$ and $\{\dlp^n_i\}_{i\in \N}$ be as in Theorem~\ref{emb-thm:main}. 
	Let $h_\Delta$ be as above and let $\{\Gamma_i\}_{i\in \N }$ be independent $\CLE_6$'s which  are also independent of $h_\Delta$.
	Let $(\wh v^n_1, \wh v^n_2,\wh v^n_3)\defeq(a^n,b^n,c^n)$. Let $\wh z_1$, $\wh z_2$, and $\wh z_3$ be equal to $(1,0,0)$, $(0,1,0)$, and $(0,0,1)$, respectively.
	Conditioning on $(\Map^n,\dlp^n_1,\dlp^n_2,\dots)$, let $\wh v^n_4$ (resp., $\{v^n_i: i\in \N\}$) be vertices of $\p \Map^n$ (resp., $\dsk^n$) which are sampled uniformly and independently at random. 
	Conditioning on $(h_\Delta,\Gamma)$, 
	let $\wh z_4$ (resp., $\{z_i:i\in \N\}$) be boundary (resp., interior) points of $\Delta$ which are sampled independently  from the measure $\xi_\Delta$ (resp.,  $\mu_\Delta$).	
	Then there exists a coupling such that for each  $m\in\N$, almost surely the following convergence holds in the GHPUL topology with $m+4$ marked points:
	\begin{equation}\label{eq:couple-GHPU}
		\lim_{n\to \infty}(\cM^n,\dlp^n_i,\wh v^n_1,\dots,\wh v^n_4,v^n_1,\dots,v^n_m)=
		(\Delta,h_\Delta, \Gamma_i,\wh z_1,\dots,\wh z_4,z_1,\dots,z_m),\qquad \textrm{for each }i\in \N.
	\end{equation}
	\label{emb-prop:GHPU-marked}
\end{lemma}
\begin{proof}
	By Skorokhod embedding theorem, it suffices to show that the convergence in~\eqref{eq:couple-GHPU} holds in law for a fixed $m$.
	The  convergence  of  $\wh v^n_1,\dots,\wh v^n_4$  follows from the uniform convergence of the boundary curve, and the convergence of $v_1^n,\dots,v^n_m$ follows from the convergence of $\mu^n$. The gives the desired convergence in law.
\end{proof}

Throughout this section we work under  a coupling as  described in Lemma \ref{emb-prop:GHPU-marked}.
We will prove that $(d^n_\Delta,\mu^n_\Delta,\xi^n_\Delta)$ converge to $(d_\Delta,\mu_\Delta,\xi_\Delta)$ in probability, which implies Theorem \ref{emb-thm:Cardy}.

First we will argue that for each fixed $i\in\N$, as $n\rta\infty$,
\eqb
\Cdy^n(v^n_i)\rta z_i\,\,\text{in probability for the Euclidean metric on $\ol\Delta$.}
\label{eq21}
\eqe
Since the total mass of $\mu^n_\Delta$ converge to that of $\mu_\Delta$ in probability and $\Cdy^n(v^n_i)$ (resp., $z_i$) has the law of a vertex (resp., point) sampled  according to $\mu^n_\Delta$  (resp., $\mu_\Delta$), \eqref{eq21} implies that $\mu^n_\Delta$ converge to $\mu_\Delta$ in probability. 

We fix $i\in \N$ while proving \eqref{eq21}.
For $j\in \N$, let $E^{j,n}_{1}\defeq\{\omega^n_j\in E_{a^n}(v^n_i)\}$, namely, $E^{j,n}_{1}$ is the event $E_{a^n}(v^n_i)$ in Definition~\ref{emb-def:cardy} for $\omega^n_j$. 
The dependence of $E^{j,n}_1$ on $i$ is dropped in the notation since $i$ is fixed.
Similarly, let $E^{j,n}_{2}\defeq\{\omega^n_j\in E_{b^n}(v^n_i)\}$ and $E^{j,n}_{3}\defeq\{\omega^n_j\in E_{c^n}(v^n_i)\}$. 
Let $E^{j}_{1},E^{j}_2,E^{j}_3$ be the continuum analogs of  $E^{j,n}_{1},E^{j,n}_2,E^{j,n}_3$ defined in terms of $z_i\in\Delta$ and the CLE$_6$ $\Gamma_j$. 
We describe $E^{j}_1$ precisely following \cite[Sections 7.9]{bhs-site-perc};  $E^{j}_2$ and $E^{j}_3$ can be defined similarly by permuting the indices.
Let $\eta$ be the interface of $\Gamma^j$ on $(\Delta, \wh z_3, \wh z_2)$ as defined in Lemma~\ref{lem:iteration}.
Then 
\begin{equation}\label{eq:event}
	\textrm{$E^j_1$ is the event that $z_i$ is strictly on the same side of $\eta$ as $\wh z_1$}.
\end{equation}
To be precise, the event $E^j_1$ occurs if and only if there is a path in $\Delta$ connecting $z_i$ and $\wh z_1$ which does not intersect $\sle$. 
By \cite[Proposition 6.7]{ghs-metric-peano} (which builds on \cite[Theorem 8.7]{bhs-site-perc}) the following convergence holds in probability
\eqb
(\1_{E^{j,n}_{1}},\1_{E^{j,n}_2},\1_{E^{j,n}_3})
\rta 
(\1_{E^{j}_{1}},\1_{E^{j}_2},\1_{E^{j}_3}),\qquad
j=1,\dots,k.
\label{emb-eq:conv1}
\eqe
\xncomment{Here we use the following basic measure theoretic fact. 
	\begin{lemma}\label{lem:law-prob}
		Suppose $(X_n,Y_n)_{n\ge 1}$ is a sequence of random variables taking values in 
		a complete separable metric space. Suppose $X_n\rta X$ in probability, $(X_n,Y_n)\rta(X,Y)$ in law, $Y_n$  is measurable with respect to $X_n$ for $n\ge 1$, and $Y$ is measurable with respect to $X$.
		Then $Y_n \rta Y$ in probability.
	\end{lemma}	
	See~\cite[Lemma~4.5]{ss-level} for  Lemma~\ref{lem:law-prob} with $X_n\equiv X$, which also works with little modification  if $X_n\rta X$ a.s.\ instead. 
	The  case where $X_n\rta X$ in probability 	follows  from the a.s.\ case by extracting subsequences. 
	
	For~\eqref{emb-eq:conv1}, we let $X_n=(\cM^n,\dlp^n_j,\wh v^n_1,\dots,\wh v^n_4,v^n_i)$ and $X=(\Delta,h_\Delta, \Gamma_j,\wh z_1,\dots,\wh z_4,z_i)$.
	We  let $Y_n$ and $Y$ be the two sides of~\eqref{emb-eq:conv1}. In the continuum, the interfaces and the crossing events are measurable with respect to  
	$(\Gamma_j,\wh z_1,\dots,\wh z_4,z_i)$; see the end of~\cite[Section 5]{camia-newman-sle6}.
}

It follows \xncomment{from~\eqref{emb-eq:conv1}} that for any fixed $k$ by choosing $n$ (depending on $k$ and $\zeta$) sufficiently large, we have with probability at least $1-\zeta$ that
\eqb
\frac 1k \sum_{j=1}^{k} (\1_{E^{j,n}_{1}},\1_{E^{j,n}_2},\1_{E^{j,n}_3}) 
=
\frac 1k \sum_{j=1}^{k} (\1_{E^{j}_{1}},\1_{E^{j}_2},\1_{E^{j}_3}).
\label{emb-eq:conv1b}
\eqe
\xncomment{By the Azuma–Hoeffding inequality for Bernoulli random variables, the following statement holds uniformly in $n$.
	For $k$ sufficiently large depending on $\zeta$, with probability at least $1-\zeta$,}
\eqb
\Big|
\frac 1k \sum_{j=1}^{k} (\1_{E^{j,n}_{1}},\1_{E^{j,n}_2},\1_{E^{j,n}_3}) 
- 
\big(\Ber_{\dsk^n} [E_{a^n}(v^n_i)] + \Ber_{\dsk^n}[E_{b^n}(v^n_i)] + \Ber_{\dsk^n}[E_{c^n}(v^n_i)] \big)
\Big|<\zeta
\label{emb-eq:conv2a}
\eqe
and
\eqb
\Big|
\frac 1k \sum_{j=1}^{k} (\1_{E^{j}_{1}},\1_{E^{j}_2},\1_{E^{j}_3}) 
- 
\big(\P[ E^1_1 ],\P[ E^1_2],\P[ E^1_3]\big)
\Big|
<\zeta.
\label{emb-eq:conv2b}
\eqe
Since  $\P[ E^1_1 ]+\P[ E^1_2]+\P[ E^1_3]=1$ by Theorem~\ref{thm:Smirnov}, on the event that \eqref{emb-eq:conv1b}, \eqref{emb-eq:conv2a}, and \eqref{emb-eq:conv2b} are satisfied,
\begin{equation}\label{eq:sum-typical}
	|\Ber_{\dsk^n} [E_{a^n}(v^n_i)] + \Ber_{\dsk^n}[E_{b^n}(v^n_i)] + \Ber_{\dsk^n}[E_{c^n}(v^n_i)]-1|<2\zeta.
\end{equation}
Combining this with \eqref{emb-eq:conv1b}, \eqref{emb-eq:conv2a}, and \eqref{emb-eq:conv2b} and the definition of the Cardy embedding, we get that with probability at least $1-3\zeta$, for all sufficiently large $n$ (depending only on $\zeta$),
\eqb
|\Cdy^n(v^n_i)-z_i| < 50\zeta.
\label{eq19}
\eqe
Since $\zeta$ was arbitrary, we obtain \eqref{eq21}, which concludes the proof that $\mu^n_\Delta\rta\mu_{\Delta}$ in probability.

We prove that $\xi^n_\Delta\rta\xi_{\Delta}$ in probability by a very similar argument. As above, it is sufficient to show that $\Cdy^n(\wh v^n_4)\rta \wh z_4$ in probability for the Euclidean metric as $n\rta\infty$. Again the result follows by applying \cite{bhs-site-perc}, which give convergence in probability of the three crossing events $\wh E^{j}_{1}, \wh E^{j}_2, \wh E^{j}_3$ (now defined with $\wh v^n_4$ instead of $v^n_i$). Note that the convergence result for $\wh E^{j}_{1}, \wh E^{j}_2, \wh E^{j}_3$ in \cite[Theorem~8.7]{bhs-site-perc} is stated for the case where the four boundary points have deterministic distances along the boundary from the root, rather than being sampled uniformly and independently at random, but the proof in \cite[Theorem~8.7]{bhs-site-perc} is identical for the case of random points.

We now establish a modulus of continuity estimate for the Cardy embedding. 
\begin{proposition}\label{prop:max}	
	\eqb
	\lim_{r\to 0}\sup_{u,v\in \cV(\dsk^n)\,:\,d^n(u,v)< r } |\Ber_{\dsk^n} [E_{a^n}(u)] -\Ber_{\dsk^n} [E_{a^n}(v)] |=0.
	\label{emb-eq:ber}
	\eqe
	Then same holds with  $E_{a^n}$ replaced by $E_{b^n}$ and $E_{c^n}$. 
\end{proposition}
Before proving Proposition~\ref{prop:max}, we first recall the notion of percolation interface  following \cite{ghs-metric-peano}.
Let $M$ be a triangulation of a polygon and let $e$ and $e'$ be two distinct edges on $M$. 
Recall that $(e,e')$ denotes the counterclockwise arc on $\p M$ from $e$ to $e'$.
The \notion{$(e,e')$-boundary condition}  for a site percolation on $M$ is the coloring of $\p M$ where vertices on $(e,e')$ (resp., $(e',e)$) are blue  (resp., red).  
Given a site percolation  $\omega_M$ on $M$, regardless of its own boundary condition, if we impose the $(e,e')$-boundary condition to it, then there is a unique edge path (recall Section~\ref{subsec:pre}) on $M$ from $e$ to $e'$, such that each edge on the path has a red vertex on its left side and a blue vertex on its right side. We call this path the \notion{percolation interface} of $\omega_M$ on $(M,e,e')$. Note that this percolation interface only depends on the coloring of the inner vertices.
\begin{proof}[Proof of Proposition~\ref{prop:max}]
	Given a percolation interface $\sle^n$ on $(M^n,c^n,b^n)$ of a site percolation on $M^n$, we call the segment between the last time  $\sle^n$ visits the counterclockwise arc $(c^n,a^n)$ and the first time $\sle^n$ visit the counterclockwise arc $(a^n,b^n)$ the \emph{middle segment} of $\sle^n$.  Here visits means passing through an edge with an endpoint on the arc.		
	Recall that $d^n$ is  the graph distance on $\dsk^n$ rescaled by $(3n/4)^{-1/4}$. 
	Given a $d^n$-metric ball $B$ on $M^n$, let $E^n(B)$ be the event that the middle segment of $\eta^n$ is
	passing though $B$.  Let $X^n(r)\defeq \max_{B}\{  \Ber_{M^n}(E^n(B)) \}$,  where $B$ ranges over all such  balls of radius $r$. We claim that 
	\begin{equation}\label{eq:max}
		\lim_{r\to 0} \limsup_{n\to \infty} X^n(r) =0.
	\end{equation}
	Let us first explain that \eqref{emb-eq:ber} follows from~\eqref{eq:max}. 
	Let $\omega^n$ be sampled from $\Ber_{\dsk^n}$ and $\eta^n$ be its percolation interface on $(M^n,c^n,b^n)$.
	It is elementary to check that the discrete analog of \eqref{eq:event} can be used to characterize the crossing events in terms of $\eta^n$; see e.g.\ \cite[Section 6.8]{bhs-site-perc}.
	As a consequence, given $u,v\in \cV(\dsk^n)$ and $r>0$ such that $d^n(u,v)< r$, if $E_{a^n}(u)\xor E_{a^n}(v)$ occurs
	then the middle segment of $\eta^n$ must cross the $d^n$-ball centered at $u$ of radius $r$.
	Therefore~\eqref{eq:max} implies~\eqref{emb-eq:ber}. 
	
	We prove \eqref{eq:max} by contradiction. Let $\sle_i$ be the interface of $\Gamma_i$ on $(\Delta,\wh z_3,\wh z_2)$ as defined in Lemma~\ref{lem:iteration}.
	We define the middle segment of $\sle_i$ to be the segment between the last time  $\sle_i$ visits the counterclockwise arc $(\wh z_3,\wh z_1)$ and the first time $\sle_i$ visit the   counterclockwise  arc $(\wh z_1,\wh z_2)$. Let 	$\omega^n_i$ be the site percolation corresponding to $\dlp^n_i$ in Lemma~\ref{emb-prop:GHPU-marked}. 
	If \eqref{eq:max} does not hold, then there exists $\zeta>0$ and  a sequence $r_n\to 0$ such that for each $n$, with probability at least $\zeta$,  
	there exists a $d^n$-ball $B$ of radius $r_n$ such that $E^n(B)$ occurs for the  each of $\omega^n_i$ $(1\le i\le 10)$.   
	\xncomment{As explained below~\eqref{emb-eq:conv1}, the discrete interfaces converge in probability in the coupling of Lemma~\ref{emb-prop:GHPU-marked}.}
	Since $r_n\to 0$, sending $n\to\infty$, we see that in this coupling with positive probability the middle segments of $\eta_i$ for $1\le i\le 10$ share a common point on $\ol\Delta$. This is not possible because $\SLE_6$ has dimension $7/4$. \xncomment{More precisely, by the one-point estimate by Beffara~\cite[Proposition 4]{beffara-dim} } the probability that an $\SLE_6$ passes through a ball of Euclidean radius $s$ decays like $s^{\xncomment{2-7/4}}$ uniformly over all balls bounded away from the corners of $\Delta$, and $(2-7/4)\cdot 10>2$. 
\end{proof}

\begin{proposition}\label{prop:max2}	
	$$
	\lim_{n\to \infty}\max_{v\in \cV(\Map^n)}    \left|\Ber_{\dsk^n} [E_{a^n}(v)] + \Ber_{\dsk^n}[E_{b^n}(v)] + \Ber_{\dsk^n}[E_{c^n}(v)]-1 \right| =0\textrm{ in probability}.
	$$ 
\end{proposition}
\begin{proof}
	Since $\mu_\Delta$ almost surely assigns positive mass to any open set of $\Delta$, $\{z_i: i\in \N\}$ is dense in $\ol\Delta$ for both Euclidean and the $d_\Delta$-metric. 
	\xncomment{(Recall that the $d_\Delta$-metric and the Euclidean metric induce the same topology as they are a.s.\ bi-H\"older with respect to each other~\cite[Theorem 2]{lqg-tbm2}.)}
	Since we are under the coupling in Lemma \ref{emb-prop:GHPU-marked}, where the convergence  is almost sure,  
	we have that $\lim_{n\to0}\sup_{v\in\cV(\Map^n)} \inf_{i\in\N} d^n(v,v_i^n)=0$  in probability.
	Proposition~\ref{prop:max2}  now follows from this observation,~\eqref{eq:sum-typical}, and Proposition~\ref{prop:max}.
\end{proof}

To conclude the proof of Theorem \ref{emb-thm:Cardy} we must show that $d^n_\Delta$ converge in probability to $d_{\Delta}$. 
For $\ol z\in\ol \Delta$ and $\zeta>0$, let $B(\ol z,\zeta)$ denote the Euclidean ball of radius $\zeta$ centered at $\ol z$.
\xncomment{As explained in the proof of Proposition~\ref{prop:max2}, since $\mu_\Delta$ almost surely has full support on $\Delta$, 
	$\{z_i: i\in \N\}$ is dense in $\ol\Delta$, for both the Euclidean metric and the $d_\Delta$-metric.} 
Combined with \eqref{eq21},  we have
\eqb
\liminf_{n\rta\infty}\P[{\textrm{for each }\ol z\in \ol\Delta},  \exists v\in\cV(\Map^n) \text{\,\,such\,\,that\,\,} \Cdy^n(v)\in B(\ol z,\zeta) ]\geq 1-\zeta.
\label{eq20}
\eqe
Therefore 
\eqb
\sup_{x\in\ol\Delta}|\Cdy^n(\frk v(x))-x|\rta 0 \qquad \textrm{in probability as }n\to \infty.
\label{eq18}
\eqe
On the other hand, we have
\eqb
\begin{split}
	\sup_{x,y\in \ol \Delta} |d^n_\Delta(x,y)-d_{\Delta}(x,y)|
	\leq\,
	& \sup_{x,y\in \ol \Delta}  |d^n(\frk v(x),\frk v(y))-
	d_{\Delta}(\Cdy^n (\frk v(x)),\Cdy^n (\frk v(y)))|\\
	&+  \sup_{x,y\in \ol \Delta} |d_{\Delta}(\Cdy^n (\frk v(x)),\Cdy^n (\frk v(y)))-d_{\Delta}(x,y)|,
\end{split}
\label{eq14}
\eqe
\xncomment{Since $d_\Delta$ induces the Euclidean topology,} the second term on the right side  of \eqref{eq14}  converges to 0 by~\eqref{eq18}.
Therefore to get the convergence of $d^n_\Delta$ it suffices to show that
\eqb
\begin{split}
	\lim_{n\to\infty}	\sup_{v',v''\in \cV(\dsk^n)} | d^n (v',v'')-d_\Delta (\Cdy^n (v')&,\Cdy^n(v'' )) |=0.
\end{split}
\label{eq11}
\eqe

For any $\zeta>0$, by Propositions~\ref{prop:max}  and~\ref{prop:max2},  
we can choose $\rho>0$ (depending only on $\zeta$) sufficiently small, such that for all sufficiently large $n$ (depending on $\zeta$), the following holds with probability at least $1-\zeta$,
\eqb
\sup_{v,u\in \cV(\dsk^n)\,:\,d^n(u,v)< \rho } |\Cdy^n(u)-\Cdy^n(v)| < \zeta.
\label{emb-eq:cdy}
\eqe
In the coupling in Lemma \ref{emb-prop:GHPU-marked},  $\lim_{n\to  \infty} d^n (v^n_i,v^n_j)=d_{\Delta}(z_i,z_j)$ a.s.\ for  each $i,j\in\N$.
Since $d_\Delta$ is continuous relative to the Euclidean metric, an application of the triangle inequality and \eqref{eq21} gives
\begin{equation}\label{eq:metric-typical}
	\lim_{n\to\infty} | d^n (v^n_i,v^n_j)-d_\Delta (\Cdy^n (v^n_i),\Cdy^n(v^n_j )) |=0	\quad \textrm{in probability}.
\end{equation}	
Combing \eqref{emb-eq:cdy} and \eqref{eq:metric-typical} and using  the density of $\{z_i:i \in \N \}$ in $\ol\Delta$ for $d_\Delta$, we get~\eqref{eq11}.

\begin{proof}[Proof of Theorem \ref{thm:crossings}]
	Recall the proof of the convergence of $\xi_\Delta^n$. The argument there implies that $\Cdy^n(\wh v^n_4)$ converge to $\Cdy^n(\wh z_4)$ in probability.
	Now conditioning on the event that $\wh v^4_n$ falls on the arc $(c^n,a^n)$ and on the event that  $\wh z_4$ falls into the counterclockwise arc on $\p \Delta$ from $(0,0,1)$ to $(1,0,0)$, we obtain Theorem~\ref{thm:crossings}. 
\end{proof}

\section{The quantum pivotal measure of CLE$_6$}\label{emb-sec:pivot} 
We recall the setting of \eqref{eq:KPZ}. Namely,  let  $\bh$ be as in Definition~\ref{def:bh}
so that $(\D,\bh,1)/\mathord\sim$ is a unit boundary length $\sqrt{8/3}$-LQG disk (Definition~\ref{def:disk}). 
Let $\Gamma$ be a $\CLE_6$ on $\D$ with monochromatic blue boundary condition (Definition~\ref{def:CLE})  which is  independent of $\bh$.
Fix $\eps>0$.  Let  $\cP_\eps$ be the set of $\eps$-pivotal points of $(\bh,\Gamma)$ as in Definition~\ref{def:eps-piv}. 
The measure $\nu^\eps_{\bh,\Gamma}$ on $\cP_\eps$  was introduced in \cite[Section~7]{bhs-site-perc}  based on the theory of mating of trees \cite{wedges},  and we will review its definition in~Section~\ref{subsec:pivot-measure}.  Let $\pivm^\eps$ be the renormalized scaling limit of  $e^{\bh/\sqrt6}d^2z$ restricted to the discrete analog of $\cP_\eps$.
We have described the discrete setting above \eqref{eq:KPZ} and will describe $\pivm^\eps$ precisely in Definition~\ref{def:pivm}. 
We now restate~\eqref{eq:KPZ} as a proposition.
\begin{proposition}\label{prop:quantum-Mink}
	In the setting right above,  there exists a deterministic constant $\constp>0$ such that
	for each fixed $\eps>0$, we have $\nu^\eps_{\bh,\Gamma}=\constp\pivm^\eps$ a.s.
\end{proposition}

We will recall the mating-of-trees theory for $\SLE_6$ on $\sqrt{8/3}$-LQG surfaces in Section~\ref{subsec:natural-time}.
In Section~\ref{subsec:pivot-measure}, we give a definition of $\nu^\eps_{\bh,\Gamma}$ which is a slight reformulation of the one in \cite{bhs-site-perc}. 
The bulk of this section, Section~\ref{subsec:Mink}, is devoted to the proof of a local version of Proposition~\ref{prop:quantum-Mink}, namely     Proposition~\ref{prop:pivot-bichordal}.
As we will show in Lemma~\ref{lem:cover}, the set $\cP_\eps$ can be covered by the points of intersection of the so-called 2-$\SLE_6$ as defined below. 
\begin{definition}\label{def:bi-chordal}
	Suppose $Q$ is a simply connected domain with simple piecewise smooth boundary and $a,b,c,d$ are four distinct boundary points ordered counterclockwise.  
	Let $\sle^{ad}_Q$ be a chordal $\SLE_6$ on $(Q,a,d)$ conditioned on not hitting the counterclockwise boundary arc $\p_{b,c} Q$ from $b$ to $c$. 
	Conditioned on $\eta^{ad}$, let $Q'$ be  the component of $Q\setminus\sle_Q^{ad}$ whose boundary contains $\p_{b,c}Q$, and let $\sle_Q^{cb}$ 
	 be a chordal $\SLE_6$ on $(Q',c,b)$. 
	 We call $(\sle_Q^{ad},\sle_Q^{cb})$ a 2-$\SLE_6$  on $(Q, a, b,c,d)$.
\end{definition} 
Proposition~\ref{prop:pivot-bichordal} is the variant of  Proposition~\ref{prop:quantum-Mink} with $\sle_Q^{ad}\cap\sle_Q^{cb}$ in place of $\cP_\eps$.
Combined with  the covering lemma (i.e.\ Lemma~\ref{lem:cover}), this will give Proposition~\ref{prop:quantum-Mink}. 
We will explain this part  in Section~\ref{sec:prop:quantum-Mink}.
The reader may skip Sections~\ref{subsec:natural-time} to \ref{subsec:Mink} and proceed directly to Section~\ref{sec:GPS}  
if she is willing to accept Proposition~\ref{prop:quantum-Mink} without a proof.
\subsection{Mating-of-trees theory for SLE$_6$ on $\sqrt{8/3}$-LQG surfaces}\label{subsec:natural-time}

The definition of $\nu^\eps_{\bh,\Gamma}$ and the proof of Proposition~\ref{prop:quantum-Mink} both rely on the mating-of-trees theory for $\SLE_6$ on $\sqrt{8/3}$-LQG surfaces. The general theory is built in the foundational paper \cite{wedges}. It is further developed in \cite{gwynne-miller-sle6} and revisited in~\cite[Section~7]{bhs-site-perc}.
In this subsection we review what is needed for Proposition~\ref{prop:quantum-Mink}. See \cite{ghs-mot-survey} for a thorough survey.

\subsubsection{Quantum wedges and disks}\label{subsub:surface}
We start by recalling the definition of a family of LQG surfaces which plays a key role in the mating-of-trees theory, namely the quantum wedges  \cite{shef-zipper,wedges}. Recall the notation $\cS=\R\times (0,\pi)$ for the horizontal strip.
\begin{definition}[Quantum wedge]\label{def:thick-wedge}
	Fix  $W>4/3$ and $a>0$ such that  $W=4/3+\sqrt{8/3}a$ \cite[Table~1.1]{wedges}.
	Let $(X_t)_{t\in \R}$ be such that
	\begin{compactitem}
		\item $(X_t)_{t\ge 0}\eqD (B_{2t}-at)_{t\ge 0}$, where $B_t$ is a standard linear Brownian motion starting at 0,
		\item $(X_{-t})_{t\ge  0}$ has the law of $(B_{2t}+at)_{t\ge  0}$ conditioned to be positive, and
		\item $(X_{-t})_{t\ge 0}$ and $(X_t)_{t\ge 0}$ are independent.
	\end{compactitem}
	Let  $h^1(t+si)=X_t$ for each $t+si\in \cS$.  
	Let $h^2$ be the random distribution on $\cS$ independent of $X$ whose law is the lateral component of the free-boundary GFF on $\cS$.
	Set $h=h^1+h^2$.
	Then the  law of the $\sqrt{8/3}$-LQG surface $(\cS,h,+\infty,-\infty)/\mathord\sim$ 
	is called the \notion{$W$-quantum wedge}.\footnote{In \cite{wedges} quantum wedges are parametrized in six different ways. See \cite[Table~1.1]{wedges} for their relations. Our choice in Definition~\ref{def:thick-wedge} is called parametrization by \notion{weight}. The notion of $\alpha$-quantum wedge in \cite{wedges} is different from the one in Definition~\ref{def:thick-wedge} since $\alpha$ refers to the log singularity parameter, while our $W$ refers to the weight. These are related by $W=\gamma(\gamma/2+{\mathbf Q}-\alpha)$, where $\gamma=\sqrt{8/3}$ and ${\mathbf Q}=5/\sqrt6$.}
	
	If in the above definition, the law of $X$ is such that $(X_t)_{t\ge 0}\eqD (B_{2t})_{t\ge 0}$ conditioned to be negative, and $(X_{-t})_{t\ge  0}$ has the law of $(B_{2t})_{t\ge  0}$, then the law of the $\sqrt{8/3}$-LQG surface $(\cS,h,+\infty,-\infty)/\mathord\sim$ 	is called the \notion{$4/3$-quantum wedge}.
\end{definition}

\begin{remark}\label{rmk:wedge}
	By \cite[Proposition 4.7]{wedges}, quantum wedges have the following symmetry.
	If $(D,\hw,a,b)/\mathord\sim$ is a $W$-quantum wedge,
	then $(D,\hw+c,a,b)/\mathord\sim\, \overset{d}{=} \, (D,\hw,a,b)/\mathord\sim $
	for each deterministic  $c\in\R$.  
	
	By \cite[Proposition 1.7]{shef-zipper}, the $2$-quantum wedge has an additional symmetry.  If 
	$(D,\hw,a,b)/\mathord\sim$ is a  $2$-quantum wedge and $s>0$, let  $a_s\in D$ be
	on the left boundary of $(D,a,b)$ such that the $\xi_{\hw}$-length of the left boundary  of $(D,a,a_s)$
	equals $s$. 
	Then $(D,\hw,a_s,b)/\mathord\sim$ has the law of a $2$-quantum wedge. 
\end{remark}

The following representative of a quantum wedge (i.e., a representative of the equivalence class defining the wedge) will be technically convenient in several of our arguments.
\begin{definition}\label{def:circle}
	Let $\cW$ be a $W$-quantum wedge for  $W\ge 4/3$ and
	let $\phi(z)\defeq e^{-z}$ be a conformal map from $\cS$ to $\bbH$.
	Suppose $\hw$ is the random distribution on $\bbH$ such that $\cW=(\bbH,\hw,0,\infty)/\mathord \sim$ and, moreover,   $\hw\circ\phi  + {\mathbf Q}\log |\phi'|$ has the law of  $h$ in Definition \ref{def:thick-wedge}.
	Then we say that $(\bbH,\hw,0,\infty)$ is the \notion{circle average embedding} of $\cW$.
\end{definition}
Existence and uniqueness of the circle average embedding is clear from Definition~\ref{def:thick-wedge}.

In order to state the mating-of-trees theorem, we  need to extend our definition of the $\sqrt{8/3}$-LQG disk to allow arbitrary boundary length.
\begin{definition}\label{def:disk2}
	Recall the notions in Section~\ref{subsec:lqg}. Suppose $D$ is a simply connected $C^0$ domain, $a$ is a point on $\p D$, and $h$ is a free Liouville field on $D$. Define $L\defeq \xi_h(\p D)$.  Recall from Convention~\ref{conv:DF} that $\gamma=\sqrt{8/3}$. If $(D, h-2\gamma^{-1}\log L,a)/\mathord\sim$ is independent of $L$ and has the law of a $\sqrt{8/3}$-LQG disk with unit boundary length (see Definition~\ref{def:disk}),
	then we say that $(D,h,a)/\mathord\sim$ is a \notion{$\sqrt{8/3}$-LQG disk} and call $L$ the \notion{boundary length} of the disk.
\end{definition}

\subsubsection{Mating-of-trees theory for SLE$_6$ on a 2-quantum wedge}
Recall notions in Section~\ref{subsec:SLE}. 
Given $(D,a,b)\in\CD$, let $\eta$ be an  $\SLE_6$ on $(D,a,b)$.
Let $\hw$ be a random distribution on $D$ such that $\cW:=(D,\hw,a,b)/\mathord\sim$ is a 2-quantum wedge. 
A set $\cB\subset D$ is called a \notion{bubble} of $\eta$ if it is a connected component of $D\setminus \eta$.  Let $t_\cB=\sup\{t\ge 0:  \cB\subset D_t \}$. We call $x_\cB:=\eta(t_\cB)$ the \notion{root} of $\cB$.   
 By \cite[Theorem~1.18 and Corollary~1.19]{wedges}, we have the following parametrization of $\eta$.
\begin{proposition}
\label{def:natural-time}
Let $X$  be a  L\'evy processes with L\'evy measure $\frac{3}{4\sqrt{\pi}}|x|^{-5/2} \1_{x<0}\,dx$ and let $J(X)=\{t\geq 0\,:\, X_t-X_{t^-}\neq 0 \}$ be the set of times at which $X$ makes a jump. 
	Conditioned on $X$, for each $t\in J(X)$ sample an independent $\sqrt{8/3}$-LQG disk $\cS_t$ with boundary length equal to the jump size $X_t-X_{t^-}$ of $X$ at time $t$, and set $\cE=\{ (t,\cS_t)\,:\,t\in J(X) \}$.
	%Conditioning on $X$, sample a collection $\cE=(\cS_t\,:\,t\in J(X) )$ of independent $\sqrt{8/3}$-LQG disks $\cS_t$ indexed by $J(X)$ such that the boundary length of $\cS_t$ is equal to the jump size $X_t-X_{t^-}$ of $X$ at time $t$.
In the setting of the previous paragraph, there exists a unique parametrization of $\eta$ such that  the following holds.
Let $J^{\op{L}}\subset\R_+$ be the set of times $t$ at which a bubble $\cB_t$ is cut off by $\eta$ on its left side, let $\cS^{\op{L}}_t=(\cB_t,h|_\cB,x_\cB)\mathord/\mathord\sim$ %which are cut off at these times 
denote the LQG surface which is cut off at time $t$, and set 
$\cE^{\op L}=\{ (t,\cS^{\op{L}}_t)\,:\, t\in J^{\op{L}} \}$. %independent  %Let $\cE^{\op L}=( \cS^{\op{L}}_t\,:\, t\in J^{\op{L}} )$ denote the set of LQG surfaces $\cS^{\op{L}}_t=(\cB_t,h|_\cB,x_\cB)\mathord/\mathord\sim$ which are cut off at these times, indexed by $J^{\op{L}}$.
Define $\cE^{\op R}$ in the same way with left replaced by right. Then $\cE^{\op L}$ and $\cE^{\op R}$ %\xncomment{$(J^{\op{L}},\cE^{\op L})$ and $(J^{\op{R}},\cE^{\op R})$} 
are independent, and each of them have the same law as $\cE^{\op R}$. %\xncomment{$(J(X),\cE)$}. 
We call this parametrization the \notion{quantum natural parametrization of $\eta$ under $\hw$}.\footnote{In fact, the quantum natural parametrization in \cite{wedges} is defined only up to a multiplicative constant, which we fix in this paper by specifying the L\'evy measure of $X$.}
\end{proposition}
\begin{proposition}\label{prop:mating0} 
	Let $(D,\hw,a,b,\eta)$ be as in Proposition~\ref{def:natural-time}, with $\eta$ having the quantum natural parametrization.
	For a fixed  $t>0$, conditioning on $Z^{\op w}|_{[0,t]}$, 
	the conditional law of $\{(\cB,h, x_\cB)/\mathord\sim:  \textrm{ $\cB$ is a bubble with }t_\cB\le t \}$ is that of independent $\sqrt{8/3}$-LQG disks with given boundary length, which are also conditionally independent of $(D_t, h, \eta(t),b,\eta)/\mathord\sim$. Furthermore, the conditional law of $(D_t, h, \eta(t),b,\eta)/\mathord\sim$ equals the law of $(D, h, a,b,\eta)/\mathord\sim$, where $D_t$ is the connected component of $D \setminus \eta([0,t])$ whose boundary contains $b$.
\end{proposition}
 
By the quantum zipper theory of Sheffield~\cite{shef-zipper}, given a variant of $\SLE_{8/3}$ coupled with an independent free Liouville field on the same domain, one can unambiguously define a notion of \notion{quantum length} measure on the $\SLE_{8/3}$-type curve, as an extension of the $\sqrt{8/3}$-LQG boundary measure.
For example,  in Proposition~\ref{def:natural-time}, let $U$ be either $D_t$ or a bubble of $\eta$. Given a segment $V$ of $\p U$,  since $h|_{U}$ is either a quantum wedge or a $\sqrt{8/3}$-LQG disk, 
the mass of $V$ under the $\sqrt{8/3}$-LQG boundary measure of $h|_U$ is well defined, which we call the \emph{quantum length} of $V$. (Recall by SLE duality that $\p U$ is either a variant of $\SLE_{8/3}$ or part of $\p D$.) In the rest of Section~\ref{emb-sec:pivot} there are a few other occasions where we consider the quantum length along  $\SLE_{8/3}$-type curves. At each place, locally the SLE$_{8/3}$ curve cuts the domain into two subdomains with the curve lying on their border. The field restricted to the two subdomains are both free Liouville fields, each  of which induces a notion of quantum length for the curve using the $\sqrt{8/3}$-LQG boundary measure.
The highly nontrivial fact established in~\cite{shef-zipper} is that the two notions agree. See Proposition~\ref{prop:zipper83} for such an instance.

\begin{figure}
	\centering
	\includegraphics{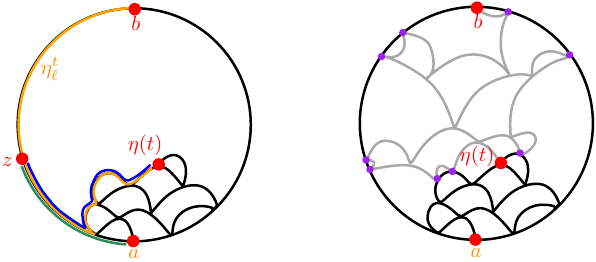}
	\caption{{\bf Left}: Illustration of the boundary length process $Z^{\op w}$ as defined in Proposition \ref{prop:Z}. We have that $Z^{\op w}_t$ is equal to the quantum length of the blue curve minus the quantum length of the green curve. {\bf Right}: The boundary touching measure of $\eta$ at time $t$ (Definition \ref{def:EBT-wedge}) is supported on the purple points.}
	\label{fig:Z}
\end{figure}

The key observable in the mating-of-trees theory is the so-called boundary length process. 
The next proposition   follows from  \cite[Corollary~1.19]{wedges}. See Figure~\ref{fig:Z} for an illustration.
\begin{proposition}\label{prop:Z}
Suppose we are in the setting of Proposition~\ref{def:natural-time}.  Set $L^{\op w}_0=R^{\op w}_0=0$.  For $t>0$, let $\eta^t_{\markl}$ and $\eta^t_{\op r}$ be the left and right, respectively, boundary of $(D_t,\eta(t),b)$.
Let $z$ be a point on $\eta^t_{\markl}\cap \bdy D$. 
Let $L^{\op w}_t$ be the quantum length  of
 the clockwise arc from $\eta(t)$ to $z$ on $\bdy D_t$
minus the quantum length of the clockwise arc from $0$ to $z$ on $\bdy D$. (It is clear that the value of $L^{\op w}_t$ does not depend on the choice of $z$.) Define $R^{\op w}_t$ similarly with $z$ on $\eta^t_{\op r}\cap \bdy D$ instead and with counterclockwise instead of clockwise.
Then $L^{\op w}$ and $R^{\op w}$ are independent and have the same distribution as $X$ in Proposition~\ref{def:natural-time}.
We call $Z^{\op w}=(L^{\op w},R^{\op w})$ the \notion{boundary length process} of $(D,\hw,a,b,\eta)$. 
\end{proposition}

The process  $L^{\op w}$ (resp., $R^{\op w}$) has a downward jump at time $t$ if and only if $t=t_\cB$ for some bubble $\cB$ to the left (resp., right) of $\eta$. Moreover, the size of the jump equals the quantum length of $\bdy \cB$.  
By \eqref{eq:conf-inv}, $Z^{\op w}$ is a.s.\ determined by the $\sqrt{8/3}$-LQG surface $(D,\hw,a,b,\eta)/\mathord\sim$.

\subsubsection{Mating-of-trees theory for SLE$_6$ on $\sqrt{8/3}$-quantum disks}

We now introduce  the quantum natural parametrization for $\SLE_6$ on a $\sqrt{8/3}$-LQG disk following \cite{gwynne-miller-sle6}.
Given constants  $\ell, r>0$, let  $(D,a,b)\in \CD$  and let $h$ be a random distribution on $D$ such that 
$(D,h,a)/\mathord\sim$ is a $\sqrt{8/3}$-LQG disk with boundary length $\ell+r$ and the right boundary length of $(D,a,b)$ equals $r$.  Let $\eta$ be a chordal $\SLE_6$ on $(D,a,b)$ independent of $h$. 
We can define the \emph{boundary length process} $Z^{\op d}=(L^{\op d},R^{\op d})$ of $(D,h,a,b,\eta)$  in the same way as $Z^{\op w}=(L^{\op w},R^{\op w})$ in Proposition~\ref{prop:Z}.
 It is easy to see that $L_t+\ell$ and $R_t+r$ measure  the quantum length of the left and right, respectively, boundary  $(D_t,\eta(t),b)$.
\begin{proposition}[\cite{gwynne-miller-sle6}] \label{def:natural-time2}
	In the setting right above,  there exists a unique parametrization of $\eta$, defined on some random interval  $[0,\sigma]$, such that the law of $Z^{\op d}=(L^{\op d},R^{\op d})$ can be characterized as follows.	Let $Z^{\op w}=(L^{\op w}, R^{\op w})$ be as in Proposition~\ref{prop:Z} and let $\sigma^{\op w}=\inf\{t\ge 0:  L^{\op w}(t)\le -\ell \; \textrm{or} \; R^{\op w}(t)\le -r \}$.  Then for each fixed $t>0$, the law of $Z^{\op d}|_{[0,t]} \cdot \1_{t<\sigma}$ is absolutely continuous with respect to $Z^{\op w}|_{[0,t]}  \cdot \1_{t<\sigma^{\op w}}$ with Radon-Nikodym derivative proportional to $(L^{\op w}(t) +R^{\op w}(t)+\ell+r )^{-5/2} \1_{t<\sigma^{\op w}}$. Moreover, $\lim\limits_{t\to \sigma}Z^{\op d}(t)=(-\ell,-r)$ almost surely.
	
	We call this parametrization the  \notion{quantum natural parametrization of $\eta$ under $h$}.
\end{proposition}
Intuitively, the law of $Z^{\op d}$ is the conditional law of $Z^{\op w}$ until exiting $(-\ell,\infty)\times (-r,\infty)$,
conditioning on exiting  at $(-\ell,-r)$.

The following proposition is the disk variant of Proposition~\ref{prop:mating0}.
\begin{proposition}[\cite{gwynne-miller-sle6}]\label{prop:mating}
	Let $(D,h,a,b,\eta)$ be as in Proposition~\ref{def:natural-time2}, with $\eta$ having the quantum natural parametrization.
	For a fixed  $t>0$, conditioning on $Z^{\op d}|_{[0,t]}$ and the event $D_t\neq \emptyset$, 
		the conditional law of $\{(\cB,h, x_\cB)/\mathord\sim:  \textrm{ $\cB$ is a bubble with }t_\cB\le t \}$ is that of independent $\sqrt{8/3}$-LQG disks with given boundary length, which are also conditionally independent of $(D_t, h, \eta(t),b,\eta)/\mathord\sim$. Furthermore, the conditional law of $(D_t, h, \eta(t),b,\eta)/\mathord\sim$ equals the law of $(D, h, a,b,\eta)/\mathord\sim$ with $(\ell, r)$ replaced by $(L_t+\ell,R_t+r)$.  
\end{proposition}

\subsection{$\sqrt{8/3}$-LQG pivotal measure as a local time}\label{subsec:pivot-measure}
In this section we provide a construction of the $\eps$-pivotal measure using the mating-of-trees theory we reviewed in Section~\ref{subsec:natural-time}. Our construction differs from the one in \cite[Section~7]{bhs-site-perc} since we rely heavily on the iterative construction of $\CLE_6$ (Lemma~\ref{lem:iteration}). However, as explained in Remark~\ref{rmk:equiv}, 
the two constructions produce  the same pivotal measure up to a multiplicative constant.

We will  rely on a natural way of constructing measures supported  on fractals.
\begin{definition}[Occupation measure]\label{def:Mink}
	Fix a positive integer $n$ and a compact set $A\subset \R^n$.  For $r>0$, let $A_r=\{z\in\C: |z-x|\le r\textrm{ for some } x\in A\}$. For $d\in (0,n]$, let $\Mink^r_{A,d}$ be the measure given by $r^{d-n}$ times  Lebesgue measure restricted to $A_r$. If the limit $\Mink_A=\lim_{r\to 0}\Mink^r_{A,d}$ exists  for the weak topology on the set of Borel measures and has finite and positive 
	total mass, we call $\Mink_A$ the \notion{$d$-occupation measure} of $A$.
\end{definition}
It is clear that there is at most one $d$ such that the $d$-occupation measure of $A$  exists. 
If $\Mink_A$ exists, then $\frk m_A(\R^n)$ is the so-called $d$-dimensional \emph{Minkowski content} of $A$.

We now recall some standard facts from fluctuation theory for L\'evy processes and stable subordinators which can be found in \cite{ladder,bertoin-sub}.  
For each $\beta\in (0,1)$, a L\'evy process $(\tau_t)_{t\geq 0}$ is called a $\beta$-stable subordinator if $\tau$ is a.s.\ increasing and $\tau_{at}\overset{d}{=}a^{1/\beta}\tau_t$ for each $a>0$. 
The closure $\cR_\tau$ of ${\{\tau_t\,:\, t\ge 0 \}}$ is called the \emph{range} of $\tau$. 
Let $m_\tau$ be the pushforward of Lebesgue measure on $[0,\infty)$ by $\tau$, so that $m_\tau$ is a measure supported on $\cR_\tau$. 
We call $m_\tau$ the \notion{local time} on $\cR_\tau$. We will  rely crucially on the occupation measure interpretation  of local time.
\begin{lemma}\label{lem:Mink-zero}
For a $\beta$-stable subordinator $(\tau_t)_{t\geq 0}$, there exists a deterministic constant $c_\beta>0$ such that almost surely the $\beta$-occupation measure  $\Mink_{\cR_\tau}$ of $R_\tau$ is well-defined, and
\begin{equation}\label{eq:local}
m_\tau([0,t])=c_\beta \Mink_{\cR_\tau}([0,t]) \qquad \textrm{for all }t>0.
\end{equation}
\end{lemma}
\begin{proof}
This follows by combining  e.g.\ \cite[Proposition 10]{py97} and \cite[Theorem 2.2]{lp93}, as explained in \cite[Section 13.4.2]{lf13}.
\end{proof}

\begin{lemma}\label{lem:Levy-sub}
	Let $X $ be as in Proposition~\ref{def:natural-time}.
	Then there exists a  $1/3$-subordinator $\tau$ such that $\cR_\tau =\{t\ge 0: X(t)=\inf_{s\in[0,t]} X(s) \}$. Moreover,  let $H_s=-X(\tau_s)$ for $s>0$. Then  $H$ is a $1/2$-stable subordinator and  almost surely $m_H$ equals the pushforward of $m_\tau$ under $-X$.
\end{lemma}
\begin{proof}
	The existence of $\tau$ and the law of $H$ can be found in \cite[Section~6]{ladder}, where $(X,H)$ is called the \emph{ladder process}.
	The fact that $m_H=(-X)_* m_\tau$ a.s.\ follows by definition.
\end{proof}
The following definition is the starting point of  the construction of $\nu^\eps_{\bh,\Gamma}$.
\begin{definition}\label{def:EBT-wedge}
	Let $(D,a,b,\hw,\eta)$ and $Z^{\op w}=(L^{\op w},R^{\op w})$ be as in Propositions~\ref{def:natural-time} and~\ref{prop:Z}, where $\eta$ has the quantum natural parametrization and $Z^{\op w}$ is the boundary length process.
	Let  $m_{\markl}$ and $m_{\op r}$ be defined as $m_\tau$ in Lemma~\ref{lem:Levy-sub} 
	with $L^{\op w}$ and $R^{\op w}$, respectively, in place of $X$, so that $m_{\markl}$ (resp., $m_{\op r}$) is a measure
	supported on  the set of times at which $L^{\op w}$ (resp., $R^{\op w}$) reach a running infimum. Let $\nu^0_\eta\defeq\eta_*m_{\markl}+\eta_*m_{\op r}$,
	which by the definition of $Z^{\op w}$ is a measure supported on $\eta\cap \bdy D$.
	For each $t>0$,  let  $\nu^t_\eta$ be defined as $\nu^0_\eta$ with $D$, $\hw$, $a$, and $\eta$ replaced by $D_t$, $\hw|_{D_t}$, $\eta(t)$, and $\eta|_{[t,\infty)}$, respectively.
	We call $\nu^t_\eta$  the \notion{boundary touching measure}  of $\eta$ at time $t$.
\end{definition}

For each $t\ge 0$, the measure $\nu^t_{\eta}$ is supported on  $\eta([t,\infty))\cap \bdy D_t$.
We now show that $\nu^t_\eta$ is determined by the set $\eta[t,\infty)\cap \bdy D_t$  and the quantum length measure on $\partial D_t$.
\begin{lemma}\label{lem:bdy}
	Let $(D,\hw,a,b,\eta)$ be as in Proposition~\ref{def:natural-time}. Let $c_{1/2}$ be as in \eqref{eq:local} with $\beta=1/2$.
	For a fixed $t\ge 0$, let $\eta^t_{\markl}$ and $\eta^t_{\op r}$ be the left and right, respectively, boundary of $(D_t,\eta(t),b)$, parametrized by quantum length
	starting from $\eta^t_{\markl}(0)=\eta^t_{\op r}(0)=\eta(t)$.
	Then the $\frac12$-occupation measure of \(\{s\ge 0: \eta^t_{\markl}(s)\in \eta([t,\infty))\cap \p D_t\}  \) on $[0,\infty)$  a.s.\ exists, which we denote by $m^t_{\markl}$. 
	 We can define $m^t_{\op r}$ in the same way  with $\eta^t_{\markl}$ replaced by $\eta^t_{\op r}$. 
	Then 
\begin{equation}\label{eq:bdy}
	\nu^{t}_{\eta}= c_{1/2}(\eta^t_{\markl})_* m^t_{\markl} + c_{1/2}(\eta^t_{\op r})_* m^t_{\op r}\qquad \textrm{a.s.}
\end{equation}
\end{lemma}
\begin{proof}
We only prove the case when $t=0$ since the general case follows from the stationarity in Proposition~\ref{prop:mating0}. 
Since $\eta^0_{\markl}$ is parametrized by its quantum length,  we have $\eta(u)=\eta^0_{\markl} (-L^{\op w}(u))$ for each $u\in \{t\ge 0: L^{\op w}(t)=\inf_{s\in[0,t]} L^{\op w} (s) \}$.
By Lemmas~\ref{lem:Mink-zero} and~\ref{lem:Levy-sub}, the measures $m^0_{\markl}$ and $m^0_{\op r}$ are well defined. By \eqref{eq:local}  and Lemma~\ref{lem:Levy-sub}, $(-L^{\op w})_*m_{\markl}=c_{1/2}m^0_{\markl}$ a.s., hence $(\eta^0_{\markl})_*(-L^{\op w})_*m_{\markl}=c_{1/2}(\eta^0_{\markl})_*m^0_{\markl}$. Furthermore, restricted to the support of $m_{\markl}$, we have $\eta=\eta^0_{\markl} \circ (-L^{\op w})$, hence $\eta_*m_{\markl}=c_{1/2}(\eta^0_{\markl})_* m^0_{\markl} $ a.s. Similarly, we have $\eta_*m_{\op r}=c_{1/2}(\eta^0_{\op r})_* m^0_{\op r} $ a.s.  
Therefore $\nu^{0}_{\eta}=c_{1/2}(\eta^0_{\markl})_* m^0_{\markl} + c_{1/2}(\eta^0_{\op r})_* m^0_{\op r}$ a.s. 
This proves Lemma~\ref{lem:bdy} for $t=0$.
\end{proof}
By the relationship between $Z^{\op d}$ and $Z^{\op w}$, we can define the boundary touching measure for 
an SLE$_6$-decorated $\sqrt{8/3}$-LQG disk in the exact same way as in Lemma~\ref{lem:bdy} via~\eqref{eq:bdy}.
\begin{definition}\label{def:bdy-touch}
	Let $(D,h,a,b, \eta)$, $\sigma$, and $Z^{\op d}=(L^{\op d}, R^{\op d})$ be as in Proposition~\ref{def:natural-time2}, so that $\eta$ has the quantum natural parametrization.
	For each $t\ge 0$, on the event $\{\sigma>t \}$, let $\dbl_{\eta,t}\defeq\eta([t,\sigma])\cap \bdy D_t$.  Let $\nu^t_\eta$ be the measure supported on $\dbl_{\eta,t}$ 
	defined in the same way as in Lemma~\ref{lem:bdy} in terms of $\eta^t_{\markl}$, $\eta^t_{\op r}$, and $\eta$ via~\eqref{eq:bdy}. We call  $\nu^t_\eta$ the \notion{boundary touching measure} of $\eta$  at time $t$.
	The  countable collection of measures $\{\nu^t_{\eta}\}_{t\in [0,\sigma)\cap \Q}$ 
	extends to a measure $\nu_{\eta}$ on the union of their supports, which we call the \notion{extended boundary touching} (EBT) measure of $\eta$ for $(D,h)$. 
\end{definition}
Given $(D,a,b)\in\CD$, let $\eta$ be an $\SLE_6$ on $(D,a,b)$ and define
\begin{equation}\label{eq:dbl}
\dbl_\eta:=\{p\in D\,:\, \exists s\neq t\textrm{ such  that }\eta(s)=\eta(t)=p \}   \quad\textrm{and} \quad\dbl_{\eta,D}:=\dbl_\eta\cup (\eta\cap \bdy D).
\end{equation}
Then $\nu_\eta$ is supported on $\dbl_{\eta,D} $ by definition. We remark (although this fact will not be used) that the measure $\nu_{\eta}$ is a.s.\ not locally finite and a.s.\ assigns infinite measure to any open ball intersecting $\eta$. However, $\nu_{\eta}$ is a.s.\ $\sigma$-finite and is a.s.\ finite when restricted to the set $\cP_\eps$ of $\eps$-pivotals (see Definition~\ref{def:eps-piv} below).

\begin{figure}
	\centering
	\includegraphics[scale=1]{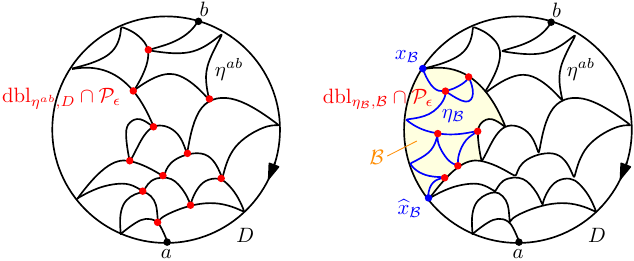}
	\caption {Illustration of the construction of the pivotal measure $\nu^\eps_{h,\Gamma}$ given in Definition \ref{def:eps-piv}. The left figure illustrates the construction for a monochromatic domain $D$ (Step 1), while the right figure considers the case of a dichromatic bubble $\cB$ (Step 2). The $\eps$-pivotal points which are captured in each step are shown in red. Note that points of intersection between an SLE$_6$ and $\partial D$ are not $\eps$-pivotal points, while in later iterations points of intersection between an SLE$_6$ interface and the boundary of some monochromatic bubble $\cB'$ could be $\eps$-pivotal points.} 
	\label{fig:pivmeasure-iteration}
\end{figure}

Now we are ready to define the measure $\mpivep$ for $(D,h,a)$, where $\nu^\ep_{\bh, \Gamma}$ in Proposition~\ref{prop:quantum-Mink} is the special case when $(D,h,a)=(\D,\bh,1)$.  See Figure~\ref{fig:pivmeasure-iteration} for an illustration.
\begin{definition}\label{def:eps-piv}
	Let $D$ be a Jordan domain and let $(D,h,a)$ be a $\sqrt{8/3}$-LQG disk with boundary length $L$. Let $\Gamma$ be a $\CLE_6$ on $D$ independent of $h$ with monochromatic blue boundary condition.
	Let $\cP_\eps$ be the set of $\eps$-pivotal points of $(h,\Gamma)$.
	The \notion{$\sqrt{8/3}$-LQG pivotal measure $\mpivep$}  on $\cP_\eps$ is the measure supported on $\cP_\eps$ which can be constructed as follows.
	\begin{enumerate}
		\item[Step 1] \label{item:dbl}
		Let $b\in\bdy D$ be such that 
		the left boundary of $(D,a,b)$ has quantum length $L/2$. Let $\Gamma_a^b$ and $\eta^{ab}$ be determined by $\Gamma$ as  in Lemma~\ref{lem:iteration}. 
		Set $\mpivep=\nu_{\eta^{ab}}$ on $\cP_\eps \cap \dbl_{\eta^{ab},D}$  where $\nu_{\eta^{ab}}$ is the EBT measure (see Definition~\ref{def:bdy-touch}) of $\eta^{ab}$ for $(D,h)$.
		\item[Step 2] \label{item:dich}
		Recalling notations in the paragraph above Lemma~\ref{lem:iteration}, for each dichromatic bubble $\cB$
		set $\mpivep=\nu_{\eta_\cB}$ on $\cP_\eps\cap \dbl_{\eta_\cB,\cB}$  where $\nu_{\eta_\cB}$ is the EBT measure of $\eta_\cB$ for $(\cB,h|_\cB)$. Here, although the domain $\cB$ itself is random, Definition~\ref{def:bdy-touch} trivially extends to $(\cB,h|_\cB,\eta_\cB)$.
	\end{enumerate}
	Given a connected component $\cB'$ of $D\setminus \Gamma_a^b$, which is a monochromatic bubble, let  $a'$ be the last point on $\p \cB'$ visited by $\eta^{ab}$ or one of the $\eta_\cB$'s 
	with $\cB$ being a dichromatic bubble.
	Namely, if $\p \cB'$ does not intersect  any of the $\eta_\cB$'s, then $a'$ is the last point on $\p \cB'$ visited by $\eta^{ab}$.  
	If $\p \cB'$ intersects an $\eta_\cB$, then $a'$ is  the last point visited by this $\eta_\cB$.
	We define the measure $\mpivep$ on $\cB'\cup\p \cB'$ by repeating Steps~1 and~2 on $(\cB',h|_{\cB'},a',\Gamma|_{\cB'})$ and then iterate.\footnote{Note that $\cP_\eps\cap \op \D=\emptyset$ but with positive probability $\cP_\eps\cap \p \cB'\neq \emptyset$, in which case $\mpivep(\p \cB')$ is non-trivial.}
\end{definition}

The fact that $\mpivep$  in Definition~\ref{def:eps-piv} is well-defined requires some justification. Let $\dbl_t$ be the support of $\nu^t_{\eta^{ab}}$. As explained in \cite[Lemma~7.9]{bhs-site-perc}, there exists a finite set $\cT$ such that $\cP_\eps\cap \dbl_{\eta^{ab},D}\subset \cup_{t\in \cT} \dbl_t$.
Therefore $\mpivep$ restricted to $\cP_\eps\cap \dbl_{\eta^{ab},D}$ is a finite Borel measure.
In Step 2, there are finitely many dichromatic bubbles with $\cP_\eps\cap \dbl_{\eta_\cB,\cB}\neq\emptyset$. On each such bubble, the same consideration shows that $\mpivep$ restricted to $\cP_\eps\cap \dbl_{\eta_\cB,\cB}$ is a finite Borel measure.  
Note that a component $U$ of $D\setminus \Gamma^a_b$ with $\mu_h$-mass smaller than $\eps$ has no intersection with $\cP_\eps$ since it would be either a point of intersection between two loops contained in $\ol U$ or a double point of a loop contained in $\ol U$. Using this and the local finiteness of $\Gamma$ in Lemma~\ref{lem:property}, we get that the iteration a.s.\ exhausts $\cP_\eps$ in finitely many steps.
%Since a component of $D\setminus \Gamma^a_b$ with $\mu_h$-mass smaller than $\eps$ have no intersection with $\cP_\eps$, by the local finiteness of $\Gamma$ in Lemma~\ref{lem:property} the iteration a.s.\ exhausts $\cP_\eps$ in finitely many steps.
By the no-triple-points property of $\Gamma$ in Lemma~\ref{lem:property}, the subsets of $\cP_\eps$  on which we define $\mpivep$ in different iterative steps are all disjoint.
In particular, our definition of $\mpivep$ has no inconsistency in different steps. Moreover,  $\mpivep$ is almost surely a finite Borel measure on $D$.

\begin{remark}[Equivalent definitions of quantum pivotal measure]
	\label{rmk:equiv}
	We now explain the equivalence between  $\nu^\eps_{h,\Gamma}$ in Definition~\ref{def:eps-piv} and
	the $\eps$-LQG pivotal measure defined in \cite[Section 7]{bhs-site-perc}. 
	The latter measure is denoted by $\mathbf{\nu}_\eps$ in \cite{bhs-site-perc}, and we adopt the same notation here.
	We do not provide  the detailed construction  in \cite{bhs-site-perc},   but only point out how one can check the equivalence.
	If we do not employ Lemma~\ref{lem:bdy} but  only use the notations in Lemma~\ref{def:EBT-wedge} to describe Definitions~\ref{def:bdy-touch}~and~\ref{def:eps-piv}, then  restricted to $\cP_\eps\cap \dbl_{\eta^{ab},D}$ as in Step~1 in Definition~\ref{def:eps-piv}, our description of  $\mpivep$ is identical to that of $c_{\op p}^{-1}\nu_\eps$ in \cite[Section 7]{bhs-site-perc}, 
	 with $c_{\op p}$ in \eqref{eq:piv-conv}.
	This multiplicative constant is needed because  the normalization of local time in \cite{bhs-site-perc} is chosen such that $\nu^n_{\eps,\pv}\rta\nu_\eps$. 
	Recall $\eta_\cB$, $\eta^\cB$, and $\ell_\cB$ as defined in the paragraph above Lemma~\ref{lem:iteration}.  
	In the notation of \cite[Section~7.5]{bhs-site-perc},  $\eta_\cB$ and $\eta^\cB$ are the so-called past and future, respectively, segments of the loop $\ell_\cB$. This observation together with a further bookkeeping inspection of \cite[Section 7.7]{bhs-site-perc} implies that $\nu_{\eps}=c_{\op p}\mpivep$ on $\cP_\eps\cap \dbl_{\eta_\cB,\cB}$ as in Step~2 in Definition~\ref{def:eps-piv}. By iteration, one can check that $\nu_\eps=c_{\op p}\mpivep$. 
\end{remark}

\subsection{$\sqrt{8/3}$-LQG pivotal  measure as a quantum occupation measure}\label{subsec:Mink}
The main result of this section is Proposition \ref{prop:pivot-bichordal}, which is a  preliminary version of Proposition \ref{prop:quantum-Mink}.
In Sections \ref{subsub:zipper} and \ref{subsub:gmc-occ} we provide the necessary background and basic results on quantum zippers and GMC over occupations measures, respectively. This allows us to prove a first variant of Proposition \ref{prop:quantum-Mink} in Section~\ref{subsub:cut},  where the $\eps$-pivotal points are replaced by the points of intersection between two SLE$_{8/3}$-like curves (see Lemma \ref{lem:pre}). In Section \ref{subsub:wedge} we prove Proposition \ref{prop:pivot-bichordal} by linking to the setting of Section~\ref{subsub:cut}.
 
\subsubsection{SLE with force points, $2/3$-quantum wedges, and quantum zippers}\label{subsub:zipper}
We start by recalling a generalization of $\SLE_\kappa$ called $\SLE_\kappa(\rho_{\markl};\rho_{\op r})$, where $\SLE_\kappa$ is the special case $\SLE_\kappa(0;0)$.
Consider tuples of the form  $(D,a,b;v_{\markl},v_{\op r})$, where $(D,a,b)\in \CD$,  
and $v_{\markl}$ (resp., $v_{\op r}$) is a point on the left (resp., right) boundary of $(D,a,b)$. The points $v_{\markl}$ and $v_{\op r}$ are allowed to be equal to $a$, in which case we will denote them by $a^-$ and $a^+$. 
Given $\kappa>0,\rho_{\markl}>-2$, and $\rho_{\op r}>-2$, the (chordal) $\SLE_\kappa(\rho_{\markl};\rho_{\markl})$  on $(D,a,b; v_{\markl}, v_{\op r})$ 
is a probability measure  on  non-self-crossing curves on $D\cup \bdy D$ from $a$ to $b$ modulo increasing reparametrization. 
Away from $\bdy D$, an $\SLE_\kappa(\rho_{\markl}; \rho_{\op r} )$ curve looks locally  like $\SLE_\kappa$ in the sense that it has the same a.s.\ properties. The points $v_{\markl}$ and $v_{\op r}$ are called the \notion{force points}.
The parameter $\rho_{\markl}$ (resp., $\rho_{\op r}$) is called the \notion{weight} of $v_{\markl}$ (resp., $v_{\op r}$), and governs the behavior of the  curve when it approaches the left (resp., right) boundary.  
An $\SLE_\kappa(\rho_{\markl};\rho_{\op r})$ curve a.s.\ does not touch the left (resp., right) boundary of 
$(D,a,b)$ except for the ending points if and only if 
\begin{equation}\label{eq:touch}
\textrm{$\rho_{\markl}$ (resp., $\rho_{\op r}$) is at least $\kappa/2 -2$. }
\end{equation}
The $\SLE_\kappa(\rho_{\markl};\rho_{\op r})$ has conformal invariance and domain Markov properties similar to those in Section~\ref{subsec:SLE}, 
with the two additional marked points  taken into account when applying conformal maps.  See  \cite{ig1,ig2,wedges,lsw-restriction,dubedat-duality,zhan-duality1} for more background on $\SLE_\kappa(\rho_{\markl};\rho_{\op r})$.  In the rest of the paper the force points are always assumed to be located at $a^-$ and $a^+$ when we refer to $\SLE_\kappa(\rho_{\markl};\rho_{\op r})$   on $(D,a,b)$.

Let $\eta$ be an $\SLE_\kappa(\rho_{\markl};\rho_{\op r})$ on $(D,a,b)$ for $\kappa>4$. The  \emph{left (resp., right) boundary} of $\eta$ is the curve starting at $a$ and ending at $b$ which consists of the points on $\eta$ which are either on the left (resp.\ right) boundary of $(D,a,b)$ or can be connected to the left (resp., right) boundary of $(D,a,b)$ by a curve which does not intersect $\partial D$ or $\eta$, except possibly at the end-points. Here is a precise variant of the aforementioned $\SLE$ duality, see e.g.\  \cite{dubedat-duality,zhan-duality1,ig1}.
\begin{proposition}\label{prop:duality}
For $\rho_\ell,\rho_{\op r}>-1$, let $\eta$ be an $\SLE_6(\rho_\ell; \rho_{\op r})$ on $(D,a,b)$. 
Let $\eta_{\markl}$ and $\eta_{\op r}$ be its left and right boundary, respectively.
Then $\eta_{\markl}$ is an $\SLE_{8/3}(\frac 23 \rho_\ell-\frac 43; \frac 23 \rho_{\op r}-\frac 23)$ on $(D,a,b)$. 
If $\rho_{\op r}\ge 0$ so that $\eta_{\markl}$ does not touch the right boundary of $(D,a,b)$ by~\eqref{eq:touch},
conditioning on $\sle_{\markl}$,  the curve $\sle_{\op r}$ is  an $\SLE_{8/3}(-\frac 43;\frac 23 \rho_{\op r}-\frac 43)$ from $a$  to $b$ on the domain bounded between $\sle_{\markl}$ and the right boundary of $(D,a,b)$, 
and $\eta$ itself is an $\SLE_6(-1;\rho_{\op r})$ on the same domain.
\end{proposition}

A crucial fact in the quantum zipper theory is the \notion{conformal removability} of $\SLE_{8/3}(\rho_{\markl}; \rho_{\op r})$ (see e.g.\ \cite[Proposition 3.16]{wedges} and \cite[Theorem 1.4]{shef-zipper}, which build on \cite{jones-smirnov-removability,schramm-sle}).
\begin{lemma}\label{lem:conf-rem}
	Let $\eta$ be an $\SLE_{8/3}(\rho_{\markl};\rho_{\op r})$ on $(D,a,b)\in\CD$ with $\rho_{\markl},\rho_{\op r}>-2$.
	Suppose $U\subset D$ is open and that $\phi :U\to\C$
	is continuous on $U$ and conformal on $U\setminus \eta$. Then $\phi$ is a.s.\ conformal on $U$. 
\end{lemma}

We will use an important variant of the quantum wedge called the  $2/3$-quantum wedge, which is an ordered collection of $\sqrt{8/3}$-LQG surfaces with two marked boundary points.
\begin{definition}[$2/3$-quantum wedge]\label{def:thin}	
Let $\cE=\{(\ell, t)\}$ be a Poisson point process on $(0,\infty)^2$ with intensity measure $\ell^{-3/2}d\ell\otimes dt$. Conditioning on $\cE$, for each $(\ell,t)\in \cE$, sample an independent $\sqrt{8/3}$-LQG disk of length $\ell$, which we denote by  $(D_t, h_t, a_t)/\mathord\sim$.
Moreover, for each $(D_t, h_t,a_t)$, sample a point $b_t$ on $\p D_t$ according to the quantum boundary measure $\xi_{h_t}$. Then  $\{ (D_t,h_t,a_t, b_t)/\mathord\sim \}$ in the increasing order of $t$ is called a  \notion{$2/3$-quantum wedge}.
\end{definition}
\begin{comment}
	By \cite[Lemma 4.19]{wedges}, the p.p.p. should be given by the zeros of  linear Brownian motion.
\end{comment}

In \cite[Section 4.4]{wedges}, the $W$-quantum wedge with $W\in(0,4/3)$ is constructed in the spirit of Definition~\ref{def:thick-wedge}. Wedges with $W\in(0,4/3)$ are called thin wedges. Just as the $2/3$-wedge, they may be described as an ordered chain of finite-volume LQG surfaces. We do not need the $W\neq 2/3$ case in this paper, and therefore omit the construction.

Let $(D,a,b)\in \CD$ 
and let $\eta'_\ell$ and $\eta'_r$ be two simple curves on $D\cup \bdy D$ from $a$ to $b$ which do not cross each other, such that $\eta'_\ell$ is between $\eta'_r$ and the left boundary of $(D,a,b)$. Let $D'\subset D\cup\bdy D$ be the open set with boundary $\eta'_\ell\cup\eta'_r$. We call $D'$ \emph{the region bounded by $\eta'_\ell $ and $\eta'_r$}.
For  each bounded connected component $\cB$ of $D'$, let $a_\cB, b_\cB\in \bdy \cB$ be the two points on the intersection of the left and right boundary of $(D,a',b')$ such that $a_\cB$ is visited before $b_\cB$ by $\eta'_\ell$ and $\eta'_r$. 
Let $\{\cB\}$ be the collection of such components ordered such that $\{a_\cB\}$ is in order of visit by $\eta'_\ell$ and $\eta'_r$. Given a distribution $h$ on $D$, we let $(D',h,a,b)/\mathord\sim:= \{ (\cB,h,a_\cB, b_\cB)/\mathord\sim \} $ be the ordered collection of $\sqrt{8/3}$-LQG surfaces with two marked boundary points.

The main fact about the $2/3$-quantum wedge  which we will use is  the following   proposition from quantum zipper theory (see  \cite{shef-zipper} and \cite[Theorem 1.2]{wedges}).  
\begin{proposition}\label{prop:zipper83}
	Let $W^\ell,W^{\op r}\in\{2/3 \}\cup[4/3,\infty)$
	and let $(\BB H,\hw,0,\infty)$ be the circle average embedding of a $(W^\ell+W^{\op r})$-quantum wedge. (Recall Definitions~\ref{def:thick-wedge} and~\ref{def:circle}).
	Let $\eta'$ be an $\SLE_{8/3}(W^\ell-2;W^{\op r}-2)$ on $(\BB H,0,\infty)$. Let $D^\ell$ (resp.\ $D^{\op r}$) be the region bounded by $\eta'$ and the left (resp.\ right) boundary of $(D,a,b)$.\footnote{In the remainder of this section we will typically use a prime ($'$) when we refer to $\SLE_{8/3}$-type curves while we use no prime when we refer to $\SLE_6$-type curves.} Then the surfaces $(D^\ell,\hw,0,\infty)/\mathord\sim$ and $(D^{\op r},\hw,0,\infty)/\mathord\sim$ are independent and have the law of quantum wedges with weight $W^\ell$ and $W^{\op r}$, respectively. Furthermore, $(D^\ell,\hw,0,\infty)/\mathord\sim$ and $(D^{\op r},\hw,0,\infty)/\mathord\sim$ almost surely determine $\hw$ (and therefore also the surface $(\BB H,\hw,0,\infty)/\mathord\sim$). Finally, the $\sqrt{8/3}$-LQG boundary measure on $\eta'$  obtained by viewing  $\eta'$ as a boundary arc of $(D^\ell,\hw)/\mathord\sim$ or $(D^{\op r},\hw)/\mathord\sim$ agree.
\end{proposition}
In Proposition~\ref{prop:zipper83}, we say that the surface $(\BB H,h,0,\infty)/\mathord\sim$  is the \emph{conformal welding} of the surfaces $(D^\ell,h,0,\infty)/\mathord\sim$ and $(D^{\op r},h,0,\infty)/\mathord\sim$. 
Let $V$ be a segment of $\eta'$. We call the mass of $V$ under the $\sqrt{8/3}$-LQG boundary measure the \emph{quantum length} of $V$. By the last assertion of Proposition~\ref{prop:zipper83}, this is unambiguously defined. 

By Propositions~\ref{prop:duality} and~\ref{prop:zipper83}, we have the following.
\begin{proposition}\label{prop:zipper6}
	Let $W^\ell,W^{\op r}\in \{2/3 \}\cup[4/3,\infty)$ and let $(\BB H,\hw,0,\infty)$ be the circle average embedding of a $W^\ell+W^{\op r}+2/3$-quantum wedge. Let $\eta$ be an $\SLE_6(\frac 32 W^\ell-1; \frac 32 W^{\op r}-1)$ on $(\BB H,0,\infty)$ which is independent of $h$.
	Let $D^\ell$, $D^{\op r}$, and $D^{\op m}$  be the regions in $D$ bounded by the left boundary of $(D,a,b)$ and the left boundary of $\eta$, the right boundary of $(D,a,b)$ and the right boundary of $\eta$, and the left and right boundaries of $\eta$, respectively.
	Then $(D^\ell,\hw,0,\infty)/\mathord\sim$, $(D^{\op m},\hw,0,\infty)/\mathord\sim$, and $(D^{\op r},\hw,0,\infty)/\mathord\sim$ are independent $\sqrt{8/3}$-LQG surfaces and they have the law of wedges of weights $W^\ell$, $2/3$, and $W^{\op r}$, respectively. Furthermore, $(D^\ell,\hw,0,\infty)/\mathord\sim$, $(D^{\op m},\hw,0,\infty)/\mathord\sim$, and $(D^{\op r},\hw,0,\infty)/\mathord\sim$ almost surely determine $\hw$ (and therefore also the surface $(\BB H,\hw,0,\infty)/\mathord\sim$).
\end{proposition}
\begin{proof}
By Proposition~\ref{prop:duality} the left boundary of $\eta$ has the law of an $\SLE_{8/3}(W^\ell-2;W^{\op r}-\frac 43)$. An application of Proposition~\ref{prop:zipper83} implies that $(D^\ell,\hw,0,\infty)/\mathord\sim$ is a $W^\ell$-quantum wedge and is independent of $(D^{\op m,\op r}, \hw,0,\infty)/\mathord\sim$, where $D^{\op m,r}$ is the interior of the closure of $D^{\op m}\cup D^{\op r}$. We conclude the proof by a second application of Propositions~\ref{prop:duality} and~\ref{prop:zipper83}, this time using that 
conditioning on $D^{\op m,\op r}$, the curve $\eta$  is  an $\SLE_6(-1; \frac 32 W^{\op r}-1)$ on $(D^{\op m,\op r},0,\infty)$. 
\end{proof}

\subsubsection{Coordinate change for GMC over occupation measures} \label{subsub:gmc-occ}
A key fact we will use in the proof of \eqref{eq:KPZ} is that the two considered measures transform in the same way under conformal coordinate changes. In this section we collect some basic facts on conformal coordinate changes of a general class of random measures.
\begin{definition}\label{def:GMC}
	Let  $h$ be a free Liouville field (Definition~\ref{def:free}) on a domain $D$ and let $\mu$ be a random finite Borel measure on $D$. 
	For each $r>0$ and $z\in \C$, let $h_r(z)$ be the average of $h$ over the circle $\{w\in \C: |w-z|=r \}$, if this circle is contained in $D$.\footnote{The process $(z,r)\mapsto h_r(z)$ is well-defined as a continuous process on $\{ (z,r)\in D\times\,:\, |z-w|>r\,\, \forall w\in\C\setminus D \}$ (see e.g.\ \cite{shef-kpz}) and is known as the \emph{circle average process}.}
	Let $h_r(z)=0$ otherwise. 
	For $\alpha>0$,  we define the measure $e^{\alpha h}\mu$ by $\lim_{r\to 0}r^{\alpha^2/2} e^{\alpha h_r} \mu$ if the limit exists almost surely in the weak topology. (Recall the convention $f\mu$ in Section~\ref{subsec:notation}).
\end{definition}
In Definition~\ref{def:GMC}, when $h$ is a Gaussian field, the measure $e^{\alpha h}\mu$ is called the \emph{Gaussian multiplicative chaos} (GMC) over $\mu$ in the literature, except that the normalization $r^{\alpha^2}$ is sometimes replaced by $\E[e^{\alpha h_r(\cdot)}]^{-1}$.
We require $\lim_{r\to 0}r^{\alpha^2/2} e^{\alpha h_r}\mu$ to exist \emph{almost surely} as $r\rta 0$, rather than considering a limit in probability (or almost surely along dyadic numbers)
as in most other literature on GMC. This will be used in Lemma~\ref{lem:scaling}.

We are interested in the coordinate change for GMC over occupation measures  (see Definition~\ref{def:Mink}) of certain SLE related fractals. We first record a preliminary deterministic fact, whose proof is left to the reader.
	\begin{lemma}\label{lem:Mink-conf}
	Let $d\in(0,2)$ and let $A$ be a compact set on $\C$ whose 
	$d$-occupation  measure $\Mink_A$ exists.  
	Let $\phi$ be a conformal map on a domain containing $A$. Then the  $d$-occupation  measure  $\Mink_{\phi(A)}$ of $\phi(A)$ exists and equals
	$|(\phi^{-1})'|^{-d} \cdot (\phi_*\Mink_A)$.
	If furthermore
	\begin{equation}\label{eq:energy}
	\textrm{$\iint_{U\times U }\frac{d \Mink_A(x) d\Mink_A(y)}{|x-y|^{d-\eps}}<\infty$  for all bounded sets $U$ 
		and $\eps \in (0,d)$,}
	\end{equation}
	then \eqref{eq:energy} still holds with $\Mink_A$ replaced by $\Mink_{\phi(A)}$. 
\end{lemma}
\begin{comment}
	\begin{proof}
	We first assume $d\in (0,1)$. Now we shift the set $A$ by the random vector $(U,U)$, where $U$ is the a uniform random variable. Then by the same argument as in Lemma~\ref{lem:avoid-bdy}, the shifted $A$ will avoid the boundary of all dyadic polygons. Now we work with $A$ under this assumption. 
	Now we can use piecewise affine map to approximate $\phi$. Since the boundary does not cause problem, the lemma holds. This covers the case $d=3/4$ which is what we need later.
	
	For $d\in[1,2)$, we can proceed similarly. In this case, we cannot guarantee that $A$ avoids all dyadic lines. But we can makes sure using random shift  that almost surely no dyadic line has a nontrivial contribution. This is also sufficient to carry the approximation scheme using piecewise affine maps. 
	\end{proof}
\end{comment}
We also record a one-dimensional variant of Lemma~\ref{lem:Mink-conf}, which will be used in the proof of Proposition~\ref{prop:pivot-bichordal}. We again leave the elementary proof to the reader. 
\begin{lemma}\label{lem:1D-measure}
Let $d\in(0,1)$ and let $A$ be a compact set on $\R$ whose 
$d$-occupation measure $\Mink_A$ exists.  
Let $\phi$ be a $C^1$ map on an interval containing $A$ such that $\phi'>0$.
Then the  $d$-occupation  measure  $\Mink_{\phi(A)}$ of $\phi(A)$ exists and equals $|(\phi^{-1})'|^{-d} \cdot (\phi_*\Mink_A)$.
\end{lemma}
\begin{comment}
	\begin{proof}
		We do a random shift of $A$ as in Lemma~\ref{lem:Mink-conf} and proceed similarly. Since $A$ has Lebesgue measure zero, the shifted $A$ will avoid all rationals and the approximation argument goes through.
	\end{proof}
\end{comment}

The following lemma guarantees the existence of GMC over an occupation measure. The lemma would have followed from e.g.\ \cite{berestycki-gmt-elementary} if we had considered convergence in probability instead of a.s.\ convergence in Definition~\ref{def:GMC}. We include its proof in the appendix.
\begin{lemma}\label{lem:GMC}
	Fix $d\in(0,2)$, $\alpha\in (0,\sqrt{d})$, and a Jordan domain $D$. Let $A$ be a compact set on $D$ whose 
	$d$-occupation measure $\Mink_A$ exists and satisfies \eqref{eq:energy}. Let $h$ be a free Liouville field on $D$.
	Then $\nu=e^{\alpha h} \Mink_A$ exists in the sense of Definition~\ref{def:GMC} and is non-atomic.
\end{lemma}
We expect that Lemma~\ref{lem:GMC} remains true for $\alpha\in [\sqrt{d}, \sqrt{2d})$, but the $\alpha\in(0,\sqrt d)$ case is more straightforward to verify by the  $L^2$ argument  and  is sufficient for our purpose.

We now formulate a coordinate change formula that is convenient for our applications.
\begin{definition}[Coordinate change]\label{def:coord-change}
	Fix $d\in(0,2)$ and a Jordan domain $D$. Define ${\mathbf Q}(\alpha,d):=\alpha/2+d/\alpha$ and let $\alpha\in (0,\sqrt{d})$ be such that ${\mathbf Q}(\alpha,d)=5/\sqrt{6}$. Consider a triple $(A,\phi, h)$ of random variables with the following properties:
	$A$ is a compact subset of $D$ whose 
	$d$-occupation  measure $\Mink_{A}$ exists and satisfies \eqref{eq:energy}, 
	$h$ is a free Liouville field on $D$ such that 
	$\nu=e^{\alpha h} \Mink_{A}$ exists in the sense of Definition~\ref{def:GMC}, and $\phi$ is a conformal map on  $D$. 
	Let 
	\eqb
	h_\phi\defeq h\circ\phi^{-1} +5/\sqrt{6}\cdot\log|(\phi^{-1})'|. 
	\label{eqn-lqg-coord}
	\eqe
	We say that \notion{coordinate change} holds for $(A,\phi, h)$  if $e^{\alpha  h_\phi}\Mink_{\phi(A)}$ exists in the sense of Definition~\ref{def:GMC} and  $e^{\alpha  h_\phi}\Mink_{\phi(A)}=\phi_*\nu$ a.s. Here $\phi_*\nu$ means the pushforward of $\nu$ under $\phi$.
\end{definition}

\begin{proposition} \label{prop:coord-change}
	Let $(A,\phi, h)$ be as in Definition~\ref{def:coord-change}. If $(\phi,A)$ is independent of $h$, then coordinate change holds for  $(A,\phi, h)$.
	\begin{comment}
	The  existence of $e^{\alpha  h_\phi}\Mink_{\phi(A)}$ is ensured by Lemma~\ref{lem:Mink-conf} and the first statement.
	\end{comment}
\end{proposition}
\begin{proof}
	The proposition follows from \cite[Proposition 2.2]{c-greater-than-1} for the case where $h$ is a GFF. (Here we  use the assumption  that $(\phi,A)$ is independent of $h$.)
	Adding a continuous function does not change the result, since the continuous function can be locally approximated by a constant. Finally, since coordinate change is an a.s.\ property, reweighting the probability measure does not change the result.
\end{proof}
\begin{remark}[KPZ]\label{rmk:KPZ}
With ${\mathbf Q}$ as in Definition~\ref{def:coord-change}, the equation ${\mathbf Q}(\alpha,d)={\mathbf Q}(\gamma,2)$ 
is a version of the KPZ formula for fractals with Euclidean dimension $d$ on a $\gamma $-LQG surface. 
Heuristically, $\alpha$ describes the magnitude of the logarithmic singularity of the field at a point $z$ 
sampled according to the $\gamma$-LQG area measure ``conditioned on $z$ being on the fractal''.  
We require ${\mathbf Q}(\alpha,d)=5/\sqrt 6$ in Definition~\ref{def:coord-change} due to Convention~\ref{conv:DF}. For the pivotal points the relevant dimension is $d=3/4$, which gives $\alpha=1/\sqrt 6$. This explains why we consider GMC  with $\alpha=1/{\sqrt6}$  in Sections~\ref{subsub:cut} and~\ref{subsub:wedge}.
\end{remark}
We will apply coordinate change to various settings where the independence in Proposition \ref{prop:coord-change} does not hold.
Lemmas~\ref{lem:scaling} and~\ref{lem:coord-change} below are what we use in those cases.
\begin{lemma}\label{lem:scaling}
	In the setting of Definition~\ref{def:coord-change}, suppose coordinate change holds for $(A, \phi,h)$. Let $C\in\R$ and $\frk s>0$ be two random numbers coupled  with $(A, \phi,h)$. 
	(Here $C,\frk s$ are \emph{not} necessarily independent of $(A,\phi,h)$.)
	Then coordinate change holds for $(A,\frk s \phi,h+C)$.
\end{lemma}
\begin{proof}
	Almost surely, for any $C\in\R$ replacing $h$ by $h+C$ changes both the measures $e^{\alpha h}\Mink_{A}$ and $e^{\alpha h_\phi}\Mink_{\phi(A)}$ by a factor of $e^{\alpha C}$. Therefore coordinate change will hold for $(A,\phi,h+C)$ if it holds for $(A,\phi,h)$.
	It remains to show that coordinate change holds for maps of the form $z\mapsto \frk s z$.
	This property holds since we required the limit in Definition \ref{def:GMC} to be almost sure (rather than e.g.\ a limit in probability or a limit along powers of 2).
\end{proof}
\begin{lemma}\label{lem:coord-change}
	Fix $W>\frac43$. Let $\hw$ be the random distribution on $\bbH$ such that $(\bbH,\hw, 0,\infty)$ is the circle average embedding of a $W$-quantum wedge (recall Definition~\ref{def:circle}).	Suppose $D$ is a Jordan domain such that  $D\cup \p D\subset\bbH$. Let $A$ be  a random compact on $D$ whose $d$-occupation  measure $\Mink_{A}$ exists and satisfies \eqref{eq:energy}. Let $\phi$ be a random conformal map on $D$. If $(A,\phi)$ is conditionally independent of $\hw|_D$ given $\hw|_{D^c}$, then coordinate change holds for $(A, \phi,\hw )$.
\end{lemma}

Lemma~\ref{lem:coord-change} is an immediate consequence of Proposition~\ref{prop:coord-change} and the following lemma.
\begin{lemma}\label{lem:Dom-Markov}
	In the setting of Lemma~\ref{lem:coord-change},
	by enlarging the probability space, $\hw|_D$ can be written as $h_D+g$, where $h_D$ is a zero-boundary GFF on $D$
	independent of $\hw|_{D^c}$ and $g$ is an almost surely continuous function on $D$. 
\end{lemma}
\begin{proof}
	We can  write $\hw=h^{\ell}+h^{\op c}$ uniquely such that $h^{\ell}$ has average zero along all circles centered at the origin and $h^{\op c}$ is radially symmetric.
	Let $\ol h^{\op c}$ be independent of $\hw$ and have the law of the radially symmetric component of a free-boundary GFF on $\bbH$.  Here we fix the additive constant for $\ol h^{\op c}$ by letting its value on $\p \D\cap \bbH$ be equal to 0. 
	Then $\ol h\defeq h^{\ell} + \ol h^{\op c}$ is a free-boundary GFF independent of $h^{\op c}$. 
	In particular, $\ol h|_{D}$ can be written as a zero-boundary GFF $h_D$ plus the harmonic extension of $\ol h|_{D^c}$. Now $h_D$ is independent of $\hw|_{D^c}$ because $h_D$ is  independent of $(\ol h, h^c)|_{D^c}$. Moreover,
	$g\defeq\hw|_D-h_D$ is a.s.\ continuous on $D$.
\end{proof}

\subsubsection{Measure equivalence I: Brownian cut points}\label{subsub:cut}
In this section we prove a first version of Proposition~\ref{prop:quantum-Mink} (see Lemma~\ref{lem:pre}), which is based on a variant of planar Brownian motion called
the \notion{Brownian excursion in the upper half-plane}.
It is defined as the planar Brownian motion starting from 0 conditioned to stay inside $\bbH$ forever. See e.g.\ \cite[Section 5.3]{lawler-book} for how to make this conditioning precise.

The following proposition extracted  from \cite{lsw-restriction} is an example of  the deep relation between planar Brownian motion and $\SLE_6$.
 \begin{proposition}\label{prop:cutpt}
 	Let $(\cB_s)_{s\geq 0}$ be a Brownian excursion in the upper half-plane. Let $\eta$ be an $\SLE_6(2;2)$ on $(\BB H,0,\infty)$.
 	Let the \emph{hull} of $\ol\cB$ (resp., $\eta$) be the closure of the set of points $z\in\BB H$ for which we can find a $t>0$ such that $z$ is disconnected by $\cB([0,t])$ (resp., $\eta([0,t])$) from infinity. 
 	Then the hulls of $\cB$ and $\eta$ have the same law.
\end{proposition} 

Let $\eta'_{\ell}$ and $\eta'_{\op r}$ denote the left and right, respectively, boundary of the $\SLE_6(2;2)$ curve $\eta$. Then the interior of the hull of $\eta$ is bounded by $\eta'_{\markl}$ and $\eta'_{\op r}$. The rest of this section is devoted to the study of the set $\cC'\defeq \eta'_{\ell}\cap\eta'_{\op r}$. 

A point $p$ on the trace of $\cB$ is called a \emph{cut point} if removing $p$ disconnects the trace.
By Proposition~\ref{prop:cutpt},  $\cC'$ has the same law as the set of cut points of $(\cB_s)_{s\geq 0}$.
The occupation measure of Brownian cut points is thoroughly studied in \cite{hlls-cut}. In particular, we have the following.
\begin{proposition}\label{prop:Mink-SLE}
	Let $U$ be a bounded domain with piecewise smooth boundary 
	satisfying  $U\Subset \bbH$ (namely, $U\cup\p U\subset \bbH$). Set $A=\cC'\cap U$.
	Then the $3/4$-occupation measure (see Definition~\ref{def:Mink})  $\Mink_A$ of $A$ exists, and for each $\eps \in (0,3/4)$,
	$\iint_{U\times U }\frac{d\Mink_A(x) \,d\Mink_A(y)}{|x-y|^{3/4-\eps}}<\infty$ a.s.
\end{proposition}
\begin{proof}
 Since $\cC'$ has the same law as the cut points of $(\cB_s)_{s\geq 0}$,  Proposition~\ref{prop:Mink-SLE} follows from \cite[Theorem~4.12]{hlls-cut}.
\end{proof}
The following fact allows us to ignore the domain boundary when considering $\cC'$. For technical convenience we focus on a particular class of domains.
A Jordan domain  $D$ with piecewise smooth boundary is called a \notion{dyadic polygon} if $\p D$ is contained in $\cup_{k\in\N} \{(x,y)\in \R^2: 2^kx\in \Z \textrm{ or  }2^ky\in\Z\}$.
\begin{lemma}\label{lem:avoid-bdy}
For  each fixed dyadic polygon $U\Subset \bbH$, we have $\cC' \cap \p U=\emptyset$ a.s.
\end{lemma}
\begin{proof}
We first prove that $\P[\cC'\cap \{z: \op{Im} z=y \} \neq \emptyset]=0$ for each $y>0$. By scaling invariance of $\cC'$,  $\P[\cC'\cap \{z: \op{Im} z=y \} \neq\emptyset]$ does not depend on $y$. By way of contradiction, assume that the probability is positive. Let $Z$ be the Lebesgue measure of the set $\{y\in(0,1): \cC'\cap \{z: \op{Im} z=y \} \neq\emptyset\}$.  Then $\E[Z]>0$, hence $\P[Z>0]>0$.
Let $A=\cC'\cap \{z: \op{Im} z\in (0,1)\}$. Using the notations in Lemma~\ref{def:Mink}, we have $m^r_{A,1}(A)\ge Z$ for each $r>0$. This contradicts  the fact that $\lim_{r\to 0}m^r_{A,1}(A)=0$ a.s.\ by Proposition~\ref{prop:Mink-SLE}. By the same argument we have $\P[\cC'\cap \{z: \op{Re} z=x \} \neq \emptyset]=0$ for each $x\in \R$. This concludes the proof.
\end{proof}
In our proof of Proposition~\ref{prop:quantum-Mink} in Section~\ref{sec:prop:quantum-Mink}, $\cP_\eps$ will be covered by a finite  union of pieces that look like $\cC'\cap U$.
By Proposition~\ref{prop:Mink-SLE} and Lemma~\ref{lem:avoid-bdy}, there exists a non-atomic measure $\Mink'$ supported on  $\cC'$ such that for each fixed dyadic polygon $U\Subset\bbH$, the $3/4$-occupation measure  of  $\cC'\cap U$ a.s.\ equals $\Mink'|_U$. For more general domains, we only need the following. 
\begin{lemma}\label{lem:finite-exp}
For each bounded set $V\subset \bbH$  we have  $\E[\Mink'(V)]<\infty$.
\end{lemma}
\begin{proof}
This is an immediate consequence of the estimate for $G^{\mathrm{cut}}_\bbH$ in \cite[Theorem 4.12]{hlls-cut}.
\end{proof}

Let $h'$ be a random distribution on $(\bbH,0,\infty)$ independent of $\eta'_{\ell}$ and $\eta'_{\op r}$ such that $(\bbH,h',0,\infty)/\mathord\sim$ is a $14/3$-quantum wedge.  
Given a dyadic polygon $U\Subset \bbH$, by Lemma~\ref{lem:GMC} and Proposition~\ref{prop:Mink-SLE}, 
$e^{h' /\sqrt{6}}(\Mink'|_U)$ exists in the sense of Definition~\ref{def:GMC} and is non-atomic. 
We abuse notation and let $e^{h' /\sqrt{6}} \Mink'$ denote the random measure supported on $\cC'$ such that
for each $U$, $(e^{h' /\sqrt{6}} \Mink')|_U=e^{h' /\sqrt{6}}(\Mink'|_U)$ a.s.

Now we are ready to state the preliminary version of Proposition~\ref{prop:quantum-Mink} for $\cC'$.
\begin{lemma}\label{lem:pre}
	With notations introduced above, suppose $\eta'_\ell$ is parametrized by quantum length. 
	Then $(\eta'_\ell)^{-1}(\cC')$	 has the law of the range of a $1/2$-stable subordinator.
  	Let $\nu'$ be  the pushforward under $\eta'_\ell$ of the $1/2$-occupation measure  of $(\eta'_\ell)^{-1}(\cC')$.
  	Then $\nu'=ce^{h' /\sqrt{6}}\Mink'$ a.s.\ for some deterministic constant  $c>0$.
\end{lemma}

\begin{figure}
	\centering
	\includegraphics[scale=1]{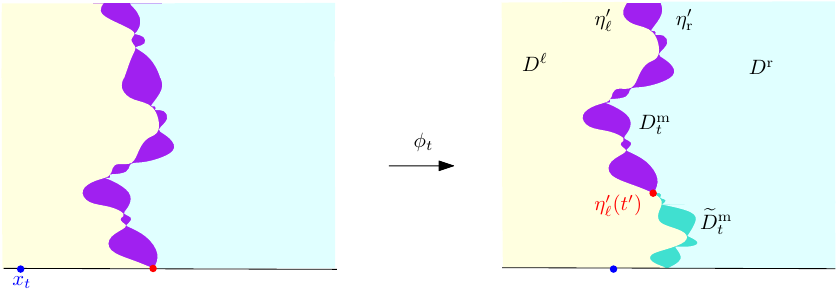}
	\caption{Illustration of the proof of Lemma~\ref{lem:pre}. 
		The green region is $\wt D^{\op m}_t$ and the purple region in the right figure is $D^{\op m}_t$. 
		The blue point $x_t$ on the left figure is such that the $h^t$-quantum boundary length of $[x_t,0]$ equals $1+t'$.} 
	\label{fig:pre}
\end{figure}
\begin{proof} 
	Using results in~\cite{wedges}, there are several ways to see that $(\eta'_\ell)^{-1}(\cC')$ can be realized as the range of a $1/2$-stable subordinator, which by Lemma~\ref{lem:Mink-zero} has $1/2$-occupation measure.
	For example, we can apply  Proposition~\ref{prop:zipper6} to the setting of Propositions~\ref{def:natural-time} to \ref{prop:Z}, which means $W^\ell=W^{\op r}=2/3$ in  Proposition~\ref{prop:zipper6}.
	Since $\eta'_{\markl}$ is parametrized by the quantum length, we see that $(\eta'_\ell)^{-1}(\cC')$ has the same law as \(\{s\ge 0: \eta^0_{\markl}(s)\in \eta([0,\infty))\cap \p D\}\) in Lemma~\ref{lem:bdy}.
		
	It remains to prove that $\nu'=ce^{h' /\sqrt{6}}\Mink'$ a.s.\ for some deterministic constant  $c>0$. The idea of our proof is to use the ergodic theorem and scaling properties to argue that that the ratio of the measures $\nu'$ and $e^{h' /\sqrt{6}}\Mink'$ is constant a.s.
	We advise the reader to look at Figure~\ref{fig:pre} while reading the rest of the proof. 
	
	Without loss of generality we assume that $(\bbH,h',0,\infty)$ is the circle average embedding of $(\bbH,h',0,\infty)/\mathord\sim$.
	Let $D^{\op m}$ be the region bounded by $\eta'_{\markl}$ and $\eta'_{\op r}$.
	Let $D^{\markl}$ and $D^{\op r}$ be the interior of the left and right, respectively, connected components of $\bbH\setminus D^{\op m}$. Let 
	$\cW_0^{\markl}=(D^{\markl},h',0,\infty,)/\mathord\sim$, $\cW_0^{\op r}=(D^{\markl},h',0,\infty)/\mathord\sim$,
	and $\cW_0^{\op m}=(D^{\op m},h',0,\infty)/\mathord\sim$. By Proposition~\ref{prop:zipper6}, $(h',\eta'_\ell,\eta'_{\op r})$ is determined by $(\cW^{\markl}_0,\cW^{\op m}_0,\cW^{\op r}_0)$.

	For each fixed $t>0$, let $t'=\inf\{s\ge 0: m'([0,s])=t   \}$, where $m'$ is the $1/2$-occupation measure  of $(\eta'_\ell)^{-1}(\cC')$.
	Let $\ol {D^{\op m}}$ be the closure of $D^{\op m}$, let $D^{\op m}_t$ be the  interior of the unbounded component of $\ol{D^{\op m}}\setminus\{\eta'_\ell(t')\}$,  
	and $\wt D^{\op m}_t$ be the closure of the bounded component of $\ol{D^{\op m}}\setminus\{\eta'_\ell(t')\}$.
	Let $\cW^{\markl}_t=(D^{\markl}, h', \eta'_\ell(t') , \infty)/\mathord\sim$, 
	$\cW^{\op r}_t=(D^{\op r}, h', \eta'_\ell(t') , \infty)/\mathord\sim$, and  
	$\cW^{\op m}_t=(D^{\op m}_{t}, h', \eta'_\ell(t'), \infty)/\mathord\sim$. By Proposition \ref{prop:zipper6}, $\cW^{\op m}_0$ is a weight $2/3$ wedge. 
		
	We claim that $\cW^{\op m}_t\eqd \cW^{\op m}_0$, i.e., $\cW^{\op m}_t$ also has the law of a weight $2/3$ wedge.
	Recall the notations of Definition~\ref{def:thin}, in particular the p.p.p.\ $\cE$ with intensity measure $\ell^{-3/2}d\ell\otimes dt$.
	Let $\cE'$ denote the collection of pairs $(\ell',s)$ such that $\ell'$ is the left boundary length of the LQG disk associated with $(\ell,s)\in\cE$. Since the two marked points of  the LQG disks constituting $\cW^{\op{m}}_0$ have the law of uniform and independent points sampled from the boundary measure (\cite[Proposition A.8]{wedges}), we get that $\cE'$ has the law of a p.p.p.\ with intensity measure $\frac 23 \ell^{-3/2}d\ell\otimes dt$. It follows that the process defined by $\tau_u:=\sum_{(\ell',s)\in\cE',s\leq u} \ell' $ is a $1/2$-stable subordinator.
	Recall that Lemma \ref{lem:Mink-zero} relates the time set and the $1/2$-occupation measure of the range of a stable subordinator. Since $m'$ is defined to be the $1/2$-occupation measure on the range of $\tau$, the lemma implies that for some deterministic constant $c_{1/2}>0$ we have $m'([0,\tau_{u}])=c_{1/2}u$ for all $u>0$ a.s. By definition, $\cW^{\op m}_t$ contains precisely the surfaces for which the point $(\ell',s)\in\cE'$ satisfies $m'([0,\tau_{s}])>t$. Using $m'([0,\tau_{u}])=c_{1/2}u$, we get that
 	$\cW^{\op m}_t$ contains precisely the surfaces for which the point $(\ell,s)\in\cE$ satisfies $s>c_{1/2}t$. By the definition of a Poisson point process $\{ (\ell,s-c_{1/2}t)\in\cE\,:\,s>c_{1/2}t \}\eqD \cE$, which implies that $\cW^{\op m}_t$ (consisting of the disks corresponding to $(\ell,s)$ with $s>c_{1/2}t$) is equal in law to $\cW^{\op m}_t$ (consisting of the disks corresponding to $(\ell,s)$ with $s>0$), i.e., $\cW^{\op m}_t\eqD \cW^{\op m}_0$. 
	
	%Then the process defined by $\tau_t:=\sum_{(\ell',s)\in\cE',s\leq t} \ell' $ is a $1/2$-stable subordinator *explain*.
	%Recall that Lemma \ref{lem:Mink-zero} relates the time set and the $1/2$-occupation measure of the range of a stable subordinator. Applying this lemma we get that 
	%Therefore, by Lemma \ref{lem:Mink-zero}, $\cW^{\op m}_t$ corresponds precisely to the surfaces for which the point $(\ell,s)\in\cE$ satisfies $s>c't$  for some deterministic constant $c'>0$. Hence $\cW^{\op m}_t$ also has the law of a weight $2/3$ wedge as claimed. } \xncomment{[Add three more sentences of justification.]}
	
	Since $\cW^{\markl}_0$ and $\cW^{\op r}_0$ are $2$-quantum wedges independent of $\cW^{\op m}_0$ and $t'$ is determined by $\cW^{\op m}_0$, 
	we see that $t'$  is independent of $\cW^{\markl}_0$ and $\cW^{\op r}_0$. Therefore, by Remark~\ref{rmk:wedge}, 
	$(\cW^{\markl}_t,\cW^{\op r}_t) \eqd (\cW^{\markl}_0,  \cW^{\op r}_0)$, so by Proposition~\ref{prop:zipper6} $(\bbH\setminus\wt D^{\op m}_t,h, \eta'_\ell(t'),\infty)/\mathord\sim$ is a $14/3$-quantum wedge.
	Let $\cW^{\op m}_0\setminus \cW^{\op m}_t$ be the collection of LQG surfaces in $\cW^{\op m}_0$ but not in $\cW^{\op m}_t$, ordered in the same way as in $\cW^{\op m}_0$.
	Then $\cW^{\op m}_0\setminus \cW^{\op m}_t$ and $(\cW^{\markl}_t,\cW^{\op m}_t,\cW^{\op r}_t)$ are independent. 
	
	Let $\phi_t:\bbH\to\bbH\setminus \wt D^{\op m}_t$ be the conformal map 
	such that $h^t\defeq h'\circ \phi_t  +{\mathbf Q}\log|\phi_t'|$ 
	has the same law as $h'$. Namely, $(\bbH,h^t,0,\infty)$ is  the circle average embedding of  $(\bbH\setminus \wt D^{\op m}_t,h, \eta'_\ell(t'),\infty)/\mathord\sim$. Then
	the set $\phi_t^{-1}(\cC')$, the field $h^t$, and $\cW^{\op m}_0\setminus \cW^{\op m}_t$
	are independent. 
	
	For a dyadic polygon $U\Subset\bbH$, set $A=\phi_t^{-1}(\cC')\cap U$. 
	We claim that $\phi_t$ can be written as $\frk s\phi$, where $\frk s$ is a random positive scaling constant and $\phi$ is determined by $h^t|_{U^c}$ and $\cW^{\op m}_0\setminus \cW^{\op m}_t$. We postpone the proof of this claim and proceed to conclude the proof of Lemma~\ref{lem:pre}.
	By Lemma~\ref{lem:coord-change} and this claim, the coordinate change in Definition~\ref{def:coord-change} holds for $(A, \phi, h^t)$. By Lemma~\ref{lem:scaling}, the same coordinate change holds for $(A, \phi_t, h^t)$.
	Let $X_t$ be the $e^{h'/\sqrt{6}}\Mink'$-mass of $\wt D^{\op m}_t$, which is almost surely finite by Lemma~\ref{lem:finite-exp}.
	For a fixed $s>0$, let $D^{\op m}_{t,s}$ be the closure of $\wt D^{\op m}_{t+s}\setminus \wt D^{\op m}_t$ so that  $X_{t+s}-X_s$ equals the $e^{h'/\sqrt{6}}\Mink'$-mass of $D^{\op m}_{t,s}$.
	Varying $U$ we see that $X_{t+s}-X_t$ equals the $e^{h^t/\sqrt{6}}\Mink'_t$-mass of $\phi_t^{-1}(D^{\op m}_{t,s})$ a.s., where 
	$\Mink'_t$ and $e^{h^t/\sqrt{6}} \Mink'_t$ are defined  in the same way as $\Mink'$ and $e^{h'/\sqrt{6}}\Mink'$ with $\phi^{-1}_t(D^{\op m}_{t})$ and $ h^t$ in place of $D^{\op m}$ and $h'$.
	Therefore, the process $(X_{t+s})_{s\ge 0}$ is determined by  $(\cW^{\markl}_t,\cW^{\op m}_t, \cW^{\op r}_t)$ in the same way as $(X_{s})_{s\ge 0}$ is determined by $(\cW^{\markl}_0,\cW^{\op m}_0, \cW^{\op r}_0)$, thus $(X_s)_{s\ge 0}$ has stationary increments. 
	
	By adding constants to $h'$ and using Remark~\ref{rmk:wedge} and \eqref{eq:scaling}, we see that
	the law of $X_t/t$ does not depend on $t$.  For $M\in (0,\infty)$, let $Y^M_i=(X_{i}-X_{i-1})\wedge M$ for $i\in\N$.  
	Then by ergodic theorem,
	$\lim_{n\to\infty}n^{-1}\sum_{i=1}^{n}Y^M_i$ exists almost surely. 
	We realize $D^{\op m}$ as the hull of a Brownian excursion $\cB$ independent of $h'$.
	Then  the limit belongs to the $\sigma$-algebra of $h'$ and $\cB'$ restricted to $\bbH\setminus(R\D)$.
	Taking $R\to\infty$, the tail triviality of $(h',\cB)$ yields that 
	$$\lim_{n\to\infty}n^{-1}\sum_{i=1}^{n}Y^M_i=\E[Y^M_1]=\E[X_1\wedge M]\quad  \textrm{a.s.}$$ 
	On the other hand,   since $n^{-1}\sum_{i=1}^{n}Y^M_i\le n^{-1}X_n$ and $n^{-1}X_n\eqd X_1$, we have $\P[X_1\ge \E[X_1\wedge M]]=1$. Letting $M\to\infty$, we get $X_1=\E[X_1]<\infty$ a.s.  Therefore $X_t=\E[X_1]t$ a.s.\ for all $t\ge 0$. This proves Lemma~\ref{lem:pre} with $c=\E[X_1]^{-1}\in (0,\infty)$.  
	
	It remains to prove the above mentioned claim that $\phi_t=\frk s \phi$. We can let $\frk s$ be such that the quantum length of $[-1,0]$ with respect  to the field 
	$h_{\frk s} (\cdot)\defeq h'(\frk s \cdot )+{\mathbf Q}\log \frk s$ equals $1$. Let $\phi=\frk s^{-1}\phi_t$ so that  $h^t=h_{\frk s}\circ \phi  +{\mathbf Q}\log|\phi'|$. Let $x_t=\phi^{-1}(-1).$ Then the quantum length of $[x_t,0]$ with respect to $h^t$ equals $t'+1$, which means that $x_t$ is determined by $h^t|_{U^c}$ and $\cW^{\op m}_0\setminus \cW^{\op m}_t$.
	Conditioning on $h^t|_{U^c}$ and $\cW^{\op m}_0\setminus \cW^{\op m}_t$, let $\wh\phi$ be a conditionally independent sample of $\phi$. It suffices to show that $\phi=\wh \phi$ a.s. 
	Note that the surface $(\bbH,h',0,\infty)/\mathord\sim$ can be obtained by identifying boundary arcs of the surfaces $(\bbH,h^t,0,\infty)/\mathord\sim$ and $\cW^{\op m}_0\setminus \cW^{\op m}_t$
	according to the quantum length. This defines a bijective map $\psi:\bbH\to\bbH$ such that $\wh\phi=\psi\circ \phi$ (in particular, $\psi$ is conformal on the image of $\phi$, which equals  $\bbH\setminus \big(\frk s^{-1}(\wt D_t^{\op m}\cup\p \wt D_t^{\op m} )\big)$), $\psi$ is conformal inside $\frk s^{-1}\wt D^{\op m}_t$, and $\psi$ is continuous everywhere. 
	By the conformal removability of $\frk s^{-1}(\eta'_\ell\cup\eta'_{\op r})$ (Lemma~\ref{lem:conf-rem}), $\psi$  is conformal  on the entire $\bbH$.\footnote{The way we apply conformal removability first appeared in the proof of \cite[Theorems 1.3 and 1.4]{shef-zipper}. See also \cite[Theorem 1.4]{wedges}.}
	Since $\psi(\infty)=\infty$, $\psi(0)=0$,  and $\phi(x_t)=\wh \phi(x_t)=-1$, we have that $\psi$ is the identity and hence $\phi=\wh \phi$ a.s.
\end{proof}

\subsubsection{Measure equivalence II:  intersections of bi-chordal SLE$_6$}
\label{subsub:wedge}
Recall the  setting of  Definition~\ref{def:bi-chordal}. In this section we formulate and prove a variant of  Proposition~\ref{prop:quantum-Mink} with $\sle_Q^{ad}\cap\sle_Q^{cb}$ in place of $\cP_\eps$, namely, Proposition~\ref{prop:pivot-bichordal} below. We will use Lemma \ref{lem:pre} and that the two measures considered in that lemma transform in the same manner when we add a continuous function to the field.
We first introduce a degenerate version of 2-$\SLE_6$ with an extra scaling invariance.
\begin{definition}\label{def:degen}
Let $\eta_1$ be an $\op{SLE}_6(0;2)$ on $(\bbH,0,\infty)$.	Let $\bbH'$ denote the component of $\bbH \setminus\eta_1$ whose boundary contains $(0,\infty)$.
Conditioned on $\eta_1$, let $\eta_2$ be an $\op{SLE}_6$ on $(\bbH',0,\infty)$.    
\end{definition}

\begin{figure}
	\centering
	\includegraphics[scale=1]{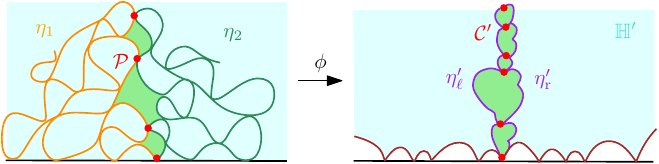}
	\caption{Illustration of the statement and proof of Lemma \ref{lem:cut}.   }
	\label{fig:pivot-bichordal1}
\end{figure}

\begin{remark}\label{rmk:degen}
To see why Definition~\ref{def:degen} gives a degenerate notion of $2$-$\SLE_6$, let $\tau\defeq\inf\{ t\geq 0\,:\, \op{Im}\eta(t)=1\}$.  	
Let $\ol\sle_i$ be the reversal of $\sle_i$ for $i=1,2$.  Let $\ol\tau$ be the first time such that  the unbounded component $\wh Q$ of  $\bbH\setminus(\sle_1([0,\tau])\cup \ol \sle_1([0, \ol\tau]))$ can be conformally mapped to $Q$ with $(\sle_1(\tau), 0,\infty, \ol \sle_1( \ol \tau))$ mapped to $(a,b,c,d)$.  It is argued in \cite[Lemma~4.4]{hlls-pivot} that $\P[\ol{\tau}<\infty]>0$ and moreover, on the event $E=\{\ol{\tau}<\infty\}$, the remainder of $\sle_1$ has the law of an $\SLE_6$ conditioned not to hit the real line.
Denote the conformal map from $\wh Q$ to $Q$ by $\psi$. (See Figure~\ref{fig:pivot-bichordal1}.) By the choice of $\tau$ and $\ol{\tau}$, the image of the remainder of $\sle_1$ under $\psi$ has the law of a chordal $\SLE_6$ on $(Q,a,d)$ conditioning on avoiding $\p_{b,c}Q$. Therefore the  image of the remainder of $\sle_1$ and $\ol\sle_2$ under $\psi$, as a pair of curves, 
have the law of $(\sle_Q^{ad},\sle_Q^{cb})$.
\end{remark}

We first prove the variant of Proposition~\ref{prop:pivot-bichordal} in the degenerate  case.   
\begin{lemma}\label{lem:cut}  
	Let $(\eta_1,\eta_2)$ be as in Definition~\ref{def:degen}.
	Let $\hw$ be a field independent of $(\eta_1,\eta_2)$ such that $(\BB H,\hw,0,\infty)$ is the circle-average embedding of a $10/3$-quantum wedge. 
	Let $\cP\defeq\sle_1\cap\sle_2$.  Then Proposition~\ref{prop:Mink-SLE} and Lemma~\ref{lem:avoid-bdy} hold with $\cP$ in place of $\cC'$ so that we can define the measures
	$\Mink_\cP$ and $e^{\hw/\sqrt{6}}\Mink_{\cP}$ in the same way as $\Mink'$ and $e^{\hw/\sqrt{6}}\Mink'$ in Lemma~\ref{lem:pre}.		
	Let $\eta_1^{\op r}:[0,\infty)\to \BB H\cup\partial \BB H$ be the right boundary of $\eta_1$ (recall Proposition~\ref{prop:duality}) parametrized by quantum length, 
	starting from $\eta_1^{\op r}(0)=0$. Then $(\eta_1^{\op r})^{-1} (\cP)$ has the law of the range of a $1/2$-stable subordinator.
	Moreover, $\nu =c e^{\hw/\sqrt{6}}\Mink_{\cP}$, where  $\nu$ is the pushforward of the $1/2$-occupation measure  of $ (\eta_1^{\op r})^{-1} (\cP)$ and $c$  is as in Lemma~\ref{lem:pre}.
\end{lemma}
\begin{proof}
	Consider two $2/3$-quantum wedges $\cW_1$ and $\cW_2$ which are independent of each other and of $\hw$. 
	Recall Lemma \ref{prop:zipper83}. Let $\cW'$ be the $14/3$-quantum wedge obtained by 
	conformally welding $\cW_1$, $\cW\defeq(\BB H,\hw,0,\infty)/\mathord\sim$, and $\cW_2$, such that $\cW_1$ (resp., $\cW_2$) is to the left (resp., right) of $\cW$. 
	Let $(\BB H,h',0,\infty)$ be the circle average embedding of $\cW'$. Let $\bbH'\subset \bbH$ be such that $\cW=(\bbH',h'|_{\bbH'},0,\infty)/\mathord\sim$ and
	let $\phi:\bbH\to \bbH'$ be the conformal map  such that $\hw=h'\circ\phi+{\mathbf Q}\log|\phi'|$ on $\bbH$.  See Figure~\ref{fig:pivot-bichordal1} for an illustration.
	 
	 Let $\eta_2^{\markl}$ be the left boundary of $\eta_2$.   Applying Proposition~\ref{prop:zipper6} twice we see that $\eta_1^{\op r}$  and  $\eta_2^{\markl}$  cut $\cW$ into three independents quantum wedges of weight $4/3$, $2/3$, and $4/3$, respectively.  Let $\eta'_{\markl}=\phi \circ \eta_1^{\op r}$ and $\eta'_{\op r} =\phi\circ\eta_2^{\markl}$. 
	 Then $\eta'_{\op r}$  and  $\eta'_{\markl}$  cut $\cW'$ into three independents quantum wedges of weights $2$, $2/3$, and $2$, respectively.  Namely, Lemma~\ref{lem:pre} applies to $(h',\eta'_{\markl},\eta'_{\op r})$ defined here. Let $\cC'=\phi (\cP)=\eta'_{\markl}\cap\eta'_{\op r}$. Then $(\eta_1^{\op r})^{-1} (\cP)=(\eta'_{\markl})^{-1}(\cC')$ has the law of the range of a $1/2$-stable subordinator. By Lemma~\ref{lem:Mink-conf}, Proposition~\ref{prop:Mink-SLE} holds with $\cP$ in place of $\cC'$. Moreover, note that the argument for Lemma~\ref{lem:avoid-bdy}  still applies if $\cC'$ is replaced by $\cP$ since it is scaling invariant.
	 
	 Define $\nu'$, $\Mink'$, and $e^{h'/\sqrt{6}}\Mink'$  as in Lemma~\ref{lem:pre}. Then $\nu=\phi_*\nu'$. To conclude our proof, we must show $e^{\hw/\sqrt{6}}\Mink_{\cP}=\phi_*(e^{h'/\sqrt{6}}\Mink')$.
	 It is sufficient to show that the coordinate change in Definition~\ref{def:coord-change} applies to $(\cP\cap U, \phi, \hw)$ for each dyadic polygon $U\Subset\bbH$. 
	 Recall that in Lemma~\ref{lem:pre} the same is proved for $(A, \phi_t, h^t)$ based on the theory of quantum zippers in Section~\ref{subsub:zipper}, as well as  Lemmas~\ref{lem:scaling} and~\ref{lem:coord-change}. A similar argument applies to $(\cP\cap U, \phi, \hw)$, where we need the conformal removability of $\p \bbH'$.
	 We leave the details to the reader.\qedhere
\end{proof}

\begin{figure} 
	\centering
	\includegraphics[scale=1]{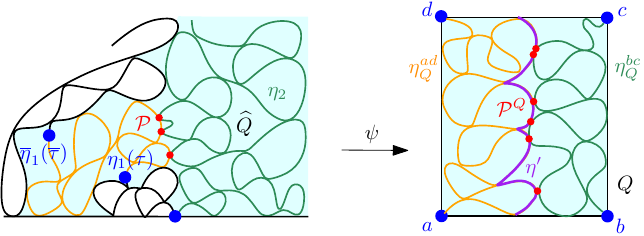}
	\caption{Illustration of Definition~\ref{def:degen} and the proof of Proposition \ref{prop:pivot-bichordal}. The domain $Q'$ (not indicated in the figure) is the subset of $Q$ to the right of $\eta'$.}
	\label{fig:pivot-bichordal2}
\end{figure}

	In the rest of this section, let $\eta_1$, $\eta_2$, $\ol\eta_2$, $ \wh Q$, $\psi$, $E$,  $(Q, a, b,c,d)$ be as in Definition~\ref{def:degen} and Remark~\ref{rmk:degen}. Moreover, we condition on the positive probability event $E$.  We identify the image under $\psi$ of the part of $(\eta_1, \ol\eta_2)$ inside $\wh Q$  as $(\sle^{ad}_Q,\sle^{cb}_Q)$ in Definition~\ref{def:bi-chordal}.  Let $\hw$ and $\cP$ be as in Lemma~\ref{lem:cut}. Let $\wt h$   be the field on $Q$ such that $(Q,\wt h)\sim_\psi (\wh Q,\hw|_{\wh Q})$.  Also recall $Q'$ from Definition~\ref{def:bi-chordal}. Proposition~\ref{prop:pivot-bichordal} follows from Lemma~\ref{lem:cut} and the  observation below.
\begin{lemma}\label{lem:Weyl}
	Let $h$ be a free Liouville field (Definition \ref{def:free}) which is independent of $(\eta^{ad}_Q,\eta^{cb}_Q)$. 
	Then we can enlarge the probability space generated by $(h,\eta^{ad}_Q,\eta^{cb}_Q)$ to a bigger probability space $(\Omega,\cF,\P)$ satisfying the following properties.
	There exists a random continuous function  $g$ measurable with respect to $(\Omega,\cF)$ and a probability measure $\wt \P$  such that 
	the $\wt\P$-law  of $h-g$ is that of $\wt h$ defined right above and $\P$ is absolutely continuous with respect to $\wt \P$.
	
	We abuse notation and set $\wt h\defeq h-g$. 
	Let $\eta'\defeq \p Q'\cap \eta^{ad}_Q$,
	and	let $\xi'_h$ and $\xi'_{\wt h}$ be the quantum length measure on $\eta'$ induced by $h$ and $\wt h$, respectively. Then $\P$-almost surely $\xi'_h=e^{2g/\sqrt{6}}  \xi'_{\wt h}$.
\end{lemma}
\begin{proof}
	To prove the first assertion, we first assume that $h$ which is a zero-boundary Gaussian free field (GFF) on $Q$.
	By Definition \ref{def:thick-wedge}, $\hw$ can be written as the sum of a free-boundary GFF and a continuous function.
	Note that $\hw$ is independent of $\wh Q$.  By the conformal invariance and domain Markov property of GFF,  there exists a coupling of a random continuous function $g$ 
	with $h$ such that $\wt h\eqd h-g$. Setting  $\wt\P=\P$ gives the first assertion in this case.  The general case follows from  the definition of a free Liouville field and the fact that a free-boundary GFF can be decomposed as a zero-boundary GFF plus a harmonic function.
	
	To prove the second assertion, let $ \phi_{Q'}: \bbH\to Q'$ be a conformal map. Let $h_1$ and $h_2$ be two random distributions on $\bbH$ such that 
	$(Q',\wt h|_{Q'})\sim_{\phi_{Q'}} (\bbH,h_1)$ and $(Q',h|_{Q'})\sim_{\phi_{Q'}} (\bbH,h_2)$. Let  $f=h_1-h_2$. 
	Then $f$ is continuous on $ \phi^{-1}_{Q'}(\eta')$. It is clear that $h_1$ is a free Liouville field, hence so is $h_2$.
	Restricted to $\phi^{-1}_{Q'}(\eta')$, we have  $e^{2h_1/\sqrt6}dx=e^{2f/\sqrt6} e^{2h_2/\sqrt6}dx$ a.s. This concludes the proof.
\end{proof}

\begin{comment}
This lemma has two purposes. Firstly, it is an input for the following Proposition~\ref{prop:pivot-bichordal}. 
Moreover, it implicit shows that  LQG boundary measure can be defined because $(D',h)$ equivalent to  a Liouville field on $\bbH$. 
\end{comment}

\begin{proposition}\label{prop:pivot-bichordal}
	Let $\cP^{Q}=\sle^{ad}_Q\cap\sle^{cb}_Q$. Let $h$ be a free Liouville field (Definition \ref{def:free}) which is independent of $(\eta^{ad}_Q,\eta^{cb}_Q)$. 
	Let $\cI=(\eta')^{-1} (\cP^Q)$ where $\eta'=\p Q'\cap \eta^{ad}_Q$ as in Lemma~\ref{lem:Weyl},  parametrized by the quantum length induced by $h$. 
	Then almost surely the following hold.
	\begin{enumerate}
		\item\label{item:Haus} 
		The $3/4$-occupation measure of $\cP^Q$ exists, which we denote by $\Mink_Q$.\footnote{The existence of $\Mink_Q$ is also proved in  \cite[Proposition~1.8]{hlls-pivot}. We include the proof here for completeness.}
		The measure $e^{h/\sqrt{6}}\Mink_Q$ exists as in Definition~\ref{def:GMC}.
		The $1/2$-occupation measure of $\cI$ exists. Let $\nu_{\cI}$ denote the pushforward of this measure by $\eta'$. 
		\item\label{item:const}   $\nu_{\cI}=ce^{h/\sqrt{6}}\Mink_Q$ with $c$ as in Lemma~\ref{lem:pre}. 
	\end{enumerate}
\end{proposition}
\begin{proof}
	Note that there almost surely exists a dyadic polygon $U\Subset\bbH$ such that $\psi^{-1}(\cP^Q)\subset U$. 
	Since $\cP\cap \p U=\emptyset$ a.s.\ in Lemma~\ref{lem:cut}, the existence of $\Mink_\cP$ in Lemma~\ref{lem:cut} combined with Lemma~\ref{lem:Mink-conf} implies that $\Mink_Q$ exists, and $\iint_{Q\times Q }\frac{d\Mink_Q(x)\, d\Mink_Q(y)}{|x-y|^{\frac34-\eps}}<\infty$ a.s.\ for $\eps\in (0,\frac34)$. Therefore $e^{h/\sqrt{6}}\Mink_Q$ exists.
	
	Let $\wt \eta'$ be $\eta'$ reparameterized by the quantum length induced by $\wt h$ from Lemma~\ref{lem:Weyl}.
	Then Lemma~\ref{lem:cut} implies that the $1/2$-occupation measure $\Mink_{\wt \cI}$ of $\wt \cI\defeq (\wt\eta')^{-1}(\cP^Q)$ exists. 
	Let $\nu_{\wt \cI}\defeq(\wt\eta')_*\Mink_{\wt \cI}$ denote the pushforward of $\Mink_{\wt \cI}$ by $\wt \eta'$. 
	Then Lemmas~\ref{lem:coord-change} and~\ref{lem:cut} further imply that $\nu_{\wt \cI}=c e^{\wt h/\sqrt{6}}\Mink_Q$ with $c$ as in Lemma~\ref{lem:pre}.   
	When we apply Lemma~\ref{lem:coord-change} here, we use in particular that $\psi$ and $\cP$ are independent of $\hw$.
	
	By Lemma~\ref{lem:Weyl},  it suffices to prove Proposition~\ref{prop:pivot-bichordal} in the case when $h$ has the form $\wt h+g$ where $g:Q\to\R$ is a random continuous function coupled with $\wt h$ in an arbitrary manner. 
	Note that  $\cP^Q\subset \wt \eta'(I)$ for some closed interval $I$. Without loss of generality, we assume that  $I=[0,A]$ for some $A>0$.
	Recall that $\eta'$ is parametrized according to the quantum length measure induced by $\wt h+g$. By Lemma~\ref{lem:Weyl},  $\eta'(s(t))=\wt \eta'(t)$ for each $t\in[0,A]$, where 
	\begin{equation}\label{eq:s}
	s(t) = \int_{0}^{t} e^{\sqrt{2/3}g(\wt \eta'(u))}\,du \;\textrm{for }t\in[0,A].
	\end{equation}
	Set $B\defeq s(A)$.
	Since $s:[0,A]\to [0,B]$ is a $C^1$ function with $s'>0$, and $s(\wt \cI)=\cI$, by Lemma~\ref{lem:1D-measure},  the  $1/2$-occupation  measure  $\Mink_\cI$ of $\cI$ exists and equals $|(s^{-1})'|^{-1/2} \cdot (s_*\Mink_{\wt \cI})$. By~\eqref{eq:s}, for each $x\in \eta'([0,B])$, we have that $|(s^{-1})'( (\eta')^{-1}(x))|^{-1/2}=(e^{-\sqrt{2/3}g(x)})^{-1/2} =e^{g(x)/\sqrt 6}$.
	Therefore
	$$(\eta')_*\Mink_\cI=e^{g/\sqrt 6} ((\eta')_*(s_*\Mink_{\wt \cI}) )=e^{g/\sqrt 6} ((\wt \eta')_*     \Mink_{\wt \cI})=e^{g/\sqrt 6} \nu_{\wt \cI}=ce^{g/\sqrt 6} e^{\wt h/\sqrt6}\Mink_Q=ce^{h/\sqrt6} \Mink_Q.$$
	Now $\nu_\cI=(\eta')_*\Mink_\cI=ce^{h/\sqrt6} \Mink_Q$ as desired.
\end{proof}

\section{Liouville dynamical percolation}\label{sec:GPS}

In this section we prove Lemmas~\ref{lem:eps-LDP} and~\ref{emb-prop:LDP}. This concludes the proof of Theorem~\ref{emb-thm:main}.
Lemma~\ref{lem:eps-LDP} is a relatively easy consequence of \eqref{eq:KPZ} and an ingredient (Proposition~\ref{prop:flip}) from \cite{ghs-metric-peano} and \cite{bhs-site-perc}.
For Lemma~\ref{emb-prop:LDP}, neither the convergence nor the ergodicity seems easy to access from random planar maps and  mating-of-trees perspective. 
To prove this lemma, we use the Liouville dynamical percolation introduced in \cite{ghss18}.
We review this object in Sections~\ref{subsec:quad0} and~\ref{subsec:LDP} and prove Lemma~\ref{emb-prop:LDP} in Section~\ref{subsec:lemma}, 
with certain ingredients supplied in later subsections.

We will use the following notions and conventions.  $\CLE_6$ in this section will be assumed to have  monochromatic blue boundary condition; see Definition~\ref{def:CLE}.
Given a finite measure $\mu$,  if $z$ is sampled from $\mu$ normalized to be a probability measure, we will simply say that $z$ is sampled from $\mu$. 
For a metric space $(X,d)$, recall that a process taking values in $X$ is called \emph{c\`adl\`ag} if it is right-continuous and has left limits everywhere.
In this section we will often consider convergence of c\`adl\`ag  processes in the Skorokhod topology. 
For functions $f_j:I_j\to X$ defined on bounded intervals $I_j\subset\R$ for $j=1,2$,
this topology is generated by the following metric 
$$
d_{\op{Sk}}(f_1,f_2) \defeq
\inf_\phi \sup_{t\in I_1} 
\Big(d\big( f_1(t),f_2( \phi(t) )  \big)
+|t-\phi(t)|\Big),
$$
where the infimum is taken over all increasing bijections $\phi:I_1\to I_2$. If $f_1$ and $f_2$ are defined on $[0,\infty)$, then we define $d_{\op{Sk}}$ similarly; more precisely,
$$
d_{\op{Sk}}(f_1,f_2) \defeq
\sum_{k=1}^\infty
\inf_{\phi} \sup_{t\vee\phi(t)\in [0,2^k]} 
2^{-k}\wedge\Big(d \big( f_1(t),f_2( \phi(t) )  \big)
+|t-\phi(t)|\Big),
$$
where the infimum is taken over all increasing bijections $\phi:[0,\infty)\to [0,\infty)$.

\subsection{Quad-crossing space}\label{subsec:quad0}
We start by recalling  a metric space due to Schramm and Smirnov~\cite{ss-planar-perc} as a method of describing the scaling limit of planar percolation other than loop ensembles. We will omit the detailed construction of the metric and only review materials necessary for this paper.

A {\it quad}  is 
a homeomorphism $Q$ from $[0,1]^2$  into $\C$,
where two homeomorphisms $Q_1$ and $Q_2$ are identified as the same quad if $Q_1([0,1]^2)=Q_2([0,1]^2)$, and $Q_1(z)=Q_2(z)$ for $z\in\{ (0,0),(0,1),(1,0),(1,1) \}$.
Let
\begin{align*}
\partial_1 Q\defeq Q(\{0\}\times [0,1]), \quad\partial_2 Q\defeq Q( [0,1]\times \{0\}),\\
\partial_3 Q\defeq Q(\{1\}\times [0,1]), \quad \partial_4 Q\defeq Q([0,1]\times \{1\}).
\end{align*}
A {\it crossing} of a quad $Q$ is a closed set in $\C$ containing a connected closed subset of $Q([0,1]^2)$ that intersects both $\partial_1 Q$ and $\partial_3 Q$. 
A natural partial order $\le$ can be defined on $\cQ_D$ by saying that  $Q_1 \le Q_2$ if and only if every crossing of $Q_1$ is also a crossing of $Q_2$.

Let $D$ be a bounded  domain. Let $\cQ_D$ denote the space  of all quads satisfying $Q([0,1]^2)\subset D$. 
We say that a subset $S\subseteq \cQ_D$ is \notion{hereditary} if, whenever $Q\in S$ and $Q'\in \cQ_D$ satisfies $Q'\leq Q$, we have $Q'\in S$.  
We call a  closed hereditary subset of $\cQ_D$ a {\it quad-crossing configuration} on $D$ and denote the space of quad-crossing configurations by $\cH(D)$. For $\omega\in\cH(D)$ we may identify it with a function $\omega:\cQ_D\to \{0,1\}$  such  that $\omega^{-1}(1)$ is closed in $\cQ_D$ and such that for any $Q_1,Q_2$ with  $Q_1\le Q_2$ and $\omega(Q_1)=1$, we have $\omega(Q_2)=1$. 
(Here we abuse notation and let $\omega$ denote both the element of $\cH(D)$ and the function from $\cQ_D$ to $\{0,1 \}$.)
By \cite{{ss-planar-perc}}, $\cH(D)$ can be endowed with a metric $d_\cH$ such that $(\cH(D),d_\cH)$ is a compact separable metric space. 
For each $Q\in \cQ_D$, the function $\omega\mapsto \omega(Q)$ is measurable with respect to the Borel $\sigma$-algebra of $(\cH(D),d_\cH)$. 
Moreover, there exists  a countable set $\{Q_n\}_{n\in \N}\subset \cQ_D$ such  that $Q_n([0,1]^2)$ has piecewise smooth boundary and 
\begin{equation}\label{eq:generator}
\textrm{ $\{\omega(Q_n)\}_{n\in \N}$ generates the Borel $\sigma$-algebra of $(\cH,d_\cH)$.}
\end{equation}

We now focus on the setting relevant to the remainder of the paper.
For $\mesh>0$, let $\omega^\mesh$  be a site percolation on $\D^\mesh$ (see the paragraph above Theorem~\ref{emb-thm:CLE} for the definition).
For each $Q\in \cQ_{\D}$, let $\omega^\mesh(Q)=1$ if and only if the union of all red hexagons on the dual lattice of $\D^\mesh$ gives a crossing of $Q$. This identifies $\omega^\mesh$ with an element in $\cH(\D)$. 
If $\omega^\mesh$ is sampled from Bernoulli-$\frac12$ site percolation,  then  $\omega^\mesh$  converges in law to a random variable $\omega$  in $\cH(\D)$ for the $d_\cH$-metric \cite{camia-newman-sle6,gps-pivotal}.  
Let $\ol \cQ_\D$ be the collection of quads such that $Q([0,1]^2) \subset \D\cup\p \D$. For each $Q\in \ol\cQ_{\D}$ we can still define $\omega^\mesh(Q)$ as before. 
In this section, we use the following lemma to extend $\omega$ from $\cQ_{\D}$ to $\ol\cQ_{\D}$.
\begin{lemma}\label{lem:bdy-quad}
	Almost surely $\omega$ admits an extension to $\ol\cQ_{ \D}$ such that for each fixed $Q\in \ol\cQ_\D$ $\lim_{n\to \infty}\omega(Q_n)=\omega(Q)$ in probability  where $Q_n$ is obtained by restricting $Q$ to $[n^{-1}, 1-n^{-1}]^2$.
	Suppose we are in a coupling such that $\lim_{\mesh\to 0}\omega^\mesh=\omega$ almost surely as elements in $\cH(\D)$.
	Then $\lim_{\mesh\to0}\omega^\mesh(Q)=\omega(Q)$ in probability for each fixed $Q\in \ol\cQ_{\D}$.
\end{lemma}	
\begin{proof}
	Suppose $\wh{\omega}^\mesh$ is defined as $\omega^\mesh$ with $2\D$ in place of $\D$. We further require that $\wh \omega^\mesh$ converge almost surely as elements in $\cH(2\D)$ and that $\omega^\mesh$ is obtained by restricting $\wh\omega^\mesh$ to $\D$. Let $\wh\omega=\lim_{\mesh\to 0} \wh\omega^\mesh$ in the $d_\cH$-metric. 
	By \cite[Lemma A.1]{ss-planar-perc},
	$\limsup_{\mesh\to0}\P[\wh\omega^\mesh(Q)\neq \wh\omega^\mesh(Q_n) ]=o_n(1)$. 
	By \cite[Corollary~5.2]{ss-planar-perc}, $\lim_{\mesh\to0}\omega^\mesh(Q)=\omega(Q)$ in probability for each fixed $Q\in \ol\cQ_{\D}$. 
	Therefore $\wh \omega$ restricted to $\ol\cQ_{\D}$ is the desired extension of $\omega$ as described in Lemma~\ref{lem:bdy-quad}. 
\end{proof}

\subsection{Liouville dynamical percolation}\label{subsec:LDP}

We first specify the setting under which we will prove Lemmas~\ref{lem:eps-LDP} and~\ref{emb-prop:LDP} in Section~\ref{subsec:lemma}.
Let $\gamma=\sqrt{8/3}$,  ${\bf Q}=5/\sqrt{6}$, and $a={\bf Q}-\gamma=1/\sqrt{6}$. 
We consider  a probability space  $(\Omega,\cF,\P)$ with random variables  $X_t,h^1,h^2,\hs$  whose law are as described in Definition~\ref{def:disk}.
Namely,   $(X_t)_{t\ge 0}$ has the law of $B_{2t}-at$, where $B_t$ is a standard Brownian motion, $(X_{-t})_{t\ge 0}$ is independent of $(X_t)_{t\ge 0}$, and $(X_{-t})_{t\ge 0}$ has the law of $B_{2t}-at$ conditioned on being negative. Moreover, $\hs=h^1+h^2$, where $h^1(z)=X_t$ for each $z\in \cS$ and $t\in \R$ with $\op{Re} z=t$. 
Finally, $h^2$ is  independent of $X_t$ with the law of the lateral component of the free-boundary GFF on $\cS$. 
Let $\Pd$ be the probability measure obtained from normalizing 
$e^{-\gamma M/4}\xi_\hs(\partial \cS)^{1/2} d\P$, where $M=\sup_{t\in \R} X_t$. (Recall from \eqref{eq:Boltzmann} that $({\bf Q}-\gamma)M=\gamma M/4$ and $4/\gamma^2 -1=1/2$.)
Let $\hd\defeq\hs-2\gamma^{-1}\log\xi_\hs(\partial\cS)$ so that under the $\Pd$-measure  $\hd$ is the field of a unit boundary length $\sqrt{8/3}$-LQG disk by definition.
Now let $\phi : \D\to \cS$ be the conformal map in Definition~\ref{def:bh}. Let $\bh$ be the field as  in Definition~\ref{def:bh}, i.e.,  
$\bh=\hd\circ\phi+{\bf Q}\log|\phi'|$. 

Let $\fh=\hs\circ \phi+{\bf Q}\log |\phi'|$. Then the fields $\fh$ and $\bh$ are related by a shift:
\eqb 
\bh=\fh-2\gamma^{-1}\log\xi_\fh(\partial\D).
\label{eq:fields}
\eqe
We are mainly interested in $\bh$  because under the $\Pd$-measure,  it is the field considered in  Lemmas~\ref{lem:eps-LDP} and~\ref{emb-prop:LDP}.
However, most technical work 
in this section will be done for $\fh$ instead because of the following lemma.

\begin{lemma}\label{lem:g}
	In the setting above,  $\fh$ can be written as $\Phi+g$, where the $\P$-law of $\Phi$ is a free boundary GFF as in Theorem~\ref{thm:LDP}
	and $g$ is a random continuous function on $\D$. Moreover,
	\begin{equation}\label{eq:g}
	g(z)\le {\bf Q}\log |\phi'(z)|-a|\op{Re}\phi(z)|\qquad \textrm{for all }z\in \D. 
	\end{equation}
\end{lemma}
\begin{proof}
	Let $h^{\op f}$ be the free boundary GFF on $\cS$ with average 0 along $i[0,\pi]$.  
	In the definition of $\hs$ in Definition~\ref{def:disk}, if the law of $X_t$ were set to be the two-sided Brownian motion $(B_{2t})_{t\in \R}$ without drift or conditioning, then the law of $\hs$ would be given by $h^{\op f}$. Since there exists a coupling of  $(B_{2t})_{t\geq 0}$ and $(X_t)_{t\geq 0}$ such that $X_t=B_{2t}-at$ for $t\ge 0$ and $X_{t}\le B_{-2t}+at$ for all $t\le 0$, we can couple $h^{\op f}$ and $\hs$ on the same probability space such that 
	\begin{enumerate}
		\item the lateral component of $h^{\op f}$ (see the paragraph above Definition~\ref{def:disk}) equals $h^2$;
		\item $\hs=h^{\op f}-a\op{Re} z$ on $\cS\cap \{z: \op{Re} z\ge 0 \}$; 
		\item $\hs\le h^{\op f}+a\op{Re} z$ on $\cS\cap \{z: \op{Re} z<0 \}$.
	\end{enumerate}
	Since $\fh=\hs\circ \phi+{\bf Q}\log |\phi'|$, taking $\Phi=h^{\op f}\circ \phi$ and $g=\fh-\Phi$ and using that $\phi$ maps $i[-1,1]$ to $[0,i\pi]$, we obtain~\eqref{eq:g}.
\end{proof}
The following immediate corollary of Lemma~\ref{lem:g} will be useful in Sections~\ref{subsec:pivot} and~\ref{subsec:stability}.
\begin{corollary}\label{cor:g}
	For $\fh$ and $\Phi$ in Lemma~\ref{lem:g}, given any $r\in (0,1)$, there exists a deterministic constant $c_r$ such that $\fh\le \Phi +c_r$ on $r\D\defeq \{z: |z|<r\}$. 
\end{corollary}

Now we review Liouville dynamical percolation in the setting specified above. Let 
\[
\xncomment{\mu'_\fh:= e^{\fh/\sqrt6} \,d^2z  = \lim_{\eps\to 0}\eps^{\frac1{12}} e^{\fh/\sqrt6} \,d^2z}
\]
be as defined in Definition~\ref{def:GMC} with $\alpha=1/\sqrt6$. 
Fix $\mesh>0$ and consider the lattice $\D^\mesh$. 
For each vertex  $v$ on $\D^\mesh$, let $\mu'_\fh(v)$ be the $\mu'_\fh$-mass of the hexagon on the dual lattice of $\D^\mesh$ corresponding to $v$.  
Let $\arm^\mesh_4(\mesh,r)$ be the probability of that Bernoulli-$\frac12$ site percolation on $\mesh\Tg$ possesses four disjoint monochromatic paths of alternating color from the origin to the boundary of the box $[-r,r]^2$.

Now we enlarge the probability space $(\Omega,\cF,\P)$ to contain random variables defined as follows.
For $\mesh>0$, let $\omega^\mesh_0$ be an instance of Bernoulli-$\frac12$ site percolation on $\D^\mesh$ with monochromatic blue boundary condition.
We assume that the loop ensembles corresponding to $\omega^\mesh_0$ converge $\P$-almost surely (see Theorem~\ref{emb-thm:CLE}).
We further require $\fh$ and $\{\omega^\mesh_0\}_{\mesh>0}$ to be independent under $\P$.
Consider a clock for each inner vertex of $\D^\mesh$ such that conditioning on $(\fh,\omega^\mesh_0)$, these  are independent exponential clocks  with rate $\mu'_\fh(v)\arm^\mesh_4(\mesh,1)^{-1}$. Namely, the  set of times when the clock at $v$ rings is a Poisson process on $(0,\infty)$ of intensity $\mu'_\fh(v)\arm^\mesh_4(\mesh,1)^{-1}$.
Now we define a dynamic on the space of site percolation configurations on $\D^\mesh$  as follows.
Letting the initial coloring be $\omega^\mesh_0$, when the clock rings at an inner vertex $v$, we flip the color at $v$. 
This defines a stationary process $(\omega^\mesh_t)_{t\ge 0}$, which by Section~\ref{subsec:quad0} can be viewed as taking values in $\cH(\D)$.
We call $(\omega^\mesh_t)_{t\ge 0}$ the \emph{discrete Liouville dynamical percolation} (LDP) on $\D^\mesh$ driven by $e^{\fh/\sqrt6}$.
We will use the following key input from \cite{ghss18}. 
\begin{theorem}\label{thm:LDP}
	There exists a probability space $(\Omega,\cF,\P)$  with random variables $\fh $,  $\{(\omega^\mesh_t)_{t\ge 0}: \mesh\in(0,1)\}$, and $(\omega_t)_{t\ge 0}$ satisfying the following.
	\begin{itemize}
		\item The joint law of $\fh $ and $\{(\omega^\mesh_t)_{t\ge 0}: \mesh\in(0,1)\}$ is as described right above.
		\item $(\omega_t)_{t\ge 0}$ is a stationary  process  taking values on $\cH(\D)$ with  following mixing property. 
		For any two events $A$ and $B$ in the Borel $\sigma$-algebra of $(\cH(\D),d_\cH)$,
		\(\lim_{t\to\infty}\P[\1_{\omega_0\in A}  \1_{\omega_t\in B}\mid \fh ] =\P[A]\P[B]\) almost surely.
		\item For each $r\in (0,1)$ and $t\ge 0$, let $\omega^\mesh_t|_{r\D}$ (resp., $\omega_t|_{r\D}$) be $\omega^\mesh_t$ (resp., $\omega_t$) restricted to $\cQ_{r\D}$, where $r\D\defeq \{z\in\C: |z|<r \}$.
		Then for each $r\in(0,1)$, $(\omega_t|_{r\D})_{t\ge 0}$ is a c\`adl\`ag process and $\lim_{\mesh\to 0}(\omega^\mesh_t|_{r\D})_{t\geq 0}=(\omega_t|_{r\D})_{t\geq 0}$ in probability  in the Skorokhod topology.
	\end{itemize}
\end{theorem}
\begin{proof}
	Note that $\Phi$ in Lemma~\ref{lem:g} under the probability measure $\P$
	is a Gaussian field on $r\D$ with kernel of the form $-\log|x-y|+C(x,y)$, where $C(\cdot,\cdot)$ is continuous up to the boundary of $r\D$.
	Therefore, if $g$ were equal to $0$ in Lemma~\ref{lem:g} so that $\fh=\Phi$, 
	Theorem~\ref{thm:LDP} would  fall into the framework of \cite{ghss18}. 
	The third assertion of Theorem~\ref{thm:LDP} would follow from \cite[Theorem~1.3]{ghss18}. 
	For the second assertion, if $A,B$  are in the Borel $\sigma$-algebra of $(\cH(r\D),d_\cH)$, then the second assertion would follow from  \cite[Theorem~1.4]{ghss18}. 
	Since the Borel $\sigma$-algebra of $(\cH(\D),d_\cH)$ is the minimal $\sigma$-algebra containing the Borel $\sigma$-algebra of $(\cH(r\D),d_\cH)$ for all $r\in(0,1)$,
	we would have the second assertion of Theorem~\ref{thm:LDP} without the constraint to $r\D$.
	
	Now, although $g\neq 0$, since $g$ is uniformly bounded from above and below on $r\D$, as explained in \cite[Remark~1.6]{ghss18},
	the non-quantitative results of \cite[Theorems~1.3 and~1.4]{ghss18} still hold and give Theorem~\ref{thm:LDP}.
\end{proof}

We call $(\omega_t)_{t\ge 0}$ the \emph{continuous  Liouville dynamical percolation} driven by $e^{\fh/\sqrt6}$. The boundary condition of $(\omega^\mesh_t)_{t\ge0}$ is irrelevant for Theorem~\ref{thm:LDP}. We impose the monochromatic boundary condition and restrict the update of colors only to inner vertices in order to mimic the dynamic
$(\Map^n,\omega^{n}_t)_{t\ge 0}$ in Section~\ref{subsub:DP}.

\subsection{Proof of Lemmas~\ref{lem:eps-LDP} and~\ref{emb-prop:LDP}}
\label{subsec:lemma}

\xncomment{In this section we prove Lemmas~\ref{lem:eps-LDP} and~\ref{emb-prop:LDP} with a few ingredients whose proofs are postponed to later subsections. The idea of the proof of Lemma~\ref{lem:eps-LDP} is to mimic the definition of $(\cM^n,\ol\dlp^{\eps,n}_i)_{i\in \N }$ in Lemma~\ref{lem:eps-LDP} using a cut-off version of the discrete LDP and show that their transition kernels as Markov processes  are identical in the continuum. Once Lemma~\ref{lem:eps-LDP} is proved this way, the desired ergodicity in Lemma~\ref{emb-prop:LDP} follows from the corresponding ergodicity of the continuum LDP from Theorem~\ref{thm:LDP}.}
We will consider a probability space $(\Omega,\cF,\P)$ satisfying the properties described in Theorem~\ref{thm:LDP}. Let $\bh$ be defined as in~\eqref{eq:fields} and let $\Pd$ be as above \eqref{eq:fields}, so that the $\Pd$-law of $\bh$ is as in Lemmas~\ref{lem:eps-LDP} and~\ref{emb-prop:LDP}.

Fix a site percolation configuration $\omega$ on $\D^\mesh$ with monochromatic blue boundary condition.
Let $\Gamma(\omega)$ be the loop ensemble of $\omega$.
Given $\ell\in \Gamma(\omega)$, by  our convention in Section~\ref{subsec:pre},  $\ell$ is viewed as  an edge path on the triangulation $\D^\mesh$.
Given each edge $e$ in $\ell$, let $e^*$ be its \emph{dual} edge obtained by rotating $e$ around its midpoint by 90 degrees. 
The collection of such dual edges forms an oriented simple loop, where the orientation is such that the red vertex of each edge $e$ is on the left side. 
We call the domain bounded by this simple loop the \notion{region} enclosed by $\ell$.
Given $\ell\in\Gamma(\omega)$, similarly as in Definition~\ref{def:eps}, 
we call the $\mu_{\bh}$-mass of  the region enclosed by $\ell$   the \notion{$\mu_{\bh}$-area of $\ell$}.
Given an inner vertex $v$ of $\D^\mesh$, let $\omega_v$ be the coloring of $\cV(\D^\mesh)$
such that for each $v'\in \cV(\D^\mesh)$,  $\omega_{v}(v')=\omega(v')$ if and only if $v'\neq v$.
Let $\cL_v$ be the symmetric difference between $\Gamma(\omega)$ and $\Gamma(\omega_v)$.
For $\eps>0$, we call $v$ an \notion{$\eps$-pivotal point} of $(\bh,\omega)$ if there are at least three loops in $\cL_v$ with $\mu_{\bh}$-area at least $\eps$.

For $\eps>0$, let  $(\omega^{\eps,\mesh}_t)_{t\ge 0}$ be the following  modification of the discrete LDP $(\omega^\mesh_t)_{t\ge 0}$ on $\D^\mesh$ driven by $e^{\fh/\sqrt6}$: 
when the clock at an inner vertex $v$ rings at time $t$, the color of $v$ is flipped if and only if $v$ is  an $\eps$-pivotal point of  $(\bh,\omega^{\eps,\mesh}_{t^{-}})$.
Note that $(\omega^{\eps,\mesh}_t)_{t\ge 0}$ is defined similarly as $(\ol\omega^{\eps,n}_t)_{t\ge 0}$ in Lemmas~\ref{lem:stationary}, i.e., by rejecting updates of vertices which are not $\eps$-pivotal. 
\xncomment{Let $\Gamma^{\eps,\mesh}_t=\Gamma(\omega^{\eps,\mesh}_t)$  for each $t\geq 0$. Then $(\Gamma^{\eps,\mesh}_t)_{t\ge 0}$ is the lattice analog of
	$(\cM^n,\ol\dlp^{\eps,n}_t)_{t\ge 0}$. Our next lemma shows that  $(\Gamma^{\eps,\mesh}_t)_{t\ge 0}$  converges in law for each fixed $\eps>0$.}

\begin{lemma}\label{prop:eps-LDP-law}
	\xncomment{In the setting of Theorem~\ref{thm:LDP}, for each $\eps>0$ let $(\omega^{\eps,\mesh}_t)_{t\geq 0}$  and  $(\Gamma^{\eps,\mesh}_t)_{t\ge 0}$  be defined as above.}
	Then there exists a process $(\Gamma^{\eps}_t)_{t\ge 0}$ coupled with $\fh$ such that
	$(\fh, \Gamma^{\eps,\mesh}_t)_{t\ge 0}$ converge in law to $(\fh, \Gamma^{\eps}_t)_{t\ge 0}$ as $\mesh\rta 0$  in the Skorokhod topology as c\`adl\`ag  processes taking values in $H^{-1}(\D)\times\cL(\D)$.
	Conditioned on $\fh$, $(\Gamma^\eps_t)_{t\ge 0}$ is a stationary Markov process, where the conditional law of $\Gamma^\eps_0$ is that of a $\CLE_6$ on $\D$. 
	Moreover, almost surely $(\Gamma^\eps_t)_{t\ge 0}$ either stays constant or has infinitely many jumps. In the latter case, it has finitely many jumps in any finite interval.
\end{lemma}

\xncomment{We will prove Lemma~\ref{prop:eps-LDP-law}  in Section~\ref{subsec:Markov}.} 
The proof will also provide a recipe for sampling  $(\fh,\Gamma^\eps_t)_{t\ge 0}$ without referring to the lattice approximation. 
Before describing it in Lemma~\ref{lem:limit-law}, we give a purely continuum description of the limiting pivotal measures involved.
Given a subset $S$ of $\mesh\D^\mesh$ and a measure $\mu$ on $\C$, by \emph{$\mu$ restricted to $S$}, 
we mean $\mu$ restricted to the union of hexagons in the dual lattice whose vertex is in $S$.   
\xncomment{We also recall the definition of $\eps$-pivotal points (Definition~\ref{def:eps}) for a $\CLE_6$ coupled with an independent Liouville field.} 
\begin{lemma}\label{lem:pmeasure}
	There exists a constant $c'>0$ such that the following holds. In the setting of Theorem~\ref{thm:LDP}, for each $\eps>0$, $\arm^\mesh_4(\mesh,1)^{-1}$ times  Lebesgue measure restricted to the set of $\eps$-pivotal points of $(\bh,\omega^\mesh_0)$ converge to  a measure $\Mink_\eps$ in probability. Moreover,
	there exists a random set  $\cA\subset \D$ measurable with respect to $(\bh,\Gamma_0)$ such that  $\Mink_\eps=(c'\Mink_\cA)|_{\cP_\eps}$, where $\Mink_\cA$ is the $3/4$-occupation measure of $\cA$ and $\cP_\eps$ is the $\eps$-pivotal points of $(\bh,\Gamma_0)$. 
\end{lemma}
We will prove Lemma~\ref{lem:pmeasure} in  Section~\ref{subsub:pivm}, where we will see  that \xncomment{the set} $\cA$ can be chosen to be the \xncomment{so-called} $\rho$-important points (Definition~\ref{def:rho-imp}) of $\Gamma_0$ for small enough $\rho$. In fact, $\Mink_\eps$ is $c'$ times the $3/4$-occupation measure of $\cP_\eps$ but we omit the proof of this fact since we do not need it. %(we do not need this fact  so we omit its proof.)

\xncomment{Given Lemma~\ref{lem:pmeasure}, let \(\cM^\eps_{\bh,\Gamma_0} :=(c'e^{\bh/\sqrt 6}\Mink_\cA)|_{\cP_\eps}\).
	Since $\cA$ is measurable with respect to $(\bh, \Gamma_0)$, so is the measure $\cM^\eps_{\bh,\Gamma_0}$.
	Recall the measure $\pivm^\eps$ from Proposition~\ref{prop:quantum-Mink}, where the law of $(\bh,\Gamma)$ is the same as that of $(\bh, \Gamma_0)$ considered here, and  the precise definition  of $\pivm^\eps$ was postponed  to this section. In fact, we will simply define $\pivm^\eps$  by applying the measurable function on $H^{-1}(\D)\times\cL(\D)$ defining  $(\bh,\Gamma_0)$ to $(\bh ,\Gamma)$ instead. This will be made  precise as Definition~\ref{def:pivm} in Section~\ref{sec:prop:quantum-Mink} after the set $\cA$ is described more concretely. Given this definition,  the content of Proposition~\ref{prop:quantum-Mink} is that 
	$\nu^\eps_{\bh,\Gamma}=\constp\pivm^\eps$ a.s.\ where $\pivm^\eps$ is the quantum natural measure on the $\eps$-pivotal points of $(\bh, \Gamma)$
	constructed from the mating-of-trees theory. We will  conclude the proof of  Proposition~\ref{prop:quantum-Mink} in Section~\ref{sec:prop:quantum-Mink}.
	We now use $\cM^\eps_{\bh,\Gamma_0}$ to describe the Markov process $(\Gamma^\eps_t)_{t \ge 0}$.}

\begin{lemma}\label{lem:limit-law}
	The law of $(\fh,\Gamma^\eps_t)_{t \ge 0}$  in Lemma~\ref{prop:eps-LDP-law} can be described as follows.
	Conditioning on $(\D,\fh,\Gamma^\eps_0)$,  an exponential clock rings with rate $(\xi_\fh(\p \D))^{1/2}\cM^\eps_{\bh,\Gamma_0}(\D)$.
	Here we make the convention that an exponential clock with rate 0 never rings. Once the clock rings, sample an $\eps$-pivotal point
	$\bz$ from $\cM^\eps_{\bh,\Gamma_0}$.  The process jumps to the loop ensemble obtained from $\Gamma^\eps_0$ (i.e.\ $\Gamma_0$) by flipping the color at $\bz$. (Recall the notion of color flipping for $\CLE_6$  above Definition~\ref{def:eps}.)
	The remaining jumps in the process, are sampled iteratively.  
\end{lemma}
Since $(\Gamma^\eps_t)_{t\ge 0}$ is stationary and  has finitely many jumps in any finite interval by Lemma~\ref{prop:eps-LDP-law}, 
\xncomment{modulo a probability zero event,  $\cM^\eps_{\bh,\Gamma_t}$  is well-defined simultaneously for all $(\bh, \Gamma^\eps_t)$.} 
Therefore the iterative sampling in Lemma~\ref{lem:limit-law} makes sense.

\begin{comment}
[Why the detail of the coupling does not matter.] Because both in the second statement of Proposition~\ref{prop:eps-LDP} and in Lemma~\ref{lem:limit-law}, it only concerns the law of $(\bh, \Gamma^\eps_t)$. Therefore, it only depends on  the law of the weak limit of $(\fh, \Gamma^\eps_t)_{t\ge 0}$. This allows us to say without looking at the detail of the coupling what is the law of $(\fh, \Gamma^\eps_t)_{t\ge 0}$. Then in the continuum,once we have the explicit description of the law, how it changes under reweighting of $\fh$ is clear.
\end{comment}

Recall the constants
$c_{\op p}$ in \eqref{eq:piv-conv} (see also Proposition~\ref{prop:flip} below) where  \xncomment{ $c_{\op p}\pivm^\eps$ 
	describes the scaling limit of  the discrete pivotal measure on random planar maps}. Recall $\constp$ in Proposition~\ref{prop:quantum-Mink} 
\xncomment{such that $\nu^\eps_{\bh,\Gamma}=\constp\pivm^\eps$ a.s.}
In the setting of Lemmas~\ref{prop:eps-LDP-law} and~\ref{lem:limit-law}, let 
\begin{equation}\label{eq:def-cz}
\ol \Gamma^\eps_t \defeq\Gamma^{\eps}_{\constp c_{\op p}t\xi_\fh(\p\D)^{-1/2}} \quad \textrm{for each }t\ge 0.
\end{equation}
\xncomment{Since $\nu^\eps_{\bh,\Gamma}=\constp\pivm^\eps$,} by  Lemma~\ref{lem:limit-law},
conditioning on $(\fh,\Gamma_0)$, the first time at which the process  
$(\ol \Gamma^\eps_t)_{t\ge 0}$ jumps has the law of an exponential random variable
with rate $c_{\op p}\nu^\eps_{\bh,\Gamma_0}(\D)$, where $\nu^\eps_{\bh,\Gamma_0}$ is as $\nu^\eps_{\bh,\Gamma}$ in Proposition~\ref{prop:quantum-Mink} with $\Gamma_0$ in place of $\Gamma$.

\xncomment{Recall that $\Pd$ is the probability measure obtained from a reweighing of the probability measure $\P$ in Theorem~\ref{thm:LDP}
as above \eqref{eq:fields}, so that the $\Pd$-law of $\bh$ is as in Lemmas~\ref{lem:eps-LDP} and~\ref{emb-prop:LDP}.	
} Let $(Y^\eps_{t})_{t\ge 0}$ be a sample of  $\big(\D, \bh,\ol \Gamma^\eps_t \big)_{t\ge 0}$ according to its $\Pd$-law, 
where $\big(\D, \bh,\ol \Gamma^\eps_t \big)$ is viewed as a random variable in $\BM^{\GHPUL}$ as in Remark~\ref{rmk:CLE}.
\xncomment{We will prove Lemma~\ref{lem:eps-LDP} by showing that $(Y^\eps_{t})_{t\ge 0}$ is the scaling  limit of 
	$(\cM^n,\ol\dlp^{\eps,n}_t)_{t\ge 0 }$. The  following lemma is the only input from random planar maps that  we need for this proof.}
\begin{lemma}\label{lem:two-jump} 
	Fix $\eps>0$. Let $S^n=(S^n_t)_{t\geq 0}$ be the Markov process $(\cM^n,\ol\dlp^{\eps,n}_t)_{t\geq 0}$ in Lemma~\ref{lem:eps-LDP} and let $(Y^\eps_t)_{t\geq 0}$ be as above. 
	For $i\in\N$, 
	let $\tau^n_i$ and $\tau_i$ be the $i$th time that $S^n_t$  and $Y^\eps_t$, respectively, jump. If no jump occurs we set all the jumping times to be $\infty$.
	Then $(S^n_{\tau^n_1},S^n_{\tau^n_2},\tau^n_1,\tau^n_2)$ and the event $\{\tau^n_1<\infty\}$ jointly converge  in law to  $(Y^\eps_{\tau_1},Y^\eps_{\tau_2},\tau_1,\tau_2)$ and $\{\tau_1<\infty\}$.
\end{lemma}

We postpone the proof  of Lemma~\ref{lem:two-jump} to Section~\ref{subsec:Markov} and proceed to the proof of Lemma~\ref{lem:eps-LDP}.
\begin{proof}[Proof of Lemma~\ref{lem:eps-LDP}]
	Suppose we are in the setting of Lemma~\ref{lem:two-jump}. 
	By Lemma~\ref{lem:two-jump},  $S^n|_{[0,\tau^n_2)}$ converges to $Y^\eps|_{[0,\tau_2)}$ in the Skorokhod topology.  
	Given $s>0$, let $\tau^{s,n}_i$ be defined in the same way as $\tau^n_i$ with $(S^n_t)_{t\ge 0}$ replaced by $(S^{s,n}_t)_{t\ge 0}:=(S^n_{t+s})_{t\ge 0}$. 
	Let $\Q_+$ be the set of positive rationals.
	Then at least along a subsequence of $\N$, there is a coupling of $(S^n)_{n\in\N}$ and a family of processes $\{(Y^{\eps,s}_t)_{t\ge 0}:s\in \Q_+\}$ such that for each $s\in\Q_+$,  it holds that $S^{s,n}|_{[0,\tau^{s,n}_2)}$ converges to $Y^{\eps,s}|_{[0,\tau^{s}_2)}$ a.s.\ in the Skorokhod topology, where each $(Y^{\eps,s}_t)_{t\ge 0}$ has the same law as $(Y^{\eps}_t)_{t\ge 0}$ above. Given a rational $s\in (\tau_1,\tau_2)$,  for $n$ large enough  $\tau^{s,n}_i+s=\tau^{n}_{i+1}$ for all $i\in\N$. In particular, $S^n|_{[s,\tau^{n}_3)}= S^{s,n}|_{[0,\tau^{s,n}_2)}$.
	This implies that in our coupling along the chosen subsequence $S^n|_{[0,\tau^{n}_3)}$ converges almost surely in the Skorokhod topology and the law of the limiting object is given by $Y^\eps|_{[0,\tau_3)}$.
	Therefore $S^n|_{[0,\tau^{n}_3)}$ converges  in law to $Y^\eps|_{[0,\tau_3)}$ in  the Skorokhod topology, without passing to a subsequence.
	By induction, the same convergence  holds with $\tau^n_3,\tau_3$ replaced by $\tau^n_i,\tau_i$ for any $i\in\{4,5,\dots \}$.  
	By Lemma~\ref{prop:eps-LDP-law}, $\lim_{i\to\infty}\tau_i=\infty$ a.s.  Therefore $(S^n_t)_{t\geq 0}$ converges to $(Y^\eps_t)_{t\geq 0}$ in the Skorokhod topology. 
	
	Since every c\`adl\`ag function has countably many discontinuous points and $(Y^\eps_t)_{t\geq 0}$ is stationary, for each fixed $t\ge 0$, $Y^\eps$ is almost surely continuous at $t$. This gives Lemma~\ref{lem:eps-LDP}.\qedhere
\end{proof} 

\xncomment{Although the convergence in Lemma~\ref{prop:eps-LDP-law} is only  in law,  the following proposition, which we will prove in Section~\ref{subsec:stability},  upgrades it to convergence in probability. This will be important to the proof of Lemma~\ref{emb-prop:LDP}.}
\begin{proposition}\label{prop:eps-LDP-proba}
	There exists a probability space $(\Omega,\cF,\P)$ satisfying Theorem~\ref{thm:LDP}  and Lemma~\ref{prop:eps-LDP-law} 
	such that for each $\eps>0$, $(\Gamma^{\eps,\mesh}_t)_{t\ge 0}$ converge in probability  as $\mesh\to0$.
\end{proposition}
For $\mesh>0$, let $\omega^\mesh$ be the Bernoulli-$\frac12$ site percolation on $\D^\mesh$ with monochromatic blue boundary condition.
Let $\Gamma^\mesh\defeq\Gamma(\omega^\mesh)$. As explained in \cite{gps-pivotal},  $\omega^\mesh$ and $\Gamma^\mesh$ jointly converge in law.
Suppose $(\omega,\Gamma)$ is a sample from the limiting joint law. Then the quad crossing configuration  $\omega$ is a.s.\ determined by $\Gamma$ \cite{camia-newman-sle6,gps-pivotal}.
In Section~\ref{subsec:quad}  we prove the inverse measurability statement conjectured in \cite{ss-planar-perc}. 
\begin{theorem}\label{thm:quad}
	$\Gamma$ is almost surely determined by $\omega$.
\end{theorem}

From now on we work on the probability space $(\Omega,\cF,\P)$ in  Proposition~\ref{prop:eps-LDP-proba} and 
let $(\Gamma^{\eps}_t)_{t\ge 0}$ be the in-probability limit of $(\Gamma^{\eps,\mesh}_t)_{t\ge 0}$  as $\mesh\to0$. This way, $(\Gamma^\eps_t)_{t\ge 0}$ for different $\eps$'s in Lemma~\ref{prop:eps-LDP-law} are coupled together.
To prove Lemma~\ref{emb-prop:LDP},  we would like to take the $\eps\to0$ limit of $(\Gamma^\eps_t)_{t\ge 0}$. 
However, this convergence is hard to establish directly in $\cL(\D)$.
Theorem~\ref{thm:quad} allows us to reduce  Lemma~\ref{emb-prop:LDP} to the following proposition on quad-crossing elements.
\begin{proposition}\label{prop:quad-conv}
	For each $\eps>0$ and $t\ge 0$, let $\omega^\eps_t \defeq \omega(\Gamma^\eps_t)$ be the element of $\cH(\D)$ corresponding to $\Gamma^\eps_t$.
	\xncomment{Recall $(\omega_t)_{t\ge 0}$ in Theorem~\ref{thm:LDP}.}
	Then for each $r\in(0,1)$, $\lim_{\eps\to0} (\omega^\eps_t|_{r\D})_{t\ge 0}=(\omega_t|_{r\D})_{t\ge 0}$ 
	in probability  in the Skorokhod topology as c\`adl\`ag processes in $\cH(r\D)$, where $\omega^\eps_t|_{r\D}$ is $\omega^\eps_t$ restricted to $\cQ_{r\D}$.
\end{proposition}
The proof of Proposition~\ref{prop:quad-conv} will be given in Section~\ref{subsec:stability}.

\begin{proof}[Proof of Lemma~\ref{emb-prop:LDP}]
	Recall $(\ol \Gamma^\eps_t)_{t\geq 0}$  in~\eqref{eq:def-cz} defined in terms of $(\Gamma^\eps_t)_{t\ge 0}$.
	\xncomment{Since the $\P^{\op d}$-law of $(\D,\bh,\ol \Gamma^\eps_i)_{i\in \N}$ equals the law of $(Y^\eps_i)_{i\ge 0}$ in Lemma~\ref{emb-prop:LDP},}
	it suffices to show that under $\Pd$, as $\eps\to 0$, $(\ol \Gamma^\eps_i)_{i\in \N}$ converge to an ergodic sequence.
	
	For each $t\ge 0$, let $\ol\omega^\eps_t\defeq \omega^\eps_{\constp c_{\op p}t\xi_\fh(\p\D)^{-1/2}}$ be the element in $\cH(\D)$ corresponding to  $\ol\Gamma^\eps_t$.
	Let $\ol\omega_t\defeq \omega_{\constp c_{\op p}t\xi_\fh(\p\D)^{-1/2}}$.
	Restricted to $r\D$, both $(\ol\omega^\eps_t)_{t\ge 0}$
	and $(\ol\omega_t)_{t\ge 0}$ are stationary c\`adl\`ag processes. As in the last paragraph in the proof Lemma~\ref{lem:eps-LDP}, for each fixed $t\ge 0$, 
	Proposition~\ref{prop:quad-conv} implies that $\lim_{\eps\to 0}\ol\omega^\eps_t|_{r\D}=\ol\omega_t|_{r\D}$ in probability. Varying $r$ we see that $\lim_{\eps\to 0}\ol\omega^\eps_t=\ol\omega_t$ in probability.
	
	In light of Theorem~\ref{thm:quad}, for each fixed $t\ge 0$, $\ol\omega_t$ a.s.\ determines an instance of $\CLE_6$ on $\D$, which we denote by $\ol\Gamma_t$.
	Since  $\lim_{\eps\to0}(\ol \Gamma^\eps_t, \ol\omega^\eps_t)=(\ol\Gamma,\ol\omega)$ in law,
	Theorem~\ref{thm:quad} implies that $\lim_{\eps\to0}\ol\Gamma^\eps_t=\ol\Gamma_t$ in probability under $\P$. 
	\xncomment{Here we use again the  measure theoretic fact Lemma~\ref{lem:law-prob} which upgrades joint converge in law to convergence in probability given measurability.}
	
	By absolutely continuity,  $\lim_{\eps\to0}\ol\Gamma^\eps_t=\ol\Gamma_t$ in probability under $\P^{\op d}$.
	By~\eqref{eq:generator}, the mixing property for $(\omega_t)_{t\ge 0}$ in Theorem~\ref{thm:LDP} also holds for $\ol\omega_t$,
	under both $\P$ and $\P^{\op d}$. In particular $(\ol\omega_i)_{i\in\N}$ is ergodic under $\P^{\op d}$. By Theorem~\ref{thm:quad}, $(\ol \Gamma_i)_{i\in \N}$ is ergodic under $\P^{\op d}$ as well.
\end{proof} 

In the rest of Section~\ref{sec:GPS}, we first
prove Proposition~\ref{prop:quantum-Mink}, Lemma~\ref{lem:pmeasure}, and Theorem~\ref{thm:quad},
and provide tools on percolation without dynamics in Sections~\ref{subsec:pivot} to \ref{subsec:quad}.
Then in Sections~\ref{subsec:Markov} and~\ref{subsec:stability} we study the various dynamics considered in Section~\ref{subsec:lemma} and prove  Lemmas~\ref{prop:eps-LDP-law}, \ref{lem:limit-law}, and~\ref{lem:two-jump} and Propositions~\ref{prop:eps-LDP-proba} and~\ref{prop:quad-conv}.

\subsection{Lattice approximation of the pivotal measure}\label{subsec:pivot}
In this section we introduce a cutoff on the set of pivotal points. The cutoff is different from the one we use when defining $\eps$-pivotal points, and we call the set of macroscopic pivotal points for the new cutoff $\rho$-important points.
The concept of $\rho$-important points has also been used in \cite{gps-dynamic,ghss18} (see the beginning of Section~\ref{subsec:stability} for further discussion).
Although lacking a natural connection to random planar maps, this cutoff is more amenable for technical analysis.

Recall that $\D^\mesh$ is a subset of the rescaled triangular lattice $\mesh\Tg$ approximating $\D$. Throughout this subsection $\omega^\mesh$ denotes a sample of Bernoulli-$\frac12$  site percolation on $\D^\mesh$ for $\mesh>0$.  Moreover, $\{\omega^\mesh\}_{\mesh>0}$ are coupled such that the loop-ensembles $\Gamma^\delta\defeq \Gamma(\omega^\mesh)$ converge to a $\CLE_6$ $\Gamma$ in $\cL(\D)$ almost surely (see Theorem~\ref{emb-thm:CLE}). 
We parametrize loops in $ \Gamma$ and $\Gamma^\mesh$  such that when listed in decreasing order according to the (Euclidean) area of the enclosed region, the $k$th  loop converges a.s.\ in the uniform topology for each $k\in \N$.
We enlarge our coupling to include a sample of $\fh$, hence $\bh$, as in Lemma~\ref{lem:g}, which is independent of $\{\omega^\mesh\}_{\mesh>0}$. Let $\nu_\mesh$ be  the renormalized weighted counting measure on $\D^\mesh$ where each vertex $x$  is assigned mass $\mu'_\fh(x)\arm^\mesh_4(\mesh,1)^{-1}$, where we recall from above Theorem~\ref{thm:LDP} that $\mu'_\fh(x)$ is the mass assigned by $e^{\fh/\sqrt{6}}d^2z$ to the hexagon associated with $x$ on the dual lattice of $\D^\mesh$. Note that the law of $\{\omega^\mesh_0\}_{\mesh>0}$ and $\fh$ in Theorem~\ref{thm:LDP}  satisfies the description of the law of $\{\omega^\mesh\}_{\mesh>0}$ and $\fh$ in this section.
\begin{comment}
In this subsection we parametrize loops in $\Gamma$ and $\Gamma^\mesh$. But the only results we used in the next Section are Lemmas~\ref{lem:flip2},~\ref{lem:UI},~\ref{lem:cover2} and~\ref{lem:cover3}, none of which relies on the particular parametrization.
\end{comment}

This subsection is organized as follows.  In Section~\ref{subsub:dis-bich}, we recall  some results from \cite{hlls-pivot} concerning $2$-$\SLE_6$. 
Then we introduce $\rho$-important points and prove its basic properties in Sections~\ref{subsub:rho-important}---\ref{subsubsec:flip-lattice},  and establish its relation with $\eps$-pivotal points in Section~\ref{subsub:mutual}.
Finally, we prove Lemma~\ref{lem:pmeasure} in Section~\ref{subsub:pivm}. 
We encourage the reader to skip the technical proofs in the first reading but keep in mind the definitions  and results for later  applications.

\subsubsection{Percolation interfaces and the discrete analog of 2-$\SLE_6$}
\label{subsub:dis-bich}
In this subsection we recall that for a quad $(Q,a,b,c,d)$, a certain pair of interfaces connecting $a,b,c,d$ pairwise in $Q$ is given by a 2-SLE$_6$ (Lemma \ref{lem:quad-chod0}) and that the points of intersection between these interfaces have a well-defined $3/4$-occupation measure which describes the scaling limit of the corresponding discrete measure for percolation on the triangular lattice (Proposition \ref{prop:dis-2SLE}).

Suppose $U\subset\D$ is a Jordan domain. 
For $x\in \p U$, let $x^\mesh$ be the  edge on $ \p U^\mesh$ closest to $x$ (if there is a tie, choose one arbitrarily).  
We always assume that $\mesh$ is small enough such that $a^\mesh\neq b^\mesh$. 
Let $\sle^{ab}_{U,\mesh}$ be the percolation interface of $\omega^\mesh$  (see the definition below Proposition~\ref{prop:max})
on $(U^\mesh,a^\mesh,b^\mesh)$. Since the triangular lattice is canonically embedded in $\C$, we identify each edge with its dual edge on the hexagonal lattice so that  $\sle^{ab}_{U,\mesh}$ and loops in $\Gamma^\mesh$ are simple curves.

As proved in \cite[Section~5]{camia-newman-sle6}, in our coupling, for a fixed $(U,a,b)$, $\sle^{ab}_{U,\mesh}$ converges in probability to a chordal $\SLE_6$ on $(U,a,b)$ which we denote by $\sle^{ab}_U$. Moreover, $\sle^{ab}_U$ is a.s.\ determined by $\Gamma$ in an explicit way. 
We call $\sle^{ab}_U$ the \notion{interface} of $\Gamma$ on $(U,a,b)$. 
In particular, when $U=\D$, then  $\sle^{ab}_{U}$ is the interface of $\Gamma$ on $(\D, a,b)$ as defined in Lemma~\ref{lem:iteration}.

Given a quad $Q$, we call $Q((0,1)^2)$ the \emph{domain} of $Q$. Abusing notation, we denote the domain of $Q$ by $Q$ for simplicity. 
Let $a,b,c,d$ be $Q(0,0), Q(1,0) , Q(1,1), Q(0,1)$, respectively.

Recall the notions in Lemma~\ref{lem:bdy-quad}. Suppose $Q \subset \ol\cQ_\D$ and $\p Q$ is piecewise smooth. 
Recall the notation $\p_{a,b} D$ in Section~\ref{subsec:notation}. 
Let $E$ be the event that $\sle^{ac}_{Q}$ hits $\p_{b,d} Q $ at a point on $\p_{c,d} Q$. 
As explained in~\cite[Section  1.2]{hlls-pivot}, we have the following. See the left part of Figure \ref{fig:Pijij} for an illustration of the event $E$ and see the right part of Figure \ref{fig:pivot-bichordal2} for an illustration of the $2$-$\SLE_6$ $(\sle^{ad}_{Q},\sle^{cb}_{Q})$.
\begin{lemma}\label{lem:quad-chod0}
	The event $E$ equals $\{\omega(Q)=0\}$ a.s., where $\omega$ is viewed as an element in $\cH_\D$.
	Moreover, the conditional joint law of $(\sle^{ad}_{Q},\sle^{cb}_{Q})$ given $E$ is a $2$-$\SLE_6$ (see Definition~\ref{def:bi-chordal}).
\end{lemma}

Let $\cP^Q=\sle^{ad}_Q\cap\sle^{cb}_Q$ on $E$ and $\cP^Q=\sle^{ab}_Q\cap\sle^{cd}_Q$ on $\neg E$ (i.e.\ the complement of $E$).  
Let $\sle^{ad}_{Q,\mesh}\cap\sle^{cb}_{Q,\mesh}$ be the set of vertices such that $v\in \sle^{ad}_{Q,\mesh}\cap\sle^{cb}_{Q,\mesh}$
if both $\sle^{ad}_{Q,\mesh}$  and $\sle^{cb}_{Q,\mesh}$ traverse an edge with $v$ as an endpoint.
Let $E_\mesh\defeq\{{\omega^\mesh(Q)=0}\}$ and $\cP^Q_\mesh$ be defined in a similar way as $\cP^Q$.
As explained in \cite[Section~1.2]{hlls-pivot}, $\cP^Q_\mesh$ is the set of pivotal points for the crossing event $E_\mesh$.
The following result is extracted from  \cite[Theorem~1.7, Proposition~1.8, Theorem~1.9]{hlls-pivot}.
\begin{proposition}[\cite{hlls-pivot}]
	\label{prop:dis-2SLE}
	The $3/4$-occupation measure  $\Mink_Q$ of $\cP^Q$ exists a.s.
	Moreover, 
	$\alpha_4^\mesh(\mesh,1)^{-1}$ times Lebesgue   measure restricted to  $\cP^Q_\mesh$ (recall this notion from above Lemma~\ref{lem:pmeasure})
	converge to   $c' \Mink_Q$ in probability, where $c'>0$ is a deterministic  constant not depending on $Q$.
\end{proposition}

\subsubsection{$A$-important points and $\rho$-important points}\label{subsub:rho-important}
In this subsection we introduce two sets of pivotal points: $A$-important points and $\rho$-important points. Furthermore, we argue in Lemma \ref{lem:cover} that the $A$-important points can be written as the disjoint union of sets $\cP^{Q}\cap\cB$ for finitely many quads $Q$ and a fixed square $\cB$.

Let $\cB$ be a square of side length $\rho$ for some $\rho>0$ and let $\wt \cB$ be the square of side length $3\rho$ centered around $\cB$. 
Let $A= A_\cB\defeq \wt\cB \setminus(\cB\cup\p \cB)$. For $\cB\cap \D\neq \emptyset$,
let $\Gamma^{A}\defeq\{\ell\in \Gamma: \ell \cap \cB \neq \emptyset\textrm{ and } \ell \cap (\C\setminus \wt \cB)\neq \emptyset  \}$.
By local finiteness of $\CLE_6$ (see Section~\ref{subsec:SLE}),  $\Gamma^A$ contains finitely many loops a.s.
Given $\ell,\ell'\in \Gamma^A$, if $\ell\neq \ell'$,
let $\cP^A(\ell,\ell')\defeq\ell\cap \ell'\cap \cB$, and if $\ell=\ell'$, let 
$$
\cP^A(\ell,\ell')\defeq\{z\in \cB:  \textrm{$\ell\setminus\{z\}$ has two connected components, each of which intersects $\C\setminus\wt \cB$}
\}.
$$ 

Let $\cP^A\defeq \cup_{(\ell,\ell') \in \Gamma^A\times \Gamma^A} \cP^A(\ell,\ell')$. 
A point $z$ is called \notion{$A$-important} for $\Gamma$ if and only if $z\in \cP^A$.
A vertex $v$ on $\cB\cap \D^\mesh$ is called {\bf $A$-important} for $\omega^\mesh$ if and only if there are four  arms from $v$ to $\p \wt\cB$ 
with alternating colors.  Here an  arm refers to a  connected monochromatic path. Let $\cP^A_\mesh$ be the set of $A$-important points for $\omega^\mesh$.

The following lemma  says that  $A$-important points for $\Gamma$ and $\omega^\mesh$ are  covered by finitely many sets of the form $\cP^Q$ and $\cP^Q_\mesh$ from Section~\ref{subsub:dis-bich}, respectively.
\begin{lemma}\label{lem:cover}
	Let $\cB$ be a square of side length $\rho$ for some $\rho>0$ such that $\cB\cap \D\neq \emptyset$ and let $A=A_\cB$.
	Let $\cC$ be a countable dense subset of $\p \wt \cB$. 
	Then, almost surely there exist $\mesh_0>0$ and   quads $Q_1,\dots Q_n$ with domain equal to $\wt \cB\cap \D$ and  marked points contained in $\cC$,  
	such that $\cP^A$ is the disjoint union of $ \{\cP^{Q_i}\cap \cB\}_{1\le i\le n}$, and $\cP^A_\mesh$ is the disjoint union of $\{\cP^{Q_i}_\mesh\cap \cB \}_{1\le i\le n}$ for $\mesh\in (0,\mesh_0)$.
\end{lemma}

\begin{proof}
	For $\ell\in \Gamma^\mesh$, let $\op V(\ell)$ be the set of vertices which are endpoints of edges traversed  by $\ell$. 
	Let  $\Gamma^{\mesh,A}=\{\ell\in \Gamma^\mesh: \op V(\ell) \cap \cB \neq \emptyset\textrm{ and } \op V(\ell) \cap (\C\setminus \wt \cB)\neq \emptyset  \}$. Then $\cP^A_\mesh \subset \cup_{\ell\in\Gamma^{\mesh,A}} \op V(\ell)$.
	We write $\Gamma^A$ and  $\Gamma^{\mesh,A}$ as $\{\ell^1,\dots, \ell^K \}$ and $\{\ell^1_\mesh,\dots, \ell^{K_\mesh}_\mesh \}$, respectively, 
	where loops are listed by decreasing enclosed Euclidean area.
	By the definition of our coupling and the way $\Gamma^A$ and  $\Gamma^{\mesh,A}$ are parametrized, 
	almost surely $\lim_{\mesh\to0 }K_\mesh=K$ and $\lim_{\mesh\to0}\ell^i_\mesh\to \ell^i$   in the uniform topology, for all $1\le i\le K_\mesh$. 
	For each $1\le i\le K$, let $(s^{i,1},t^{i,1}), \dots (s^{i,m^i},t^{i,m^i})$ be the list of intervals of the form $\{(s,t): \ell^i(s),\ell^i(t)\in \p \wt\cB, \ell^i((s,t)) \subset \wt\cB, \ell^i([s,t])\cap \p\cB\neq \emptyset \}$ ordered by increasing left end-point.   
	Since $\ell^i$ is a continuous closed curve, we have $m^i<\infty$ a.s. Let $(s^{i,1}_\mesh,t^{i,1}_\mesh), \dots (s^{i, m^i_\mesh}_\mesh,t^{i,m^i_\mesh}_\mesh)$ be defined similarly for $\Gamma^\mesh$. Define $\ell^{i,j}_\mesh:=\ell^i_\mesh|_{[s^{i,j}_\mesh, t^{i,j}_\mesh]}$ and   $\ell^{i,j}:=\ell^i|_{[s^{i,j}, t^{i,j}]}$.
	Then almost surely $m^i_\mesh\to m^i$ and $\ell^{i,j}_\mesh\rta\ell^{i,j}$   for all $1\le i\le K$ and $1\le j\le m^i$.  
	This convergence follows from the fact that $\SLE_6$ a.s.\ crosses a (fixed) smooth curve upon hitting it. (See e.g. \cite[Lemma~2.2]{hlls-pivot}).

	For $1\le i,i'\le K, 1\le j\le  m^i, 1\le j'\le m^{i'}$ such that $(i,j)\neq(i',j')$, let $\cP^A(i,j;i',j')=\ell^{i,j}([s^{i,j}, t^{i,j}])\cap \ell^{i',j'}([s^{i',j'}, t^{i',j'}])$.
	Let $\op V(\ell^i_\mesh([s^{i,j}_\mesh, t^{i,j}_\mesh]))$ be the vertex set defined as $\op V(\ell)$ above with $\ell$ replaced by $\ell^i_\mesh([s^{i,j}_\mesh, t^{i,j}_\mesh]$, 
	and  let $\cP^A_\mesh(i,j;i',j')=\op V(\ell^i_\mesh([s^{i,j}_\mesh, t^{i,j}_\mesh]))\cap \op V(\ell^{i'}_\mesh([s^{i',j'}_\mesh, t^{i',j'}_\mesh]))$.
	By the no-triple-point property of $\CLE_6$ (see Section~\ref{subsec:SLE}), the sets $\cP^A(i,j;i',j')$ are disjoint.
	Therefore $\cP^A$ is the disjoint union of $\cP^A(i,j;i',j')\cap \cB$
	for all $(i,j)\neq (i',j')$. A similar statement holds for $\cP^A_\mesh$ for small enough $\mesh$. 
	
	\begin{figure}
		\centering
		\includegraphics{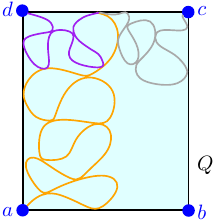}\qquad\qquad
		\includegraphics[scale=1]{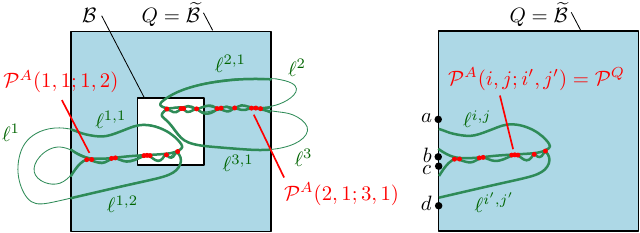}
		\caption{
			{\bf Left}: Illustration of the event $E$. We have that $\sle^{ac}_{Q}$ is the concatenation of the orange curve and the gray curve, while $\sle^{ad}_{Q}$ is the concatenation of the orange curve and the purple curve.
			{\bf Middle}: Illustration of objects defined in the proof of Lemma \ref{lem:cover}. In the case shown, we have  $\wt \cB\subset\D$ so that $Q=\wt \cB$. The annulus $A=A_{\cB}$ is shown in blue. The disk $\D$ is not drawn.
			{\bf Right}: The points $a,b,c,d\in\cC\subset\partial\wt{\cB}$ are chosen such that $\mathcal P^A(i,j;i',j')=\mathcal P^Q$.	
	}
		\label{fig:Pijij}
	\end{figure}

	For $(i,j)\neq (i',j')$ such that $\cP^A(i,j;i',j')\neq\emptyset$, 
	by the parity property of $\CLE_6$ (Lemma \ref{lem:property}), we may assume $\ell^{i,j}(s^{i,j})$, $\ell^{i,j}(t^{i,j})$, $\ell^{i',j'}(s^{i',j'})$, $\ell^{i',j'}(t^{i',j'})$ are in  cyclic  order on $\p \wt \cB$, 
	either counterclockwise or clockwise.  We focus on the former case since the latter case can be treated similarly. 
	Let $Q$ be a quad with domain $\wt\cB\cap \D$ and marked points $a,b,c,d$ in $\cC$ that are to be determined.
	Choose $a,b,c,d\in \cC$ counterclockwise aligned  such that $\p_{\ell^{i,j}(s^{i,j}),\ell^{i,j}(t^{i,j})} Q \subset \p_{a,b}Q$ and  $\p_{\ell^{i',j'}(s^{i',j'}), \ell^{i',j'}(t^{i',j'})} Q\subset \p_{c,d}Q$. For $a$, $b$, $c$, $d$ sufficiently close to $\ell^{i,j}(s^{i,j})$, $\ell^{i,j}(t^{i,j})$, $\ell^{i',j'}(s^{i',j'})$, $\ell^{i',j'}(t^{i',j'})$, respectively,  we have $\cP^A(i,j;i',j')=\cP^Q$.
	For small enough $\mesh$, we also have $\cP^A_\mesh(i,j;i',j')=\cP_\mesh^Q$. This concludes the proof.\qedhere
	\begin{comment}
	To see the existence of $Q$, we may consider the $\CLE_6$ in $\wt\cB$ that is the scaling limit of the loop ensemble of $\omega^\mesh|_{\wt\cB}$. Then the $\ell^{i,j}$ must be part of a boundary touching loop in this smaller $\CLE_6$. By looking at this $\CLE_6$, we may find $a,b,c,d$ as desired. 
	\end{comment}
\end{proof} 
We say that $\cB$ is a \emph{square on $\rho\Z^2$} if it is a square of side length $\rho$ such that all four vertices lie on $\rho\Z^2$.
\begin{definition}\label{def:rho-imp}
	For each $\rho>0$  
	let $$\cP^\rho_\mesh\defeq\bigcup_\cB \cP^{A_\cB}_\mesh\textrm{ for each $\delta>0$,} \qquad \textrm{and} \qquad\cP^\rho\defeq\bigcup_\cB \cP^{A_\cB},$$ 
	where the union is over all squares $\cB$ on $\rho\Z^2$ with $\cB\cap \D\neq \emptyset$. 
	Points  in $\cP^\rho_\mesh$ and $\cP^\rho$ are called \notion{$\rho$-important points} of $\omega^\mesh$ and $\Gamma$, respectively.
\end{definition}

\subsubsection{Scaling limit of discrete pivotal measures}\label{subsub:scaling-piv} 
We now gather some facts concerning the scaling limit of measures on $\cP^A_\mesh$ and $\cP^\rho_\mesh$. We will see in Proposition \ref{prop:piv-measure} and Lemma \ref{lem:UI} that various measures (both Euclidean and quantum) defined on these points converge to their continuum counterparts defined in terms of the $3/4$-occupation measure. Recall that $\nu_\mesh$ is the measure on $\D^\mesh$ where each vertex $x$  is assigned mass $\mu'_\fh(x)\arm^\mesh_4(\mesh,1)^{-1}$ where  $\mu'_\fh(x)$ is the mass assigned by $e^{\fh/\sqrt{6}}d^2z$ to the hexagon associated with $x$ on the dual lattice of $\D^\mesh$.
\begin{proposition}	\label{prop:piv-measure}
	In the setting of  Lemma~\ref{lem:cover},  the $3/4$-occupation measure of $\cP^A$  exists a.s., which we denote  by $\Mink^A$.
	Let $m^A_\mesh$ be $\alpha_4^\mesh(\alpha,1)^{-1}$ times Lebesgue measure restricted to $\cP^A_\mesh$. 
	Let $\nu^A_\mesh$  be the measure $\nu_\mesh$ restricted to $\cP^A_\mesh$.
	Then $\lim_{\mesh\to 0} m^A_\mesh=c' \Mink^A$ and $\lim_{\mesh\to 0} \nu^A_\mesh=c'e^{\fh/\sqrt6}\,\Mink^A$ in probability in the weak topology, 
	where $c'$ is as in Proposition~\ref{prop:dis-2SLE}.
	If $A\Subset\D$ (i.e.\ $A\cup \p A\subset \D$), then $\lim_{\mesh\to 0} m^A_\mesh(\D)=c'\Mink^A(\D)$  and 
	$\lim_{\mesh\to 0} \nu^A_\mesh(\D)=c'\int_{\D}e^{\fh/\sqrt6}\, \Mink^A$ in $L^2$.
\end{proposition}
\begin{proof}
	We obtain the existence of $\Mink^A$ and the convergence of $m^A_\mesh$ in probability from Proposition~\ref{prop:dis-2SLE} and  Lemma~\ref{lem:cover}.   
	If  $A\Subset\D$, the $L^2$ convergence of $m^A_\mesh(\D)$  follows from the moment bounds of $m^A_\mesh(\D)$ given in \cite[Lemma~4.5]{gps-pivotal}. 
	
	Recall $\fh=\Phi+g$ as in Lemma~\ref{lem:g}. If $g$ were equal 0, then by \cite[Propositions~A.1 and A.2]{ghss18},  $\lim_{\mesh\to 0} \nu^A_\mesh=c'e^{\fh/\sqrt6}\,\Mink^A$ in probability,
	and if $A\Subset\D$ then  $\lim_{\mesh\to 0} \nu^A_\mesh(\D)=c'\int_{\D}e^{\fh/\sqrt6}\, \Mink^A$ in $L^2$. Although $g\neq 0$, Corollary~\ref{cor:g} yields the same conclusion.
\end{proof}

Let $\nu^\rho_\mesh$ be the restriction of $\nu_\mesh$ to $\cP^\rho_\mesh$. The next lemma concerns the scaling limit of $\nu^\rho_\mesh$. Both in the proof of the lemma and later in this section we will use the quasi-multiplicativity of $\alpha^\mesh_4$ (see e.g.\ \cite{smirnov-werner-percolation}), namely that for some constant $c>0$ and $\delta\leq r_1\leq r_2\leq r_3$, 
\eqb
c\alpha^\mesh_4(r_1,r_2)\alpha^\mesh_4(r_2,r_3) 
\leq \alpha^\mesh_4(r_1,r_3)
\leq \alpha^\mesh_4(r_1,r_2)\alpha^\mesh_4(r_2,r_3). 
\label{eq:quasi-mult}
\eqe
\begin{lemma}\label{lem:UI}
	Fix $\rho>0$. 
	The $3/4$-occupation measure of $\cP^\rho$  exists a.s. We denote this measure by $\Mink^\rho$.	
	Then   $\lim_{\delta\rta 0}\nu^\rho_\mesh=c'e^{\fh/\sqrt 6} \Mink^\rho$ in probability, where $c'$ is the constant in Proposition~\ref{prop:dis-2SLE}.
	Moreover,  $\lim_{\mesh\to0}\nu_\mesh^\rho(\D)=c'\int_\D e^{\fh/\sqrt6}\Mink^\rho$ in $L^1$.
\end{lemma}
\begin{proof}
	Since the sets $\cP^{A_\cB}_\mesh\subset\cB$ are disjoint for distinct squares $\cB$ on $\rho\Z^2$, the existence of $\Mink^\rho$ and the convergence in probability in Lemma~\ref{lem:UI} follows from Proposition~\ref{prop:piv-measure}. 
	It remains to  prove the $L^1$ convergence of $\nu_\mesh^\rho(\D)$.	For $k\in\N$, set  $r\defeq 1-e^{-k}/2$.	
	By Proposition~\ref{prop:piv-measure}, for each $k\in\N$ and $\rho>0$, 
	$\lim_{\mesh\to0}\nu_\mesh^\rho(r\D)=c'\int_{r\D} e^{\fh/\sqrt6}\Mink^\rho$ in $L^2$.
	\begin{comment}
	Here also $A_\cB$ could touch the boundary,  we can smaller annulus to cover it.
	\end{comment}
	It suffices to prove
	\begin{equation}\label{eq:first-moment}
	\lim_{k\to \infty}\limsup_{\mesh\to0}\E[ \nu^\rho_\mesh(\D\setminus r\D) ]=0.
	\end{equation} 
	For each  $x\in \D^\mesh$, let $E_x$ be event that $x$ is $\rho$-important. 
	Recall that $\nu_\mesh(x)=\mu'_\fh(x)\arm^\mesh_4(\mesh,1)^{-1}$, where $\mu'_\fh(x)$ is the $\mu'_\fh$-mass of the hexagon corresponding to $x$ 
	in the dual lattice. 
	Therefore 
	\[\E[\nu_\mesh(x)\1_{E_x}]=\P[E_x]\arm^\mesh_4(\mesh,1)^{-1}\E[\mu'_\fh(x)].\] 
	
	For $r_2>r_1>0$, let $\wt\alpha^\mesh_4(r_1,r_2)$ be the probability that Bernoulli-$\frac12$ site percolation on $\bbH^\mesh$ has four alternating arms in the semi-annulus 
	$(r_2\D\cap \bbH)\setminus r_1\D $.  We claim that 
	\[\P[E_x]\le C \alpha^\mesh_4(\delta, 1-|x|) \wt\alpha^\mesh_4(1-|x|,\rho),\]
	where $C$ is a constant not depending on $\delta,r,\rho$. It is sufficient to consider the case $|x|>0.9$ and $\rho\in( 10(1-|x|),0.1 )$  since this implies the general case. Let $E'_x$ denote the event that $x$ has four arms to distance $0.5(1-|x|)$, and let $E''_x$ denote the event that there are four alternating arms in the annulus of radii $3(1-|x|)$ and $0.7\rho$ centered at $x/|x|$. Then $E_x\subset E'_x\cap E''_x$, $\P[E'_x]\leq C'\alpha^\mesh_4(\delta, 1-|x|)$, and $\P[E''_x]\leq C'' \wt\alpha^\mesh_4(1-|x|,\rho)$ for constants $C',C''>0$. The claim follows from this and independence of $E'_x$ and $E''_x$:
		$$
		\P[E_x]\le \P[E'_x]\P[E''_x] \le C \alpha^\mesh_4(\delta, 1-|x|) \wt\alpha^\mesh_4(1-|x|,\rho).
		$$
	
	From here on we use $C_\rho$ to denote a constant only depending on $\rho$ that can vary from  place to place.
	By \cite{smirnov-werner-percolation}, the half-plane four-arm exponent is $10/3$ while the plane alternating  four-arm exponent is $5/4$. By this and quasi-multiplicativity \eqref{eq:quasi-mult}, we have the following for some $C>0$
	\eqbn
	\alpha^\mesh_4(\delta, 1-|x|) \arm^\mesh_4(\mesh,1)^{-1} \leq C\alpha^\mesh_4(1-|x|,1)^{-1}=(1-|x|)^{-5/4+o(1)},
	\eqen
and moreover
$\wt\alpha^\mesh_4(1-|x|,\rho)=(1-|x|)^{10/3+o(1)}$. Therefore
	\begin{equation}\label{eq:UI}
		\begin{split}
			\E[\nu_\mesh(x)\1_{E_x}]
			&\le  C\alpha^\mesh_4(\delta, 1-|x|) \wt\alpha^\mesh_4(1-|x|,\rho)\arm^\mesh_4(\mesh,1)^{-1}\E[\mu'_\fh(x)]
			\leq C_\rho (1-|x|)^2 \E[\mu'_\fh(x)].
		\end{split}
	\end{equation}
	Here we have $(1-|x|)^2$ because $2<\frac{10}3-\frac54$.
	
	Let $\phi:\D\to\cS$ be as in Lemma~\ref{lem:g}, i.e., it is the conformal map from $\D$ to $\cS$ satisfying $\phi(0)=\pi i/2$ and $\phi(1)=+\infty$.
	Let $\wt\cB_n=\phi^{-1}([n,n+1] \times (0,\pi))$ where $\phi$ is as in Lemma~\ref{lem:g}. 
	For $n\ge k$, define
	\[A^+_n\defeq \{z\in \D: \op{Re} z\ge 0, \; 1-|z|\in  (e^{-n-1}/2, e^{-n}/2)  , \phi(z)\in [0,k] \times(0,\pi) \}.\]  
	Recall that $a=Q-\gamma=1/\sqrt{6}$ in Lemma~\ref{lem:g}. 
	Since $e^{(B_{2t}-at)/\sqrt6}$ (with $B$ as in Definition \ref{def:disk}) is a martingale,  the value of $\E[\mu'_\fh(\wt\cB_n)]$ does not depend on $n\in \N$.   
	By  \eqref{eq:UI}, for $n\ge k$ we have 
	\[
	\E[\nu^\rho_\mesh(A^+_n)]\le C_\rho e^{-2n}  \E[\mu'_\fh(A^+_n)]\le C_\rho e^{-2n}  \E[\mu'_\fh(\phi^{-1}([0,k]))]\le C_\rho k e^{-2n} .
	\]
	By the definition of $\phi$  we have $e^{\phi(z)}=i(1+z)/(1-z)$ for each $z\in \D$. 
	Therefore  $1-|z|\le 2e^{-n}$ for all  $n\in\N$ and $z\in \wt \cB_n$.
	Now by \eqref{eq:UI}, $\E[\nu^\rho_\mesh(\wt \cB_n)]\le C_\rho e^{-2n}$ for all $n\ge k$. 
	Since $(\D\setminus r\D) \cap \{z:\op{Re} z\ge 0 \}\subset \cup_{n\ge k} (A^+_n\cup \wt \cB_n)$, \eqref{eq:first-moment} holds with $(\D\setminus r\D) \cap \{z:\op{Re} z\ge 0 \}$ in place of $\D\setminus r\D$. 
	
	For the remaining part of $\D\setminus r\D$, we recall from Definition~\ref{def:disk}  that $(X_{-t})_{t>0}$ has the law of $B_{2t}-at$ conditioned to stay negative, which is stochastically dominated by the unconditional law of $B_{2t}-at$. Therefore \eqref{eq:first-moment} holds with $(\D\setminus r\D) \cap \{z:\op{Re} z<0 \}$ in place of $\D\setminus r\D$.
\end{proof}

\subsubsection{Convergence of loop ensemble after flipping a $\rho$-important point}\label{subsubsec:flip-lattice}
The main result of this subsection is the following lemma, which gives convergence of the loop ensemble after flipping a $\rho$-important point. See Figure~\ref{fig:flip} for an illustration of the proof idea.
\begin{lemma}\label{lem:flip3}
	Let $\rho>0$. Suppose $z^\mesh$ and $z$ are random points such that $z^\mesh\in\cP^\rho_\mesh$, $z\in \cP^\rho$,  and $\lim_{\mesh\to 0}z^\mesh=z$ in probability. 
	Let  $\wh \Gamma^\mesh$ and $\wh \Gamma$ be the loop ensembles obtained after flipping the color of $z^\mesh$ and $z$ for $\Gamma^\mesh$ and $\Gamma$, respectively. 
	Then  $\lim_{\mesh\to 0} \wh \Gamma^\mesh=\wh\Gamma$ in probability  in $\cL(\D)$. 
\end{lemma}
\begin{proof}
	Let $\cB$ be the box on $\rho\Z^2\cap \D$ such that $z\in \cB$, let $A\defeq A_\cB$, and fix a small $r_0>0$. 
	We retain the notations in the proof of Lemma~\ref{lem:cover}, including the parametrizations of loops in $\Gamma^\mesh,\Gamma$.
	Then $z$ must belong to some $\cP^A(i,j;i',j')$. Since $\lim_{\mesh\to0}z^\mesh=z$ a.s., $z^\mesh \in \cP^A_\mesh(i,j;i',j')$ with probability 
	$1-o_\mesh(1)$.
	Here and below the implicit constant in $o_\mesh(1)$ may depend on $\rho$ and $r_0$ but is independent of all other parameters.
	From now on whenever we declare an event $E_\mesh$ to have probability $1-o_\mesh(1)$, we will work on $E_\mesh$ thereafter without explicitly mentioning it.
	Without loss of generality, assume $\omega^\mesh(z^\mesh)$ is blue. %Fix a small $r_0>0$ and 
	Let $B(z^\mesh,r_0)$  be the Euclidean ball of radius $r_0$ centered at $z^\mesh$.
	Let $\ell_\mesh$  be the segment of $\ell^{i,j}_\mesh$ from $s^{i,j}_\mesh$ until the first edge that has  $z^\mesh$ as an endpoint, excluding this edge.
	Let $\ol \ell_\mesh$  be the segment of the time reversal of $\ell^{i,j}_\mesh$ from $t^{i,j}_\mesh$ to  the first edge  that has  $z^\mesh$ as an endpoint, excluding this edge.
	Define $(\ell'_\mesh,\ol \ell_\mesh')$  in the same way as $({\ell}_\mesh,\ol\ell_\mesh)$ with $\ell^{i',j'}_\mesh$ in place of $\ell^{i,j}_\mesh$.
	Since  the alternating five-arm exponent for Bernoulli-$\frac12$ site percolation on $\Tg$ is strictly smaller than the four-arm exponent \cite{smirnov-werner-percolation}, 
	with probability $1-o_\mesh(1)$,
	after the color of $z^\mesh$ is flipped to red, we have that $\ell_\mesh$, an edge path $\ell''_\mesh$ contained in $B(z^\mesh,r_0)$, and $\ol \ell_\mesh'$ form a segment of a loop in $\wh\Gamma^\mesh$.
	The same statement holds for $\ol{\ell}_\mesh$, an edge path $\ol\ell''_\mesh$, and $\ell'_\mesh$.
	The two segments $\ell''_\mesh$ and $\ol\ell''_\mesh$ trace 
	small red clusters of $\omega^\mesh$ in $B(z^\mesh,r_0)$ which have a vertex adjacent to $z^\mesh$  but have no vertex in $\op V(\ell^{i,j}_\mesh) \cup  \op V (\ell^{i',j'}_\mesh)$. See Figure~\ref{fig:flip} for an illustration.
	
	\begin{figure}
		\centering
		\includegraphics[scale=1]{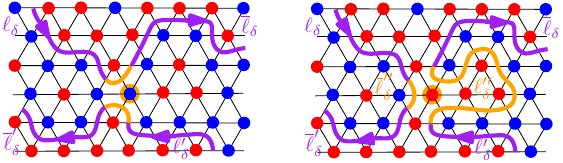}
		\caption{Illustration of the proof of Lemma~\ref{lem:flip3}. The left (resp., right) figure shows two of the percolation interfaces before (resp., after) the color of the 
			$\rho$-important point $z^\delta$ 
			(marked in orange) has been flipped from blue to red. We show that the percolation interfaces after the flip converge by using that the orange paths $\ell''_\delta$ and $\ol\ell''_\delta$ have diameter $o_\delta(1)$ with probability $1-o_\delta(1)$.}\label{fig:flip}
	\end{figure}
	
	Let $\Gamma^\mesh(r_0)=\{\gamma^\mesh\in \Gamma^\mesh: \gamma^\mesh \not\subset B(z^\mesh,r_0), \ell^{i,j}_\mesh\not\subset \gamma^\mesh, \ell^{i',j'}_\mesh \not\subset \gamma^\mesh\}$. 
	By the no-triple-point property (see Section~\ref{subsec:SLE})  of $\CLE_6$, with probability $1-o_\mesh(1)$, 
	$z^\mesh\notin \op V(\gamma^\mesh)$ for any loop $\gamma^\mesh\in\Gamma^\mesh(r_0)$.
	Therefore $\Gamma^\mesh(r_0) \subset \wh\Gamma^\mesh$.
	On the other hand, with probability $1-o_\mesh(1)$,
	$\ell^{i,j}_{\delta}([s^{i,j}_\mesh,t^{i,j}_\mesh])\setminus (\ell_\mesh\cup\ol{\ell}_\mesh)$ 
	and $\ell^{i',j'}_{\delta}([s^{i',j'}_\mesh,t^{i',j'}_\mesh])\setminus (\ell'_\mesh\cup\ol \ell_\mesh')$ are contained in $B(z^\mesh,r_0)$; this follows by symmetry in red and blue and the exact same argument as we used above to argue that $\ell''_\delta,\ol\ell''_\delta\subset B(z^\delta,r_0)$ with probability $1-o_\delta(1)$. By the convergence $\Gamma^\delta\rta\Gamma$ a.s., the segments $\ell_\mesh$, $\ol{\ell}_\mesh$, $\ell'_\mesh$ and $\ol{\ell}'_\mesh$ converge a.s., and we denote the limits by
	%Let
	$\ell$, $\ol{\ell}$, $\ell'$, and $\ol{\ell}'$, respectively. %be the $\mesh\to 0$ limit  of  $\ell_\mesh$, $\ol{\ell}_\mesh$, $\ell'_\mesh$ and $\ol{\ell}'_\mesh$.
	In the continuum, the loop ensemble $\wh \Gamma$ is obtained from $\Gamma$ by concatenating  $\ell$ with $\ol\ell'$, and $\ol{\ell}$ with $\ell'$, while keeping other loops unchanged. 
	Therefore, there is vanishing function $o_{r_0}(1)$ such that for any fixed $r_0>0$, with probability $1-o_\mesh(1)$, $\BB d^{\op L}_d(\wh\Gamma,\wh \Gamma^\mesh)\le o_{r_0}(1)$.
	This concludes the proof.
\end{proof}  

\begin{remark}\label{rmk:flip3}
	Lemma~\ref{lem:flip3} remains true if the assumption that $\Gamma^\mesh\to \Gamma$ almost surely is weakened to convergence in probability. 
	This observation will be used in Section~\ref{subsec:stability}.
\end{remark}

\subsubsection{Mutual inclusion of $\eps$-pivotal points and $\rho$-important points}\label{subsub:mutual}
Recall $\bh$ from \eqref{eq:fields} and the notion of $\eps$-pivotal point for $(\bh,\Gamma)$ and $(\bh, \omega^\mesh)$  in Section~\ref{subsec:lemma}.
The next three lemmas give certain mutual inclusion relations of $\eps$-pivotal points and $\rho$-important points, allowing us to study the former through the latter.
\begin{lemma}\label{lem:cover2}
	Fix $\eps>0$.
	There almost surely exists $b>0$ such that $\mu_{\bh}(\cB)<\eps$ for all squares with side length less than $b$. Let $b^\eps$ be the supremum of all such $b$'s and set $\rho^\eps=0.01 b^\eps$.  
	Then each $\eps$-pivotal point  of $(\bh,\Gamma)$ (resp., $(\bh, \omega^\mesh)$)  is $\rho$-important for $\Gamma$ (resp., $\Gamma^\mesh$) for $\rho\in(0,\rho^\eps)$ and $\mesh\in (0,\rho)$. 
\end{lemma} 
\begin{proof}
	Since $\mu_{\bh}$ is a.s.\ non-atomic, we obtain the existence of $b$ with the desired property.
	Given $\rho\in (0,\rho^\eps)$ and an $\eps$-pivotal point $z$ for $(\bh,\Gamma)$, 
	let $\cB$ be a box of $\rho\Z^2$ such that $z\in \cB$. Set $A\defeq A_\cB$.
	Recall $\Gamma^A$ in the proof of Lemma~\ref{lem:cover}.
	If $z\in \ell\cap \ell'$ for some  distinct loops $\ell,\ell'\in\Gamma$, then we must have $\ell,\ell' \in \Gamma^{A}$. Similarly, if $z$ is a double point on some $\ell\in\Gamma$, then  the two new loops $\ell'$ and $\ell''$ which we get after flipping the color of $z$ must intersect both boundaries of $A$. Therefore $z$ is $\rho$-important for $\Gamma$. The statement for $\omega^\mesh$ follows from the same argument.
\end{proof}
\begin{lemma}\label{lem:rho-is-eps}
	Fix $\rho>0$. There almost surely exists  $\eps'>0$ such that $\cP^\rho\subset\cP_{\eps'}$. 
\end{lemma}
\begin{proof}
	Recall  the setting of Lemma~\ref{lem:cover} and its proof. %It  suffices to prove that $\cP^A(i,j;i',j')\subset \cP_{\eps'}$ for small enough $\eps'$.	
	It suffices to prove that for sufficiently small $\eps'$, $\cP^A(i,j;i',j')\subset \cP_{\eps'}$ for all $i,j,i',j'$.
	%If $\cP^A(i,j;i',j')\neq \emptyset$, 
	%consider 
	Suppose $i,j,i',j'$ are such that $\cP^A(i,j;i',j')\neq \emptyset$. Consider
	the segment of $\ell^{i,j}$ starting from $\ell^{i,j}(s^{i,j})$ until the first time  when it hits $\ell^{i',j'}$. Then the complement of this segment in $\D$ contains a countable collection of components with clockwise boundary orientation. Let $\eps_1(i,j;i'j')$ be the largest $\mu_{\bh}$-area of components in this collection.  Let $\eps_2(i,j;i'j')$ be similarly defined with counterclockwise in place of clockwise. 
	We  define $\eps_3(i,j;i'j'),\eps_4(i,j;i'j')$ in the same way as $\eps_1(i,j;i'j'),\eps_2(i,j;i'j')$ where we trace $\ell^{i,j}$ in the reverse direction until it hits $\ell^{i',j'}$.	
	Define $\eps_k(i,j;i'j')$ with $k=5,6,7,8$ in the same way where the roles of $\ell^{i,j}$ and $\ell^{i',j'}$ are swapped. 
	Let \(E_\eps=\{\eps_k(i,j;i'j') > \eps \textrm{ for }k=1,\dots,8,  \textrm{ if }\cP^A(i,j;i'j')\neq \emptyset\}.\)
	On the event $E_\eps$, if $v$ is an $A$-important point for $\Gamma$, there exists a $\cP^A(i,j;i'j') $ containing $v$. For each loop $\ell\in \cL_v$, 
	one of the eight types of bubbles in the definition of $\eps_k(i,j;i'j')$ must be contained in the region enclosed by $\ell$. Therefore $v\in\cP_\eps$.
	%Since $\eps_k(i,j;j'j')>0$ a.s.\ for all $1\le k\le 8$ and all $i,j,i',j'$ and $\cP^A(i,j;i',j')\neq \emptyset$, this concludes the proof.
	Since $\eps_k(i,j;j'j')>0$ a.s.\ for all $k,i,j,i',j'$ such that $\cP^A(i,j;i',j')\neq \emptyset$, and since $\cP^A(i,j;i',j')= \emptyset$ except for finitely many $i,j,i',j'$, this concludes the proof.
\end{proof}	
The next lemma says that for a vertex $v$ away from the boundary and for $s$ sufficiently small, conditioning on $v$ being $s$-important, it is $\eps$-pivotal with probability $1-o_\eps(1)$.
\begin{lemma}\label{lem:cover3}
	Let $r\in (0,1)$. For each $s\in(0,0.1(1-r))$ and $\zeta'\in(0,1)$, there exists $\eps>0$ and $\delta_0>0$ only depending on $s$, $\zeta'$, $r$ such that
	\[
	\P\left[\textrm{$v$ is not $\eps$-pivotal for $(\bh,\omega^\mesh)$} \mid \textrm{$v$ is $s$-important for $\omega^\mesh$}\right]\le \zeta', \,\, \forall\;\delta\in(0,\delta_0)\textrm{ and } v\in \D^\mesh\cap r\D.
	\]
\end{lemma}
\begin{proof}
	For $\bz\in r\D$, let $\cB_\bz$ be the square of side length $s$ centered at $\bz$ and set the annulus $A=A_{\cB_\bz}$. Consider the set of pairs $(\omega,v)$ where $\omega$ is a site percolation configuration on $\D^\mesh$ with monochromatic boundary condition and $v$ is an $A$-important point.
	Suppose $( \omega^\mesh,\bv^\mesh)$ is uniformly chosen  from this set. Here we use the same symbol $\omega^\mesh$ as in Lemma~\ref{lem:cover3}, although the law of $\omega^\mesh$ here is not uniform.
	One way to sample $(\omega^\mesh,\bv^\mesh)$ is the following. 
	First sample a Bernoulli-$\frac12$ site percolation $\omega^\mesh$ on $\D^\mesh$ with monochromatic boundary condition. Then reweight the law of $\omega^\mesh$ by the number of $A$-important points.
	Finally, conditioning on $\omega^\mesh$, sample the point $\bv^\mesh$ according to the uniform measure on $A$-important points of $\omega^\mesh$.

	Let $\Gamma^\mesh=\Gamma(\omega^\mesh)$ be the associated loop ensemble. 
	By Proposition~\ref{prop:piv-measure},  $(\Gamma^\mesh,\bv^\mesh)$ jointly converge to a pair $(\Gamma,\bv)$ that can be sampled as follows. First sample a $\CLE_6$ $\Gamma$ in $\D$.  Define $\cP^A$ as in Lemma~\ref{lem:cover}. Then reweight the law of $\bd\Gamma$ by $\Mink^A(\D)$,  where $\Mink^A$  is the $3/4$-occupation measure of  $\cP^A$.
	Note that this is well-defined since the measure we reweight by has finite expectation by Proposition~\ref{prop:piv-measure}.
	Finally, conditioning on $\Gamma$, sample the point $\bv$ according to $\Mink^A$. 
	By the Skorokhod representation theorem we may assume that the convergence above holds almost surely. We enlarge the sample space by considering an independent sample of the field $\bh$ from~\eqref{eq:fields}. Denote this probability measure by $\wh \P$. 
	
	Recall $E_\eps$ in the proof of Lemma~\ref{lem:rho-is-eps}. From that proof, we see that on the event $E_\eps$, each $A$-important point is $\eps$-pivotal for $(\bh,\Gamma)$.
	Moreover, $\lim_{\eps\to 0}\wh\P[E_\eps]=1$.
	Let $E^\mesh_\eps$ be the exact analog of $E_\eps$ defined for $\omega^\mesh$. By the scaling limit result, for each $\zeta>0$, there exist  $\eps>0$ and $\delta_0>0$ small enough only depending on $s,\zeta$ such that for each $\delta\in(0,\delta_0)$, 
	on the event $E^\mesh_\eps$ every $A$-important point for $\omega^\mesh$ is $\eps$-pivotal for $(\bh,\omega^\mesh)$, and moreover,
	\begin{equation}\label{eq:bv}
	\wh \P[E^\mesh_\eps]>1-\zeta,  \qquad \forall \;\delta\in(0,\delta_0).
	\end{equation}

	Now let us sample $(\omega^\mesh,\bv^\mesh)$ in another way. We first sample $\bv^\mesh$ according to its marginal law. Then we sample the Bernoulli-$\frac12$ site percolation $\omega^\mesh$ on $\D^\mesh$ conditioned on the event $F^\mesh_s$ that $\bv^\mesh$ is $A$-important. 
	Let $\neg E^\mesh_\eps$ be the complement of $E^\mesh_\eps$. For the choice of $\zeta,\eps$ in \eqref{eq:bv}, 
	\begin{equation}\label{eq:cond}
	\P\left[ \neg E^\mesh_\eps \mid F^\mesh_s \right]
	=\wh\P[\neg E^\delta_\eps]
	\le \zeta \qquad \forall \;\delta\in(0,\delta_0).
	\end{equation}
	For each $v\in \D^\mesh\cap \cB_{\bz}$, 
	$\alpha_4^\mesh(\mesh,10s)\le \P[ \textrm{$v$ is $s$-important for $\omega^\mesh$}]\le \alpha_4^\mesh(\mesh,s).$ 
	By the quasi-multiplicativity of $\alpha^\mesh_4(\cdot,\cdot)$ \eqref{eq:quasi-mult},
	there is a constant $C>0$ not depending on $\bz,s$ such that
	\begin{equation}\label{eq:mono2}
	\P[F^\mesh_s]\le C\P[ \textrm{$v$ is $s$-important for $\omega^\mesh$}], \quad \forall\;\delta\in(0,0.1)\textrm{ and } v\in \D^\mesh\cap \cB_{\bz}.
	\end{equation}
	
	If $v\in \D^\mesh\cap \cB_{\bz}$ is $s$-important for $\omega^\mesh$, then $v$ must be $A$-important for $\omega^\mesh$. On the event $E^\mesh_\eps$, we further have that $v$ is $\eps$-pivotal for $(\bh,\omega^\mesh)$. Therefore
	\[ 
	\P\left[\textrm{  $v$ is not $\eps$-pivotal for $\omega^\mesh$ while $v$ is $s$-important for $\omega^\mesh$}\right]\le  
	\P[\neg E^\mesh_\eps,\; F^\mesh_s],\quad \forall\; \mesh\in(0,\delta_0).
	\]
	By \eqref{eq:cond} and \eqref{eq:mono2}, for small enough $\zeta$ the upper bound in Lemma~\ref{lem:cover3} holds for $v\in \D^\mesh\cap \cB_{\bz}$. 
	We can choose finitely many $\bd z_i$'s such that $\cB_{\bd z_i}$   cover $r\D$. 
	This concludes the proof of Lemma~\ref{lem:cover3}.
\end{proof}

\subsubsection{Measures on $\eps$-pivotal points and the proof of Lemma~\ref{lem:pmeasure}}\label{subsub:pivm}
In Section \ref{subsub:scaling-piv} we proved that the natural Euclidean and quantum measure defined on the $\rho$-important points converge in the scaling limit to their continuum counterparts. In this subsection we prove a similar result for the $\eps$-pivotals, and %(building on Section \ref{subsub:mutual}) 
	we prove that the measures on the $\eps$-pivotals can be obtained via a restriction of the measures on the $\rho$-important points. This allows us to conclude the proof of Lemma~\ref{lem:pmeasure}.
\begin{proposition}\label{prop:pivm}
	Fix $\eps>0$. As $\mesh\to0$, $\arm^\mesh_4(\mesh,1)^{-1}$ times the Lebesgue measure restricted to $\cP^\mesh_\eps$ converge to  a measure $\Mink_\eps$ in probability.  
	The restriction of $\nu_\mesh$ to $\cP^\mesh_\eps$ converge to a measure  $\cM^\eps(\fh,\Gamma)$ in probability.
	Recall the constant $c'>0$ in Proposition~\ref{prop:dis-2SLE} and $\rho^\eps$ in Lemma \ref{lem:cover2}. 
	For each fixed $u\in(0,1)$, almost surely
	\begin{equation}\label{eq:restriction}
	\Mink_\eps=c'\Mink^\rho|_{\cP_\eps},\qquad \textrm{and}\qquad		\cM^\eps(\fh,\Gamma)=(c'e^{\fh/\sqrt 6} \Mink^\rho)|_{\cP_\eps} \textrm{ with } \rho=u\rho^\eps.
	\end{equation}
\end{proposition}
\begin{proof}
	Conditioning on $\fh$ and $\omega^\mesh$, let $z^\mesh$ be sampled uniformly from $\cP^\rho_\mesh$. 
	By Proposition~\ref{prop:piv-measure}, we can assume that $z^\mesh$ converge almost surely to a random point $z\in \cP^\rho$. Moreover, conditioning on $(\fh,\Gamma)$, the conditional law of $z$ is $(\Mink^\rho(\D))^{-1}\Mink^\rho$. Let $A(z^\mesh,\eps)$ (resp. $A(z, \eps)$)  be the event that  $z^\mesh$ (resp., $z$) is $\eps$-pivotal for $(\bh, \omega^\mesh)$ (resp., $(\bh,\Gamma)$).
	We claim that if $A(z,\eps)$ occurs, then almost surely there exists $\eps'>\eps$ such that $A(z,\eps')$ occurs. In fact, if $\ell\in \Gamma$ is chosen in a manner independent of $\bh$, then it is clear from the definition of GMC that $\mu_\bh(\reg(\ell))$ is a non-atomic random variable.  Therefore $\eps\notin \{ \mu_\bh(\reg(\ell)): \ell\in\Gamma \}$ a.s., which proves the claim.
	
	If $A(z,\eps)$ occurs, due to the existence of $\eps'>\eps$ above, Lemma~\ref{lem:flip3} implies that  $A(z^\mesh,\eps)$ occurs for sufficiently small $\mesh$.
	If $A(z,\eps)$ does not occur, again by Lemma~\ref{lem:flip3}, $A(z^\mesh,\eps)$ does not occur for sufficiently small $\mesh$. Therefore $\lim_{\mesh\to 0}\1_{A(z^\mesh,\eps)} = \1_{A(z, \eps)}$.
	Hence for any bounded continuous function $f:\C\to\R$, we  have 
	\begin{equation}\label{eq:eps-piv-conv}
	\lim_{\mesh\to 0}\E[ f(z^\mesh) \1_{A(z^\mesh,\eps)} \mid (\bh, \omega^\mesh,\Gamma) ]=\E[ f(z) \1_{A(z,\eps)} \mid (\bh,\omega^\mesh,\Gamma)] \qquad a.s.
	\end{equation}
	Since $\lim_{\mesh\to 0}\Mink^\rho_\mesh(\D)=\Mink^\rho(\D)$, $\arm^\mesh_4(\mesh,1)^{-1}$ times Lebesgue measure restricted to $\cP^\mesh_\eps$ converge to some limiting measure $\Mink_\eps$ in probability, and we have $\Mink_\eps=c'\Mink^\rho|_{\cP_\eps}$ a.s. Therefore  $\Mink_\eps=c'\Mink^\rho|_{\cP_\eps}$ a.s.
	
	The results concerning $\nu^\mesh,\cM^\eps(\fh,\Gamma)$, and $e^{\fh/\sqrt 6} \Mink^\rho$ follow from the exact same argument, where we assume that $z^\mesh$ is sampled according to  $\nu^\mesh|_{\cP^\rho_\mesh}$ and invoke Lemma~\ref{lem:UI} instead Proposition~\ref{prop:piv-measure}.
\end{proof}

\begin{proof}[Proof of Lemma~\ref{lem:pmeasure}]
	The coupling of $(\omega^\mesh,\bh, \Gamma_0)$ in Lemma~\ref{lem:pmeasure} is exactly as $(\omega^\mesh,\bh, \Gamma)$ in Proposition~\ref{prop:pivm}. 
	Now Lemma~\ref{lem:pmeasure} follows from Proposition~\ref{prop:pivm}.
	Moreover, the set $\cA$ in Lemma~\ref{lem:pmeasure} can be taken to be $\cP^\rho$ for small enough $\rho$.
\end{proof}

\subsection{Proof of Proposition~\ref{prop:quantum-Mink}}
\label{sec:prop:quantum-Mink}

We now conclude the proof of Proposition~\ref{prop:quantum-Mink} using results of the previous subsection. We first provide a precise definition of the measure  $\pivm^\eps$ in Proposition \ref{prop:quantum-Mink}.
\begin{definition}\label{def:pivm}
	Fix $\eps>0$. Recall  $\bh, \Gamma$  in Proposition~\ref{prop:quantum-Mink} and  $c'$ in Proposition~\ref{prop:dis-2SLE}.
	Let $\rho^\eps$ be defined as in Lemma \ref{lem:cover2} in terms of $\bh$.
	We set $ \pivm^\eps\defeq(c'e^{\bh/\sqrt{6}}\Mink^\rho)|_{\cP_\eps}$, where $\rho=0.5\rho^\eps$.
\end{definition}

Recall from Section~\ref{subsub:LDP} that $\Mink$ is the renormalized scaling limit of Lebesgue measure restricted to macroscopic pivotal points. To be more precise, we define the measure $\Mink$ to be the unique measure such that $\Mink|_{\cP^\rho}= c'\Mink^\rho$ for each $\rho>0$. 
Similarly, we define $e^{\bh/\sqrt 6}\Mink$ to be the unique measure such that 
$(e^{\bh/\sqrt 6}\Mink)|_{\cP^\rho}= e^{\bh/\sqrt 6}\Mink^\rho$	for each $\rho>0$. (Note that $\Mink$ itself is not locally finite as a Borel measure on $\R^2$ so we cannot directly define  a GMC on it.) In this sense, we may write $\pivm^\eps=(e^{\bh/\sqrt 6}\Mink)|_{\cP_\eps}$ as we did above~\eqref{eq:KPZ}.

Recall the definition of $\nu_\mesh$ and $\cP^\mesh_\eps$ from Section~\ref{subsec:pivot}. 
By Proposition~\ref{prop:pivm} and \eqref{eq:fields}, the measure $\pivm^\eps$ can be obtained as  the renormalized scaling limit of  $e^{\bh/\sqrt6}d^2z$ restricted to $\cP^\mesh_\eps$ (viewed as a collection of hexagons). Moreover, $\cM^\eps(\fh,\Gamma)$ from \eqref{eq:restriction} equals  $\xi_\fh(\p\D)^{\frac12}\pivm^\eps$ almost surely.
\begin{proof}[Proof of Proposition~\ref{prop:quantum-Mink}]
	Given $\rho$ in  Definition~\ref{def:pivm}, by Lemma~\ref{lem:cover}, we can find quads  $Q_1,\cdots Q_n$ such that  $\cP^\rho=\cup_{i=1}^n\cP^{Q_i}$ and the sets $\cP^{Q_i}$ are disjoint.
	By Lemma~\ref{lem:rho-is-eps}, we can find $\eps'\in (0,\eps)$  small enough such that $\cP^\rho\subset\cP_{\eps'}$.
	In Proposition~\ref{prop:pivot-bichordal}, let $h=\bh$  and $Q=Q_i$ for some $1\le i\le n$. By Definition~\ref{def:eps-piv}, $\nu_\cI$ from Proposition~\ref{prop:pivot-bichordal} then agrees with $\nu^{\eps'}_{\bh,\Gamma}|_{\cP^Q}$. Therefore
	$e^{\bh/\sqrt{6}}\Mink_Q=c\nu^{\eps'}_{\bh,\Gamma}|_{\cP_Q}$ with $c$ as in Proposition~\ref{prop:pivot-bichordal} (and Lemma~\ref{lem:pre}). Therefore $e^{\bh/\sqrt{6}}\Mink^\rho=c\nu^{\eps'}_{\bh,\Gamma}|_{\cP^\rho}$. 
	By Definition~\ref{def:eps-piv}, $\nu^{\eps}_{\bh,\Gamma}=\nu^{\eps'}_{\bh,\Gamma}|_{\cP_\eps}$.  Therefore $(c')^{-1} \pivm^{\eps}=c\nu^{\eps}_{\bh,\Gamma}$, so  
	$\pivm^\eps=\constp\nu^{\eps}_{\bh,\Gamma}$ for $\constp=cc'$.
\end{proof}

\subsection{The quad-crossing configuration determines the CLE$_6$}\label{subsec:quad}
By the iterative construction of $\CLE_6$ in Lemma~\ref{lem:iteration}, 
Theorem~\ref{thm:quad}  can be deduced from the following single interface variant.
\begin{proposition}\label{prop:quad-chod}
	In the setting of Theorem~\ref{thm:quad}, let $\eta$ be the interface of $\Gamma$ on $(\D,-i,i)$. Then $\eta$ is a.s.\ determined by $\omega$.
\end{proposition}
\begin{proof}[Proof of Theorem~\ref{thm:quad} given Proposition~\ref{prop:quad-chod}]
	Let $a=-i$ and $b=i$.  By Proposition~\ref{prop:quad-chod}, $\eta^{ab}$ is a.s.\ determined by $\omega$.
	Let $\cB$ be a dichromatic bubble of $\eta^{ab}$. Recall $x_\cB$, $\wh x_\cB$ and $\eta_\cB$ as defined above Lemma~\ref{lem:iteration}.
	Let $\phi:\cB\to\D$ be a conformal map with $\phi(x_\cB)=-i$ and $\phi(\wh x_\cB)=i$.
	Let $\phi_*\omega\in\cH(\D)$ be defined by $\phi_*\omega(Q)=\omega(\phi^{-1}\circ Q)$ for each $Q\in \cQ_\D$.
	Then $(\phi_*\omega, \phi\circ\eta_\cB)\overset{d}{=} (\omega,\eta^{ab})$, where  $\phi\circ\eta_\cB$ and $\eta^{ab}$ are viewed as curves modulo increasing reparametrization.
	Therefore $\phi\circ\eta_\cB$ is a.s.\ determined by $\phi_*\omega$,  hence $\eta_\cB$ is a.s.\ determined by $\omega$. 
	Therefore $\omega$ a.s.\ determine $\Gamma^b_a$.
	In light of Lemma~\ref{lem:iteration}, Theorem~\ref{thm:quad} follows by iterating this argument.
\end{proof}

\begin{figure}
	\centering
	\includegraphics[scale=1]{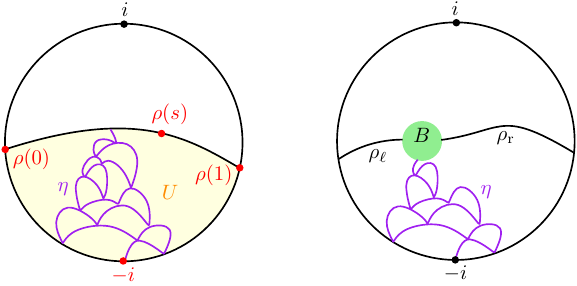}
	\caption{Illustration of the proof of Lemma \ref{prop:quad-chod}. {\bf Left}: The quad-crossing configuration $\omega$ determines whether the quad $(U,a, \rho(1),\rho(s), \rho(0))$ (in light yellow, with marked points in red) is crossed, and therefore whether $\eta'$ hits $\rho([0,s])$ or $\rho([s,1])$ first. {\bf Right}: Illustration of the event $E(B,\rho_{\markl},\rho_{\op r})$. By varying $\rho_{\markl}$ and $\rho_{\op r}$ we can determine whether $\eta'\cap B=\emptyset$.}
\end{figure}
It remains to prove Proposition~\ref{prop:quad-chod}.
In the following proof, given a quad $Q$, we write $Q=(U,a,b,c,d)$ if $Q((0,1)^2)=U$ and the four marked points are $a,b,c,d$ in counterclockwise order from $Q(0,0)=a$.
\begin{proof}[Proof of Proposition \ref{prop:quad-chod}]
	We first argue that the range of $\sle$ is determined by $\omega$.
	Let $\rho:[0,1]\to \D\cup\p\D$ be a simple smooth curve such that $\rho(0)$ and $\rho(1)$ are on the left and right boundary of $(\D,-i,i)$ (not including endpoints), respectively,  and  
	$\rho((0,1))\subset \D$.  Let $\tau=\inf\{t: \sle(t) \in\rho \}$.
	Let $U$ be the connected component of $\D\setminus \rho$ whose boundary contains $-i$.
	Then for each fixed  $s\in(0,1)$,  we claim that it is a.s.\ the case that 
	\[
	\xncomment{\textrm{$\eta(\tau)\in \rho([0,s])$ if and only if $\omega(Q)=1$ with $Q=(U,-i, \rho(1),\rho(s), \rho(0))$.}} 
	\]
	Given this claim, since $Q\in \cQ_{\ol \cD}$, by Lemma~\ref{lem:bdy-quad},  $\sle(\tau)$ is a.s.\ determined by~$\omega$.
	
	\xncomment{To prove the claim above, we write $a=-i$, $b=i$, and $c=\rho(s)$.  Recall the discussion about percolation interfaces and crossing events above Lemma \ref{lem:quad-chod0}. For $\mesh>0$, let $\omega^\mesh$ be the Bernoulli-$\frac12$ site percolation on $\D^\mesh$ with monochromatic blue boundary condition.
	Let $\sle^{ab}_{\mesh}$ be the percolation interface of $\omega^\mesh$ on $(\D^\mesh,a^\mesh,b^\mesh)$. 
		Let $\sle^{ac}_{U,\mesh}$ be the percolation interface of $\omega^\mesh$ on $(U^\mesh,a^\mesh, c^\mesh)$. Then
		$\sle^{ab}_{\mesh}, \sle^{ac}_{U,\mesh}, \omega^\mesh, \Gamma(\omega^\mesh)$ jointly converge. Denoting the joint limit by $(\sle^{ab}, \sle^{ac}_{U}, \omega, \Gamma)$,
		the joint law of $(\omega,\sle^{ab})$ is the same as that of $(\omega,\eta)$ in  Proposition \ref{prop:quad-chod}.
		Moreover, before hitting $\rho([0,1])$, the curves $\sle^{ab}$ and $\sle_U^{ac}$ a.s.\ coincide. 
		Therefore it is a.s.\ the case that  $\eta(\tau)\in \rho([0,s])$ if and only if  $\sle^{ac}_{U}$ hits $\rho([0,1])$ at a point on $\rho([0,s])$. 
		The latter event equals $\{\omega(Q)=1\}$ a.s.\ by Lemma \ref{lem:quad-chod0}.  This proves the claim above.}

	Let $B$ be a ball contained in $\D$. For $\bullet\in\{\markl,\op{r}\}$, let $\rho_\bullet:[0,1]\to \D\cup\p\D$ be a simple smooth curve such that $\rho_\bullet(0)\in\bdy B$, $\rho_\bullet(t)\in \D\setminus \bdy B$, and $\rho_\bullet(1)$ is on the left (resp., right) boundary of $(D,-i,i)$ when $\bullet$ equals $\markl$ (resp., $\op{r}$). Furthermore, we require $\rho_{\markl}\cap \rho_{\op r}=\emptyset$.
	By the previous paragraph the location where $\sle$ hits $B\cup \rho_{\markl} \cup \rho_{\op r}$ is a.s.\ determined by $\omega$. In particular, the event $E(B,\rho_{\markl},\rho_{\op r})$ that $\sle$ hits $B$ before $\rho_{\markl}\cup\rho_{\op r}$ is a.s.\ determined by $\omega$. Note that $\sle\cap B\neq \emptyset$ if and only if there exists $\rho_{\markl},\rho_{\op r}$  such that $E(B,\rho_{\markl},\rho_{\op r})$ occurs. Furthermore, if $E(B,\rho_{\markl},\rho_{\op r})$ occurs for some $\rho_{\markl},\rho_{\op r}$, then it holds a.s.\ that $E(B,\rho_{\markl},\rho_{\op r})$ occurs for $\rho_{\markl},\rho_{\op r}$ chosen from some countable set.
	This implies that the event $\sle\cap B\neq \emptyset$ is a.s.\ determined by $\omega$. Therefore, the range of $\sle$ is determined by $\omega$. 
	
	Now recall $\rho,U,\tau$ as defined above.  
	Since $\eta([0,\tau])$ is the intersection of the range of the percolation interfaces of $\Gamma$ on $(U,-i,\rho(0))$ and $(U,-i,\rho(1))$, by the previous paragraph $\sle([0,\tau])$ is a.s.\ determined by $\omega$. We  assume that $\psi(\eta)$ is parameterized by its half plane capacity, where $\psi(z)=\frac{z+i}{1+iz}$ maps $(\D,-i,i)$ to $(\bbH,0,\infty)$.
	Then for a fixed $t>0$, the event $\{\sle([0,t])\subset U\}=\{\tau>t\}$ is a.s.\ determined by $\omega$.
	Using the inclusion-exclusion principle and varying $U$,  we see that $\sle([0,t])$ is a.s.\ determined by $\omega$, hence $\eta$ is a.s.\ determined by $\omega$.
\end{proof}
\xncomment{\begin{remark}
		Based on~\cite[Section~2.4]{gps-pivotal}, it was essentially known  to Garban, Pete and Shcramm that the quad crossing configuration 
		determines the range of $\SLE_6$. 
		However, since the range of $\SLE_6$ does not determine its order~\cite{msw-non-simple}, 
		those authors considered Proposition~\ref{prop:quad-chod} as an open question; see Question~2.14 there. 
		The novel step in our proof as compared to \cite[Section~2.4]{gps-pivotal} is the observation that due to the target invariance property of $\SLE_6$ we can determine the range in every domain, which allows us to recover the ordering.  
	\end{remark}}
	
	\subsection{Proof of Lemmas \ref{prop:eps-LDP-law}, \ref{lem:limit-law}, and \ref{lem:two-jump}}\label{subsec:Markov} 
	\xncomment{In this section, we first prove in Section~\ref{subsub:flip-lattice} the existence of a probability space $(\Omega,\cF,\P)$ where Theorem~\ref{thm:LDP} holds and furthermore, the Poisson point process corresponding to the updates in the discrete LDP converge in a strong sense. 
		In Section~\ref{subsub:Markov}, we put $\fh$ and $(\omega^{\eps,\mesh}_t)_{t\geq 0}$ into the framework of continuous time finite-state Markov chains.
		Then in Section~\ref{subsub:proof}, we show that  $(\Omega,\cF,\P)$ from Section~\ref{subsub:flip-lattice} satisfies  Lemma~\ref{prop:eps-LDP-law}, which 
		asserts that $(\Gamma^{\eps,\mesh}_t)_{t\ge 0}$ converge in law to a process $(\Gamma^\eps_t)_{t\ge 0}$. Moreover, we prove Lemma~\ref{lem:limit-law} which describes the law of $(\fh,\Gamma^\eps_t)_{t\ge 0}$.  Finally, in Section~\ref{subsub:flip-map}, we prove Lemma~\ref{lem:two-jump} which gives  convergence of the $\eps$-dynamics on the planar map until the second jump.
	}

	\subsubsection{Assumptions on $(\Omega,\cF,\P)$}\label{subsub:flip-lattice}
	Let $(\Omega,\cF,\P)$ be a probability space satisfying Theorem~\ref{thm:LDP}.
	Recall  that $\lim_{\mesh\to0}\Gamma^{\eps,\mesh}_0=\Gamma^\eps_0$ a.s.
	Let $(\omega^\mesh,\Gamma^\mesh,\Gamma)\defeq (\omega^\mesh_0,\Gamma^{\eps,\mesh}_0,\Gamma^{\eps}_0)$
	so that $(\omega^\mesh,\Gamma^\mesh,\Gamma,\fh)$ satisfies the conditions in Section~\ref{subsec:pivot}.
	Recall $\nu_\mesh$ at the beginning of Section~\ref{subsec:pivot}. 
	Let $\nu^\rho_\mesh$ and  $\Mink^\rho$ be as in Lemma \ref{lem:UI}.
	Conditioning on $\fh$, the ringing locations and times for the clocks in the discrete LDP $(\omega^\mesh_t)_{t\ge 0}$ is a Poisson point process (p.p.p.) with intensity $\nu_\mesh\otimes dt$, which we denote by $\op{PPP}_\mesh$. 
	If we only look at updates in $\cP^\rho_\mesh$, namely, $\rho$-important points of $\omega^\mesh$, then we get a p.p.p.\ with intensity $\nu^\rho_\mesh\otimes dt$, which we denote by $\op{PPP}^\rho_\mesh$.
	For $t>0$ and $x\in \D^\mesh$, we write $(x,t)\in \op{PPP}_\mesh$ if the clock at $x$ rings at time $t$. The same convention applies to other p.p.p.'s.
	In the rest of this section we further require that the probability space  $(\Omega,\cF,\P)$ satisfies the property in the following lemma.
	\begin{lemma}\label{lem:setting}
		There exists  $(\Omega,\cF,\P)$  satisfying both Theorem~\ref{thm:LDP} and the following condition.
		For each fixed $\rho>0$, $\op{PPP}^\rho_\mesh$
		converge almost surely to a p.p.p.\  $\op{PPP}^\rho$ with intensity $c' e^{\fh/\sqrt6} \Mink^\rho\xncomment{\otimes dt}$ in the following sense. 
		For each $T>0$,  as $\mesh\to0$, 
		$\{(x,t) \in \op{PPP}^\rho_\mesh: t\in [0,T] \}$ converge to $\{(x,t) \in \op{PPP}^\rho: t\in [0,T] \}$ almost surely.
	\end{lemma}
	\begin{proof}
		Let $(\Omega,\cF,\P)$ be a probability space satisfying Theorem~\ref{thm:LDP}. In particular, $(\omega^\mesh,\Gamma^\mesh,\Gamma,\fh)$ satisfies the conditions in Section~\ref{subsec:pivot}.
		Fix $k\in\N$ and set $s=10^{-k}$. By Lemma~\ref{lem:UI}, $\lim_{\mesh\to0}\nu^{s}_\mesh=c' e^{\fh/\sqrt{6}} \Mink^{s}$ in probability. 
		By \cite[Lemma~7.5 and Corollary~7.6]{gps-dynamic}, we can find  a coupling of $(\omega^\mesh,\Gamma^\mesh,\fh)$ and $\op{PPP}_\mesh$ such that $\op{PPP}^{s}_\mesh$ converge almost surely to a p.p.p.\  $\op{PPP}^{s}$ with intensity $c' e^{\fh/\sqrt6} \Mink^{s}\xncomment{\otimes dt}$ in the sense specified in Lemma~\ref{lem:setting}.
		By Definition~\ref{def:rho-imp} and elementary geometric considerations, for each $\rho\ge 10s$, we have $\cP^{\rho} \subset \cP^{s}$, and  $\cP^{\rho}_\mesh \subset \cP^{s}_\mesh$ for small enough $\mesh$. 	
		Fix $T>0$. By Lemma~\ref{lem:cover}, 
		there almost surely exists $\rho'\in (s,\rho)$ and $\rho''>\rho$ sufficiently close to $\rho$,  such that for each 
		$(x,t) \in \op{PPP}^{s}$ with $t\in [0,T]$,  if $x\in \cP^{\rho}$, then $x\in \cP^{\rho''}$, otherwise,  $x\notin \cP^{\rho'}$.
		By the convergence of loops, $\{(x,t) \in \op{PPP}^{\rho}_\mesh: t\in [0,T] \}$ converge to $\{(x,t) \in \op{PPP}^{\rho}: t\in [0,T] \}$ almost surely.
		In particular, the convergence holds for $\rho=10^{-k+1}$.
		
		By the Skorokhod embedding theorem, we can further require $(\Omega,\cF,\P)$  to be such that
		$\op{PPP}^s_\mesh$ 	converge to $\op{PPP}^s$ a.s.\ for $s\in \{10^{-k}:k\in\N\}$.  In such a  coupling, for a fixed $\rho>0$, by considering $s=10^{-k}$ with $\rho\ge 10s$ and repeating the argument in the previous paragraph, we see that $\op{PPP}^\rho_\mesh$ converge  to  $\op{PPP}^\rho$ a.s. This concludes the proof.
	\end{proof}

	\begin{comment}
	Even if we only require rational $\rho$'s in Lemma~\ref{lem:setting}, we would still have to keep the entire first paragraph. Because we need to deal with different values of $\rho$ at the same time. 
	We would only skip the last sentence of the second paragraph.
	\end{comment}
	
	\subsubsection{A continuous time  Markov chain}\label{subsub:Markov}
	To prove Lemmas~\ref{prop:eps-LDP-law} and~\ref{lem:limit-law},
	we put $\fh$ and $(\omega^{\eps,\mesh}_t)_{t\geq 0}$ into the framework of continuous time finite-state Markov chains.
	Let $\cS^\mesh$ be the space of site percolation configurations of $\D^\mesh$ with monochromatic blue boundary condition.
	Then conditioning on $\fh$,  $(\omega^{\eps,\mesh}_t)_{t\geq 0}$ is  a  continuous time  Markov chain on the  state space $\cS^\mesh$ whose initial distribution is the uniform measure. 
	Let $Q_\fh:=(q_{ij})_{i,j\in \cS^\mesh}$ be the transition rate matrix of $(\omega^{\eps,\mesh}_t)_{t\geq 0}$. 
	For any two distinct states $i$ and $j$ in $\cS^\mesh$, if
	\begin{enumerate}
		\item the colorings of $i,j$ only differ at one vertex $v\in \D^\mesh$, and
		\item $v$ is an $\eps$-pivotal point for $i$, or, equivalently, for $j$,
	\end{enumerate} 
	then $q_{ij}=q_{ji}=\mu'_\fh(v)\arm^\mesh_4(\mesh,1)^{-1}$.  Otherwise, $q_{ij}=0$.
	Since $Q_\fh$ is symmetric, the uniform measure on $\cS^\mesh$ is a stationary distribution. Namely, $(\omega^{\eps,\mesh}_t)_{t\geq 0}$ is stationary conditioning on $\fh$.
	
	For each state $i\in \cS^\mesh$, let $N^{\eps}_\mesh(i)\defeq \sum  \mu'_\fh(v)\arm^\mesh_4(\mesh,1)^{-1}$, where the summation ranges over $\eps$-pivotal points of $(\bh,i)$. 
	Let $\cS^\mesh_+\defeq \{i\in \cS^\mesh: N^{\eps}_\mesh(i)>0\}$. If $\omega^{\eps,\mesh}_0\notin \cS^\mesh_+$, then $\omega^{\eps,\mesh}_t=\omega^{\eps,\mesh}_0$ for all $t\ge 0$.
	On the event $\omega^{\eps,\mesh}_0\in \cS^\mesh_+$, the process $(\omega^{\eps,\mesh}_t)_{t\ge 0}$  evolves as a  stationary Markov chain on $\cS^\mesh_+$.
	Let $(J^{\eps,\mesh}_k)_{k\in\N}$ be the \notion{discrete  skeleton} of  $(\omega^{\eps,\mesh} _t)_{t\geq 0}$.
	Namely, on the event $\omega^{\eps,\mesh}_0\in \cS^\mesh_+$, $(J^{\eps,\mesh}_k)_{k\in\N}$ is the discrete-time Markov chain on $\cS^\mesh_+$  keeping track of the jumps of $(\omega^{\eps,\mesh} _t)_{t\geq 0}$.  If $\omega^{\eps,\mesh}_0\notin \cS^\mesh_+$, then $J^{\eps,\mesh}_k=\omega^{\eps,\mesh}_0$ for each $k\in\N$.
	
	Conditioning on $\fh$, we can sample $(\omega^{\eps,\mesh}_t)_{t\geq 0}$ in a two-step procedure:
	\begin{enumerate}
		\item Run $(J^{\eps,\mesh}_k)_{k\in\N}$ with its  $\P$-law (conditioning on $\fh$).
		\item Conditioning on $\fh$ and $(J^{\eps,\mesh}_k)_{k\in\N}$, the time spent in each state $J^{\eps,\mesh}_k$ is an independent exponential random variable with rate $N^{\eps}_\mesh(J^{\eps,\mesh}_k)$.
	\end{enumerate}
	Let $P_\fh$ be the transition matrix of $(J^{\eps,\mesh}_k)_{k\in\N}$ conditioning on $\fh$.
	It  is elementary  to see that the uniform measure on $\cS^\mesh$ reweighted by $N^\eps_\mesh$ is a stationary measure for $P_{\fh}$. 
	In other words, 
	define $\cN^\eps_\mesh\defeq N^{\eps}_\mesh( \omega^{\eps,\mesh}_0 )$.
	Then $(J^{\eps,\mesh}_k)_{k\in\N}$ is stationary under the probability measure  obtained by normalizing $\cN^{\eps}_\mesh d\P$.
	
	\begin{comment}
	Whenever we talk about $\eps$-pivotal point,  we always have to use $\bh$, even if the Poisson clocks are driven by $\mu_\fh$.
	\end{comment}

 	\subsubsection{Proof of Lemmas~\ref{prop:eps-LDP-law} and~\ref{lem:limit-law}}\label{subsub:proof}
 	We now prove  Lemma~\ref{prop:eps-LDP-law} and~\ref{lem:limit-law}.
 	\xncomment{Given the Markov chain description of $(\omega^{\eps,\mesh}_t)_{t\geq 0}$, hence $(\Gamma^{\eps,\mesh}_t)_{t\geq 0}$, in Section~\ref{subsub:Markov}, the convergence in law for $(\Gamma^{\eps,\mesh}_t)_{t\geq 0}$ asserted in Lemma~\ref{prop:eps-LDP-law} to the desired limit desired described in Lemma~\ref{lem:limit-law}  is quite straightforward. We just need to show that the exponential clocks and the skeleton of 
 		$(\Gamma^{\eps,\mesh}_t)_{t\geq 0}$	converges to the desired distribution. 
 		For the skeleton convergence, we will use the fact that under the reweighted measure  
 		$\cN^{\eps}_\mesh d\P$ at the end of Section~\ref{subsub:Markov},
 		the skeleton $(J^{\eps,\mesh}_k)_{k\in\N}$  is stationary. If $\cN^{\eps}_\mesh$ converge in $L^1$, 
 		then given all the work done in Section~\ref{subsec:pivot}, this convergence would be trivial. 
 		The main difficulty that we face in this final step is that $\cN^{\eps}_\mesh$ may not converge in $L^1$, at least this is not clear to us. 
 		Most of the technical work in this subsection is devoted to circumventing  this issue.}
 	We start by	some basic limiting properties of $\cN^{\eps}_\mesh$. 
 	\begin{lemma}\label{lem:one-step}
 		Recall $\cM^\eps(\fh,\Gamma)$ in Proposition~\ref{prop:pivm}.
 		Fix $\eps>0$.  
 		Let $\cN^\eps$ be the $\cM^\eps(\fh,\Gamma)$-mass of $\eps$-pivotal points of $(\bh,\Gamma)$.
 		Then $\lim_{\mesh\to 0} \cN^\eps_\mesh=\cN^\eps$ and $\lim_{\mesh\to 0} \1_{\cN^\eps_\mesh=0}=\1_{\cN^\eps=0}$ in probability.
 	\end{lemma}
 	\begin{proof}
 		Note that $\cN^{\eps}_\mesh$ is the total  $\nu_\mesh$-mass of the $\eps$-pivotal points of $(\bh,\omega^\mesh)$. Setting $f\equiv 1$ in~\eqref{eq:eps-piv-conv},
 		we get $\lim_{\mesh\to 0} \cN^\eps_\mesh=\cN^\eps$. For the second assertion,  recall $\op{PPP}^\rho_\mesh$ and $\op{PPP}^\rho$ in Lemma~\ref{lem:setting}. For each fixed $T$ and $\rho>0$, we query whether points in $\{x: (x,t)\in \op{PPP}^\rho_\mesh, t\in [0,T]\}$ are $\eps$-pivotal for $(\bh,\omega^\mesh)$. By Lemma~\ref{lem:flip3} the answer converges to its counterpart  for $\op{PPP}^\rho$. Sending $T\to\infty$ and $\rho\rta 0$ we conclude.
 	\end{proof}

 	\xncomment{Recall $\rho^\eps$ from Lemma~\ref{lem:cover2}. Namely, $\rho^\eps>0$ is a random number such that each $\eps$-pivotal point  of $(\bh,\Gamma)$ (resp., $(\bh, \omega^\mesh)$)  is $\rho$-important for $\Gamma$ (resp., $\Gamma^\mesh$) for $\rho\in(0,\rho^\eps)$ and $\mesh\in (0,\rho)$.}
 	The following variant of Lemma~\ref{lem:flip3} is immediate from Lemma~\ref{lem:cover2}.
 	\begin{lemma}\label{lem:jump-prob}
 		Let $\tau^\mesh\defeq\inf\{t>0:\omega^{\eps,\mesh}_t\neq \omega^\mesh_0 \}$   be the first time $(\omega^{\eps,\mesh}_t)_{t\ge 0}$  jumps. 
 		Let  $\wh \Gamma^{\mesh}\defeq \Gamma^{\eps,\mesh}_{\tau^\mesh}$ 
 		if $\tau^\mesh<\infty$ and $\wh \Gamma^{\mesh}\defeq \Gamma^{\eps,\mesh}_{0}$ if $\tau^\mesh=\infty$.
 		Then the limit $\wh\Gamma=\lim_{\mesh\to 0}\wh\Gamma^{\mesh}$  exists in probability for the $\cL(\D)$-metric.
 	\end{lemma}
 	\begin{proof}
 		The case $\cN^\eps\neq 0$ is immediate, so we focus on the event that $\cN^\eps\neq 0$. By Lemma~\ref{lem:one-step}, on this event for small enough $\delta$ we have $\cN^\eps_\mesh>0$, hence $\tau^\mesh<\infty$.
 		Define $\wh\omega^\mesh\defeq\omega^{\eps,\mesh}_{\tau^\mesh}$, so that $\wh \Gamma^{\mesh}=\Gamma(\wh \omega^\mesh)$. 
 		Let $\bz^\mesh\in\D^\mesh$ be such that $ \wh\omega^\mesh(\bz^\mesh )\neq \omega^\mesh(\bz^\mesh)$. \xncomment{Then $\bz^\mesh$ must be $\eps$-pivotal for  $(\bh, \omega^\mesh)$. Therefore $(\bz^\mesh, \tau^\mesh) \in \op{PPP}^\rho_\mesh$ for $\rho\in(0,\rho^\eps)$ and $\mesh\in (0,\rho)$. 
 			By the almost sure convergence of $\op{PPP}^\rho_\mesh$ for arbitrary $\rho$ from Lemma~\ref{lem:setting},we see that
 			$\bz^\mesh$ converge almost surely to a random point $\bz\in\D$.  Now the convergence $\wh\Gamma=\lim_{\mesh\to 0}\wh\Gamma^{\mesh}$  follows from Lemma~\ref{lem:one-step}.}
 		%	By definition, $\bz^\mesh$ is sampled according to the restriction of $\nu_\mesh$ to $\cP^\mesh_\eps$. 
 		%	Therefore $\bz $ must be  sampled from the measure $\cM^\eps(\fh,\Gamma)$ in Proposition~\ref{prop:pivm}, which is the limit of 
 		%	the restriction of $\nu_\mesh$ to $\cP^\mesh_\eps$.
 	\end{proof}
 	For a fixed $\rho>0$,  let $\P^{\rho}= \P[\rho<\rho^\eps]^{-1} \1_{\rho<\rho^\eps} \P$ with $\rho^\eps$ in Lemma~\ref{lem:cover2}.	
 	We introduce $\P^\rho$ because \xncomment{under $\P^\rho$ every $\eps$-pivotal point of also $\rho$-important and  we do have convergence of $\cN^\eps_\mesh$  in $L^1$.}
 	\begin{lemma}\label{lem:UI2}
 		$\cN^\eps_\mesh$ converge to $\cN^\eps$ in $L^1$ under $\P^{\rho}$.
 	\end{lemma}
 	\begin{proof}
 		Since $\cN^\eps_\mesh\1_{\rho<\rho^\eps}\le \nu^\rho_\mesh(\D)$, by Lemma~\ref{lem:UI}, $\{\cN^\eps_\mesh\}_{\mesh>0}$ is uniformly integrable under $\P^{\rho}$. Since $\cN^\eps_\mesh$ converge to $\cN^\eps$ in $\P^{\rho}$-probability by Lemma~\ref{lem:one-step}, we have Lemma~\ref{lem:UI2}.
 	\end{proof}

 	Let $\wt \P^{\rho}_\mesh$ be the probability measure obtained by normalizing $\cN^\eps_\mesh \P^{\rho}$.  Let $\wt \P^\rho$ be the probability measure obtained by normalizing $\cN^\eps\P^{\rho}$. \xncomment{We first prove the variant of Lemmas \ref{prop:eps-LDP-law} and \ref{lem:limit-law} (i.e., Lemma~\ref{lem:J-conv}) under these truncated and reweighted measures. Then we use Lemma~\ref{lem:UI2} to remove the truncation and reweighting.

 		\begin{lemma}\label{lem:J-conv}
 			Fix $\eps>0$ and $\rho>0$.
 			We can enlarge the sample space $(\Omega,\cF,\wt\P^\rho)$ to admit a process $(\Gamma^{\eps}_t)_{t\ge 0}$ 
 			such that the $\wt \P^{\rho}_\mesh$-law of $(\fh, \Gamma^{\eps,\mesh}_t)_{t\ge 0}$ weakly converge to 
 			the $\wt\P^\rho$-law of $(\fh, \Gamma^{\eps}_t)_{t\ge 0}$ as $\mesh\rta 0$. Moreover,
 			the law of $(\fh,\Gamma^\eps_t)_{t \ge 0}$ can be described as follows.
 			Conditioned on $\fh$, $(\Gamma^\eps_t)_{t\ge 0}$ is a stationary Markov process, where $\Gamma^\eps_0$ is a $\CLE_6$ on $\D$. 
 			Moreover, almost surely $(\Gamma^\eps_t)_{t\ge 0}$ has infinitely many jumps but only has finitely many in any finite interval.
 			Conditioning on $(\D,\fh,\Gamma^\eps_0)$,  an exponential clock rings with rate $(\xi_\fh(\p \D))^{1/2}\cM^\eps_{\bh,\Gamma^\eps_0}(\D)$.
 			Once the clock rings, sample an $\eps$-pivotal point
 			$\bz$ from $\cM^\eps_{\bh,\Gamma^\eps_0}$.  The process $(\Gamma^\eps_t)$
 			jumps to the loop ensemble obtained from $\Gamma^\eps_0$ by flipping the color at $\bz$.
 			The remaining jumps in the process are sampled iteratively.
 		\end{lemma}}

\xncomment{	
	 \begin{proof}
 	Recall that  $(J^{\eps,\mesh}_k)_{k\ge 0}$ is the discrete  skeleton of  $(\omega^{\eps,\mesh} _t)_{t\geq 0}$.	In this proof we abuse notation and identify $J^{\eps,\mesh}_k$ with its loop ensemble $\Gamma(J^{\eps,\mesh}_k)$. 
 	Since $\wt\P^\rho_\mesh[\cN^\eps_\mesh>0]=1$ due to the reweighting of $\cN^\eps_\mesh$ in $\wt\P^\rho_\mesh$,  the skeleton chain is stationary and  a.s.\ non-trivial (i.e., it is not constant in time) %(i.e.\ it does not stay unchanged.)
 	 We will carry out the proof of Lemma~\ref{lem:J-conv} in three  steps.
 	In Step 1, we show that in any coupling where $(\fh, J^{\eps,\mesh}_0)$ converge in probability, we can enlarge the coupling such that $J^{\eps,\mesh}_1$ converge in probability as well. In Step 2 we use the stationarity of $(J^{\eps,\mesh}_k)_{k\ge 0}$ to inductively show that  there exists a coupling where $(\fh, J^{\eps,\mesh}_0,\cdots, J^{\eps,\mesh}_n)$ converge in probability for all $n\ge 0$. In Step 3 we deal with the convergence  of 
 	the waiting time between jumps and their local finiteness in the continuum. Putting these three steps together we get the desired description of the limiting process $(\Gamma^\eps_t)_{t\ge 0}$.

 				\bigskip
 				
 				\noindent{\bf Step 1.  From $J^{\eps,\mesh}_0$ to $J^{\eps,\mesh}_1$.} 
 				Let $J^\eps_0=\Gamma_0$. Since $(\fh, J^{\eps,\mesh}_0)$  converge almost surely to  $(\fh,J^\eps_0)$ under  $\P$, the same convergence   holds  under $\P_\rho$.
 				By  the $L^1$  convergence of the Radon-Nykodim derivative in Lemma~\ref{lem:UI2}, the $\wt\P^{\rho}_\mesh$-law of $(\fh, J^{\eps,\mesh}_0)$  
 				converge to  the $\wt\P^{\rho}$-law of $(\fh, J^\eps_0)$. Recall the measure $\nu_\mesh$ from Section~\ref{subsec:pivot}  where each vertex $x$  is assigned mass $\mu'_\fh(x)\arm^\mesh_4(\mesh,1)^{-1}$. Let $\nu^{\eps,\mesh}_0$ be $\nu_\mesh$ restricted to the set of $\eps$-pivotal points for $(\fh, J^{\eps,\mesh}_0)$.
 				Then by Proposition~\ref{prop:pivm}, $\nu^{\eps,\mesh}_0$ converge in probability to the measure $\cM^\eps(0,J^\eps_0)=\cM^\eps(\fh,\Gamma)$. 
 				Again by Lemma~\ref{lem:UI2}, the $\wt\P^{\rho}_\mesh$-law of $(\fh, J^{\eps,\mesh}_0, \nu^{\eps,\mesh}_0)$  
 				converge to  the $\wt\P^{\rho}$-law of $(\fh, J^\eps_0,\cM^\eps(\fh,\Gamma))$.

 				Let $\fh^\mesh=\fh$ for each $\delta>0$. Suppose we are under an arbitrary coupling of  the $\wt\P^{\rho}_\mesh$-law of $(\fh^\mesh, J^{\eps,\mesh}_0)$  and the $\wt\P^{\rho}$-law of $(\fh, J^\eps_0)$ such that  $(\fh^\mesh, J^{\eps,\mesh}_0)$ converge to $(\fh, J^\eps_0)$ in probability. 
 				We claim that  the sample space of this coupling  can be enlarged to become a  coupling of the $\wt\P^{\rho}_\mesh$-law of $(\fh^\mesh, J^{\eps,\mesh}_0, J^{\eps,\mesh}_1)$ where  $J^{\eps,\mesh}_1$ converge in probability as well. 
 				
 				Note that $\cM^\eps(\fh,\Gamma)$ can be viewed as a measurable function from $H^{-1}(\D)\times \cL(\D)$ to the space of Borel measures on $\D$, which is well defined modulo a $\P$-probability zero event.  
 				Therefore by Lemma~\ref{lem:law-prob},	under the  coupling in the previous paragraph, $\nu^{\eps,\mesh}_0$ converge to $\cM^\eps(\fh,J^\eps_0)$ in probability
 				so we can enlarge the sample space for this coupling to contain random points $\bz^\mesh_0$ and $\bz_0$
 				sampled according to  $\nu^{\eps,\mesh}_0$ and  $\cM^\eps(\fh,J^\eps_0)$, respectively, such that $\bz^\mesh_0\rta \bz_0$ in probability. 
 				By the definition of $\wt\P^{\rho}_\mesh$ and $\wt\P^{\rho}$, almost surely each $\eps$-pivotal point is $\rho$-important. Therefore 
 				$\bz^\mesh_0$ (resp. $\bz$) is $\rho$-important for $J^{\eps,\mesh}_0$ (resp. $J^\eps_0$). 
 				Let $J^{\eps,\mesh}_1$ be the loop ensemble obtained from $J^{\eps,\mesh}_0$ by flipping the color at $\bz^\mesh_1$.  
				By Lemma~\ref{lem:flip3}  and Remark~\ref{rmk:flip3}, $J^{\eps,\mesh}_1$ converge in probability to the loop ensemble  
 				$J^{\eps}_1$ obtained by flipping the color of  $J^{\eps}_0$ at  $\bz_0$. 
 				Although  Lemma~\ref{lem:flip3} is proved under $\P$, given the $L^1$-convergence  from Lemma~\ref{lem:UI2}, 
 				the proof of Lemma~\ref{lem:flip3} works in this setting.  		
 				This gives a  coupling of the $\wt\P^{\rho}_\mesh$-law of $(\fh^\mesh, J^{\eps,\mesh}_0, J^{\eps,\mesh}_1)$ that converge in probability.
 				\bigskip
 				
 				\noindent{\bf Step 2. Convergence of  the full skeleton.}  
 				We perform the following induction. Suppose for some $ n\ge 0$ there exists a coupling of  the $\wt\P^{\rho}_\mesh$-law of $(\fh^\mesh, J^{\eps,\mesh}_0,\cdots ,J^{\eps,\mesh}_n)$ such that the convergence holds in probability. We write $(\fh^0, J^{\eps}_0,\cdots ,J^{\eps}_n)$ as the in-probability limit. 
 				By the two-step sampling procedure for $\{J^{\eps,\mesh}_k\}_{k\in \N}$ at the end of Section \ref{subsub:Markov}.
 				we see that the $\wt\P^{\rho}_\mesh$-conditional law of $(J^{\eps,\mesh}_k)_{k\ge 0}$ given $\fh^\mesh$ is a stationary Markov chain. 
 				In particular, under this coupling the law of $(\fh^\mesh, J^{\eps,\mesh}_n)$ equals  the $\wt\P^{\rho}_\mesh$-law of $(\fh^\mesh, J^{\eps,\mesh}_0)$ and 
 				the law of $(\fh^0,J^\eps_n)$ equals  the $\wt\P^{\rho}$-law of $(\fh, J^{\eps}_0)$. 
 				Then we can enlarge the sample space of this coupling to admit $\bz^\mesh_n$ and $ \bz_n$, such that $\bz^\mesh_n\rta \bz_n$ in probability, 
 				and $\bz^\mesh_n$ and $\bz_n$ are sampled in the same way as $\bz^\mesh_0$ and $\bz_0$ in Step 1 with $J^{\eps,\mesh}_n, J^{\eps}_n$ in place of 
 				$J^{\eps,\mesh}_0, J^{\eps}_0$.
 				Define $J^{\eps,\mesh}_{n+1}$ and  $J^{\eps}_{n+1}$  in the same way as $J^{\eps,\mesh}_{1}$ and  $J^{\eps}_{1}$ by flipping colors. 
 				Then by the same argument as in Step 1, $J^{\eps,\mesh}_{n+1}$ converge to  $J^{\eps}_{n+1}$ in probability.
 				We also note that conditioning on $(\fh^0, J^{\eps}_0,\cdots ,J^{\eps}_n)$,  the point $\bz_n$ is sampled according to $\cM^\eps(\fh^0,J^\eps_n)$,  
 				and  $J^\eps_{n+1}$ is obtained from $J^\eps_n$ by flipping color at $\bz_n$.

 				By the above induction, we see that for each $n\ge 0$,  the	$\wt\P^{\rho}_\mesh$-law of $(\fh^\mesh, J^{\eps,\mesh}_0,\cdots ,J^{\eps,\mesh}_n)$	
 				weakly converge. Since the law of the field equals the $\wt \P^\rho$-law of $\fh$,
 				we can enlarge the sample space of $(\Omega, \cF,\wt \P^\rho)$ to admit random variables $(J^{\eps}_0,\cdots ,J^{\eps}_n)$  such that 
 				the weak limit is the $\wt \P^\rho$-law of 	
 				$(\fh, J^{\eps}_0,\cdots ,J^{\eps}_n)$. Also from the above induction, we see that conditioning on  $\fh$, 
 				the conditional law of $(J^{\eps}_k)_{k\in \N}$ is a Markov chain  whose transition kernel  is as 
 				described in Lemma~\ref{lem:J-conv}, i.e., first sample an $\eps$-pivotal point and then flip the color.
 				
 				\bigskip

 				\noindent{\bf Step 3. Convergence of jumping times and their local finiteness.}  Recall the enlarged sample space $(\Omega,\cF,\wt \P^\rho)$
 				at the end of Step 2. We further enlarge it as follows. Let $\cM_k$ be the total mass of $\cM^\eps(\fh,J^\eps_k)$ for $k\ge 0$.
 				Conditioned on  $(\fh, J^\eps_k)_{k\ge 0}$,  sample an independent sequence $(\theta_i)_{i\ge 1}$
 				such that $\tau_k$ has the law of the ringing time of an exponential clock with rate  $\cM_{i-1}$ for each $i\ge 1$. 	  
 				Let  $\cN^{\eps,\mesh}_k$ be defined as $\cN^\eps_\mesh$ with $(\fh,J^{\eps,\mesh}_k)$ in place of 
 				$(\fh, \Gamma^{\eps,\mesh}_0)$. Then $\cN^{\eps,\mesh}_k$ convergence in law  to $(\cM_k)_{k\ge 0}$. 
 				Moreover, the convergence holds jointly with $(\fh, J^\eps_k)_{k\ge 0}$. 
 				For $k\ge 1$, let $\tau^\mesh_k$ be the $k$-th jumping time of $(\Gamma^{\eps,\mesh}_t)_{t \ge 0}$ and  $\tau_k=\sum_{i=1}^k \theta_i$.
 				Then  by the two-step sampling procedure at the end of Section \ref{subsub:Markov} and the conclusion of Step~2,
 				the $\wt \P^\rho_\mesh$-law of $\fh$, $ (J^{\eps,\mesh}_k)_{k\ge 0}$, and $(\tau^{\mesh}_k)_{k\ge 1}$ converge to 
 				the $\wt \P^\rho$-law of $\fh$, $ (J^{\eps}_k)_{k\ge 0}$, and $(\tau_k)_{k\ge 1}$.
 				
 				Since $(\fh,J^\eps_k)_{k\ge 0}$ is stationary, so is $(\cM_k)_{k\ge 0}$. On the other hand $\wt\P^\rho[\cM_k\in(0,\infty)]=1$ due to the $\cN^\eps$-reweighting in $\wt \P^\rho$. 
 				Now the ergodic theorem yields that
 				\begin{equation}\label{eq:M-infinity}
 				\sum_{i=1}^{\infty}\cM^{-1}_i=\infty,\quad \wt\P^\rho\textrm{-a.s.}
 				\end{equation}
 				Therefore, $\sum_{i=1}^{\infty}\theta_i=\infty$ almost surely under  $\wt\P^\rho$.
 				Now we define $\Gamma^\eps_t= J^\eps_k$ if $t\in [\tau_k,\tau_{k+1})$.  
 				Then $(\Gamma^\eps_t)_{t\ge 0}$  makes finitely many jumps in any bounded interval $\wt\P^\rho$-a.s. 
 				Moreover,  the $\wt\P^\rho$-law of $(\fh, \Gamma^{\eps}_t)_{t\ge 0}$  is the weak limit of the $\wt \P^{\rho}_\mesh$-law of $(\fh, \Gamma^{\eps,\mesh}_t)_{t\ge 0}$. 
 				Recall from Section~\ref{sec:prop:quantum-Mink} that $\cM^\eps(\fh,\Gamma)=\xi_\fh(\p\D)^{\frac12}\pivm^\eps$. 
 				Therefore the rate of $(\theta_k)_{k\ge 1}$ given $\fh $ and $(J^{\eps}_k)_{k\ge 0}$ is exactly as described in Lemma~\ref{lem:J-conv}.
 				Combined with the description of the law of $(\fh, J^{\eps}_0,\cdots ,J^{\eps}_n)$ from Step 2, we see that 
 				the $\wt\P^\rho$-law of $(\fh, \Gamma^{\eps}_t)_{t\ge 0}$ is as described in Lemma~\ref{lem:J-conv}.
 			\end{proof}
 		}

 		\begin{proof}[Proof of  Lemmas \ref{prop:eps-LDP-law} and \ref{lem:limit-law}.]
 			\xncomment{It suffices  to show that in the setting of Lemma~\ref{lem:J-conv},  
 				for a fixed $\rho>0$, as $\mesh\to 0$, the $\P^{\rho}$-law of $(\fh,\Gamma^{\eps,\mesh}_t)$  weakly converge 
 				to the $\P^\rho$-law of  $ (\fh,\Gamma^\eps_t)_{t\ge 0}$. 
 				Once  this is done, Lemma~\ref{lem:J-conv} shows that the $\P^\rho$-law of  $ (\fh,\Gamma^\eps_t)_{t\ge 0}$ is as described in  Lemma~\ref{lem:limit-law} with the further condition $\{\rho^\eps>\rho\}$.
 				So sending $\rho\to0$ we will conclude the proof of  Lemmas \ref{prop:eps-LDP-law} and \ref{lem:limit-law}.}

 			For $n\in\N$, define $g_n(x)=n^2x\1_{x<n^{-1}}+x^{-1}\1_{x>n^{-1}}$ for $x\in [0,\infty)$. 
 			Let $f$ be a bounded continuous function on the space of $H^{-1}(\D)\times \cL(\D)$-valued processes on $[0,\infty)$ under Skorokhod topology.
 			Let $\E^{\wt \P^\rho}$ be the expectation with respect to $\wt\P^\rho$. Define $\E^{\P^\rho}$ and $ \E^{\wt \P^\rho_\mesh}$ similarly. 
 			\xncomment{For a fixed $n$, we see that $g_n$ is a bounded continuous function on $[0,\infty)$. Therefore }
 			\[
 			\lim_{\mesh\to 0} \E^{\wt \P^\rho_\mesh}[f((\fh,\Gamma^{\eps,\delta}_t)_{t\ge 0}) g_n(\cN^{\eps}_\mesh)]=\E^{\wt \P^\rho}[f((\fh,\Gamma^\eps_t)_{t\ge 0}) g_n(\cN^\eps)]\quad \textrm{for each }n\in\N.
 			\] 
 			\xncomment{Since $\cN^\eps_\mesh\P^\rho=\E^{\P^\rho}[\cN^\eps_\mesh] \wt \P^\rho_\mesh$, $\cN^\eps\P^\rho=\E^{\P^\rho}[\cN^\eps] \wt \P^\rho$, and 
 				\(\lim_{\mesh\to0 }\E^{\P^\rho}[\cN^\eps_\mesh]=\E^{\P^\rho}[\cN^\eps]\),   for each $n\in \N$ we have
 				\begin{align}
 				\lim_{\mesh\to 0} \E^{\P^\rho}[f((\fh,\Gamma^{\eps,\delta}_t)_{t\ge 0}) g_n(\cN^\eps_\mesh)\cN^\eps_\mesh]&= \lim_{\mesh\to 0}\E^{\P^\rho}[\cN^\eps_\mesh]
 				\E^{\wt \P^\rho_\mesh}[f((\fh,\Gamma^{\eps,\delta}_t)_{t\ge 0}) g_n(\cN^{\eps}_\mesh)]  \nonumber\\
 				&=\E^{\P^\rho}[\cN^\eps] \E^{\wt \P^\rho}[f((\fh,\Gamma^\eps_t)_{t\ge 0}) g_n(\cN^\eps)]
 				=\E^{\P^\rho}[f((\fh,\Gamma^\eps_t)_{t\ge 0}) g_n(\cN^\eps) \cN^\eps]. \label{eq:lim1}
 				\end{align}}
 			Since $0\le g_n(x)x\leq 1$ for all $x>0$,
 			\begin{equation}\label{eq:lim2}
 			\left|\E^{\P^\rho}[f((\fh,\Gamma^{\eps,\delta}_t)_{t\ge 0}) g_n(\cN^\eps_\mesh)\cN^\eps_\mesh]-\E^{\P^\rho}[f((\fh,\Gamma^{\eps,\delta}_t)_{t\ge 0})\1_{\cN^\eps_\mesh>0}]\right|\le \|f\|_\infty \P^\rho[0<\cN^\eps_\mesh<n^{-1}].
 			\end{equation}
 			Moreover, \eqref{eq:lim2} remains true if $\Gamma^{\eps,\delta}_t$ and $\cN^\eps_\mesh$ are replaced by $\Gamma^{\eps}_t$ and $\cN^\eps$, respectively.  Therefore  
 			\[
 			\lim_{n\to\infty} \E^{\P^\rho}[f((\fh,\Gamma^\eps_t)_{t\ge 0}) g_n(\cN^\eps) \cN^\eps] = \E^{\P^\rho}[f((\fh,\Gamma^\eps_t)_{t\ge 0})  \1_{\cN^\eps>0}].
 			\]
 			This combined with  \eqref{eq:lim1} and~\eqref{eq:lim2} gives that
 			\begin{align}
 			\lim_{\mesh\to 0}\E^{\P^\rho}[f((\fh,\Gamma^{\eps,\delta}_t)_{t\ge 0})\1_{\cN^\eps_\mesh>0}]= \E^{\P^\rho}[f((\fh,\Gamma^\eps_t)_{t\ge 0})  \1_{\cN^\eps>0}]. \label{eq:non-zero}
 			\end{align}
 			On the event that $\cN^\eps_\mesh=0$, we have $\Gamma^{\eps,\delta}_t=\Gamma^{\eps,\delta}_0$ for $t>0$. Combined with Lemma~\ref{lem:one-step} and~\eqref{eq:non-zero}, the $\P^{\rho}$-law of $(\fh,\Gamma^{\eps,\mesh}_t)$  weakly converge as $\mesh\to 0$
 			\xncomment{to the $\P^\rho$-law of  $ (\fh,\Gamma^\eps_t)_{t\ge 0}$ as desired}. 
 		\end{proof}

 		\xncomment{The content of Lemmas \ref{prop:eps-LDP-law} and \ref{lem:limit-law} 
 			is the existence and a description of the weak limit of the $\wt \P^{\rho}_\mesh$-law of $(\fh, \Gamma^{\eps,\mesh}_t)_{t\ge 0}$. 
 			The enlargement of the sample space in the statement of Lemma~\ref{prop:eps-LDP-law}  is not essential.  	
 			Since the marginal law of the field stays the same, for  convenience we just use the field $\fh$ to generate a sample of the limiting dynamic. On the other hand, in  the  proof of
 			Proposition~\ref{prop:eps-LDP-proba} we will show  that in a sample space satisfying Lemma~\ref{lem:setting}, the convergence already holds in probability. 
 			Proving this  requires additional input from LDP which we will provide in Section~\ref{subsec:stability}.}

 		\subsubsection{Convergence after the first flip: planar map case} \label{subsub:flip-map}
 		
 		Now we turn our attention to Lemma~\ref{lem:two-jump}.
 		Suppose we are in the setting of Lemma~\ref{lem:jump-prob} and the proof of Lemmas \ref{prop:eps-LDP-law} and~\ref{lem:limit-law} in Section~\ref{subsub:proof}.	
 		Let $\Piv_\eps$ (resp., $\wh \Piv_\eps$) be the set of $\eps$-pivotal points of $(\bh,\Gamma)$ (resp., $(\bh,\wh\Gamma)$). Let $\Piv^\mesh_\eps$ and $\wh\Piv^\mesh_\eps$ be their counterpart for $(\bh,\omega^\mesh)$ and $(\bh,\wh\omega^\mesh)$, respectively.  
 		The following lemma is extracted from \cite[Section~8.3]{bhs-site-perc}.\footnote{With $c_{\op p}$ as in Proposition~\ref{prop:flip}, $c_{\op p}\nu^{\eps}_{\bh, \wh \Gamma}$ is the measure  $\wh \nu_{\op{D},\eps}$  in \cite[Section~8.3]{bhs-site-perc}.}
 		\begin{lemma}\label{lem:hatnu}
 			For $\eps>0$, there exists a measure  $\nu^{\eps}_{\bh, \wh \Gamma}$ supported on $\wh \Piv_\eps$ such that for each fixed $\eps'>0$, $\nu^{\eps}_{\bh, \wh \Gamma}=\nu^{\eps'}_{\bh,\Gamma}$ on $\Piv_{\eps'}\cap \wh \Piv_\eps$ a.s.
 		\end{lemma}
 		Since $\cup_{\eps>0}\Piv_\eps=\cup_{\eps>0}\wh \Piv_\eps$ almost surely, Lemma~\ref{lem:hatnu} characterizes $\wh \nu^\eps_{\bh,\Gamma}$ modulo a probability-zero event.
 		Recall that in Section~\ref{subsub:proof} we view $\cM^\eps(\fh,\Gamma)$ as a measurable function from $H^{-1}(\D)\times \cL(\D)$ to the spaces of Borel measures on $\D$.
 		In particular, the measure $\cM^\eps(\fh,\wh\Gamma)$ is well-defined and supported on $\wh\Piv_\eps$. 
 		In light of  Definition~\ref{def:pivm} and the discussion below it, we set $\cM^\eps_{\bh,\wh \Gamma}\defeq \xi_\fh(\D)^{-\frac12}\cM^\eps(\fh,\wh\Gamma)$.

 		\begin{lemma}\label{eq:piv-equi}
 			$\nu^{\eps}_{\bh, \wh \Gamma}= \constp\cM^\eps_{\bh,\wh \Gamma}$ a.s.\ with the constant $\constp$ as in Proposition~\ref{prop:quantum-Mink}.
 		\end{lemma}
 		\begin{proof}
 			By a reweighting consideration as in  the proof of Lemmas \ref{prop:eps-LDP-law} and~\ref{lem:limit-law} in Section~\ref{subsub:proof},
 			the $\nu_\mesh$-measure restricted to $\wh\Piv^\mesh_\eps$ converge in probability 
 			to $\cM^\eps(\fh,\wh\Gamma)$.
 			For $\eps,\eps'>0$,  
 			the $\nu_\mesh$-measure restricted to $\Piv^\mesh_{\eps'}\cap \wh \Piv^\mesh_\eps$ converge in probability  to both
 			$\cM^{\eps'}(\fh,\Gamma)|_{\Piv_{\eps'}\cap \wh \Piv_\eps}$ and $\cM^{\eps}(\fh,\wh\Gamma)|_{\Piv_{\eps'}\cap \wh \Piv_\eps}$.
 			The first convergence  can be shown by the same argument as in Proposition~\ref{prop:pivm}.
 			The second convergence follows from the first one and the stationarity of $\{J^{\eps,\mesh}_k\}_{k\in\N}$ under the measure $\wt\P^\rho_\mesh$ for arbitrarily small $\rho>0$.
 			Therefore $\cM^{\eps'}(\fh,\Gamma)=\cM^{\eps}(\fh,\wh\Gamma)$ on $\Piv_{\eps'}\cap \wh \Piv_\eps$. Lemma~\ref{eq:piv-equi} now follows from 
 			\eqref{eq:fields}, Definition~\ref{def:pivm} and Proposition~\ref{prop:quantum-Mink}.
 		\end{proof}
 		
 		Now let us consider $(\fh,\Gamma,\wh \Gamma)$ on the probability space $(\Omega,\cF,\Pd)$, where $(\D,\bh,1)$ is a $\sqrt{8/3}$-LQG disk.
 		For $n\in\N$, let $(\cM^n,\dlp^n)$ be as in Theorem~\ref{emb-thm:upper}.
 		Let $\bz^n$ be a uniformly sampled $\eps$-pivotal point of $\dlp^n$ and let $\wh \dlp^n$  be the loop ensemble obtained by flipping the color of $\bz^n$. 
 		Let $\nu^\eps_n$ and $\wh \nu^\eps_n$ be $n^{-1/4}$ times the counting measure of $\eps$-pivotal points of  $\dlp^n$ and $\wh \dlp^n$,  respectively. 
 		We view $(\cM^n,\dlp^n,\wh \dlp^n,\nu^\eps_n,  \wh \nu^\eps_n)$ as a metric space decorated with one boundary curve,  two loop ensembles, and three measures. In the continuum, similarly as $(\D,\bh,\Gamma)$ in Remark~\ref{rmk:CLE},  
 		we view $(\D, \bh,\Gamma,\wh\Gamma,\nu^\eps_{\bh,\Gamma}, \nu^{\eps}_{\bh, \wh \Gamma})$ 
 		as a metric space with the same kind of decorations. We straightforwardly extend the GHPUL distance in Section~\ref{subsec:pre} to this setting. 
 		\xncomment{With these notations, the following proposition is a restatement of~\cite[Proposition~6.4]{ghs-metric-peano}.}
 		\begin{proposition}\label{prop:flip}
 			In the setting right above, there exists a constant $c_{\op p}>0$ satisfying the following.
 			For each $\eps>0$,  there exists a coupling of $(\cM^n,\dlp^n,\wh \dlp^n)_{n\in\N}$ and $(\bh, \Gamma,\wh \Gamma)$ such that almost surely $(\cM^n,\dlp^n,\wh \dlp^n,\nu^\eps_n,  \wh \nu^\eps_n)$ converge to  $(\D, \bh,\Gamma,\wh\Gamma,c_{\op p}\nu^\eps_{\bh,\Gamma}, c_{\op p}\nu^{\eps}_{\bh, \wh \Gamma})$
 			in  the GHPUL topology.   
 		\end{proposition}
 		\xncomment{In the original statement of~\cite[Proposition~6.4]{ghs-metric-peano}, the limiting measures 
 			$c_{\op p}\nu^\eps_{\bh,\Gamma}$ and $ c_{\op p}\nu^{\eps}_{\bh, \wh \Gamma}$ are as defined using the terminology of~\cite{bhs-site-perc}; see~\cite[Proposition 8.12]{bhs-site-perc}.
 			As explained in Remark~\ref{rmk:equiv} and Lemmas~\ref{lem:hatnu} and~\ref{eq:piv-equi}, these measures agree with the ones considered in our paper.
 			Therefore Proposition~\ref{prop:flip} is equivalent to~\cite[Proposition~6.4]{ghs-metric-peano}.}
 		
 		\begin{proof}[Proof of Lemma~\ref{lem:two-jump}]
 			By Proposition \ref{prop:flip}, $\lim_{n\to\infty}\nu^\eps_n(\cM^n)=c_{\op p}\nu^\eps_{\bh,\Gamma}(\D)=\constp c_{\op p}\pivm^\eps(\D)$  and $\lim_{n\to\infty}\wh \nu^\eps_n(\cM^n)=c_{\op p}\nu^\eps_{\bh,\wh\Gamma}(\D)=\constp c_{\op p}\cM^\eps_{\bh,\wh \Gamma}(\D)$, by Proposition~\ref{prop:quantum-Mink} and Lemma~\ref{eq:piv-equi}.
 			The two-step sampling procedure in Section~\ref{subsub:Markov} applies to $(\cM^n,\dlp^{\eps,n}_t)_{t\ge 0}$. 
 			Recall the definition of $(Y^\eps_t)_{t\ge 0}$.   Lemma~\ref{lem:two-jump}  follows from the sampling recipe for $(Y^\eps_t)_{t\ge 0}$ prescribed by  Lemma~\ref{lem:limit-law}.
 		\end{proof}

 \subsection{Stability of the cutoff and proof of Propositions~\ref{prop:eps-LDP-proba} and \ref{prop:quad-conv}}\label{subsec:stability}
 In this section, we identify a site percolation configuration on $\D^\mesh$ with an element in $\cH(\D)$ (see Section~\ref{subsec:LDP}) as needed.
 We will first show that $(\Omega,\cF,\P)$ in Section~\ref{subsub:flip-lattice} satisfies Proposition~\ref{prop:eps-LDP-proba}, and then prove Propositions~\ref{prop:quad-conv}.
 
 Our proofs rely on some stability results established in \cite{gps-dynamic,ghss18}, asserting that  the importance of a vertex is rather stable in time.
 Before stating them formally, we point out that our definition of $\rho$-important pivotal points is slightly different from the definition in \cite{gps-dynamic,ghss18}. 
 In \cite{gps-dynamic} $\rho$-important pivotal points are defined in terms of how far the alternating four arms starting at the pivotal point can reach. 
 For a square $\cB$, recall the annulus $A=A_\cB$ in Section~\ref{subsub:rho-important}. Our notion of $A$-important point agrees with the one in \cite{gps-dynamic,ghss18} as long as $A\subset\D$. There is a small deviation in definition when $A\cap\p D\neq \emptyset$, but this is irrelevant as the results we will use from \cite{gps-dynamic,ghss18} are about $\rho$-important points in $r\D$ with $r\in(0,1)$.
 In this case, as explained in \cite[Section~4.7]{gps-pivotal}, these two notions of $\rho$-importance are effectively equivalent.
 In particular, the results we will be relying on hold for both notions.
 
 Having the clarification above, the following stability result is an immediate consequence of \cite[Lemma 3.7]{ghss18} and \cite[Proposition 3.9]{gps-dynamic}. 
 \xncomment{Intuitively, it says that with probability $1-o_\rho(1)$ uniform in $\mesh$,  the influence to the quad-crossing configuration of the updates in  the dynamics $(\omega^\mesh_t)_{t\ge 0}$ is captured by 
 	the updates on the $\rho$-important points.}

 \begin{lemma}
 	\label{lem:stability}
 	Fix $T>0$ and $r\in(0,1)$. Let $ X_{\mesh}$ be the set of vertices on $\D^\mesh$ which are updated for the dynamics $(\omega^\mesh_t)_{t\in[0,T]}$.  
 	Let $\Omega_\mesh$ be the set of percolation configurations $\omega'$ on $\D^\mesh$ such that $\omega'(v)=\omega^\mesh_0(v)$ for all $v\notin X_{\mesh}$.  
 	Let $\cP^\rho_\mesh$ be the set of $\rho$-important points for $\omega^\mesh_0$.  Given  $\omega',\omega''\in\cH(\D)$, let $d_r(\omega',\omega'')$ be the $d_\cH$-distance of the restriction of $\omega'$ and $\omega''$ to $\cQ_{r\D}$. For all $\zeta\in (0,1)$, there exists constants $\rho_1>0$ and $\mesh_0>0$ depending only on $r$, $T$, and $\zeta$ such that for all $\rho\in (0,\rho_1)$ and $\mesh \in (0,\mesh_0)$,
 	\begin{equation*}
 	\BB P\Big[\max\{d_r(\omega', \omega'')\,:\,  \omega'(v)=\omega''(v)\text{\,\,for\,\,} v\in \cP^\rho_\mesh\;\textrm{and}\; \omega',\omega'' \in \Omega_\mesh \} >\zeta\Big] <\zeta.
 	\end{equation*}
 \end{lemma}

 We also need the following variant of stability which is also essentially from \cite{ghss18}. 
 \xncomment{In terms of notations in  Lemma~\ref{lem:stability}, it says that for each $\rho>0$,  with probability $1-o_s(1)$ uniform in $\mesh$, each $\rho$-important point in $\cP^\rho_\mesh$ that is updated during $[0,T]$ remains $s$-important with respect to all configurations in $\Omega_\mesh$.} 
 \begin{lemma}\label{lem:Z}
 	In the setting of Lemma~\ref{lem:stability}, with $r\in(0,1)$ and $\rho, T>0$ fixed, let
 	\begin{align*}
 	Z_\mesh(v)&\defeq\inf\{\rho'>0:  \exists \omega'\in\Omega_\mesh\;\textrm{such\,that}\; v\;\textrm{is $\rho'$-important for}\; \omega' \} \textrm{ for } v\in r\D;\\
 	N_\mesh(\rho, s)&\defeq\#\{ v\in \cP^\rho_\mesh\cap X_\mesh\cap r\D: Z_\mesh(v)\le s\}\textrm{ for } s>0, \textrm{ where $\#$ means the cardinality.}
 	\end{align*} 
 	Then for all $\zeta\in (0,1)$,  there exist constants 
 	$s>0$ and $\mesh_0>0$ depending only on $\rho, r, T,\zeta$, such that \(\BB P[N_\mesh(\rho,s) =0 ] >1-\zeta\) for all $\mesh \in (0,\mesh_0)$.
 \end{lemma}
 \begin{proof}
 	By \cite[Lemma 3.5]{ghss18}, there  exists an almost surely finite random number $C(\fh,T)$, such that for every $\mesh,s,\rho$ satisfying
 	$2\,\mesh<s<2^4\,s<\rho\le 1$ and every vertex $v\in \D^\mesh\cap r\D$,
 	\begin{equation*}
 	\BB P\Big[v\in \cP^\rho_\mesh, Z_\mesh(v)\le s \mid  \fh\Big] \le C(\fh, T)s^\beta \alpha^\mesh_4(\mesh,\rho),
 	\end{equation*}	
 	where $\beta>0$ is a constant and $\alpha_4^\mesh(\mesh,\cdot)$ is defined as above Theorem~\ref{thm:LDP}. Therefore
 	\begin{align*}
 	&\E[N_\mesh(\rho,s)\mid \fh] = \sum_{v\in \D^\mesh\cap r\D} \P\left[ v\in \cP^\rho_\mesh\cap X_\mesh, Z_\mesh(v)\le s \mid  \fh\right]\\
 	&\quad\le C(\fh, T)s^\beta \alpha^\mesh_4(\mesh,\rho) \E[\#(X_\mesh\cap r\D)\mid \fh] 
 	\le \sum_{v\in \D^\mesh\cap r\D} 
 	C(\fh,T)s^\beta \alpha^\mesh_4(\mesh,\rho)\cdot T\mu'_{\fh} (v) \alpha_4^\mesh(\mesh,1)^{-1}.
 	\end{align*}
 	Here we recall that $\mu'_\fh(v)$ is the $\mu'_\fh$-mass of the hexagon corresponding to $v$ in the dual lattice of $\D^\mesh$.
 	By  the quasi-multiplicativity of $\alpha_4^\mesh(\cdot,\cdot)$ (see e.g.\ \cite{smirnov-werner-percolation}),  $\alpha_4^\mesh(\mesh,1)^{-1}\alpha^\mesh_4(\mesh,\rho)\le c\rho^{5/4}$, so $\alpha_4^\mesh(\mesh,1)^{-1}\alpha^\mesh_4(\mesh,\rho)$ is upper bounded by a constant $\wh c$ only depending on $\rho$. 
 	Therefore \(\E[N_\mesh(\rho,s)\mid \fh] \le  \wh cT\mu'_\fh(\D)C(\fh,T)s^\beta.\)
 	Now Lemma~\ref{lem:Z} follows from Markov's inequality.
 \end{proof}
 
 \xncomment{We are now ready to prove Proposition~\ref{prop:eps-LDP-proba}, which upgrades the convergence in law of $(\Gamma^{\eps,\mesh}_t)_{t\ge 0}$ from Lemma~\ref{prop:eps-LDP-law} to convergence in probability simultaneously for all $\eps>0$. In the proof we will retain the setting from Section~\ref{subsub:proof}, where we proved Lemmas \ref{prop:eps-LDP-law} and \ref{lem:limit-law}.
 	The key idea  to use the stability results above to  prove that with probability $1-o_s(1)$,  
 	if a jump of  $(\omega^{\eps,\mesh}_t)_{t\ge 0}$ occurs at location $x$ and time $t$, where $x$ is $\eps$-pivotal for $(\bh, \omega^{\eps,\mesh}_t)$,  then $x$ must be $s$-important for $\omega^{\eps,\mesh}_0$; see~\eqref{eq:eps-is-rho} below.
 }
 
 \begin{proof}[Proof of Proposition~\ref{prop:eps-LDP-proba}]
 	We claim that for $(\Omega,\cF,\P)$ as in Lemma~\ref{lem:setting}, $(\Gamma^{\eps,\mesh}_t)_{t\ge 0}$ converge in probability rather than just  in law.
 	Fix ${\rho_0}>0$. Let $\P^{\rho_0}$, $\wt\P^{\rho_0}_\mesh$, and $\wt \P^{\rho_0}$ be defined as $\P^{\rho}$, $\wt\P^{\rho}_\mesh$, and $\wt \P^{\rho}$ in Section~\ref{subsub:proof} with $\rho=\rho_0$.	We denote a jump of $(\Gamma^{\eps,\mesh}_t)_{t\ge 0}$ by $(x,t)$ where  $t$ is the jumping time and $x$ is the pivotal point being flipped at $t$.
 	For each $s>0$ and $T>0$, let $E^s_\mesh(T)$ be the event that for each jump $(x,t)$ of $(\Gamma^{\eps,\mesh}_t)_{t\ge 0}$ with $t\leq T$, if  $x$ is $\eps$-pivotal for $(\bh, \omega^{\eps,\mesh}_t)$ then $x$ is $s$-important for $\omega^{\eps,\mesh}_0$.
 	We claim that for all $\zeta\in(0,1)$, there exist $\mesh_0>0$ and $s>0$ only depending on $\zeta,T$ such that 
 	\begin{equation}\label{eq:eps-is-rho}
 	\wt\P^{\rho_0}[E^s_\mesh(T)]>1-\zeta\qquad\textrm{for all }\mesh\in(0,\mesh_0).
 	\end{equation}
 	
 	We first explain why \eqref{eq:eps-is-rho} is sufficient to conclude the proof. Let $\tau_k^\mesh$ denote the time of the $k$th jump of $(\Gamma^{\eps,\delta}_t)_{t\geq 0}$. 
 	By  Lemma~\ref{lem:jump-prob}, $(\Gamma^{\eps,\delta}_t)_{t\in[0,\tau_2^\mesh)}$ converge in $\wt\P^{\rho_0}$-probability. Let us write
 	$\wh \Gamma=\lim_{\mesh\to 0}\Gamma^{\eps,\mesh}_{\tau^\mesh_1}$ as in Lemma~\ref{lem:jump-prob}.
 	Let $z^\mesh$ be such that $(z^\mesh,\tau_2^\mesh)$ is the 2nd jump of $(\Gamma^{\eps,\mesh}_t)_{t\ge 0}$. 
 	By \eqref{eq:eps-is-rho} and the convergence of PPP$^{s}_\mesh$ for each $s$ (see Lemma~\ref{lem:setting}),  
 	$z^\mesh$ converge in $\wt\P^{\rho_0}$-probability to a point $z$.  On the other hand, the $\wt\P^{\rho_0}$-law of $\Gamma$ and $\wh\Gamma$
 	are the same. Observe that Lemma \ref{lem:flip3} applies to $(\Gamma^{\eps,\mesh}_{\tau^\mesh_1},z^\mesh)$  under $\wt\P^{\rho_0}$ by absolute continuity. 
 	Therefore $(\Gamma^{\eps,\delta}_t)_{t\in[\tau_1^\mesh,\tau_3^\mesh)}$ converge in $\wt\P^{\rho_0}$-probability. (Since $\wh \Gamma=\lim_{\mesh\to 0}\Gamma^{\eps,\mesh}_{\tau^\mesh_1}$, we in fact need Remark~\ref{rmk:flip3} here.) We can repeat the same argument to get 
 	$(\Gamma^{\eps,\delta}_t)_{t\in[\tau^\mesh_k,\tau_{k+2}^\mesh)}$ converge in $\wt\P^{\rho_0}$-probability for each $k\ge 1$.
 	This gives the  convergence in $\wt\P^{\rho_0}$-probability of $(\Gamma^{\eps,\delta}_t)_{t\geq 0}$. 
 	Therefore the same convergence holds under $\P^{\rho_0}$ if we further condition on  $\{\cN^\eps\neq 0\}$.
 	On the event $\cN^\eps=0$, the dynamic is trivial. We conclude that $(\Gamma^{\eps,\delta}_t)_{t\geq 0}$ converge in $\P^{\rho_0}$-probability.
 	Sending $\rho_0\to 0$ gives the desired convergence in Proposition~\ref{prop:eps-LDP-proba}.
 	
 	It remains to prove~\eqref{eq:eps-is-rho}.
 	We first argue that $\wt \P^{\rho_0}$ and $\wt \P^{\rho_0}_\delta$ are close in total variational distance when $\delta$ is small. 
 	For any event $E\in\cF$, we have that $\left|\E^{\rho_0}[\cN^\eps_\mesh\1_{E}] -\E^{\rho_0}[\cN^\eps\1_{E}] \right|\le\E^{\rho_0}[|\cN^\eps_\mesh-\cN^\eps|]$,
 	where $\E^{\rho_0}$ is the expectation corresponding to $\P^{\rho_0}$. 
 	By Lemma~\ref{lem:UI2}, there exists a function $\zeta^{\rho_0}(\delta)$ not depending on $E$ such that $\lim_{\delta\to0}\zeta^{\rho_0}(\delta)=0$ and 
 	\begin{equation}\label{eq:Error}
 	\Big|\wt\P^{\rho_0}_\mesh[E] -\wt\P^{\rho_0}[E]\Big|\le \zeta^{\rho_0}(\delta) \qquad\textrm{for all }E\in\cF.
 	\end{equation}		
 	
 	We now fix $K\in\N$ large enough and $\mesh_0>0$ small enough such that
 	$\wt\P_\mesh^{\rho_0}[\tau^\mesh_{K+1}>T]>1-0.1\zeta$ for $\mesh\in(0,\mesh_0)$.
 	Let $G^1_\mesh(r)$ be the event that if $(x,\tau^\mesh_k)$ is a jump of $(\Gamma^{\eps,\mesh}_t)_{t\ge 0}$ for $1\le k\le K$, then $x\in r\D$.
 	By possibly shrinking $\mesh_0$, we can find $r\in(0,1)$ such that $\wt\P^{\rho_0}_\mesh[G^1_\mesh(r)] \ge 1-0.1\zeta$ for $\mesh\in(0,\mesh_0)$.
 	For $1\le k\le K$ and $\rho>0$, let $G^2_\mesh(k;\rho)$ be the event that every $\eps$-pivotal point of $(\bh,\Gamma^{\eps,\mesh}_{\tau^\mesh_k})$ is $\rho$-important for $\Gamma^{\eps,\mesh}_{\tau^\mesh_k}$.  Set $G^2_\mesh(\rho)\defeq\cup_1^KG^2_\mesh(k;\rho)$.
 	Recall Lemma~\ref{lem:cover2}. By choosing $\rho$ small enough and possibly shrinking $\mesh_0$, 
 	we  can have $\min_{1\le k\le K}\wt\P_\mesh^{\rho_0}[G^2_\mesh(k;\rho)]\ge 1-0.1K^{-1}\zeta$ and hence $\wt\P_\mesh^{\rho_0}[G^2_\mesh(\rho)]\ge 1-0.1\zeta$
 	for $\mesh\in(0,\mesh_0)$.	
 	
 	For $i,j\in\{0,1,\cdots,K\}$, let $G_\mesh(i,j;\rho,s)$ be the event that every $\rho$-important point for $\omega^{\eps,\mesh}_{\tau^\mesh_i}$ is $s$-important for $\omega^{\eps,\mesh}_{\tau^\mesh_j}$, where we set $\tau^\mesh_0=0$. By Lemma \ref{lem:Z} and \eqref{eq:Error}, after possibly shrinking $\mesh_0$, we can find $s$ small enough such that $\wt\P^{\rho_0}_\mesh[G_\mesh(0,k;\rho, s)] \ge 1-K^{-1}0.01\zeta$ for each $1\le k\le K$. Since $\{\omega^{\eps,\mesh}_{\tau^\mesh_k}\}$ is reversible under $\wt\P^{\rho_0}$,
 	we have $\wt\P^{\xncomment{\rho_0}}_\mesh[G_\mesh(k,0;\rho, s)] \ge 1-K^{-1}0.01\zeta$ as well.
 	On the event  $\left(\{\tau^\mesh_{K+1}>T\}\cap G^1_\mesh(r)\cap G^2_\mesh(\rho)\right)\setminus E^s_\mesh(T)$, there exists $1\le k\le K$ such that $G_\mesh(k,0;\rho, s)$ does not occur.
 	Therefore $\wt\P^{\xncomment{\rho_0}}_\mesh[E^s_\mesh(T)]\ge 1-0.5\zeta$. By possibly shrinking $\mesh_0$ such that  $ \zeta^{\xncomment{\rho_0}}(\delta_0)<0.5\zeta$, we get \eqref{eq:eps-is-rho}  from \eqref{eq:Error}.
 	\qedhere
 \end{proof}
 
 The following lemma is the key to the proof of  Proposition~\ref{prop:quad-conv}.
 \xncomment{It says that for a fixed $\rho>0$, with probability $1-o_\eps(1)$ uniform in $\mesh$, 
 	each point in   $\cP^\rho_\mesh$ that is updated during $[0,T]$ must also be updated for the $\eps$-cutoff dynamic $(\omega^{\eps,\mesh}_t)_{t\ge 0}$.
 	The main idea is to use Lemma~\ref{lem:cover3} to show that for a fixed $s>0$, with probability $1-o_\eps(1)$ uniform in $\mesh$, if a vertex is $s$-important for  $(\omega^{\eps,\mesh}_t)_{t\ge 0}$ at the time it rings, it must be $\eps$-pivotal at the same time; see~\eqref{eq:N'}.
 	Once this is done, together with  Lemma~\ref{lem:Z} we will conclude the proof.}
 \begin{lemma}\label{lem:eps-rho}
 	In the setting of Lemma~\ref{lem:stability}, for each $\eps>0$, let $X^\eps_\mesh$ be set of vertices on $\D^\mesh$ where update occurs for the dynamic $(\omega^{\eps,\mesh}_t)_{t\in[0,T]}$. Then for all $\zeta,\rho\in(0,1)$, there exists $\eps>0$ and $\mesh_0>0$ depending only on $\rho, r, T,\zeta$ such that \(\BB P[\cP^\rho_\mesh\cap X_\mesh \subset X^\eps_\mesh]>1-\zeta\) for $\mesh\in(0,\delta_0)$.
 \end{lemma}
 \begin{proof}	
 	Suppose we are in the setting of Lemma~\ref{lem:stability} with $r\in(0,1)$ and $T>0$ fixed. For each $v\in \D^\mesh$, let $\tau_v$ be the time when the clock of $v$ rings for the first time so that $X_\mesh=\{v\in \D^\mesh: \tau_v\le T  \}$. For $s>0$ and $\eps>0$, let 
 	\[
 	N'_\mesh(s,\eps) \defeq\#\{v\in X_\mesh\cap r\D: \textrm{$v$ is $s$-important for $\omega^{\eps,\mesh}_{\tau_v}$ but not $\eps$-pivotal for $(\bh, \omega^{\eps,\mesh}_{\tau_v})$}\}.
 	\]
 	We claim that  for all $\zeta\in (0,1)$,  there exists $\eps>0$ and $\mesh_0>0$ depending only on $s, r, T,\zeta$, such that
 	\begin{equation}\label{eq:N'}
 	\BB P[N'_\mesh(s,\eps) =0 ] >1-\zeta/3\qquad \textrm{for }\mesh \in (0,\mesh_0).
 	\end{equation}	
 	
 	Given \eqref{eq:N'}, we first choose $s$ such that $\P[N_\mesh(\rho,s)=0]>1-\zeta/3$ with $N_\mesh(\rho,s)$ as defined in Lemma~\ref{lem:Z}. Then we choose $\eps$ such that 
 	$\P[N'_\mesh(s,\eps)=0]> 1-\zeta/3$. Let $E_\mesh$ be the event that  the clock at each $\rho$-important \xncomment{vertex} in $r\D$ rings at most once.
 	By a first moment calculation and possibly shrinking $\delta_0$ depending on $\zeta$, we can have $\P[E_\mesh]\ge 1-\zeta/3$  for $\delta\in(0,\delta_0)$.
 	\xncomment{More precisely, given each $v\in \D^\mesh \cap r\D$, the probability that $v$ is a $\rho$-important and the clock at $v$ rings at least twice during $[0,T]$ is of order $o(\delta^2)$. Therefore the expected number of such vertices is of order $o_\mesh(1)$. Hence by Markov's inequality $\P[E_\mesh]\ge 1-\zeta/3$  for $\delta\in(0,\delta_0)$.}
 	On $E_\mesh\cap\{N_\mesh(\rho,s)=0, \; N'_\mesh(s,\eps)=0 \} $, 
 	each $v $ in $\cP^\rho_\mesh\cap X_\mesh$ must be $s$-important for  $\omega^{\eps,\mesh}_{\tau_v}$, hence be $\eps$-pivotal for   $(\bh,\omega^{\eps,\mesh}_{\tau_v})$.
 	Therefore $v\in  X^\eps_\mesh$, which concludes the proof of Lemma~\ref{lem:eps-rho}.

 	It remains to prove~\eqref{eq:N'}.
 	Fix $v\in \D^\mesh \cap r\D$. Given a percolation configuration $\omega$ on $D^\mesh$, whether $v$ is $\eps$-pivotal  for $(\bh,\omega)$ only depends on $\omega|_{\D^\mesh\setminus\{v\}}$ and $\bh$.
 	The same statement holds for $s$-importance without involving $\bh$.
 	For $t\ge 0$, let $S^v_t(u)=\omega^{\eps,\mesh}_t(u)$ for $u\in\D^\mesh\setminus\{v\}$ and $S^v_t(v)=\omega^{\eps,\mesh}_0(v)$.
 	In other words, $(S^v_t)_{t\geq 0}$ is the same dynamics as $(\omega^{\eps,\mesh}_t)_{t\ge 0}$ except that the color of $v$ never changes.
 	Then $\tau_v$ is independent of $(S^v_t)_{t\geq 0}$. Note that $(S^v_t)_{t\geq 0}$ is still stationary. Thus $S^v_{\tau_v}$ has the same law as $\omega^\mesh_0$.
 	Fix $\zeta'\in(0,1)$ to be determined later and choose $\eps$ and $\mesh_0\in(0,0.1)$ such that Lemma~\ref{lem:cover3} holds with $s,\zeta'$ here.
 	Since $S^v_{\tau_v}$ and $\omega^{\eps,\mesh}_{\tau_v}$ agree on $\D^\mesh\setminus\{v\}$, for $\mesh\in(0,\mesh_0)$, we have 
 	\begin{align*}
 	&\P\left[ v\in X_\mesh\cap r\D, \textrm{$v$ is $s$-important for $\omega^{\eps,\mesh}_{\tau_v}$ but not $\eps$-pivotal for $(\bh,\omega^{\eps,\mesh}_{\tau_v})$}\right]\\
 	&\qquad=\P\left[\tau_v\le T, \textrm{$v$ is $s$-important for $S^v_{\tau_v}$ but not $\eps$-pivotal for $(\bh,S^v_{\tau_v})$}\right]\\
 	&\qquad=\P[\tau_v\le T]\P\left[\textrm{$v$ is $s$-important for $S^v_{\tau_v}$ but not $\eps$-pivotal for $(\bh,S^v_{\tau_v})$}\right]\\
 	&\qquad=\P[\tau_v\le T]\P\left[\textrm{$v$ is $s$-important for $\omega^\mesh_0$ but not $\eps$-pivotal for $(\bh,\omega^\mesh_0)$}\right]\\
 	&\qquad\le \P[\tau_v\le T]\P\left[\textrm{$v$ is $s$-important for $\omega^\mesh_0$}\right] \zeta'\\
 	&\qquad=\zeta'  \P[\tau_v\le T,\;\textrm{$v$ is $s$-important for $\omega^\mesh_0$}].
 	\end{align*}
 	The purpose of introducing $S_t^v$ can be seen in the third step of this equality, where we use the independence of two events.
 	By the definition of $N'_\mesh(s,\eps)$,
 	\begin{align*}
 	&\E[N'_\mesh(s,\eps)] = \sum_{v\in\D^\mesh\cap r\D} \P\left[ v\in X_\mesh\cap r\D, \textrm{$v$ is $s$-important for $\omega^{\eps,\mesh}_{\tau_v}$ but not $\eps$-pivotal for $(\bh,\omega^{\eps,\mesh}_{\tau_v})$}\right]\\
 	&\quad\le \zeta' \sum_{v\in\D^\mesh\cap r\D}  \P[\tau_v\le T,\;\textrm{$v$ is $s$-important for $\omega^\mesh_0$}]=\zeta'\E[ \#(\cP^s_\mesh\cap X_\mesh\cap r\D)]
 	\le \zeta'T\E[\nu_\mesh(\cP^s_\mesh\cap r\D)].
 	\end{align*}
 	\xncomment{By	Lemma~\ref{lem:UI}, $\lim_{\mesh\to0}\nu_\mesh^s(\D)$ exists in $L^1$, where $\nu^s_\mesh$ is the restriction of $\nu_\mesh$ to $\cP^s_\mesh$.
 		This yields that  $\max_{\delta\in(0,0.1)}\E[\nu^\mesh(\cP^s_\mesh\cap(r\D))]<\infty$.}
 	Therefore	we can choose $\zeta'$ small enough depending on $s, r, T,\zeta$ such that 
 	$\max_{\delta\in(0,0.1)}\E[N'_\mesh(s,\eps)]\le \zeta/3$.  Now  \eqref{eq:N'} follows from Markov's inequality.
 \end{proof}

 \begin{figure}
 	\centering
 	\includegraphics[scale=1]{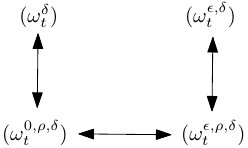} 
 	\caption{Illustration of the proof of Proposition~\ref{prop:quad-conv}. By Lemma \ref{lem:stability}, for $\rho$ sufficiently small we know that the two processes to the left (resp., right) connected by a vertical arrow are close with high probability for the metric $d_r$ at any time $t\in[0,T]$. By Lemma~\ref{lem:eps-rho}, we know that $(\orhodelta_t)_{t\in[0,T]}=(\epsrhodelta_t)_{t\in[0,T]}$ with high probability for $\eps$ sufficiently small compared to $\rho$.} 
 	\label{emb-fig:epsrho}
 \end{figure}
 
 \begin{proof}[Proof of Proposition~\ref{prop:quad-conv}]
 	We refer to Figure~\ref{emb-fig:epsrho} for an illustration of the proof. 
 	Let $(\orhodelta_t)_{t\ge0}$  be defined just as $(\omega^\mesh_t)_{t\ge 0}$  except that when the clock at a vertex $v$ rings,  we do not flip its color unless $v\in\cP_\mesh^\rho$. 
 	We define $(\epsrhodelta_t)_{t\ge 0}$ similarly with $(\omega^{\eps,\mesh}_t)_{t\ge 0}$ in place of $(\omega^\mesh_t)_{t\ge 0}$.
 	More precisely, if the clock at a vertex $v$ rings at some time $t$, the color of $v$ is flipped along the $(\epsrhodelta_t)_{t\ge 0}$ dynamic if and only if
 	$v\in \cP_\mesh^\rho$ and  $v$ is  an $\eps$-pivotal for $(\bh,\omega^{\eps,\mesh}_t)$. Recall $d_r$ in Lemma~\ref{lem:stability}.
 	For any $t\in[0,T]$, by the triangle inequality,
 	\eqbn
 	d_r( \epsodelta_t, \oodelta_t )
 	\leq
 	d_r( \epsodelta_t, \epsrhodelta_t )
 	+
 	d_r( \epsrhodelta_t, \orhodelta_t )
 	+
 	d_r( \orhodelta_t, \oodelta_t ).
 	\eqen
 	Fix $\zeta\in (0,1)$. Recall $\rho_1$ and $X_\mesh$ in Lemma \ref{lem:stability}. For $\rho\in (0,\rho_1)$, with probability at least $1-\zeta$,
 	\begin{equation*}
 	\max_{t\in[0,T]} 
 	d_r(\epsodelta_t, \epsrhodelta_t )	+ 	d_r( \orhodelta_t, \oodelta_t)\le 2\zeta.
 	\end{equation*}
 	Recall $X_\mesh^\eps$ in Lemma~\ref{lem:eps-rho}. 
 	On the event $\{\cP^\rho_\mesh\cap X_\mesh \subset X^\eps_\mesh\}$, we have $ (\orhodelta_t)_{t\in[0,T]}  = (\epsrhodelta_t)_{t\in[0,T]}$. 
 	By Lemma~\ref{lem:eps-rho}, this occurs with probability as least $1-\zeta$ if $\eps$ is small enough. 
 	For such $\eps$,
 	\eqb
 	\P\left[\max_{t\in[0,T]} d_r (\omega^{\eps,\mesh}_t,\omega^\mesh_t )>2\zeta\right]<2\zeta.
 	\label{emb-eq40}
 	\eqe
 	Sending $\mesh\rta 0$, we have \(\P\left[\max_{t\in[0,T]} d_r (\omega^{\eps}_t,\omega_t )>2\zeta\right]<2\zeta\), which  concludes the proof. 
 \end{proof}

\appendix
\section{Proof of Lemma~\ref{lem:GMC}}
	We will prove Lemma~\ref{lem:GMC} using ideas from \cite{shef-wang-lqg-coord}, where related results for the case of $\Mink_A$ equal to Lebesgue measure are proved.
	By the definition of a free Liouville field (Definition~\ref{def:free}), it is sufficient to consider the case where $h$ is a zero-boundary Gaussian free field. Let $\nu_r=r^{\alpha^2/2}e^{\alpha h_r}\Mink_A$. 
	By the argument in \cite[Section 6]{berestycki-gmt-elementary}, in order to prove that $e^{\alpha h}\Mink_A$ exists it is sufficient to prove that for a fixed set $U\Subset D$ (recall that $U\Subset D$ means $U\cup\p U\subset D$), $\nu_r(U)$ has an a.s.\ limit as $r\rta 0$.  
	
	Define $\ol h_r(z)=\alpha h_r(z)+\frac{\alpha^2}{2}\log r$. For any $s\in(0,r)$, 
	\eqb
	\begin{split}
		\E[ (\nu_r(U)-\nu_s(U))^2 ]
		=\iint_{U\times U} 
		\E\big[\big(e^{\ol h_r(z)}-e^{\ol h_s(z)}\big)
		\big(e^{\ol h_r(w)}-e^{\ol h_s(w)}\big)\big]
		\,d\Mink_A(z)\,d\Mink_A(w).
	\end{split}
	\label{eq:moment}
	\eqe
	Let $G:\ol U\times\ol U\to\R$ denote the Green's function and for $z\in U$ let $C(z;U)$ denote the conformal radius of $z$ in $U$. Recall that $\Var( h_r(z))=\log r^{-1}+\log C(z;U)$ and that $\Cov[h_r(z),h_s(w)]=G(z,w)$  if $|z-w|>r+s$. 
	Using these identities, we get that the integrand on the right side of \eqref{eq:moment} is zero when $|z-w|>2r$.
	Furthermore, $\ol h_s(z)-\ol h_r(z)\eqD a N - 0.5 a^2$ for a standard normal random variable $N$ and $a:=\alpha \sqrt{\log \frac rs}$, which gives that for some $c>0$
		\eqbn
		\E\big[\big(e^{\ol h_r(z)}-e^{\ol h_s(z)}\big)^2 \big] 
		= \E\big[e^{2\ol h_r(z)}\big] \cdot \E\big[\big(1-e^{\ol h_s(z)-\ol h_r(z)}\big)^2 \big] 
		= c \E\big[\big(1-e^{aN-0.5a^2}\big)^2\big] 
		%= c (e^{a^2}-1)
		%= c (e^{\alpha^2 \log \frac rs}-1)
		= c ((r/s)^{\alpha^2}-1).
		\eqen Therefore, for any $\wh d\in(0,d)$ and some constant $c>0$,
	\eqbn
	\begin{split}
		\E[ (\nu_r(U)&-\nu_s(U))^2 ]
		\leq 
		\iint_{U\times U,|z-w|<2r} 
		\E\big[\big(e^{\ol h_r(z)}-e^{\ol h_s(z)}\big)^2
		\big]
		\,d\Mink_A(z)\,d\Mink_A(w)\\
		&\leq c((r/s)^{\alpha^2}-1)\iint_{U\times U,|z-w|<2r} 
		\,d\Mink_A(z)\,d\Mink_A(w)
		\leq c((r/s)^{\alpha^2}-1)(2r)^{\wh d}\cdot 
		\iint_{U\times U} \frac{d\Mink_A(z)\,d\Mink_A(w)}{|z-w|^{\wh d}}.
	\end{split}
	\eqen
	The integral on the right side is finite by \eqref{eq:energy}.
	We see from this estimate that for any $N\in\N$, we have a.s.\ convergence of $\nu_r(U)$ as $r\rta 0$ along integer powers of $2^{-1/N}$. To obtain a.s.\ convergence as $r\rta 0$ (without requiring that $r$ is a power of $2^{1/N}$), we proceed similarly as in the proof of \cite[Theorem 1.1]{shef-wang-lqg-coord}, and the argument is therefore omitted. 
	
	We can find a small $\delta>0$ such that 
	$\E[\iint_{U\times U} \frac{d\nu(z)\,d\nu(w)}{|z-w|^{\delta}}]<\infty$. Therefore $\nu$ is a.s.\ non-atomic. 
	\begin{comment}
	We use a similar sandwich argument as in Menglu-Scott. We need to calculate a similar L2 norm as them, where we compare the 2^{-k/N} approximation to the approximation using supremum. We move the expectation inside the double integral. We consider an on-diagonal term and an off-diagonal term. When bounding the off-diagonal term, if z and w are the two points, then we consider balls of radius |z-w|/2 centered at each point. Given the field outside these two balls, the circle averages around z and w behave independently and according to the averages on the radius |z-w|/2 circles.
	\end{comment}

	\begin{comment}
	To conclude the proof, we need to show that the limiting measure $\nu$ is non-atomic. It is sufficient to show that for any $\delta>0$ we can find an $r>0$ such that no $r$-ball has mass larger than $\delta$ with probability at least $1-\delta$. Dividing the domain into $O(r^{-2})$ squares $S_i$ of diameter $r$, we have that for some $\wh d\in(0,d)$ and constant $C>0$,
	\eqbn
	\begin{split}
	\P[ \exists i \text{\,\,such\,\,that\,\,} \nu(S_i)>\delta ]
	&\leq \sum_i \P[\nu(S_i)>\delta]
	\leq \sum_i \delta^{-2}\E[\nu(S_i)^2]\\
	&\leq C\iint_{|z-w|<r} |z-w|^{-\wh d}\,d\Mink_A(z)\,d\Mink_A(w),
	\end{split}
	\eqen
	where the last inequality follows from \cite[Section 3]{berestycki-gmt-elementary}. The right side converges to 0 as $r\rta 0$ by \eqref{eq:energy}, which concludes the proof.
	\end{comment}


\begin{thebibliography}{GKMW18}
	\bibitem[ABA17]{ab-simple}
	L.~Addario-Berry and M.~Albenque.
	\newblock The scaling limit of random simple triangulations and random simple
	quadrangulations.
	\newblock {\em Ann. Probab.}, 45(5):2767--2825, 2017. \MR{3706731}
	
	\bibitem[Abr16]{abraham-bipartite}
	C.~Abraham.
	\newblock Rescaled bipartite planar maps converge to the {B}rownian map.
	\newblock {\em Ann. Inst. Henri Poincar\'e Probab. Stat.}, 52(2):575--595,
	2016, \arxiv{1312.5959}. \MR{3498001}
	
	\bibitem[AGMT16]{agmt16}
	D.~Ahlberg, S.~Griffiths, R.~Morris, and V.~Tassion.
	\newblock Quenched {V}oronoi percolation.
	\newblock {\em Adv. Math.}, 286:889--911, 2016. \MR{3415699}
	
%	\bibitem[AB99]{ab-random-curves}
%	M.~Aizenman and A.~Burchard.
%	\newblock H\"older regularity and dimension bounds for random curves.
%	\newblock {\em Duke Math. J.}, 99(3):419--453, 1999, \arxiv{math/9801027}.
%	\MR{1712629}
%	
	
	\bibitem[AHS19]{aasw-type2}
	M.~{Albenque}, N.~{Holden}, and X.~{Sun}.
	\newblock {Scaling limit of large triangulations of polygons}.
	\newblock {\em Electron. J. Probab.}, Vol 25, 1-43, (2020)
	
	\bibitem[AFS20]{afs-volume}
	M. Ang, H. Falconet, and X. Sun
	\newblock Volume of metric balls in Liouville quantum gravity.
	\newblock {\em Electron. J. Probab.}, Vol 25, No. 160, 1-50, (2020)
	
	
	\bibitem[Ang03]{angel-uipt}
	O.~Angel 
	\newblock Growth and percolation on the uniform infinite planar triangulation.
	\newblock {\em Geom. Funct. Anal.}, 13(5):935--974, 2003.
	
	\bibitem[{Ang}05]{angel-uihpq-perc}
	O.~{Angel}.
	\newblock {Scaling of Percolation on Infinite Planar Maps, I}.
	\newblock {\em ArXiv  e-prints}, December 2005,
	\arxiv{math/0501006}.
	
	\bibitem[AS03]{angel-schramm-uipt}
	O.~Angel and O.~Schramm.
	\newblock Uniform infinite planar triangulations.
	\newblock {\em Comm. Math. Phys.}, 241(2-3):191--213, 2003. \MR{2013797
		(2005b:60021)}
	
	\bibitem[BB07]{berger-biskup-perc-rw}
	N.~Berger and M.~Biskup.
	\newblock Quenched invariance principle for simple random walk on percolation
	clusters.
	\newblock {\em Probab. Theory Related Fields}, 137(1-2):83--120, 2007.
	\MR{2278453}
	
	\bibitem[Bef08]{beffara-dim}
	V.~Beffara.
	\newblock The dimension of the SLE curves.
	\newblock {\em 	Ann. Probab.}, 36 (2008), no. 4, 1421–1452.
	
	
	
	
	\bibitem[Bef07]{beffara-easy}
	V.~Beffara.
	\newblock Cardy's formula on the triangular lattice, the easy way.
	\newblock In {\em Universality and renormalization}, volume~50 of {\em Fields
		Inst. Commun.}, pages 39--45. Amer. Math. Soc., Providence, RI, 2007.
	\MR{2310300}
	
	\bibitem[Ber99]{bertoin-sub}
	J.~Bertoin.
	\newblock Subordinators: examples and applications.
	\newblock In {\em Lectures on probability theory and statistics
		({S}aint-{F}lour, 1997)}, volume 1717 of {\em Lecture Notes in Math.}, pages
	1--91. Springer, Berlin, 1999. \MR{1746300 (2002a:60001)}
	
	\bibitem[Ber07]{bernardi-dfs-bijection}
	O.~Bernardi.
	\newblock Bijective counting of {K}reweras walks and loopless triangulations.
	\newblock {\em J. Combin. Theory Ser. A}, 114(5):931--956, 2007.
	
	\bibitem[Ber17]{berestycki-gmt-elementary}
	N.~Berestycki.
	\newblock An elementary approach to {G}aussian multiplicative chaos.
	\newblock {\em Electron. Commun. Probab.}, 22, Paper No. 27, 2017.
	\MR{3652040}
	
	\bibitem[BHS18]{bhs-site-perc}
	O.~{Bernardi}, N.~{Holden}, and X.~{Sun}.
	\newblock {Percolation on triangulations: a bijective path to Liouville quantum
		gravity}.
	\newblock {\em ArXiv e-prints}, July 2018, \arxiv{1807.01684}. To appear in {\em Memoirs of the American Mathematical Society}.
	
	\bibitem[BJM14]{bjm-uniform}
	J.~Bettinelli, E.~Jacob, and G.~Miermont.
	\newblock The scaling limit of uniform random plane maps, {\it via} the
	{A}mbj\o rn-{B}udd bijection.
	\newblock {\em Electron. J. Probab.}, 19:no. 74, 16, 2014.
	\MR{3256874}
	
	\bibitem[BLG13]{beltran-legall-pendant}
	J.~Beltran and J.-F. Le~Gall.
	\newblock Quadrangulations with no pendant vertices.
	\newblock {\em Bernoulli}, 19(4):1150--1175, 2013.
	\MR{3102547}
	
	\bibitem[BM17]{bet-mier-disk}
	J.~Bettinelli and G.~Miermont.
	\newblock Compact {B}rownian surfaces {I}: {B}rownian disks.
	\newblock {\em Probab. Theory Related Fields}, 167(3-4):555--614, 2017. \MR{3627425}
	
	\bibitem[Car92]{cardy-formula}
	J.~L. Cardy.
	\newblock Critical percolation in finite geometries.
	\newblock {\em J. Phys. A}, 25(4):201--206, 1992.
	
	\bibitem[CK15]{curien-kortchemski-looptree-perc}
	N.~Curien and I.~Kortchemski.
	\newblock Percolation on random triangulations and stable looptrees.
	\newblock {\em Probab. Theory Related Fields}, 163(1-2):303--337, 2015. \MR{3405619}
	
	\bibitem[CN06]{camia-newman-sle6}
	F.~Camia and C.~M. Newman.
	\newblock Two-dimensional critical percolation: the full scaling limit.
	\newblock {\em Comm. Math. Phys.}, 268(1):1--38, 2006. \MR{2249794}
	
	\bibitem[Cur15]{curien-glimpse}
	N.~Curien.
	\newblock A glimpse of the conformal structure of random planar maps.
	\newblock {\em Comm. Math. Phys.}, 333(3):1417--1463, 2015.
	\MR{3302638}
	
	\bibitem[DDDF19]{DDDF19}
	J.~{Ding}, J.~{Dub{\'e}dat}, A.~{Dunlap}, and H.~{Falconet}.
	\newblock {Tightness of Liouville first passage percolation for $\gamma \in
		(0,2)$}.
	\newblock {\em Publications mathematiques de l'IHES},  132, 353--403, 2020.
	
	\bibitem[DKRV16]{dkrv-lqg-sphere}
	F.~David, A.~Kupiainen, R.~Rhodes, and V.~Vargas.
	\newblock Liouville quantum gravity on the {R}iemann sphere.
	\newblock {\em Comm. Math. Phys.}, 342(3):869--907, 2016.
	\MR{3465434}
	
	\bibitem[DMS14]{wedges}
	B.~{Duplantier}, J.~{Miller}, and S.~{Sheffield}.
	\newblock {Liouville quantum gravity as a mating of trees}.
	\newblock {\em ArXiv e-prints}, September 2014, \arxiv{1409.7055}.
	
	\bibitem[DS11]{shef-kpz}
	B.~Duplantier and S.~Sheffield.
	\newblock Liouville quantum gravity and {KPZ}.
	\newblock {\em Invent. Math.}, 185(2):333--393, 2011.
	\MR{2819163 (2012f:81251)}
	
	\bibitem[Dub09]{dubedat-duality}
	J.~Dub{\'e}dat.
	\newblock Duality of {S}chramm-{L}oewner evolutions.
	\newblock {\em Ann. Sci. \'Ec. Norm. Sup\'er. (4)}, 42(5):697--724, 2009.
	\MR{2571956 (2011g:60151)}
	
	\bibitem[GHPR19]{c-greater-than-1}
	E.~Gwynne, N.~Holden, J.~Pfeffer, and G.~Remy.
	\newblock Liouville quantum gravity with central charge in $(1, 25) $: a
	probabilistic approach.
	\newblock {\em Commum.\ Math.\ Phys.}, 376, 1573--1625, 2020.
	
	\bibitem[GHS19a]{ghs-metric-peano}
	E.~{Gwynne}, N.~{Holden}, and X.~{Sun}.
	\newblock {Joint scaling limit of site percolation on random triangulations in
		the metric and peanosphere sense}.
	\newblock {\em ArXiv e-prints},  May 2019, \arxiv{1905.06757}.
	
	\bibitem[GHS19b]{ghs-mot-survey}
	E.~{Gwynne}, N.~{Holden}, and X.~{Sun}.
	\newblock {Mating of trees for random planar maps and Liouville quantum
		gravity: a survey}.
	\newblock {\em ArXiv e-prints}, October 2019, \arxiv{1910.04713}.
	
	\bibitem[GHSS19]{ghss18}
	C.~{Garban}, N.~{Holden}, A.~{Sep{\'u}lveda}, and X.~{Sun}.
	\newblock {Liouville dynamical percolation}.
	\newblock {\em ArXiv e-prints},  May 2019, \arxiv{1905.06940}. To appear in Probab. Theory Related Fields.
	
	\bibitem[GM17a]{gwynne-miller-perc}
	E.~{Gwynne} and J.~{Miller}.
	\newblock {Convergence of percolation on uniform quadrangulations with boundary
		to SLE$_{6}$ on $\sqrt{8/3}$-Liouville quantum gravity}.
	\newblock {\em ArXiv e-prints}, January 2017, \arxiv{1701.05175}.
	
	\bibitem[GM17b]{gwynne-miller-uihpq}
	E.~Gwynne and J.~Miller.
	\newblock Scaling limit of the uniform infinite half-plane quadrangulation in
	the {G}romov-{H}ausdorff-{P}rokhorov-uniform topology.
	\newblock {\em Electron. J. Probab.}, 22:Paper No. 84, 47, 2017. \MR{3718712}
	
	\bibitem[GM18]{gwynne-miller-sle6}
	E.~Gwynne and J.~Miller.
	\newblock Chordal {${\rm SLE}_6$} explorations of a quantum disk.
	\newblock {\em Electron. J. Probab.}, 23:Paper No. 66, 24, 2018. \MR{3835472}
	
	\bibitem[GM19]{GM-metric}
	E.~{Gwynne} and J.~{Miller}.
	\newblock {Existence and uniqueness of the Liouville quantum gravity metric for
		$\gamma \in (0,2)$}.
	\newblock {\em Invent. Math.}, 223:213--333, 2021.
	
	\bibitem[GMS17b]{gms-tutte}
	E.~{Gwynne}, J.~{Miller}, and S.~{Sheffield}.
	\newblock {The Tutte embedding of the mated-CRT map converges to Liouville
		quantum gravit}.
	\newblock {\em Ann. Probab.}, 49(4):1677--1717, 2021.
	
	
	\bibitem[GMS18a]{gms-random-walk}
	E.~Gwynne, J.~Miller, and S.~Sheffield.
	\newblock {An invariance principle for ergodic scale-free random environments}.
	\newblock {\em ArXiv e-prints}, July 2018, \arxiv{1807.07515}.
	
	\bibitem[GMS18b]{gms-voronoi}
	E.~Gwynne, J.~Miller, and S.~Sheffield.
	\newblock {The Tutte embedding of the Poisson-Voronoi tessellation of the
		Brownian disk converges to $\{ 8/3 \}$-Liouville quantum gravity}.
	\newblock {\em Commun. Math. Phys.}, 374, 735--784, 2020.
	
	\bibitem[GP19]{gp-kpz}
	E.~{Gwynne} and J.~{Pfeffer}.
	\newblock {KPZ formulas for the Liouville quantum gravity metric}.
	\newblock {\em ArXiv e-prints}, May 2019, \arxiv{1905.11790}.
	
	\bibitem[GPS10]{gps-fourier}
	C.~Garban, G.~Pete, and O.~Schramm.
	\newblock The {F}ourier spectrum of critical percolation.
	\newblock {\em Acta Math.}, 205(1):19--104, 2010. \MR{2736153}
	
	\bibitem[GPS13]{gps-pivotal}
	C.~Garban, G.~Pete, and O.~Schramm.
	\newblock Pivotal, cluster, and interface measures for critical planar
	percolation.
	\newblock {\em J. Amer. Math. Soc.}, 26(4):939--1024, 2013. \MR{3073882}
	
	\bibitem[GPS18a]{gps-dynamic}
	C.~Garban, G.~Pete, and O.~Schramm.
	\newblock The scaling limits of near-critical and dynamical percolation.
	\newblock {\em J. Eur. Math. Soc.}, 20(5):1195--1268, 2018. \MR{3790067}
	
	\bibitem[GPS18b]{gps-mst}
	C.~Garban, G.~Pete, and O.~Schramm.
	\newblock The scaling limits of the minimal spanning tree and invasion
	percolation in the plane.
	\newblock {\em Ann. Probab.}, 46(6):3501--3557, 2018. \MR{3857861}
	
	
	\bibitem[HLLS18]{hlls-cut}
	N.~{Holden}, G.~F. {Lawler}, X.~{Li}, and X.~{Sun}.
	\newblock {Minkowski content of Brownian cut points}.
	\newblock {\em ArXiv e-prints}, March 2018, \arxiv{1803.10613}.
	
	\bibitem[HLS18]{hlls-pivot}
	N.~{Holden}, X.~{Li}, and X.~{Sun}.
	\newblock {Natural parametrization of percolation interface and pivotal
		points}.
	\newblock {\em ArXiv e-prints}, April 2018, \arxiv{1804.07286}.
	
	\bibitem[JS00]{jones-smirnov-removability}
	P.W.~{Jones} and S.K.~{Smirnov}.
	\newblock {Removability theorems for Sobolev functions and quasiconformal maps}.
	\newblock {\em Ark. Mat.}, 38(2):263--279, 2000.
	\MR{1785402}
	
	\bibitem[Khr18]{khristoforov-thesis}
	M.~Khristoforov.
	\newblock {\em Low dimensional defects in percolation model}.
	\newblock PhD thesis, University of Geneva, 2018.
	
	\bibitem[Kyp14]{ladder}
	A.~E. Kyprianou.
	\newblock {\em Fluctuations of {L}\'evy processes with applications}.
	\newblock Universitext. Springer, Heidelberg, second edition, 2014.
	\newblock Introductory lectures. \MR{3155252}
	
	\bibitem[Law05]{lawler-book}
	G.~F. Lawler.
	\newblock {\em Conformally invariant processes in the plane}, volume 114 of
	{\em Mathematical Surveys and Monographs}.
	\newblock American Mathematical Society, Providence, RI, 2005. \MR{2129588
		(2006i:60003)}
	
	\bibitem[{Le }13]{legall-uniqueness}
	J.-F. {Le Gall}.
	\newblock Uniqueness and universality of the {B}rownian map.
	\newblock {\em Ann. Probab.}, 41(4):2880--2960, 2013.
	\MR{3112934}
	
	\bibitem[LG14]{legall-icm}
	J.-F. Le~Gall.
	\newblock Random geometry on the sphere.
	\newblock In {\em Proceedings of the {I}nternational {C}ongress of
		{M}athematicians---{S}eoul 2014. {V}ol. 1}, pages 421--442. Kyung Moon Sa,
	Seoul, 2014. \MR{3728478}
	
	\bibitem[LP93]{lp93}
	M.~L. Lapidus and C.~Pomerance.
	\newblock The {R}iemann zeta-function and the one-dimensional {W}eyl-{B}erry
	conjecture for fractal drums.
	\newblock {\em Proc. London Math. Soc. (3)}, 66(1):41--69, 1993. \MR{1189091}
	
	\bibitem[LSAP94]{lsp94}
	R.~Langlands, Y.~Saint-Aubin, and P.~Pouliot.
	\newblock Conformal invariance in two-dimensional percolation.
	\newblock {\em Bull. Am. Math. Soc.}, 30:1--61, 1994.
	
	\bibitem[LSW03]{lsw-restriction}
	G.~Lawler, O.~Schramm, and W.~Werner.
	\newblock Conformal restriction: the chordal case.
	\newblock {\em J. Amer. Math. Soc.}, 16(4):917--955, 2003.
	\MR{1992830 (2004g:60130)}
	
	\bibitem[LvF13]{lf13}
	M.~L. Lapidus and M.~van Frankenhuijsen.
	\newblock {\em Fractal geometry, complex dimensions and zeta functions}.
	\newblock Springer Monographs in Mathematics. Springer, New York, second
	edition, 2013.
	\newblock Geometry and spectra of fractal strings. \MR{2977849}
	
	\bibitem[Mie13]{miermont-brownian-map}
	G.~Miermont.
	\newblock The {B}rownian map is the scaling limit of uniform random plane
	quadrangulations.
	\newblock {\em Acta Math.}, 210(2):319--401, 2013.
	\MR{3070569}
	
	\bibitem[MS20]{lqg-tbm1}
	J.~Miller and S.~Sheffield.
	\newblock Liouville quantum gravity and the {B}rownian map {I}: the {${\rm
			QLE}(8/3,0)$} metric.
	\newblock {\em Invent. Math.}, 219(1):75--152, 2020. \MR{4050102}
	
	\bibitem[MS15b]{sphere-constructions}
	J.~{Miller} and S.~{Sheffield}.
	\newblock {Liouville quantum gravity spheres as matings of finite-diameter
		trees}.
	\newblock {\em Ann. Inst. H. Poincaré Probab. Statist.} 55(3):1712--1750.
	
	\bibitem[MS16a]{lqg-tbm2}
	J.~{Miller} and S.~{Sheffield}.
	\newblock {Liouville quantum gravity and the Brownian map II: geodesics and
		continuity of the embedding}.
	\newblock {\em ArXiv e-prints}, May 2016, \arxiv{1605.03563}.
	
	\bibitem[MS16b]{lqg-tbm3}
	J.~{Miller} and S.~{Sheffield}.
	\newblock {Liouville quantum gravity and the Brownian map III: the conformal
		structure is determined}.
	\newblock {\em Probab. Theory Related Fields}, 179(3): 1183--1211, 2021.
	
	\bibitem[MS16c]{ig1}
	J.~Miller and S.~Sheffield.
	\newblock Imaginary geometry {I}: interacting {SLE}s.
	\newblock {\em Probab. Theory Related Fields}, 164(3-4):553--705, 2016.
	\MR{3477777}
	
	\bibitem[MS16d]{ig2}
	J.~Miller and S.~Sheffield.
	\newblock Imaginary geometry {II}: {R}eversibility of {$\operatorname{SLE}\sb
		\kappa(\rho\sb 1;\rho\sb 2)$} for {$\kappa\in(0,4)$}.
	\newblock {\em Ann. Probab.}, 44(3):1647--1722, 2016.
	\MR{3502592}
	
	\bibitem[MSW20]{msw-non-simple}
	J.~Miller, S.~Sheffield and W.~Werner.
	\newblock Non-simple SLE curves are not determined by their range.
	\newblock {\em Journal of the European Mathematical Society}, 22(3):669--716, 2020.
	
	
	
	\bibitem[Pol81]{polyakov-qg1}
	A.~M. Polyakov.
	\newblock Quantum geometry of bosonic strings.
	\newblock {\em Phys. Lett. B}, 103(3):207--210, 1981. \MR{623209 (84h:81093a)}
	
	\bibitem[PY97]{py97}
	J.~Pitman and M.~Yor.
	\newblock The two-parameter {P}oisson-{D}irichlet distribution derived from a
	stable subordinator.
	\newblock {\em Ann. Probab.}, 25(2):855--900, 1997. \MR{1434129}
	
	\bibitem[RS87]{rs87}
	B.~Rodin and D.~Sullivan.
	\newblock The convergence of circle packings to the {R}iemann mapping.
	\newblock {\em J. Differential Geom.}, 26(2):349--360, 1987. \MR{906396}
	
	\bibitem[RS05]{schramm-sle}
	S.~Rohde and O.~Schramm.
	\newblock Basic properties of {SLE}.
	\newblock {\em Ann. of Math. (2)}, 161(2):883--924, 2005.
	\MR{2153402 (2006f:60093)}
	
	\bibitem[RV14]{rhodes-vargas-review}
	R.~Rhodes and V.~Vargas.
	\newblock Gaussian multiplicative chaos and applications: {A} review.
	\newblock {\em Probab. Surv.}, 11:315--392, 2014.
	\MR{3274356}
	
	\bibitem[Sch97]{schaeffer-bijection}
	G.~Schaeffer.
	\newblock Bijective census and random generation of {E}ulerian planar maps with
	prescribed vertex degrees.
	\newblock {\em Electron. J. Combin.}, 4(1):Research Paper 20, 14 pp., 1997. \MR{1465581 (98g:05074)}
	
	\bibitem[Sch00]{schramm0}
	O.~Schramm.
	\newblock Scaling limits of loop-erased random walks and uniform spanning
	trees.
	\newblock {\em Israel J. Math.}, 118:221--288, 2000.
	\MR{1776084 (2001m:60227)}
	
	
	\bibitem[SS13]{ss-level}
	O.~Schramm and S.~Sheffield.
	\newblock A contour line of the continuum Gaussian free field.
	\newblock {\em Probab. Theory Related Fields}, 157 (1-2), 47--80, 2013.
	
	\bibitem[SS11]{ss-planar-perc}
	O.~Schramm and S.~Smirnov.
	\newblock On the scaling limits of planar percolation.
	\newblock In {\em Selected works of {O}ded {S}chramm. {V}olume 1, 2}, Sel.
	Works Probab. Stat., pages 1193--1247. Springer, New York, 2011,
	\arxiv{1101.5820}.
	\newblock With an appendix by Christophe Garban. \MR{2883400}
	
	\bibitem[She07]{shef-gff}
	S.~Sheffield.
	\newblock Gaussian free fields for mathematicians.
	\newblock {\em Probab. Theory Related Fields}, 139(3-4):521--541, 2007.
	\MR{2322706 (2008d:60120)}
	
	\bibitem[She16a]{shef-zipper}
	S.~Sheffield.
	\newblock Conformal weldings of random surfaces: {SLE} and the quantum gravity
	zipper.
	\newblock {\em Ann. Probab.}, 44(5):3474--3545, 2016.
	\MR{3551203}
	
	\bibitem[She16b]{shef-burger}
	S.~Sheffield.
	\newblock Quantum gravity and inventory accumulation.
	\newblock {\em Ann. Probab.}, 44(6):3804--3848, 2016.
	\MR{3572324}
	
	\bibitem[Smi01]{smirnov-cardy}
	S.~Smirnov.
	\newblock Critical percolation in the plane: conformal invariance, {C}ardy's
	formula, scaling limits.
	\newblock {\em C. R. Acad. Sci. Paris S\'er. I Math.}, 333(3):239--244, 2001.
	\MR{1851632 (2002f:60193)}
	
	
	
	\bibitem[SW01]{smirnov-werner-percolation}
	S.~Smirnov and W.~Werner.
	\newblock Critical exponents for two-dimensional percolation.
	\newblock {\em Math. Res. Lett.}, 8(5-6):729--744, 2001.
	\MR{1879816 (2003i:60173)}
	
	\bibitem[SW16]{shef-wang-lqg-coord}
	S.~{Sheffield} and M.~{Wang}.
	\newblock {Field-measure correspondence in Liouville quantum gravity almost
		surely commutes with all conformal maps simultaneously}.
	\newblock {\em ArXiv e-prints}, May 2016, \arxiv{1605.06171}.
	
	\bibitem[Zha08]{zhan-duality1}
	D.~Zhan.
	\newblock Duality of chordal {SLE}.
	\newblock {\em Invent. Math.}, 174(2):309--353, 2008.
	\MR{2439609 (2010f:60239)}
	
\end{thebibliography}
\end{document}